\documentclass[12pt]{amsart}

\usepackage{tikz-cd}
\usepackage{latexsym,amssymb,amsmath}
\usepackage{hyperref}
\usepackage{comment}
\usepackage{mathrsfs}
\usepackage{amsmath,amscd}
\usepackage[enableskew]{youngtab}
\usepackage[all]{xy}
\usepackage{empheq}

\usepackage{amsmath}

\usepackage{empheq}
\usepackage{xcolor}
\definecolor{lightgreen}{HTML}{90EE90}

\newcommand{\cocc}{\co(\cc)}

\newcommand{\sgst}{\sg^{\mtr{st}}}

\newcommand{\dltv}{\delta_V}
\newcommand{\bdlt}{{\bar{\delta}}}

\newcommand{\bpi}{\bar{\pi}}
\newcommand{\mtj}{\mtc L}
\newcommand{\mtjcd}{\mtj_{\cd}}
\newcommand{\lag}{\langle}
\newcommand{\rag}{\rangle}

\newcommand{\evx}{{\ev_X}}

 \newcommand{\tev}{\tilde{ev}}

\newcommand{\mtl}{{  \mtc L}}

 \newcommand{\vsk}{\vskip 0.15cm \noindent}

 \newcommand{\chad}{\ch_{\ad}}
 
 \newcommand{\fpcc}{\fp(\cc)}

\newtheorem{theorem}{Theorem}[section]

\newtheorem{corollary}[theorem]{Corollary}

\newcommand{\vect}{\mtr{Vec}}

\newcommand{\sent}{\mapsto}

\numberwithin{equation}{section}

\newcommand{\ccrm}{\ccr_M}

\newcommand{\qst}{q_*}
\newcommand{\qust}{q^*}\newcommand{\cec}{\mtr{  CE}(\cc)}
\newcommand{\ced}{\mtr{  CE}(\cd)}

\newcommand{\bll}{\blue}

\DeclareMathOperator{\rev}{rev}

  \newcommand\CC{{\mathbb C}} \newcommand\FF{{\mathbb F}}

\newcommand\C{\mathcal{C}}
\DeclareMathOperator{\id}{id}

\DeclareMathOperator{\im}{im}
\DeclareMathOperator{\ev}{ev}
\DeclareMathOperator{\coev}{coev}

\newcommand\groth{{G_0(A)}}

\newcommand{\card}[1]{\# #1}

\excludecomment{verlong}
\includecomment{vershort}
\excludecomment{noncompile}

\usepackage{tikz}
\usetikzlibrary{matrix}
\usepackage{combelow}
\newcommand{\ccb}{\mathcal B}
 \newcommand{\bdelta}{\Delta}
 \newcommand{\bann}{\mtr{Ann}}
\usepackage{amsmath,amsthm,amsfonts,amssymb}
\usepackage{amsmath, amsthm, amssymb, amscd, amsfonts}
\usepackage[all]{xy}
\usepackage{amssymb, amsthm, amsmath}
\usepackage{amsfonts}
\usepackage{color}
\DeclareMathOperator{\ad}{ad} 

\newcommand{\ra}{\rightarrow}

\newcommand{\Z}{\mathbb Z}

\newcommand{\ot}{\otimes}
\newcommand{\co}{\mathcal O}
\newcommand{\xra}{\xrightarrow}
\newcommand{\mtc}{\mathcal}
\newcommand{\cs}{\mtc S}
\newcommand{\onh}{On the other hand}

\newcommand{\lam}{\lambda}
\newcommand{\lb}{\label}
\newcommand{\Lam}{\Lambda}

\newcommand{\ann}{\mtr{Ann}}

\newcommand{\al}{\alpha}
\newcommand{\eps}{\epsilon}

\newcommand{\D}{\Delta}

\newcommand{\ul}{\underline}

\newcommand{\rh}{\rightharpoonup}
\newcommand{\lh}{\leftharpoonup}
\numberwithin{equation}{section}
\newtheorem{lem}[equation]{Lemma}

\theoremstyle{plain}
\newcommand{\ovr}{\overline}
\newtheorem{thm}[equation]{Theorem}
\newtheorem{prop}[equation]{Proposition}
\newtheorem{defn}[equation]{Definition}
\newtheorem{cor}[equation]{Corollary}
\newtheorem{rem}[equation]{Remark}

\newcommand{\dw}{\downarrow}
\newcommand{\uw}{\uparrow}
\newcommand{\dlt}{\delta}

\newcommand{\ch}{\chi}
\newcommand{\mtr}{\mathrm}

\numberwithin{equation}{section}
\newcommand{\ncm}{\newcommand}
\ncm{\np}{\newpage}
\ncm{\ebl}{\end{thebibliography}}
\ncm{\bbl}{\begin{thebibliography}}
\ncm{\chd}{_{ _{\ch}}}
\ncm{\ald}{_{ _{\al}}}
\newcommand{\blam}{\Lam}
\ncm{\cP}{\mathcal{P}}
\ncm{\ei}{e_i}
\ncm{\eij}{e_{i,\;j}}
\ncm{\bt}{\begin{thm}}
\ncm{\bdef}{\begin{defn}}
\ncm{\edf}{\end{defn}}
\ncm{\et}{\end{thm}}
\ncm{\bc}{\begin{cor}}
\ncm{\bl}{\begin{lem}}
\ncm{\el}{\end{lem}}
\ncm{\bpf}{\begin{proof}}
\ncm{\epf}{\end{proof}}
\ncm{\ec}{\end{cor}}
\ncm{\ord}{\mtr{ord}}
\ncm{\er}{\end{rem}}
\ncm{\br}{\begin{rem}}
\ncm{\bn}{\begin}

\ncm{\bp}{\begin{prop}}
\ncm{\ep}{\end{prop}}
\ncm{\bd}{
\begin{document}}
\ncm{\ed}{\end{document}}
\ncm{\beq}{\begin{equation}}
\ncm{\beqn}{\begin{equation*}}
\ncm{\eeq}{\end{equation}}
\ncm{\eeqn}{\end{equation*}}
\ncm{\bea}{\begin{eqnarray}}
\ncm{\eea}{\end{eqnarray}}
\ncm{\beanon}{\begin{eqnarray*}}
\ncm{\eeanon}{\end{eqnarray*}}\ncm{\ek}{\eps|_K}\ncm{\diez}{\#}
\ncm{\bwt}{\bowtie}
\ncm{\cC}{\mtc{C}}\ncm{\cc}{\mtc{C}}
\ncm{\cX}{\mtc{X}}
\ncm{\wt}{\widetilde}
\ncm{\sg}{\sigma}
\ncm{\Rep}{\mathrm{Rep}}
\DeclareMathOperator{\Irr}{Irr}
\ncm{\X}{\mathcal{X}}
\ncm{\cA}{\mathcal{A}}
\ncm{\HKer}{\mtr{HKer}}
\ncm{\LKER}{\mtr{LKer}}
\ncm{\aad}{\mtr{ad}}
\newcommand{\mbf}{\mathbb F}
\ncm{\Dr}{\mtr{D}}
\ncm{\cD}{{\mathcal{D}}}\ncm{\cd}{{\mathcal{D}}}\ncm{\ce}{{\mathcal{E}}}
\ncm{\G}{\mathcal{G}}
\ncm{\Dc}{\mtc{D}}
\ncm{\E}{\mtc{E}}
\ncm{\fp}{\mtr{FPdim}}
\ncm{\Vc}{\mtr{Vec}}
\ncm{\cK}{\mtc{K}}
\ncm{\cM}{\mtc{M}}
\ncm{\cE}{\mtc{E}}
\ncm{\cS}{\mtc{S}}

\newcommand{{\ipr}}{i'}

\DeclareMathOperator{\End}{End}
\ncm{\cop}{\mtr{cop}}
\ncm{\op}{\mtr{op}}
\ncm{\chr}{character }\ncm{\ck}{\mtc{K}}
\ncm{\bw}{\bwt}
\ncm{\hker}{\mtr{HKer}}
\ncm{\bx}{\boxtimes}
\ncm{\blue}{\textcolor[rgb]{.00, .00, 1.00}}
\ncm{\red}{\textcolor[rgb]{1.00, .00, .00}}
\ncm{\green}{\textcolor[rgb]{.50, 0.20, .90}}
\ncm{\bne}{\begin{enumerate}}
\ncm{\ene}{\end{enumerate}}
\ncm{\lker}{\mtr{LKer}}
\ncm{\md}{\medbreak}
\ncm{\rep}{\Rep}\ncm{\ind}{\mtr{ind}}
\ncm{\mdn}{\md\noindent}
\ncm{\dd}{$}
\ncm{\up}{^}
\newcommand{\tcs}{\text}
\newcommand{\mbb}{\mathbb B}
\newcommand{\vs}{\mathbb V}
\newcommand{\sth}{suppose that\;}
\newcommand\rad{\operatorname{rad}}
\newcommand{\itm}{\item}
\newcommand{\dbd}{$$}
\newcommand{\mol}{\mtr{mod}}
 \newcommand{\ro}{\rho}
\newcommand{\irr}{\mathrm{Irr}}
\newcommand{\mbc}{\mathbb C}
\newcommand{\mbs}{\mathbb S}
\newcommand{\mbz}{\mathbb Z}
\newcommand{\ct}{\mtc T}
\newcommand{\sm}{\setminus}
\newcommand{\epl}{^{+}}
\newcommand{\sbsq}{\subseteq}
\newcommand{\sbs}{\subset}
\newcommand{\cco}{\mtr{co}}
\newcommand{\cz}{\mathcal{Z}}
\newcommand{\dual}{^{*}}
\newcommand{\Gm}{\Gamma}
\ncm{\cY}{\mtc{Y}}
\newcommand\ZZ{{\mathbb Z}} 
\newcommand{\bab}{\color{DarkOrchid}{}}
\newcommand{\eab}{\normalcolor{}}
\newcommand{\subs}{\subsection}
\newcommand{\cv}{\mtc{V}}
  \newcommand{\grn}{\green}
\newcommand{\dt}{\delta}

\newcommand{\ccf}{\mathrm{ {CF}(\cc)}}
\newcommand{\cce}{\mathrm{ {CE}(\cc)}}
\newcommand{\cecc}{\mathrm{ {CE}(\cc)}}
\newcommand{\cecd}{\mathrm{ {CE}(\cd)}}
\newcommand{\kk}{\Bbbk}
\newcommand{\otL}{\ot_{L}}
\newcommand{\otl}{\ot_{L}}
\newcommand{\unpsi}{1_{\psi}}
\newcommand{\epsi}{e_{\psi}}
\newcommand{\ephi}{e_{\phi}}
\newcommand{\ech}{e_{\ch}}
\newcommand{\nleftcid}{\text{left normal  coideal subalgebra}}
\newcommand{\dimL}{\dim_{\kk}L}
\newcommand{\cl}{\mtc L}
\newcommand{\mj}{\mtc J}
\newcommand{\tl}{\tilde L}
\newcommand{\tL}{\tilde L}
\newcommand{\tpsi}{\tilde(\psi)}
\newcommand{\tmx}{\tilde{\mtc X}}
\newcommand{\zlh}{\mathrm{ZL}}
\newcommand{\ba}{\mathrm A}
\newcommand{\bv}{\mathrm V}
\newcommand{\zhopf}{\mtc{Z}_{\mtr{Hopf}}}
\newcommand{\lstar}{L^{*}}
\newcommand{\ldstar}{L^{**}}
\newcommand{\mstar}{M^{*}}
\newcommand{\mdstar}{M^{**}}
\newcommand{\lkera}{\lker_{A}}
\newcommand{\mdprime}{M''}
\newcommand{\ldprime}{L''}
\newcommand{\cm}{\mtc M}
\newcommand{\ccm}{\mathcal M}
\newcommand{\cn}{\mathcal N}
\newcommand{\ccn}{\mathcal N}
\newcommand{\rx}{\mtr{Rex}}
\newcommand{\cca}{\ca}
\newcommand{\ih}{\underline{\mtr{Hom}}}
\newcommand{\cih}{\underline{\mtr{coHom}}}
\newcommand{\hm}{\mtr{ {Hom}}}
\newcommand{\cov}{\mtr{coev}}
\newcommand{\rora}{\rho^{\mtr{ra}}}
\newcommand{\rola}{\rho^{\mtr{la}}}
\newcommand{\cx}{\mtc X}
 \newcommand{\cZ}{\cz}
 \newcommand{\ca}{\cA}
 \newcommand{\stat}{\noindent}
 \newcommand{\bfa}{{\bf A}}
 \newcommand{\unu}{\mathbf{1}}
 \newcommand{\barzu}{{\bar {  Z}(\unu)}}
 
\newcommand{\idx}{\id_X}
\newcommand{\lprime}{L'}
\newcommand{\mprime}{M'}
\newcommand{\nat}{ \mtr{{  Nat}}}
\newcommand{\ft}{\mtc F_\lam}
\newcommand{\rhau}{\rightharpoonup}
\newcommand{\lhau}{\leftharpoonup}
\newcommand{\cf}{\mathrm{ {CF}}}

\newcommand{\cfc}{\mathrm{{CF}}(\cc)}
\newcommand{\csu}{\overline{\mathfrak{  C}}}
\newcommand{\cfcc}{\mathrm{ {CF}}(\cc)}
\newcommand{\cfcd}{\mathrm{{CF}}(\cd)}
\newcommand{\cfd}{\mathrm{{CF}}(\cd)}
\newcommand{\czcc}{{\cz(\cc)}}
\newcommand{\czcd}{{\cz(\cd)}}
\newcommand{\czt}{{\cz(\cz(\cc))}}
\newcommand{\enx}{\mtr{  End}}
\newcommand{\runu}{R(\unu)}

\newcommand{\bdfn}{\bn{defn}}
\newcommand{\edfn}{\end{defn}}
\newcommand{\deltax}{\delta_X}
\newcommand{\deltav}{\delta_V}
\newcommand{\repcca}{\rep_\cc(A)}
\newcommand{\xotay}{X \ot_A Y}
\newcommand{\xoty}{X \ot Y}
\newcommand{\votw}{V \ot W}
\newcommand{\votaw}{V \ot_A W}
\newcommand{\dimax}{\dim_AX}
\newcommand{\dimccx}{\dim_\cc(X)}
\newcommand{\dimcca}{\dim_\cc(A)}
\newcommand{\dimccv}{\dim_\cc(V)}
\newcommand{\dima}{\dim_A}
\newcommand{\biga}{A}
\newcommand{\comp}{\mathbb C}
\newcommand{\tehtaa}{\theta_A}
\newcommand{\tetaa}{\theta_A}
\newcommand{\ida}{\id_A}
\newcommand{\hma}{\hm_A}
\newcommand{\hmcc}{\hm_\cc}
\newcommand{\fv}{F(V)}
\newcommand{\fw}{F(W)}
\newcommand{\ota}{\ot_A}
\newcommand{\repza}{\rep_\cc^0(A)}
\newcommand{\epsa}{\eps_A}
\newcommand{\bndefn}{\bn{defn}}
\newcommand{\edefn}{\end{defn}}
\newcommand{\bdefn}{\bn{defn}}

\newcommand{\vld}{V^{*}}
\newcommand{\vldd}{V^{**}}
\newcommand{\xld}{X^{*}}
\newcommand{\xldd}{X^{**}}
\newcommand{\yld}{Y^{*}}
\newcommand{\yldd}{Y^{**}}
\newcommand{\aldu}{A^{*}}
\newcommand{\aldd}{A^{**}}

\newcommand{\ia}{\mtr{i}_A}
\newcommand{\aota}{A\ot A}

\newcommand{\idv}{\id_V}

\newcommand{\ld}{^*}
\newcommand{\repg}{\rep(G)}

\newcommand{\thetav}{\theta_V}

\newcommand{\tta}{\theta_A}

\newcommand{\muv}{\mu_V}
\newcommand{\muw}{\mu_W}

\newcommand{\dimcc}{\dim(\cc)}
\newcommand{\chii}{\chi_i}
\newcommand{\chistar}{\ch_{i^*}}
\newcommand{\chj}{\ch_j}
\newcommand{\chm}{\ch_m}
\newcommand{\chn}{\ch_n}
\newcommand{\dimvi}{\dim(V_i)}
\newcommand{\mtcd}{Q}
\newcommand{\mtca}{\mtc A}
\newcommand{\lamcd}{\lam_\cd}
\newcommand{\fpdimcd}{\fp(\cd)}
\newcommand{\laml}{\lam_L}
\newcommand{\apm}{A//M}
\newcommand{\apl}{A//L}
\newcommand{\repapm}{\rep(\apm)}
\newcommand{\repapl}{\rep(\apl)}
\newcommand{\dimvj}{\dim(V_j)}
\newcommand{\dvi}{\dim(V_i)}
\newcommand{\dvj}{\dim(V_j)}
\newcommand{\sumjtom}{\sum_{j=0}^m}
\newcommand{\sumitom}{\sum_{i=0}^m}
\newcommand{\sij}{s_{ij}}
\newcommand{\sji}{s_{ji}}
\newcommand{\dxj}{d_j}
\newcommand{\dxi}{d_i}
\newcommand{\dimka}{\dim_{\kk}(A)}
\newcommand{\dimk}{\dim_{\kk}}
\newcommand{\blaml}{\blam_L}
\newcommand{\sumjtor}{\sum_{j=0}^r}
\newcommand{\dimkl}{\dim_{\kk}(L)}
\newcommand{\mtcjl}{\mtc J_L}
\newcommand{\vota}{ V\ot A}
\newcommand{\vi}{V_i}
\newcommand{\vj}{V_j}
\newcommand{\dimcd}{\dim(\cd)}

\newcommand{\alij}{\al_{ij}}
\newcommand{\alji}{\al_{ji}}
\newcommand{\rcc}{r_\cc}
\newcommand{\rcd}{r_\cd}
\newcommand{\clsx}{[X]}
\newcommand{\clsy}{[Y]}
\newcommand{\clsz}{[Z]}
\newcommand{\rcdp}{r_{\cd'}}
\newcommand{\sumjtorp}{\sum_{j=0}^{r'}}
\newcommand{\aljm}{\al_{jm}}
\newcommand{\aljn}{\al_{jn}}
\newcommand{\sjm}{s_{jm}}
\newcommand{\smj}{s_{mj}}
\newcommand{\snj}{s_{nj}}

\newcommand{\betaij}{\beta_{ij}}
\newcommand{\betaji}{\beta_{ji}}

 \newcommand{\ip}{i'}
\newcommand{\sumjtoprp}{\sum_{j=0}^{r'}}
\newcommand{\sumjtopr}{\sum_{j=0}^{r}}
 \newcommand{\teh}{\tilde{h}}
\newcommand{\cdp}{\cd'}
\newcommand{\xphii}{X_{\phi(i)}}
\newcommand{\inv}{^{-1}}

\newcommand{\fq}{f_Q}
\newcommand{\tr}{\mtr{tr}}
\newcommand{\rtwone}{R_{21}R}

\newcommand{\ccad}{{\cc_{\mtr{ad}}}}
\newcommand{\ccpt}{{\cc_{\mtr{pt}}}}
\newcommand{\qtr}{quasi-triangular\;}
\newcommand{\trq}{\tr_q}

\newcommand{\repal}{\mtr{Rep}(A//L)}
\newcommand{\lkeravi}{\lker_A(V_i)}
\newcommand{\lkeravj}{\lker_A(V_j)}
\newcommand{\cross}[1][1pt]{\ooalign{%
 \rule[1ex]{1ex}{#1}\cr
 \hss\rule{#1}{.7em}\hss\cr}}
\newcommand{\blml}{\blam_L} 
\newcommand{\phir}{\phi_R}
\newcommand{\kda}{{  \Phi(A)}}

\newcommand{\mtcil}{\mtc{I}_L}

\newcommand{\un}{\unu}
\newcommand{\tfl}{\mtc{T}}
\newcommand{\barzm}{\barz(M)}
\newcommand{\barzn}{\barz(N)}
\newcommand{\ccr}{\mtc R^{\cc}}
\newcommand{\ulc}{\ul{\cc}}

\newcommand{\pimx}{\pi_{M;\;X}}
\newcommand{\pinx}{\pi_{N;\;X}}
\newcommand{\acc}{{\mathrm A_\cc}}
\newcommand{\epsu}{\eps_\unu}

\newcommand{\ob}{\mtr{Obj}}
\newcommand{\obc}{\mtr{Obj(\cc)}}
\newcommand{\ccop}{\cc^{\mtr{op}}}
\newcommand{\mtf}{\mtc F_\lam}
\newcommand{\mtfi}{\mtc F^{-1}_\lam}
\newcommand{\elcd}{\ell_\cd}
\newcommand{\mcid}{\mtc I_\cd}
\newcommand{\mcidp}{\mtc I_{\cd'}}
\newcommand{\wtildelcd}{\widetilde{\elcd}}
\newcommand{\wtildelcdp}{\widetilde{\ell_{\cd'}}}
\newcommand{\cpt}{\cc_{\mtr{pt}}}
\newcommand{\barzr}{\barz_\cd}
\newcommand{\barzv}{\barz(V)}
\newcommand{\acd}{\mathrm A_\cd}
\newcommand{\czrcd}{\cz_\cc(\cd)}
\newcommand{\sml}{\Small}
\newcommand{\bs}{\blue{\Small }}
\newcommand{\yd}{Yetter-Drinfeld}

\newcommand{\sumitor}{\sum_{i=0}^r}
\newcommand{\cdop}{\cd^{\mtr{op}}}
\newcommand{\ccrev}{\cc^{\mtr{rev}}}
\newcommand{\barz}{{\bar{\mathrm Z}}}
\newcommand{\etl}{etale\;}
\newcommand{\czca}{\cz(\ca)}


\bd
\title[Coends]{Conjugacy classes and centralizers for pivotal fusion categories}

\author{Sebastian Burciu}
\address{Inst.\ of Math.\ ``Simion Stoillow" of the Romanian Academy
P.O. Box 1-764, RO-014700, Bucharest, Romania}
\email{sebastian.burciu@imar.ro}
\date{\today}
\maketitle


\begin{abstract}
A criterion for M\"uger centralizer of a fusion subcategory of a braided non-degenerate fusion category is given. Along the way we extend some identities on the space of class functions of a fusion category introduced by Shimizu in \cite{scalg}.  We also show that in a modular tensor category the product of two conjugacy class sums is a linear combination of conjugacy class sums with rational coefficients. 
\end{abstract}

\section{Introduction} One of the main tools in the study braided fusion categories is the notion of the  centralizer of a fusion subcategory that in its full generality was introduced in \cite{dgno2} (see also \cite{proclond}). Given a braided fusion category  $\cc$ and $\cd\subseteq \cc$ a fusion subcategory of $\cc$,  the centralizer $\cd'$ is defined as the full fusion subcategory $\cd'$ of $\cc$ generated by all simple objects $X$ of $\cc$ that centralize any object of $\cd$, i.e. $ c_{X, Y}c_{Y, X}=\mtr{id}_{X\ot Y}$ for all objects $Y \in \co(\cd)$.  The notion of the centralizer was used in important classification results for braided fusion categories, see for example papers  \cite{dgno2,  ENO, DGNO} and the references therein.

In general, the centralizer of a given fusion subcategory is very hard to compute.  There are only few cases in the literature where concrete formulae are known.  For instance, in the category of representations of a (twisted) Drinfeld double of an arbitrary finite group a  formula for the centralizer of an arbitrary fusion subcategory was then given in \cite{nnw}. 

In the paper \cite{new-sub} we have given a complete description for the M\"uger centralizer of a fusion subcategory of a fusion category of the type $\cc=\rep(H)$ for a factorizable semisimple Hopf algebra $H$. More precisely, we have shown that if $L$ is a left normal coideal subalgebra of $A$ and $\cd=\rep(H//L)$ is a fusion subcategory of $\rep(H)$ then 
\beq\label{fm1}
\rep(H//L)'=\rep(H//M)
\eeq where $M=f_Q((H//L)^*)$.  Here $f_Q:H^*\ra H$ is the Drinfeld map associated to $(H, R)$.
 \vsk Denoting by $\lam_L$ the idempotent cointegral of the quotient Hopf algebra $H//L$ this result can be written as
 \beq\label{fm2}
 \mtf (f_Q(\lam_{L}))=\frac{1}{\dimk(M)}\lam_{M}
 \eeq
 where $ \mtf :H^*\ra H$ is the Fourier transform associated to the Hopf algebra $H$ with $\lag \lam, 1\rag=1$.

Recall that recently Shimizu   introduced in \cite{sh-int} the notion of cointegral $\lam_\cc$ of a finite tensor category $\cc$. For a semisimple Hopf algebra $H$, and a fusion subcategory $\cd=\rep(H//L)$ of  $\cc=\rep(H)$ one has that $\lam_\cd=\lam_L$.

The main result of the paper is the following theorem generalizing Equation \eqref{fm2} to arbitrary ribbon categories.
\bt\label{main1}
Let $\cd$ be a fusion subcategory of a ribbon fusion category $\cc$ and $\lam_{\cd}$ be the associated idempotent cointegral of $\cd$. Then
\dbd
{\mtf}(f_Q(\lam_{\cd}))=\frac{\dim(\cd')}{\dimcc}\lam_{\cd'}.
\dbd
where $\mtf$ is the Fourier transform introduced by Shimizu in \cite{scalg} and $\lam_{\cd'}$ is the idempotent cointegral of $\cd'$.
\et

Recently, in \cite{scalg},  Shimizu also extended the notion of conjugacy classes for fusion categories, similarly to the conjugacy classes introduced  by Cohen and Westreich in \cite{CW2} for semisimple Hopf algebras. They generalize the notion of conjugacy classes in finite groups. In the same paper \cite{scalg} the author associated to each conjugacy class a central element called conjugacy class sum.  These class sums play the role of the sum of group elements in a conjugacy class of a finite group.

In \cite{scalg} the author asked if the results from \cite{CW4, CW6} concerning conjugacy classes can be extended from semisimple Hopf algebras to fusion categories. For example, for semisimple factorizable Hopf algebras,  Cohen and Westreich in \cite{CW2} proved that prove that a product of two conjugacy class sums is a linear combination with rational coefficients of all conjugacy class sums. In subsection \ref{pcc} we show that this result also  holds for modular tensor categories. We also extend some other results of \cite{CW2} from semisimple factorizable Hopf algebras to modular tensor categories. 

For any fusion subcategory $\cd$ of a spherical fusion category $\cc$ we exhibit an element $\ell_\cd\in \cecd$ which plays the role of the integral of a coideal subalgebra $L$ of a semisimple Hopf algebra $H$. More precisely, if $\cd=\rep(H//L)$ is a  fusion subcategory of $\cc=\rep(H)$, then $\ell_\cd=\blam_L$, the unique element of $L$ such that $l\blam_L=\eps(l)\blam_L$ for any $l\in L$, see \cite{sk}. Using this element and a Class Equation type argument in Theorem \ref{lcad} we show that fusion subcategories of prime index $p$ are in bijection with normal subgroups of the same prime index $p$ of the universal grading group $U_\cc$ of $\cc$. 

This paper is organised as follows. In Section \ref{pr} we review the basic notions of fusion categories that are needed through the paper. In Section \ref{adj} we recall the notion of central Hopf comonads and adjoint algebras associated to fusion categories from \cite{scalg}. This section also contains the main properties of the (co)integrals and Fourier transforms associated to fusion categories. In Section \ref{cfc} we prove several identities inside the ring $\cfcc$ of class functions that are needed through the rest of the paper.  In Section \ref{br} we prove the main result of the paper mentioned above.  In the last section we study conjugacy classes for modular tensor categories. In this section we also prove that a product of two conjugacy class sums is a linear combination with rational coefficients of all conjugacy class sums. In the Appendix we prove that a canonical natural transformation between certain coends is a morphism of Hopf comonads, result needed in the proof of the main result.

Throughout this paper we work over an algebraically closed field $\kk$ of  arbitrary characteristic.
\section{Preliminaries}\label{pr}
In this section we review the main properties of fusion categories that are needed through the paper. 
 For the basic theory of monoidal categories, we refer the reader to \cite{ML98} and \cite{Kas}.

For a monoidal category $\cc = (\cc,\ot, \unu)$ with tensor product 
$\ot$ and unit object $\unu$, we set $\ccop = (\ccop,\cc, \unu)$ and $\ccrev = (\cc,\ot^{\rev}, \unu)$, where $ \ot^{\rev}$ is the reversed
tensor product given by $X \ot^{\rev} Y = Y \ot X$.
A monoidal functor \cite[XI.2]{ML98} is a functor $F : \cc \ra \cd$ between monoidal
categories $\cc$ and $\cd$ endowed with a morphism \dd F_0 : \unu \ra F(\unu)$ in $\cc$ and a natural
transformation $F_2(X, Y ) : F(X) \otimes F(Y ) \ra F(X \otimes Y )$ ($X, Y \in \cocc$) satisfying a certain coherence condition. We say that a monoidal functor $(F,F_2, F_0)$ is {\it strong monoidal} 
if $F_0$ and $F_2$ are invertible. We also say that $F$ is strict if the above monoidal structures are identities. A {\it comonoidal functor} from $\cc$ to $\cd$ is the same thing as a monoidal functor from $\ccop$ to $\cdop$.

A {\it left dual object} of $X \in \cocc$ is a triple $(X', e, d)$ consisting of an object $X' \in \cocc$
and morphisms \dd e : X' \otimes X \ra \unu$ and \dd d : \unu \ra X \otimes X'$ such that
\beqn
(e \otimes id_{X'}) \circ (id_{X'} \ot d) = id_{X'}\;\;\text{ and}\;\; (id_X \ot e) \circ (d \ot id_X) = id_X.
\eeqn

We say that $\cc$ is {\it left rigid} if every object of $\cc$ has a left dual object. For each \dd X \in\cocc$, we choose a left dual object denoted by $(X^*, \ev_X, \cov_X)$ with the morphisms  $\ev_X:X^*\ot X\ra \unu$ and $\cov_X:\unu \ra X\ot X^*$ called the evaluation and the coevaluation respectively. Then the assignment $X \mapsto X^*$ gives rise to a strong monoidal functor
$(-)^* : \cc \ra \cc^{\op, \rev}$, which we call the {\it  left duality functor}. Similarly,  a {\it right dual} of $X$ as is triple $(X'', e', d')$ where $X' \in \cocc$ and the morphisms \dd e' : X \otimes X' \ra \unu$ and \dd d' : \unu \ra X '\otimes X$  are defined such that the triple $(X'', e', d')$ is a left dual of $X$ in $\cc^{rev}$. 

A rigid monoidal category $\cc$ is a monoidal category $\cc$ such that both $\cc$ and $\cc^{\rev}$ are left rigid. If $\cc$ is a rigid monoidal category, then the left duality functor $(-)^*$ of $\cc$ is an equivalence. A quasi-inverse of $(-)^*$, denoted by $^*(-)$, is called the right duality functor.

\br\label{sh15-1} As explained in \cite[Lemma 5.4]{sh15} it may be assumed that these duality functors are strict and mutually inverse one to another. Thus one can assume through the rest of the paper that:
$$(X\ot Y)^*=Y^*\ot X^*, \un^*=\unu,\;\;^*(X^*)=X.$$
\er
A  finite abelian category  is a $\kk$-linear category that is equivalent to the category $A$-mod of finite dimensional $A$-modules  for some finite-dimensional $\kk$- A {\it finite tensor category } \cite{EO} is a rigid monoidal category $\cc$ such that
\bne
\item $\cc$ is a finite abelian category,
\item The tensor product $\ot:\cc \times \cc \ra\cc$ is $\kk$-linear in each variable,
\item $\enx_\cc(\unu)\simeq\kk$ as algebras.
\ene
A fusion category \cite{ENO} is a semisimple finite tensor category.
\subsection{Braided  categories} By definition, the reversed category $
\cc^{rev}$ of a monoidal category $\cc$ is the same underlying category with the tensor product reversed. \vsk A braiding \cite[XIII.1]{Kas} of a monoidal category $(B, \ot_{\ccb}, \unu)$ is a natural isomorphism $\sg:\ot_B\ra \ot^{\rev}_{\ccb}$ satisfying the hexagon axiom. A braided monoidal category is a monoidal category endowed with a braiding.  The reversed braided category $\cc^{\rev}$ of a braided category $\cc$ is defined as the reversed monoidal category $\cc$ with the braiding given by:
$$\sg^{\rev}_{X,Y}:=\sg_{Y, X}.$$
\subsection{Monoidal centers}
For a monoidal category $\cc$, the monoidal center (or the Drinfeld center) of $\cc$ is a category $\czcc$ defined as follows: An object of $\czcc$ is a pair $(V, \sg_V)$ consisting of an object $V\in \cc$ and a natural isomorphism
$$\sg_{V, X}: V \ot  X \ra X \ot V $$
for all $X \in \cocc$, satisfying a part of the hexagon axiom. A morphism  $f:(V, \sg_V)\ra (W, \sg_W)$ in $\czcc$
 is a morphism in $\cc$  such that $(id_X\ot f) \circ\sg_{V, X}=\sg_{W, X}\circ(f\ot id_X)$ for all $X\in \cc$. The composition of morphisms is defined in an obvious way. The category $\czcc$ is in fact  a braided monoidal category, see, e.g., \cite[XIII.3]{Kas} for details.
\subsection{Pivotal tensor categories}
Recall that a pivotal structure of a rigid monoidal category $\cc$ is an isomorphism $j:\id_\cc\ra ()^{**}$ of monoidal functors. A pivotal monoidal category is a rigid monoidal category endowed with a pivotal structure. 
\br\label{sh15} As explained  in \cite{sh15} we may  assume further that in  a pivotal category one has $X^{**}=X$, i.e.  the pivotal structure is the identity.
\er

A pivotal structure a on a tensor category C is called {\it spherical} if
$\dim(V)=\dim(V^*)$ for any object $V \in \co(\cc)$. A tensor category is spherical if it is
equipped with a spherical structure.
\subsection{Ribbon categories} 
A braided category is called {\it pre-modular} if it  has a spherical structure. Equivalently, this  is a ribbon fusion category, that is, a fusion category equipped with a braiding and a twist
(also called a balanced structure), see \cite{ENO}.
\vsk
Assume $\kk=\mathbb C$. Then $\cc$ is called {\it pseudo-unitary} if $\fpcc =\dimcc$. If such is the case,
then by \cite[Proposition 8.23]{ENO}, $\cc$ admits a unique spherical structure with respect to
which the categorical dimensions of simple objects are all positive. It is called the canonical spherical structure. For this structure, the categorical dimension of an object coincides with its Frobenius-Perron dimension, i.e. $\fp(X)=\dim(X)$ for any object $X\in \co(\cc)$. If $\cc$ is a fusion category  such that  $\fpcc$ is an integer, then $\cc$ is pseudo-unitary by \cite[Proposition 8.24]{ENO}. Moreover, every full fusion
subcategory of $\cc$ is pseudo-unitary.
\subsection{Dinatural transformations, ends and coends.}
Let $\cc $ and $\cd$ be two categories and $S,T:\cc^{\op}\times \cc\ra \cd $ be two functors. A dinatural transformation $\xi:T\xra{:} S$ is a family of morphisms 
$$
\xi=\{\xi_X:T(X,X)\ra S(X,X)\}_{X\in \cocc}$$ such that:
$$S (\id_X, f) \circ \xi_X \circ T(f, \id_X) = S (f, \id_Y ) \circ \xi_Y \circ T(\id_Y , f)$$ for all morphisms $f:X\ra Y$ in $\cc$. An {\it end } of the functor $S$ is a pair $(E, p)$ consisting of an
object $E \in \co(\cd)$ (regarded as a constant functor from $\ccop\times \cc \ra \cd$ ) and a dinatural
transformation $p : E \xra{:} S$ that enjoys the following universal property: For any pair $(E', p')$ consisting of an object $E' \in \co(\cd)$ and a dinatural transformation $p' : E' \xra{:} S$,
there exists a unique morphism $\phi : E' \ra E$ in $\cd$ such that $p'_X
= p_X \circ \phi$ for all $X \in \cc$. The end of $S$ is expressed as $\int_{X\in \cc}S(X, X)$.

A {\it coend} of $T$ is a pair $(\cc, i)$ consisting of an object $\cc \in \co(\cd)$ and a dinatural transformation $i : T\xra{:} \cc$ having a similar universal property. The  coend of $T$ is  expressed as $\int^{X\in \cc}T(X, X)$.

See \cite[IX]{ML98} for the basic results on (co)ends.
\section{Adjoint algebra and internal characters}\label{adj}
\subsection{Hopf monads and comonads} Let $\cc$ be a monoidal category. Recall, \cite{bv07} that a {\bf bimonad} $T:\cc \ra \cc$ is a monad $(T, \mu, \eta)$   with {\it  comonoidal structures}
$$T_0:T(\unu)\ra \unu, \;\; T_2(V, W):T(V\ot W)\ra T(V)\ot T(W)$$
such that $\mu$ and $\eta$ are {\it comonoidal natural transformations.}  
 In this case the category of $T$-modules $_T\cc$ is a monoidal category and the forgetful functor $U:\;_T\cc\ra \cc$ is a strong monoidal functor.

If $\cc$ is a rigid monoidal category then a Hopf monad on $\cc$ is defined as a bimonad for which the category $\;_T\cc$ is also rigid. Recall that in this case $T$ is endowed with left and right bijective antipodes $S_V^l:T(T(V)^*)\ra V^*,\;\; S^r_V:T^*(T(V))\ra\;^*V$, see \cite{bv07}.

Also recall that a Hopf comonad is a comonad on $\cc$ endowed with a monoidal structure such that $T^{op}:\cc^{\op}\ra \cc^{op}$ is a Hopf monad.
\subsection{The central Hopf monad $\mtl$}
Let $\cc$ be a finite tensor category and $F:\czcc\ra \cc$ the forgetful functor. Then $F$ admits a left adjoint $L$ and $\mtl:=FL:\cc \ra \cc$ is a Hopf monad defined as a coend. As in \cite[Section 3.1]{scalg} one has that
\beq
\mtl(V)\simeq \int^{X\in \cc}X^*\ot V\ot X
\eeq
We denote by $\iota_{V; X}: X^*\ot V\ot X\ra \mtl(V)$ the universal dinatural transformation associated to the coend $\mtl(V)$.

Day and Street \cite{ds07} showed that the functor $V\mapsto \mtl(V )$ has a structure of a monad such that the category $_\mtl\cc$ of $\mtl$-modules is canonically isomorphic to the monoidal center $\czcc$.
\subsection{The central Hopf comonad $Z $}
The forgetful functor $F$ also admits a right adjoint functor $R:\cc \ra \czcc$  and $Z :=FR:\cc \ra \cc$ is a Hopf comonad defined as an end.  Indeed, following \cite[Section 2.6]{scalg} one has that 
\beq
Z(V)\simeq \int_{X\in \cc}X\ot V\ot X^*
\eeq
We denote by $\pi_{V;X}:Z (V)\ra X\ot V\ot X^*$ the universal dinatural transformation associated to the end $Z (V)$.

Recall that for a functor $T:\ccb \ra \cc$ between rigid monoidal categories, we may  define $T^!:\ccb \ra \cc$ by
$T^!(V ) = T(^*V)^*$, (see \cite[Subsection 2.6]{scalg}). Then one can take $R=L^!$ and $\pi_{V;X}=(i_{^*V; ^*X})^*$.

The Hopf comonad structure of $Z $ can be described in terms of the dinatural transformation $\pi$. Indeed, the
comultiplication $\delta : Z  \ra {Z }^2$ is the unique natural transformation such that
\beq\label{cobz}
(\id_X \otimes \pi_{V ;Y}\otimes \id_{X^*} ) \circ \pi_{Z (V );X} \circ \delta_V = \pi_{V ;X\otimes Y}
\eeq
The counit of $\eps:Z\ra \id_\cc $ is given by $\eps_V:= \pi_{V ;1}.$


\subsection{Adjoint algebra} It is known that $A:=Z(\unu)$ has the structure of central commutative algebra in $\cz(\cc)$.

The multiplication $m:A\ot A \ra A$ and the unit $u:\unu\ra A$ of the
adjoint algebra $A = Z (\unu)$ are uniquely determined by  by the universal property of the end $A$ as:
\beq \pi_{\unu;X} \circ u = \cov_X,
\eeq
\beq \pi_{\unu;X} \circ m = (\id_X \otimes \ev_X \otimes \id_{X^*} ) \circ (\pi_{\unu;X} \otimes \pi_{\unu;X})
\eeq

Moreover $\epsu:A\ra \unu$ is a morphism of algebras, see \cite{scalg}. Moreover it is well known that $A$ is a commutative algebra in the center $\czcc$. 
\subsection{Central elements}
The vector space $\cecc:= \hm_{\C}(\unu, A) $ is called {\it the set of central elements.} There is a canonical bijection:
{\Small
\beq\label{canisom}
\cce\xra{\psi} \enx(\id_{\C}),\;\;\psi(a)_X=\ro_X(a\ot \\id_X),\;\; X=\unu\ot X\xra{a\ot \\id_X}X\ot Z(\unu)\xra{\ro_X} X
\eeq
}
The inverse of this bijection is described using ends, see \cite[Proposition 5.2.5]{kl01}. 

\subsubsection{Multiplication in the center}
For $a, b \in \cecc$ , we set $a ·b :=m \circ (a \ot b)$. Then the set $\cecc$ is a monoid with respect to this operation. Moreover, the bijection in Equation \eqref{canisom} is in fact an isomorphism of monoids.

For a natural transformation $\al:S\ra T$ between two functors $S, \;T:\ccb \ra \cc$  we define
$\al^{!} : T^{!} \ra S^{!} $ by 
$$\al_V : T^!(V ) = T(^*V)^* \xra{(\al_{^*V} )^*}  S(^*V )^* = S^!(V) \;\;(V \in \co(\ccb)).$$

The antipodal operator  $\mtc S:\cecc\ra \cecc$ on $\cecc$ is induced by
\beqn
(-)^! : \enx(\id_\cc) \ra  \enx(\id_\cc),\; \xi\mapsto \xi^!,
\eeqn
via the bijection from Equation \eqref{canisom}. 

\subsection{Internal characters of pivotal tensor categories}
Recall that a pivotal structure $j$ on a tensor category $\cc$ is a tensor isomorphism $j:\id_\cc\ra ()^{**}$.  Using the pivotal structure one can construct a {\it right evaluation} as follows:
$$\widetilde{ev}_X:X\ot X^*\xra{j\ot id}X^{**}\ot X^*\xra{ev_{X^*}} \unu$$ 
Then {\it the right partial pivotal trace} of $f:A\ot X\ra B\ot X$ is defined as follows:
\beqn
\tr_{A, B}^X:A=A\ot \unu\xra{\id \ot coev_X}A\ot X \ot X^*\xra{f \ot id} B\ot X\ot X^*\xra{\\id \ot \widetilde{ev}_X} B
\eeqn
The usual {\it right pivotal trace} of an endomorphism $f :X\ra X$ is obtained as a particular case for $A=B=\unu$. In particular,  the {\it right pivotal dimension of $X$} is defined as the right trace of the identity of $X$.
\subsubsection{Internal characters and their multiplication}
Given an object $X \in \co(\cc)$ the internal character $\mtr{ch}(X)$ is defined as the  morphism \beqn
\mtr{ch}(X):=\tr^{X}_{A, \unu}(\rho_{X}):A\ra \unu.
\eeqn
Then the space $\cfcc:=\hm_\cc(A_\cc, \unu)$ is called the {\it space of class functions} of $\cc$.
 
  For two class functions $f, g\in \cfcc$ one can define a multiplication by
 $$f\star g:=f \circ Z(g) \circ \delta_{\unu}.$$ 
 Here $\delta: Z \ra  Z^2$ is the comultiplication structure of $Z$ defined in the Equation \eqref{cobz}.
 This multiplication coincides to the composition of morphisms in the co-Kleisli category of the comonad $Z$. Thus $\cfcc$ is a monoid with respect  to $\star$.
 
 By \cite[Theorem 3.10]{scalg} one has that $\mtr{ch}(X\ot Y)=\mtr{ch}(X)\mtr{ch}(Y)$ for any two objects $X$ and $Y$ of $\cc$. 
For  a finite tensor category  $\cc$ the space of class functions 
 $\cfcc$ is a finite-dimensional algebra.

By \cite[Theorem 4.1]{scalg} if $\cc$ is a finite pivotal tensor category over an algebraically closed field $\kk$ then the  set of irreducible characters $\ch(X)$ with $X \in \irr(\cc)$ is linearly independent in $\cfcc$.

Recall $R:\cc \ra \czcc$ is a right adjoint to the forgetful functor $F:\czcc \ra \cc$.
As explained in \cite[Theorem 3.8]{scalg} this adjunction  gives an isomorphism of monoids
\beq\label{adjisom}
\cfcc \xra{\cong} \mtr{End}_{\czcc}(R(\unu)),\;\; \ch\mapsto Z (\ch)\circ \delta_\unu.
\eeq

In the theory of semisimple Hopf algebras, it is well-known that the evaluation pairing between characters and central elements is non-degenerate. Shimizu has generalized this fact, see \cite{scalg},   by considering a paring
$\langle\;,\; \rangle : \cfcc \times \cecc\ra \unu$, given by $ \langle f, a\rangle  \id_{\unu}= f \circ a,$
for all $f \in \cfcc$ and $a\in \cecc$. 
\subsection{On the inclusion of class functions of a fusion subcategory}\label{cfinc}Let $\cd$ be a fusion subcategory of a given fusion category $\cc$.  For any $V\in \cc$ consider as in \cite[Sect. 4.3]{scalg} also the end 
$$\barz(V):\cd^{\op}\times \cd\ra \cc,\;\; \barz(V):=\int_{X\in \cd}  X\ot V\ot X^*.$$ 
\vsk Let $\bar{\pi}_{V;\; X}:\barz(V) \ra X\ot V\ot X^*$ be the universal dinatural maps defining this end. From the universal property of $\barz$ for any $V\in \cc$ there is a unique canonical map in $\cc$:
\beq
Z(V)\xra{q_V}\barz(V)
\eeq
\noindent such that $\bar{\pi}_{V;\; X}\circ q_V={\pi}_{V;\; X}$ for any  object $X$ of $\cd$.

In the appendix we prove that $q$ is a map of Hopf comonads. The map $q_\unu:Z(\unu)\ra \barzu$ induces two maps $$q_\unu^*:\cfcd \ra \cfcc, \;\ch\sent \ch \circ q_\unu,$$
\vsk $${q_{\unu}}_*:\cecc \ra \cecc, \;z\sent q_\unu \circ z.$$ 

Moreover, by Lemma \ref{rvq} from appendix it is known that $q_\unu^*$ is a monomorphism and ${q_\unu}_*$ is an epimorphism and both maps  are in fact  $\kk$-algebra homomorphisms.

\br Since $\cc$ is a fusion category one can give a direct description of the map $q_\unu^*$. Indeed, there is a canonical isomorphism  $ch^\cc:\mtc{G}r_\kk(\cc)\xra{\simeq} \cfcc, [X]\sent \ch^\cc(X)$ where $\mtc{G}r_\cc$ is the Grothendieck group of $\cc$,  see \cite[Example 4.4]{scalg}. For a fusion subcategory $\cd\subseteq \cc$  clearly we we have $\mtc{G}r_\kk(\cd)\hookrightarrow \mtc{G}r_\kk(\cc)$ as $\kk$-algebras. It is easy to see that this induces an inclusion $\cfcd\hookrightarrow \cfcc$ via the following commutative diagram:
{\center
\begin{tikzcd}[column sep=large,row sep=large]
\mtc{G}r_\kk(\cc) \arrow[r,"ch^\cc"] \arrow[d,"j"]&   \cfcc \arrow[d,"q^*"]
 \\ \mtc{G}r_\kk(\cd) \arrow[r,"ch^\cd"]& \cfcd.
\end{tikzcd}
}
\vsk
The commutativity of the diagram follows since for any $X\in \co(\cd)$ one has $$q^*(\ch^\cd(X))=\ch^\cd(X)\circ q=\tilde{ev}_X\circ \bpi_{\unu, X}\circ q=\tilde{ev}_X\circ  \pi_{\unu, X}=\ch^\cc(X).$$
\er
 \subsection{Unimodular tensor categories}
Let $\cc$ be a finite tensor category and $F:\czcc\ra \C$ be the forgetful functor. It is well known that $F$ has both left and right adjoints which are denoted by $L$ and respectively $R$. 

Etingof, Nikshych and Ostrik \cite{enor} introduced a distinguished invertible object $D\in \cc$ of a finite tensor category $\cc$ over $\kk$. The object $D$ is a category-theoretical analogue of the modular function (also called the {\it distinguished grouplike element}) of a finite-dimensional Hopf algebra, and therefore we say that $\cc$ is {\it unimodular} if $D\simeq \unu$.

By \cite{sh-imrn} a finite tensor category is  unimodular if and only if $R\simeq L$. It is well known that  any fusion category is unimodular, \cite{enor}.

\subsection{Integrals, cointegrals and Fourier transform for finite unimodular tensor categories}\label{ft}
Let $\cc$ be a {\it unimodular} finite tensor category and $A=Z(\unu)$ be its adjoint algebra as defined above.
 An  integral in $\C$ is a morphism $\Lambda: 1 \ra A$   in $  \C$   such that
 \beqn
 m \circ (\id_{A}\otimes \Lambda)=\eps_{1}\ot \Lam.
 \eeqn
\noindent
A  cointegral in $  \C $  is a morphism $  \lambda : A\ra 1 $  such that
\beqn
\bar Z(\lam)\circ \delta_{1}=u \otimes \lam
\eeqn
where $u:\unu \ra A$ is the unit of the algebra $A$. It is well known that the integral and cointegral of a finite unimodular tensor category are unique up to a scalar.

There is also a right action denoted by $\lh$ of $\cecc$ on $\cfcc$ given by 
\beqn
f \lh b=f \circ m \circ (b\ot \id_{A})
\eeqn
for all $f \in \cfcc$ and $b \in \cecc$. 

It is easy to check that $(f \lh b)\lh b'=f \lh (bb')$, thus $\cfcc$ becomes naturally a right $\cecc$-module. Similarly one can define a left action of $\cecc$ on $\cfcc$.

Let  $\lambda\in \cfcc$ be  a non-zero integral of $\cc$. The {\it Fourier transform} of $\cc$  associated to $\lambda$ is the linear map
\beq
\mtc F_{\lambda}:\cecc\ra \cfcc\;\;\text{given by}\;\;a \mapsto \lambda \lh \mtc S(a)
\eeq
\noindent
where $\mtc S: \cecc\ra \cecc$ is the above antipodal operator. The Fourier transform is a bijective $\kk$-linear map whose inverse is given in  \cite[Equation 5.17]{scalg}.
\subsection{The case of fusion categories}
For the rest of this section, suppose that $\cc$ is a pivotal fusion category over an algebraically closed field $\kk$. Furthermore, let $\{V_0, \dots, V_m\}$ be a complete set of representatives of isomorphism classes of simple objects with $V_0=\unu$. For $i \in \{0, \dots, m\}$, we define $i^*\in \{0, \dots, m\}$ by $V_i^*\simeq V_{i^*}$. Then $i \sent i^*$ is an involution on $\{0, \dots, m\}$.
As an object of $\cc$, the adjoint algebra decomposes as
\beq
A\simeq \bigoplus_{i=0}^r V_i\ot V_{i}^*
\eeq

Shimizu has defined in  \cite{scalg}  the elements
$$E_i:\unu \xra{coev} V_i\ot V_i^* \hookrightarrow A,\;\; \ch_i: A \xra{\pi_i} V_i\ot V_i^* \xra{\tev} \unu$$
It is easy to see that $\{E_i\}_{i=0, \dots , m}$ and $\{\ch_i\}_{i=0, \dots , m}$ are bases for $\cecc$ and $\cfcc$ respectively, such that
$$\lag\ch_i, E_j\rag=\delta_{i, j}.$$

Moreover, $E_iE_j=\delta_{i,j}$ and $\mtc S(E_i)=E_{i^*}$, where $\tilde S:\cecc \ra \cc$ is the antipodal map defined above. Note that under the isomorphism Equation \eqref{canisom} $E_i$ corresponds to the natural transformation that is identity on $V_i$ and zero on the other simple objects.

By \cite[Equation 6.8]{scalg} one has that 
\beq\label{intregform}\lam_{\cc}=\frac{1}{\dimcc}(\sum_{[V_i]\in \irr(\cc)}\dim(V_i^*)\ch_{i}).
\eeq 

\subsection {Examples from Hopf algebras}\label{ha}
Let $H$ be a semisimple Hopf algebra and let $\cc=\rep_\kk(H)$ be the fusion category of its finite dimensional representations over $\kk$.  As in \cite[Section3.7]{scalg} we identify the left center of $\cc$ with the category $^H_H\mtc{YD}$ of left-left Yetter Drinfeld modules of $H$. Then the right adjoint $R:\czcc\ra \cc$ can be written as $R(V)=H\ot V$ with $h.(l\ot v)=h_1lS(h_3)\ot h_2v$ and $\delta(h\ot v)=h_1\ot (h_2\ot v)$. Thus $Z(V)=R(V)=H\ot V$ as $H$-modules and the universal diantural maps are given by
\beq\label{endgen}
\pi_{M;\; V}:Z(V)\ra M\ot V\ot M^*,\;\; h\ot v\mapsto \sum_ih_1m_i\ot h_2v\ot m_i^*
\eeq
In particular one has that $Z (\unu)=H$  whose structure as $H$-module is given by  the left  adjoint action $h.a=h_1aS(h_2)$. The multiplication and the unit of $Z (\unu)$ are given the usual multiplication and unit of $H$. The universal dinatural maps are given by: $$\pi_{\unu;\;M}:H\ra  M\ot M^*,\;\; h\mapsto \sum_i  hm_i\ot m_i^*$$ and the comultiplication map $\delta_\unu: Z (\unu)\ra H\ot Z (\unu)$ coincide to the usual comultiplication of $\Delta:H\ra H\ot H$.

Moreover, in this case  $\cfcc=C(H)$, the character ring of $H$ and $\cecc=Z(H)$, the center of $H$. A categorical cointegral $\lam\in \cfcc$ is the same as a Hopf algebra (two-sided) cointegral $\lam\in H^*$ satisfying $\lam(h_1)h_2=\lam(h)1$ for all $h\in H$. The Fourier transform becomes as usually, $\mtf:Z(H)\ra C(H), z\sent \lam\lh S(z)$ for all $z \in Z(H)$.
\section{Some identities in the space of class functions}\label{cfc}
Throughout this subsection $\cc$ is a pivotal fusion category over an algebraically closed field $\kk$ of  arbitrary characteristic. Recall that by \cite[Lemma 6.2]{scalg} the Grothendieck ring $\mathrm{Gr}_\kk(\cc)$ is a symmetric Frobenius algebra with the trace $\tau:\mathrm{Gr}_\kk(\cc)\ra \cc$ given by 
$[X]\mapsto \dimk\hm_\cc(\unu, X)$. 

Recall also by \cite[Section 4]{scalg} that in the case of fusion categories the canonical map $\mathrm{Gr}_\kk(\cc)\simeq \cfcc, \;[X]\sent \mtr{ch}(X)$ is in fact an isomorphism of $\kk$ algebras.

We let  $\irr(\cc)=\{[V_0], \dots [V_m]\}$ be the set of isomorphism classes of simple objects of $\cc$ with $V_0=\unu$. In this case the characters $\ch_i:=\mtr{ch}(V_i)$ form a $\kk$-linear basis for $\cfcc$. We also denote by $d_i=\dimvi $ the quantum dimension of each simple object $V_i$.

It follows by \cite[Equation (6.4)]{scalg} that
\beq\label{Eq (6.4)}
\tau(\ch_i)=\delta_{i,0}=\langle \ch_i, E_0\rangle  \eeq
for all irreducible characters $\ch_i$.  This shows that 
\beq\label{recoverform}
\tau(\ch)=\langle \ch, E_0\rangle,\;\; \text{for any}\;\;\ch \in \cfcc,
\eeq since the irreducible characters $\chii$ form a basis on $\cfcc$.  Also we have that $
\;\;\tau(\chii\star \chj)=\delta_{i,\;j^*}$ for all irreducible characters $\ch_i$. Note that $E_0$, the idempotent associated to the unit object $\unu$ of $\cc$, is an integral $\blam_\cc \in \cecc$, see \cite[Lemma 6.1]{scalg}.

\subsection{Conjugacy classes of fusion categories}
Let $\cc$ be a pivotal fusion category over an algebraically closed field of arbitrary characteristic.
The global dimension of $\cc$ is defined as
 \beq
\dimcc=\sumitom \dimvi \dim(V_{i^*})\in \kk
\eeq
Recall by  \cite[Definition 9.1.]{ENO} that a pivotal finite tensor category is called {\it non-degenerate} if $\dim(\cc)\neq 0$ in $\kk$. By \cite[Theorem 6.6.]{scalg}  a pivotal fusion category $\cc$ is non-degenerate if and only if one of the following holds:

(1)$\mathrm{Gr}_\kk(\cc)$ is a semisimple algebra.

(3)$\czcc$ is a semisimple abelian category.

(4)$R(\unu) \in \czcc$ is a semisimple object.

\noindent
Recall that $R:\cc \ra \czcc$ is the right adjoint of the  forgetful functor $F:\czcc \ra \cc$.
 Note also that in this case by \cite[Theorem 6.6]{scalg} the object $R(\unu) \in \co(\czcc)$ is multiplicity-free since $\mathrm{Gr}_\kk(\cc)$ is a commutative ring.

Note that in this case if $\cd$ is a fusion subcategory of $\cc$ then $Gr_\kk(\cd)\subset \mtr{Gr}_\kk(\C)$ is also a semisimple ring since $\mtr{Gr}_\kk(\C)$ is a semisimple commutative ring. Therefore, by above, $\cd$ is also non-degenerate.

Let as above $F:\cz(\cc)\ra \cc$  be the forgetful functor with its right adjoint $R:\cc \ra \czcc$.  A {\it conjugacy class} of $\cc$ is defined as a simple subobject of $R(\unu)$ in $\czcc$. Since the monoidal center $\czcc$ is also a fusion category we can write $R(\unu)=\bigoplus_{i=0}^m\mathfrak C_i$ as a direct sum of simple objects in $\czcc$. Thus  $\mathfrak C_{0},\dots, \mathfrak C_{m}$ are the conjugacy classes of $\C$.  Since the unit object $\unu_{\czcc }$ is always a subobject of $R(\unu)$, we can assume $\mathfrak C_{0} = \unu_{\czcc }$.

Let $\tilde F_0, \tilde F_1, \dots, \tilde F_m\in \enx_{\czcc}(R(\unu))$ be the canonical projections on each of these conjugacy classes. Let also $F_0, F_1, \dots F_m$ be also the corresponding primitive idempotents of $\cfcc$ under the canonical adjunction isomorphism $\cfcc \simeq \enx_{\czcc}(R(\unu))$ from Equation \eqref{canisom}.

We define $\csu_i:={\mtf}^{-1}(F_i)\in \cecc$ to be the {\it conjugacy class sums} corresponding to the  conjugacy class $\mathfrak C_i$ where  $\lam\in \cfcc$ is  a cointegral such that $\langle \lam, u\rangle =1$.

For the rest of this section we suppose that $\cc$ is a  non-degenerate spherical fusion category with $\mtr{Gr}_\kk(\C)$  a commutative ring. 
 \bl Let $\cc$ be a spherical non-degenerate fusion category. With the above notations one has that:
 \beq\label{otb}
\langle F_i,\; \csu_j\rangle =\delta_{i,j}\dimcc\tau(F_i)
\eeq
 \el
\bpf
By \cite[Equation 6.10]{scalg} one has that:
\beq\label{610}
\mtf^{-1}(\ch_j)=\frac{\dimcc}{\dimvj}E_{j^*}
\eeq
for any irreducible character $\chj\in \cfcc$, where as above $\lam$ is an integral with $\langle \lam, u\rangle =1$. It follows from here that for any two irreducible characters $\chii, \chj\in \cfcc$ one has that 
$
\langle \ch_i,\; \mtf^{-1}(\ch_j)\rangle =\langle \ch_i,\;\frac{\dimcc}{\dimvj}E_{j^*}\rangle =\delta_{i,j^*}\frac{\dim(V_{j^*})}{\dimvj}\dimcc$. Since $\cc$ is spherical it follows that $\dim(V_{i^*})=\dimvi$ and therefore 
$$
\langle \ch_i,\; \mtf^{-1}(\ch_j)\rangle =\dimcc\tau(\ch_i\ch_j).
$$
Since $\{\chii\}_i$ form a basis for $\cfcc$ and $\tau$ is a bilinear form one can deduce that
\beq\label{comp1}
\langle \ch, \; \mtf^{-1}(\mu)\rangle =\dimcc \tau(\ch\mu).
\eeq
 for any $\ch, \mu\in \cfcc$.
\noindent In particular, for $\ch=F_i$ and $\mu=F_j$ one has that
 \beqn
\langle F_i,\; \csu_j\rangle =\langle F_i,\;\mtf^{-1}(F_j)\rangle =\dimcc \tau(F_iF_j)=\delta_{i,j}\dimcc\tau(F_i).
\eeqn
\epf
\vsk
We note that the size $|\mathfrak C^j|$ of $\mathfrak C^j$ is  a non-zero scalar of $\kk$ since it is the pivotal dimension of a simple object in a pivotal fusion category $\czcc$. By the proof of \cite[Lemma 6.10]{scalg} one has that
$$
\langle \epsu, \csu_i\rangle =\dimcc\langle F_i, \blam\rangle =\dimcc\tau(F_i)=|\mathfrak C^i|\neq 0.
$$
Define 
\beq\label{ni}
n_i:=\frac{1}{ \langle F_i, \blam\rangle }=\frac{1}{\tau(F_i)}=\frac{\dimcc}{|\mathfrak C^i|}.
\eeq
\vsk
Thus in the case of a non-degenerate spherical fusion category $\cc$, by Equation \eqref{otb}, the pair $\{F_i, \frac{n_i}{\dimcc}\csu_i\}$ constitutes  another pair of dual bases for the canonical pairing $\langle ,\;\rangle :\cfcc \times \cecc \ra \kk$.
\vsk
Since $\{\ch_i, \frac{1}{\dimvi }E_i\}$ is another pair of dual bases for the same pairing 
it follows that in $\cfcc\ot \cecc$ we have the following:
\beqn
\sumitom F_i\ot \frac{n_i}{\dimcc}\csu_i=\sumitom \ch_i\ot  \frac{1}{\dimvi }E_i.
\eeqn
\vsk Let $\lam \in \cfcc$ be an idempotent integral. Applying $\id \ot \mtf$ to the above identity one obtains
\beqn
\sumitom F_i\ot \frac{n_i}{\dimcc}F_i=\sumitom \ch_i\ot \frac{1}{\dimvi }\frac{\dim(V_{i^\star})}{\dimcc}\ch_{i^*}
\eeqn
which can be written as
\beq\label{db}
\sumitom F_i\ot n_iF_i=\sumitom \ch_i\ot  \ch_{i^*}.
\eeq
since $\dim(V_i)=\dim(V_{i^*})$, see \cite{EGNO15}.
\vsk  
Write $\ch_i=\sumjtom \al_{ij}F_j$ with $\al_{ij}\in \kk$. Applying $\mtf^{-1}$ to this equation one has that
\beq\label{einc}
\frac{\dimcc}{\dimvi }E_{i^*}=\sumjtom \al_{ij}\csu_j
\eeq
 \vsk Nothe that by \cite[Corollary 6.11]{scalg} one has 
\beq\label{613}
\alij=\langle \ch_i, g_j\rangle
\eeq
where $g_j:=\frac{\csu_j}{|\mathfrak C^j|}$. 
\bl \label{511}Let $\cc$ be a  non-degenerate spherical fusion category with $\mtr{Gr}_\kk(\C)$  a commutative ring. With the above notations one has that:
\beq\label{cjc}
\csu_i=\frac{\dimcc}{n_i}(\sumjtom \frac{1}{\dimvj }\alji E_j)
\eeq
\el
\bpf

The second orthogonality relation from \cite[Theorem 6.11]{scalg} can be written as 
\beq\label{ortrel}
\sumitom \alji\al_{j^*l}=\delta_{i, l}\frac{\dimcc}{\dim(\mathfrak C^i)}
\eeq
Then applying Equation \eqref{einc} one has
{\small
\beqn
\frac{\dimcc}{n_i}(\sumjtom \frac{1}{\dimvj }\alji E_j)=\frac{\dimcc}{n_i}(\sumjtom \frac{1}{\dimvj }\alji ((\frac{\dimvj}{\dimcc})\sum_{l=0}^m\al_{j^*l}\csu^l))=
\eeqn
\beqn
=\frac{1}{n_i}\sum_{l=0}^m(\sumjtom \alji\al_{j^*, l})\csu^l=\frac{\dimcc}{\dim(\mathfrak C^i)n_i}\csu_i=\csu_i.
\eeqn
}
\epf 
\subsection{On the central element $\ell_\cd$-general case.}\label{ellcd}
Let $\cd$ be a fusion subcategory of a non-degenerate spherical fusion category $\cc$. We exhibit a central  element $\elcd:=\mtfi(\lam_\cd)\in \cecc$ that completely determines the fusion subcategory $\cd$.

\bp 
Suppose that $\cc$ is a non-degenerate  spherical fusion category and $\lam$ is an integral with $\lag\lam, \unu\rag=1$. With the above notations:
\beq\label{idmptsums}
\elcd=\frac{\dimcc}{\dimcd}(\sum_{[V_i]\in \irr(\cd)}E_i).
\eeq
\ep
\bpf By Equation \eqref{intregform} one has that 
\beqn\lam_{\cd}=\frac{1}{\dimcd}(\sum_{[V_j]\in \irr(\cd)}\dim(V_j^*)\ch_{j}).
\eeqn Thus 
$$\elcd=\mtfi(\lam_\cd)=\frac{1}{\dimcd}(\sum_{[V_j]\in \irr(\cd)}\dim(V_j^*)\mtfi(\ch_{j})).$$
On the other hand by Equation \eqref{610} one has that:
$
\mtf^{-1}(\ch_j)=\frac{\dimcc}{\dimvj}E_{j^*}
$
for any irreducible character $\chj$. Therefore one can write that:
$$
\elcd=\frac{\dimcc}{\dimcd}(\sum_{[V_j]\in \irr(\cd)}\frac{\dim(V_j^*)}{\dim(V_j)}E_{j^*}). 
$$
Since $\dim(V_{j^*})=\dim(V_j)$ holds for spherical categories one obtains the desired formula. 
\epf 
\noindent
 Suppose as above that $\cc$ is a non-degenerate spherical fusion category and $\cd\subseteq \cc$ a fusion subcategory. By Subsection \ref{cfinc} there is an inclusion of $\kk$-algebras $q^*_\unu:\cfcd\hookrightarrow \cfcc$.  Thus there is a subset $\mtj_\cd \subseteq \{0, \dots, m\}$ such that
\beq\label{idemptint}
\lam_{\cd}=\sum_{j \in \mtj_\cd}F_j
\eeq
 since $\lam_\cd$ is an idempotent element inside $\cfcc$.   Then by applying the inverse $\mtf^{-1}$ of the Fourier transform to the above equation it follows that
\beq\label{idclss}
\ell_\cd=\sum_{j \in \mtj_\cd}\csu_j
\eeq
Since $\{F_i, \;\frac{1}{|\mathfrak C^i|}\csu^i\}$ are dual bases for the evaluation form $\cfcc\times \cecc \ra \kk$ it follows that
$$j \in \mtj_\cd \iff  \langle F_j,\;\ell_{\cd}\rangle\neq 0.$$
\bp\label{epsp}Suppose that $\cc$ is a non-degenerate fusion category and $\cd\subseteq \cc$. Then
\beq\label{eps}
\eps_\unu(\ell_\cd)=\frac{\dimcc}{\dimcd}=\sum_{j \in \mtjcd}|\mathfrak C^j|
\eeq
\ep
\bpf
Note that $\lag \eps_\unu, \;E_j \rag=\delta_{j,\;0}$ from the orthogonality relations of the characters, see Equation \eqref{ortrel}. Then the first equality follows from Equation \eqref{idmptsums}. The second equality follows from Equation \eqref{idclss}.
\epf
\subsection{On the lattice of fusion subcategories}
\bl\label{int} In any spherical non-degenerate fusion category with a commutative Grothendieck ring one has
$$\ell_{\cd \cap \ce}=\frac{\dimcd\dim(\ce)}{\dim(\cd\cap \ce)\dimcc}\ell_\cd \ell_{\ce},\;\lam_{\cd \cap \ce}=\mtf(\ell_\cd\ell_\ce),\; $$ 
\el
\bpf  Let $\cc$ be a fusion category with commutative character ring. As above one has that
\beqn
\elcd=\frac{\dimcc}{\dimcd}(\sum_{[V_i]\in \irr(\cd)}E_i).
\eeqn
Therefore 
\beqn
\ell_\cd \ell_{\ce}=\frac{\dimcc^2}{\dimcd\dim(\ce)}(\sum_{[V_i]\in \irr(\cd \cap \ce)}E_i)
\eeqn
On the other hand one has
\begin{eqnarray*}
\ell_{\cd\cap \ce}& = & \frac{\dimcc}{\dim(\cd\cap \ce)}(\sum_{[V_i]\in  \irr(\cd \cap \ce)}E_i)=\frac{\dimcc}{\dim(\cd\cap \ce)}\frac{\dimcd\dim(\ce)}{\dimcc^2}\ell_\cd \ell_{\ce}\\ &= &\frac{\dimcd\dim(\ce)}{\dim(\cd\cap \ce)\dimcc}\ell_\cd \ell_{\ce}.
\end{eqnarray*}
\epf
\bl Let $\cc$ be a fusion category and $\cd, \ce\subseteq \cc$ be two fusion subcategories. Then
\beq\label{inccomp}
\mtj_\cd\subseteq \mtj_\ce \iff \cd\supseteq \ce.
\eeq
\el
\bpf
If $\mtj_\cd\subseteq \mtj_\ce$ then by Equation \eqref{idemptint} one has $\lam_\cd\lam_\ce=\lam_\cd$. Thus $\lam_{\cd \vee \ce}=\lam_\cd.$ Then Equation \eqref{intregform} implies the statement.
\vsk The converse is obvious. If $\cd\supseteq \ce$ then $\lam_\cd\lam_\ce=\lam_\cd$ and therefore $\mtj_\cd\subseteq \mtj_\ce$.
\epf
\subsection{Fusion subcategories of prime index}
In this subsection we let $\kk=\mathbb C$ and $\cc$ be a pseudo-unitary fusion category with a commutative character ring $\cfcc$.  We let as above  $F_i$ be the primitive central idempotents of $\cfcc$. Without loss of generality we may assume that $F_0=\lam_\cc$, the idempotent cointegral of $\cc$. We denote by $\mu_i:\cfcc \ra \mathbb C$ the characters of the semisimple commutative $\CC$-algebra $\cfcc$.

By \cite[Proposition 9.5.1]{EGNO15} any pseudo-unitary fusion category admits a unique spherical structure $a_X : X\ra X^{**}$ with respect to which $d_X = \fp(X)$ for every simple object $X$. Thus in this case the regular character 
$ r_\cc=\frac{1}{\dimcc}(\sum_{X\in \irr(\cc)}\fp(X)\ch(X))$  coincides with the idempotent cointegral $\lam_\cc$.

Since $\mtr{Gr}_\CC(\cc)\simeq \cfcc$ as $\CC$-algebras, it follows that the results from \cite{bmonat} can be applied directly  inside the ring $\cfcc$ instead of $\mtr{Gr}_\CC(\cc)$. Let $L_{[X]}:\cfcc \ra \cfcc$ be the operator given by the left multiplication by $\ch(X)$ on $\cfcc$. Recall that $L_{[X]}$ has all non-negative entries with respect to the the basis of $\cfcc$ given by the characters $\ch(S)$ of all simple objects $S\in \irr(\cc)$. Moreover  the largest  eigenvalue (in absolute value) of $L_{[X]}$  is positive and it is called $\fp(X)$. Any eigenvector corresponding to this eigenvalue is called a {\it principal eigenvector}.

In \cite{bmonat} for any object $X\in \co(\cc)$ the kernel of $X$ is defined as:
$$\ker_{\cc}(X)=\{\mu_i\;|\; \ch(X)F_i=\fp(X)F_i\},$$
the set of characters $\mu_i$ corresponding to the principal eigenvectors $F_i$ of the operator $L_{[X]}$. Note that $\mu_0=\fp()\in \ker_\cc(X)$ since $F_0=\lam_\cc$.
\bl\label{vee} Let $\cc$ be a pseudo-unitary fusion category $\cc$ commutative Grothendieck ring. For any two fusion subcategories $\cd, \ce \subseteq \cc $ with the above notations one has:
$$\lam_{\cd \vee \ce}=\lam_{\cd}\lam_\ce, \; \ell_{\cd \vee \ce}=\mtf^{-1}(\lam_{\cd \vee \ce})=\sum_{j \in \mtj_\cd \cap \mtj_\ce} \csu_j,\; \mtj_{\cd \vee \ce}=\mtj_{\cd}\cap \mtj_{\ce}. $$ 
\el
\bpf
By \cite[Subsection 6.1]{scalg}, in a fusion category $\cc$ the cointegral element $\lam_\cc$ is the unique element $t \in \cfcc$ with the property that $$\ch t=\lag \ch, \; \unu\rag t,\;\; \text{for any}\;\; \ch \in \cfcc.$$
Since $\cfcc$ is commutative one has $\lam_\cd\lam_\ce=\lam_\ce\lam_\cd$ and $\lam_\cd\lam_\ce\in \cf(\cd\vee \ce)$. Moreover, since $\cc$ is pseudo-unitary note that for any object $X\in \cc$ one has that $[X]\in \co(\cd)$ if and only if $\ker_\cc(X)\subseteq \mtj_\cd$.  Therefore any $F_i$ with $i \in \mtj_\cd\cap \mtj_\ce$ is a principal eigenvector for any $L_{[X]}$ with $X\in \cd$ or $X\in \ce$.  By \cite[Proposition 3.3]{bmonat} it follows that any $F_i$ with $i \in \mtj_\cd\cap \mtj_\ce$ is a principal eigenvector for any $L_{[X]}$ with $X\in \cd\vee \ce$.  Thus $\ch(X)(\lam_\cd\lam_\ce)=\lag\ch(X),\;
\unu\rag\lam_\cd\lam_\ce$ for any $X\in \irr(\cd \vee \ce)$. Therefore $\lam_{\cd \vee \ce}=\lam_\cd\lam_\ce$. 
\vsk On the other hand note that:
$$\lam_{\cd \vee \ce}=\lam_{\cd}\lam_\ce=(\sum_{j\in \mtj_\cd}F_j)(\sum_{j\in \mtj_\ce}F_j)=\sum_{j \in \mtj_\cd\cap \mtj_\ce}F_j$$
\noindent This implies that $\mtj_{\cd \vee \ce}=\mtj_{\cd}\cap \mtj_{\ce}.$
\epf 
\bl 
Let $\cc$ be a pseudo-unitary  fusion category with a commutative character ring $\cfcc$  over $\kk=\mathbb C$. Then
\beq
\mtj_{\ccad}=\{i\;|\; |\mathfrak C_i|=1\}.
\eeq
\el
\bpf
Let $\chad=\sumitom \ch_i\ch_{i^*}$ be the character of the adjoint subcategory.  From Equation \eqref{db} it follows that $\chad=\sumitom {n_i}F_i,$
where by Equation \eqref{ni} one has  $n_i=\frac{\dimcc}{|\mathfrak C^i|}$.
Thus  $$\ker_\cc(\chad)=\{\mu_i\;|\; n_i=\dimcc\}.$$ By \cite[Proposition 3.12]{bmonat} it follows that the regular character $r_\ccad=\frac{1}{\dimcc}(\sum_{\ch_i\in \irr(\ccad)}\fp([X])[X])$ of $\ccad$ with respect to $\fp$ can be written in $\cfcc$  as follows
\beq
r_{\ccad}=\sum_{\{\mu_i\;|\; |\mathfrak C_i|=1\}}F_i.
\eeq
On the other hand since $\cc$ is pseudo-unitary it follows that $r_{\ccad}=\lam_\ccad$ and therefore $\mtj_{\ccad}=\{i\;|\; |\mathfrak C_i|=1\}$.
 \epf

\bt \label{lcad} Let $\cc$ be an integral fusion category with a commutative character ring $\cfcc$ over $\mathbb C$. 
 Then fusion subcategories $\cd$ of the smallest prime index dividing $\fp(\cc)$ are in bijection with normal subgroups of the same prime index $p$ of the universal grading group $U_\cc$ of $\cc$. Moreover,  $\cd\supseteq \ccad$ for any such fusion subcategory.
\et
\bpf Clearly the center $\czcc$ is also an integral fusion category.
By Equation \eqref{eps} one has that
$$\epsu(\ell_\cd)=p=1+\sum_{j \in \mtjcd \setminus \{0\}}|\mathfrak C^j|$$ and $|\mathfrak C^j||\dimcc$. Thus $|\mathfrak C^j|=1$ and therefore $\mtjcd \subseteq \{j\;|\; \ch_j \in \irr(\ccpt)\}=\mtj_{\ccad}$. This implies that $\cd\supseteq \ccad$. 

Moreover, by \cite[Lemma 2.5]{DGNO} fusion subcategories containing $\ccad$ are in bijection with subgroups of the universal grading group $U_\cc$. Thus $\cd=\cc(H)$ where $H\leq U_\cc$ is a subgroup of $U_\cc$ and $\cc(H)=\oplus_{h \in H}\cc_h$. Since $\fp(\cc(H))=|H|\fp(\ccad)$ it follows that $H$ is a subgroup of index $p$ of $G$. Since $p$ is also the least divisor of $G$ it follows that $H$ is a normal subgroup of $G$.
\epf

\section{Proof of the main centralizer result}\label{br}
The goal of this section is to prove Theorem \ref{main1}. In order to do that we need to recall the following preliminaries.

Recall that two objects $X, Y$ of a braided fusion category $\cc$ centralize each other if and only if the monodromy $c_{X,Y}c_{Y, X}:X\ot Y\ra X\ot Y$ is the identity map. Given two simple objects $V_i$ and $V_j$ of a  ribbon fusion category $\cc$ one can define  $$\sij:=\tr_q(c_{V_j,\; V_{i^*}}c_{V_{i^*},\; V_j}),$$
 where $\tr_q$ denotes the canonical quantum trace in  $\cc$.

Throughout  this section we consider a {\it pre-modular fusion category} $\cc$, i.e. a braided fusion category  with a spherical structure. Recall that by \cite[Prop 8.23]{ENO}  any braided weakly integral fusion category is pre-modular. As in the previous section we use the same notation   $V_0, \dots, V_m$ for the set of simple objects of $\cc$ up to isomorphism. We let also $E_0, \dots, E_m$ be their associated primitive idempotents inside $\cecc$.   Without loss of generality we assume that  $V_0=\unu$.   Then by \cite[Section 6]{scalg} $E_0$ is the idempotent integral of $\cc$ and $\eps_\unu=\ch_0$. It is also well known that with this notations one has $d_i =s_{0i}=s_{i0}$, the quantum dimension of the simple object $V_i$. 

\br\label{ext} Let $\cc$ be a pre-modular fusion category and $\cd, \ce$ be two fusion subcategories. It is well known that in a pre-modular fusion category two simple objects $V_i$ and $V_j$ centralize each other if and only if $\sij=d_id_j$. 

Extending linearly the $S$-matrix by $S(\sum_i\al_i\ch_i, \sum_j\beta_j\ch_j)=\sum_{ij}\al_i\beta_j\sij$ we see that if $
\cd'\subseteq \ce$ then $S(\ch, \mu)=d(\ch)d(\mu)$ for any $\ch \in \cfcd$ and $\mu \in \cf(\ce)$. Here $d:\cfcc \ra \kk $ is the quantum dimension homomorphism defined on the basis by $d(\ch_i)=d_i$.
\er
\subsection{Definition of the Drinfeld map} In this subsection we recall the construction of the Drinfeld map in a braided fusion category.
Consider as in Section \ref{adj} also the central Hopf monad of $\cc$ given by:
$$\mtc L(V)=\int^{X\in \cc} X^*\ot V \ot X$$
and  $\iota_{V; X}: X^*\ot V\ot X\ra \mtl(V)$ the associated universal dinatural maps. As already mentioned  is well known that $\mtc L$ is a Hopf monad and $\mtc L$ is a left adjoint of the forgetful functor $F:\czcc\ra \cc$.
\vsk
 Following \cite[Subsection 4.4]{FRG} one can define the map $T_{X, Y}:X\ot X^{*} \ra Y^{*}\ot Y$ as the classical Hopf link map.  Since this is dinatural with respect to $X$ and $Y$, by the universal properties of ends
and coends  there exists a unique  $Q:\mtc L(\unu) \ra Z(\unu)$ such that $T_{X, Y}$ factors as
 \beqn
T_{X, Y}=\pi_{\unu;\; X} \circ Q \circ \iota_{\unu;\;Y}
 \eeqn
 \vsk The map $Q$ is called the Drinfeld map associated to $\cc$.\vsk 
\vsk
Since $\mtc L(\unu)\simeq Z(\unu)^*$ this map induces a map on the ring of class functions $f_Q: \cfcc \ra \cecc$. In the case of a ribbon category $\cc$ the map $f_Q$ has the following formula (see \cite[Example 6.14]{scalg} ) 
\beq\label{exmp614}
f_Q(\ch_i)=\sumjtom \frac{\sij}{\dxj}E_j
\eeq
where as above $\sij:=\tr_q(c_{V_j,\; V_{i^*}}c_{V_{i^*},\; V_j}),$ Here $E_i \in \cecc$ represents the primitive idempotent of $\cecc$ associated to the simple object $V_i$.  
\vsk
Recall that for a pivotal fusion category we denoted by $\lam_\cc\in \cfcc$ the idempotent integral of $\cc$.

\bl \label{evu} Let $\cc$ be a ribbon fusion category and $f_Q:\cfcc\ra\cecc$ defined as above. Then
$$\eps_\unu(f_Q(\ch))=\langle\ch, \;u\rangle,$$
for any character $\ch \in \cfcc$.
\el
\bpf
As mentioned above $\eps_\unu=\ch_0$ and $\lag \ch_i, E_j\rag =\delta_{i, j}d_i$. Then the statement follows from Equation \eqref{exmp614}.
\epf
\subsection{Proof of the main result}
As explained in \cite[Section 6]{scalg} for a braided fusion category the spaces $\cfcc$ and $\cecc$ are $\kk$-semisimple commutative algebras. Thus, since $f_Q$ is an algebra map, it follows that $f_Q(q_\unu^*(\lam_\cd))$ can be written as sum of the central primitive idempotents $E_j\in \cecc$. By abuse of notations we will write $f_Q(\lam_\cd)$ instead of $f_Q(q_\unu^*(\lam_\cd))$ in the next theorem.

Now we are ready to prove the main result Theorem \ref{main1} concerning the centralizer in pre-modular fusion categories.
\bt
Let $\cd$ be a fusion subcategory of a ribbon fusion category $\cc$ and $\lam_{\cd}$ be the associated idempotent cointegral of $\cd$. Then
\dbd
\mtf(f_Q(\lam_{\cd}))=\frac{\dim(\cd')}{\dimcc}\lam_{\cd'}.
\dbd
where $\mtf $ is the Fourier transform introduced by Shimizu in \cite{scalg} and $\lam_{\cd'}$ is the idempotent cointegral of $\cd'$.
\et

\bpf
Let $\cd$ be a fusion subcategory of a ribbon tensor category $\cc$. Suppose that  for the idempotent integral $\lam_\cd$ of $\cd$ has the following writing in $\cecc$:
\beq\label{intcd}
f_Q(\lam_\cd)=\sum_{j \in \mtca_0}E_j.
\eeq 
\noindent
First we will prove that 
\beq\label{idsum}
\irr(\cd')=\{V_j\;|\;j \in\mtca_0\}.
\eeq 
For any irreducible character $\chii\in \cfcd$ one has by Equation \eqref{mtc} 
\beqn
f_Q(\ch_i)=\sumjtom \frac{\sij}{d_j }E_j.\;\;
\eeqn
Moreover, by the definition of the categorical cointegral one has $\chi_i\lam_\cd=d_i\lam_\cd$ for any $\ch_i\in \irr(\cd)$. This implies that 
$$d_if_Q(\lam_\cd)=f_Q(\ch_i\lam_\cd)=f_Q(\ch_i)f_Q(\lam_\cd).$$ Using  Equation \eqref{intcd} this gives that $\sij=\dxi\dxj$ for any $j\in \mtca_0$ and $\ch_i \in \cd$. Thus $\{\ch_j\;|\;j \in\mtca_0\}\subseteq \cd'$.
\vsk
Conversely if $\ch_j\in \cf(\cd')$ then  $S(\ch_j,\; \lam_\cd)= d_j$ by Remark \ref{ext}. Since $r_\cd=\dimcd\lam_\cd$ it follows that
\beq
\langle \ch_j, f_Q(r_\cd)\rangle =\dimvj \dimcd
\eeq 
From formula \eqref{intcd} this happens if and only if $j \in \mtca_0$.
\vsk 
As above one has that $\mtf(E_i)=\frac{d_{i}}{\dimcc}\ch_{i^*}$. Thus $$\mtf(f_Q(\lam_\cd))=\sum_{j \in \ca_0}\mtf(E_j)=\sum_{j \in \ca_0}\frac{d_{j}}{\dimcc}\ch_{j^*}=\frac{\dim(\cd')}{\dimcc}\lam_{\cd'}
$$
\epf
\bc Let $\cc$ be a ribbon category.
 Consider the particular case $\cd=\cc'$. Then  $f_Q(\lam_{\cd})=1_\cecc$, the sum of all primitive central the idempotents.
\ec 
 \bpf 
 The Corollary follows since $\cd'=(\cc')'=\cc$.
 \epf
 In the next example we show how to recover the centralizer formula \cite{new-sub}  from our main result Theorem \ref{main1}.

\bn{example}Suppose that $(H, R)$ is a semisimple quasi-triangular Hopf algebra and let $L$  be a left normal coideal subalgebra of $H$. Let $\cd:=\rep(H//L)$. Then $\lam_{\cd}=\lam_L$
is the integral of the Hopf algebra $(H//L)^*$.  Recall that $H$ is a ribbon Hopf algebra with ribbon element given by $v=u^{-1}$, the inverse of the Drinfeld element.
\vsk It was shown in \cite[Lemma 4.1]{new-sub} that $f_Q(\lam_L)=\blam_{M}$ where $M:=\phi_R((H//L)^*)$ is a normal left coideal subalgebra of $H$. Then in this case 
$\mtc F_\lam(\blam_{M})=\lam_{M}$ by \cite[Lemma 1.1]{CW10}. Note that the result $$\rep(H//L)'=\rep(H//M)$$ from \cite[Theorem 1.1]{new-sub} follows now from the fact $M=
\blam_{M}\lh H^*$, see \cite[Lemma 1.1]{CW10}.
\end{example}

\section{On the conjugacy classes of fusion categories}\label{cc}
A pre-modular category $\cc$ is called {\it modular }  if its $S$-matrix is non-degenerate. By \cite[Proposition 3.7]{dgno2} the $S$-matrix is non-degenerate if and only if $\cc'=\vect$.
\vsk
In the case of a modular tensor category $\cc$, the map $f_Q$ is a bijective map and after a relabelling of indices one may suppose that $f_Q(F_j)=E_j$ for any primitive idempotent  $E_j\in \cecc$ associated to the character $\ch_j\in \cfcc$. In this case as in \cite[Example 6.14]{scalg} one can also write:
\beq\label{mtc}
\ch_i=\sumjtom \frac{\sij}{d_j}F_j.\;\;
\eeq
\bl Let $\cc$ be a modular category. With the above notations it follows that $F_0=\lam_\cc$ and $f_Q(F_0)=\blam_\cc$ is the categorical integral of $\cc$.
\el
\bpf Note that from \cite[Lemma 6.9]{scalg} for $r=0$ it follows that $F_0=\lam_\cc$ is the idempotent integral of $\cc$. One has that $f_Q(f)f_Q(F_0)=f_Q(fF_0)=\langle f, u\rangle f_Q(F_0)$. Since $\langle f, u\rangle =\langle \epsu, f_Q(f)\rangle $ it follows that $f_Q(f)f_Q(F_0)=\langle \epsu, f_Q(f)\rangle$ for all $f \in \cfcc$. On the other hand since $\cc$ is modular it follows that $f_Q$ is surjective and therefore $f_Q(F_0)$ is an integral of $\cc$.
\epf 
\subsection{Product of conjugacy class sums in the modular case}\label{pcc}
Assume that $\cc$ is a modular tensor category and let as above $F_j$ be the primitive idempotents of $\cfcc$.
As in Section \ref{cfc} one can write $\ch_i=\sumjtom \al_{ij}F_j$ with $\al_{ij}\in \kk$. 
\bl Assume that $\cc$ is a modular tensor category. With the above notations one has that:
\beq
\alij=\frac{d_i }{d_j }\alji
\eeq
\el
\vsk
\bpf
In the modular case one has that $h_i: \cfcc \ra \kk,\;\; [X_i]\sent \frac{\sij}{d_j}$ are all the characters of the character ring $\cfcc$. Then as in \cite{scalg} it follows directly that $\alij=\frac{\sij}{d_j }$ . Thus the equality $\alij=\frac{d_i }{d_j }\alji$  follows from the symmetry of the $S$-matrix, i.e. $\sij=\sji$.
\epf
\bt\label{imfq} Let $\cc$ be a modular fusion category. Then
\beq
f_Q(\ch_i)=\frac{d_i n_i}{\dim(\cc)}\csu_i=\frac{d_i}{|\mathfrak C^i|}\csu_i.
\eeq
\et
\bpf Using Lemma \eqref{511} one has that
\beq\label{fqch}
f_Q(\ch_i)=\sumjtom \al_{ij} E_j=\sumjtom \frac{d_i }{d_j }\alji E_i=d_i \sumjtom \frac{1}{d_j }\alji E_i=\frac{d_i n_i}{\dimcc}\csu_i.
\eeq
\vsk The second equality of the theorem follows from Equation \eqref{ni}.
\epf
\subsection{$\ell_\cd$ in the pre-modular case}
Note that in the case of a  ribbon category $\cc$, the main Theorem \ref{main1} can be written as
\beq\label{from1}
\ell_{\cd'}=\frac{\dim(\cd')}{\dimcc}f_Q(\lam_\cd)
\eeq
\bt\label{new1} Let $\cc$ be a modular tensor category and $\cd$ be a fusion subcategory of $\cc$. With the above notations one has
\bne
\item $d_j^2=|\mathfrak C^j|$ for all $j$.
\item For the centralizer $\cd'$ one has that: \beqn
\ell_{\cd'}=\sum_{\{j\;|\; \ch_j \in \irr(\cd)\}}\csu_j,
\eeqn i.e. $\mtj_{\cd'}=\{j\;|\; \ch_j\in \irr(\cd)\}$.
\ene
\et
\bpf Note that $\dimcd \dim(\cd')=\dimcc$ since $\cc$ is a modular tensor category. By Equation \eqref{from1} one has 
{\small
\beqn
\ell_{\cd'}=\frac{\dim(\cd')}{\dimcc}f_Q(\lam_\cd)=\frac{\dim(\cd')}{\dimcc}(\frac{1}{\dim(\cd)}\sum_{j \in \irr{\cd}}d_{j^*}f_Q(\ch_j))=\sum_{j \in \irr{\cd}}d_{j^*}\frac{d_j}{|\mathfrak C^j|}\csu_j
\eeqn
}
On the other hand by Equation \eqref{idemptint} one has that
$$\ell_{\cd'}=\sum_{j \in \mtj_{\cd'}}\csu_j$$
This implies both statements of the theorem.
\epf
\bc Let $\cc$ be a modular tensor category. Then 
\bne
\item $|\mathfrak C^j|=1 \iff \ch_j\in \irr(\ccpt) \iff F(\mathfrak C^j)=\unu$.
\item 
$\ell_{\ccad}=\sum_{\{j\;|\; \ch_j\in \irr(\ccpt)\}}\csu_j,\; \mtj_{\ccad}=\{j\;|\; \ch_j\in \irr(\ccpt)\}$, and \;\; $\lam_{\ccad}=\sum_{\{j\;|\; \ch_j\in \irr(\ccpt)\}}F_j$.
\item $\ell_{\ccpt}=\sum_{\{j\;|\; \ch_j\in \irr(\ccad)\}}\csu_j$, and \;\; $\lam_{\ccpt}=\sum_{\{j\;|\; \ch_j\in \irr(\ccad)\}}F_j$.
\ene
\ec
\bpf The first item follows from the fact that $|\mathfrak C^j|=d_j^2$ and under the forgetful functor $F:\czcc \ra \cc$ the object $F(\mathfrak C^j)$ always contains $\unu_\cc$. The last two items follow from the fact that $\ccad'=\ccpt$ and $\ccpt'=\ccad$ in a modular tensor category.
\epf

\subsection{On the structure constants and product of conjugacy classes} 
In \cite{scalg} the author asked if results from \cite{CW4, CW6} concerning conjugacy classes can be extended from semisimple Hopf algebras to fusion categories.
\vsk Since the class sums $\csu_i$ form a basis for the center $\cfcc$ one has that
\beq
\csu_i\csu_j=\sum_{l=0}^rc^l_{ij}\csu_l
\eeq
for some scalars $c^l_{ij}\in \kk$. These scalars were called {\it  the structure constants} of conjugacy classes in \cite{CW2}. For factorizable Hopf algebras,  by \cite[Theorem 4.3]{CW2}, one has that
\beq\label{cwconj}
c^{l}_{ij}=\frac{d_{l}d_{i}}{d_{j}}N^{l}_{ij}
\eeq
where $\ch_i\ch_j=\sum_{l=0}^rN^{l}_{ij}\ch_l$. 
\vsk 
The following result follows from Equation \eqref{1fus} and extends the above mentioned   result to arbitrary modular tensor categories.  
\bp
Let $\cc$ be a modular tensor category. With the above notations we have that the conjugacy class sums form a $\mathbb Q$-subalgebra of the algebra $\cecc$ of central elements. More precisely, 
\beqn
\csu_i\csu_j=\sum_{l=0}^m\frac{d_id_j}{d_l}N^l_{ij}\csu^l
\eeqn
where $N^l_{ij}$ are the fusion coefficients of $\cc$, i.e. $\ch_i\ch_j=\sum_{l=0}^mN^l_{ij}\ch_l$.
\ep
\bpf Since $|\mathfrak C_i|=d_i^2$, Equation \eqref{fqch} implies that $f_Q(\ch_i)=\frac{1}{d_i}\csu_i$. For any two irreducible characters write $\ch_i\ch_j=\sum_{l=0}^mN^l_{ij}\ch_l$ for some non-negative integers $N^l_{ij}$.
Then one can write
\beq\label{1fus}
\frac{1}{d_id_j}\csu_i\csu_j=f_Q(\ch_i)f_Q(\ch_j)=f_Q(\ch_i\ch_j)=\sum_{l=0}^mN^l_{ij}f_Q(\ch_l)=\sum_{l=0}^m\frac{1}{d_l}N^l_{ij}\csu^l
\eeq
which proves the statement.
\epf

\section{Appendix} Let $\cc$ be any finite tensor category and $\cd$ a full subcategory of $\cc$.
In this appendix we prove the main properties of the ends 
$$
Z_\cc(M)=\int_{X\in \cc} X\ot M\ot X^*
$$ 
and 
$$
Z_\cd(M)=\int_{X\in \cd} X\ot M\ot X^*
$$  
that were used in the previous sections of the paper. For simplicity we write $Z:=Z_\cc$ and $\barz:=Z_\cd$. We will also denote by $\pi_{M;\; X}: Z(M)\ra X\ot M\ot X^*$ and respectively $\bar{\pi}_{M;\; X}: \barz(M)\ra X\ot M\ot X^*$ their universal diantural transformations.

\vsk
By  Fubini theorem for ends one has that 
$$
Z^2(M)=\int_{(X, Y)\in \cc\times \cc} (X\ot Y)\ot M \ot (X\ot Y)^*
$$ 
and 
$$
\barz^2(M)=\int_{(X, Y)\in \cd\times \cd} (X\ot Y)\ot M \ot (X\ot Y)^*
$$ 
\vsk with universal dinatural transformations:
$$\pi^2_{M; (X, Y)}=Z^{2}(M)\xra{\pi_{Z(M); X}} X\ot Z(M) \ot X^*\xra{\id\ot \pi_{M,Y}\ot \id}(X\ot Y)\ot M \ot (X\ot Y)^*$$
respectively
$$\bar{\pi}^2_{M; (X, Y)}=\barz^{2}(M)\xra{\bar{\pi}_{\barz(M); X}} X\ot \barz(M) \ot X^*\xra{\id\ot \bar{\pi}_{M,Y}\ot \id}(X\ot Y)\ot M \ot (X\ot Y)^*.$$
\vsk

Using the universal property of the ends there are unique maps 
$$\delta_M: Z(M)\ra Z^2(M), \; \eps_M:Z(M)\ra M$$
making $Z:\cc \ra \cc$ a comonad on $\cc$. Similarly, $\barz:\cc \ra \cc$ is a comonad with structure maps 
$$\bar{\delta}_M: \bar Z(M)\ra \barz^2(M), \; \bar{\eps}_M:\bar Z(M)\ra M$$
 As in subsection \ref{cfinc} there is a unique natural transformation $q:Z\ra \barz$ such that
$$\bar{\pi}_{M;\; X}\circ q_M={\pi}_{M;\; X}$$ for any  object $X$ of $\cd$. 
\vsk
By the universal property of the coend $Z(M\ot N)$, the canonical morphisms $$Z(M)\ot Z(N)\xra{\pimx\ot\pinx} X\ot M\ot X^*\ot X\ot N\ot X^*\xra{\id \ot \evx\ot \idx}X\ot (M\ot N)\ot X$$
induces a canonical map $Z^2_{M, N}:Z(M)\ot Z(N)\ra Z(M\ot N)$. Together with $Z^0:=\pi_{\unu;\;\unu}:Z(\unu)\ra \unu$ one can verify that $(Z, Z^2, Z^0)$ is a monoidal functor.
\vsk
Similarly one can define a monoidal structure $(\barz, \barz^2, \barz^0)$ on the relative coend $\barz$.

\subsection{Morphisms of monoidal functors.}
A monoidal morphism between two monoidal functors on a monoidal category $\cc$ is a morphism $q:T\ra S$ such that the diagrams \vsk
{\Tiny
\begin{tikzcd}[column sep=large,row sep=large]
 T(\unu) \arrow[rr,"q_{\unu}"] & & \bar S(\unu) 
 \\ & \unu \arrow[ul,"T_{0}", sloped] \arrow[ur,"S_{0}", sloped]
\end{tikzcd},
\begin{tikzcd}[column sep=large,row sep=large]
 T(M)\ot T(N) \arrow[rr,"q_{M}\ot q_N"] \arrow[d, "T_{M, N}^2"]  
 &  & S(M) \ot S(N) \arrow[d, "S^2_{M, N}"]
\\ T(M\ot N)  \arrow[rr,"q_{M\ot N}"]   
&  &  S(M\ot N)
\end{tikzcd}
}
 are commutative.
 \bl The natural transformation $q: Z \ra \barz$ from above is a morphism of monoidal functors.
 \el
 \bpf The  compatibility with the unit follows directly form the definition of $q_\unu$ since $\epsu=\pi_{\unu,\;\unu}$ and $\bar{\eps}_\unu=\bpi_{\unu,\;\unu}$.
 \vsk
In order to show that $q$ is a morphism of monoidal functors one needs to verify:
$$q_{M\ot N}\circ Z^2_{M, N}=\barz^2_{M, N}\circ (q_M\ot q_N):Z(M)\ot Z(N)\ra \barz(M\ot N)$$ i.e. 
 that the upper rectangle of the  following diagram is commutative:
\vsk
{\Tiny
\begin{tikzcd}[column sep=large,row sep=large]
&  Z(M)\ot Z(N) \arrow[r,"Z^2_{M, N}"] \arrow[ddl, "{\pi}_{M,\;X}\ot \pi_{N,\;X}", bend right] \arrow[d, "{q}_{M}\ot q_N"]  & Z(M\ot N) \arrow[d,"q_{M\ot N}"]  \arrow[ddr, "{\pi}_{M\ot N,\;X}", bend left]
 \\  & \barzm\ot \barzn \arrow[r,"\bar Z^2_{M, N}"]  \arrow[dl, dashed, "{\bpi}_{M,\;X}\ot \bpi_{N,\;X}", sloped]  &   \barz(M\ot N) \arrow[dr, dashed, "\bar{\pi}_{M\ot N,\;X}"]
 \\ X \ot M\ot  X^*\ot X \ot N \ot X^* \arrow[rrr,"\id_{X\ot M} \ot \ev_X\ot  \ot \id_N \ot \id _{X^*}"]  & & & X\ot  M\ot N \ot X^*
\end{tikzcd}
}
Note that the left and right rectangles of the above diagram are commutative by the definition of $q$. The bottom bent rectangle commutes by the definition of $\barz^2$ while the whole bent rectangle by the  definition of $Z^2$. It follows that the two maps from above are the universal maps corresponding to the same dinatural map $Z(M)\ot Z(N)\ra X\ot (M\ot N)\ot X^*$. By the unicity of the universal map it follows that they are equal.
\epf

\subsubsection{Definition of $q^{(2)}_M$. }

By the universal property of the end $\barz^2(M)=\int_{(X, Y)\in \cd\times \cd} X\ot Y\ot M\ot (X\ot Y)^*$ there  is a unique natural transformation $q^{(2)}: Z^2\ra \bar Z^2$ such that for any two objects $X, Y\in \co(\cd)$  the following diagram
{\Tiny
\begin{tikzcd}[column sep=large,row sep=large]
 Z^2(M) \arrow[rr,"q^{(2)}_{M}"] \arrow[dr,"{\pi}^2_{M;\;(X, Y)}"] & & \bar Z^2(M) \arrow[dl,"\bar{\pi}^2_{M;\;(X, Y)}"]
 \\ & (X\ot Y)\ot  M \ot   (X\ot Y)^* 
\end{tikzcd}
}
\vsk is commutative.

\subsubsection{Definition of $r^{(2)}_M$}

There is a dinatural map in $X \in \co(\cd)$ given by
$$Z^2(M)\xra{\pi_{Z(M);\;X}} X\ot Z(M) \ot X^*\xra{\id_X\ot q_M\ot \id_{X^*}} X\ot \barzm \ot X^*. $$
Since $\barz^2(M)=\int_{X\in \cc} X\ot \barzv \ot X^*$ it follows that there exists a unique map $$r^{(2)}_M: Z^2(M) \ra \barz^2(M)$$
\vsk such that:
{\Small
\begin{tikzcd}[column sep=large,row sep=large]
 Z^2(M) \arrow[rr,"r^{(2)}_M"] \arrow[d, "{\pi}_{Z(M);\;X}"] & & \bar Z^2(M) \arrow[d,"\bar{\pi}_{\barzm;\;X}"]
 \\  X\ot Z M \ot X^* \arrow[rr,"\id_X\ot q_M\ot \id _{X^*}"]  & & X\ot  \barzm \ot X^*.
\end{tikzcd}
}
\subsubsection{$r^{(2)}_M=q^{(2)}_M$}

\bl With the above notations $r^{(2)}_M=q^{(2)}_M$.
\el
\bpf
By the definition of $q^{(2)}_M$ this we have to show that 
\vsk
{\Small
\begin{tikzcd}[column sep=large,row sep=large]
 Z^2(M) \arrow[rr,"r^{(2)}_M"] \arrow[dr,"{\pi}^2_{M, (X, Y)}", sloped] & & \bar Z^2(M) \arrow[dl,"\bar{\pi}^2_{M, (X, Y)}", sloped]
 \\ & (X\ot Y)\ot  M \ot   (X\ot Y)^* .
\end{tikzcd}
}

\noindent Expanding $\pi^2$ and $\bar{\pi}^2$ we have to show that:
\vsk
{\Small
\begin{tikzcd}[column sep=large,row sep=large]
 Z^2(M) 
 \arrow[rr,"r^{(2)}_M"] \arrow[d,"{\pi}^2_{Z(M);\;X}"] 
 & 
 & 
 \bar Z^2(M) 
 \arrow[d,"\bar{\pi}^2_{Z(M);\; X}"]
 \\ 
 X\ot Z(M) \ot X^*  
 \arrow[dr,"\idx \ot {\pi}^2_{M;\; Y}\ot \id_{X^*}", sloped] \arrow[rr, , dashed, "\idx \ot q_M\ot \id_{X^*}"] 
 & 
 & 
 X\ot \barzv \ot  X^*  \arrow[dl,"\idx \ot \bar{\pi}^2_{M, Y}\ot \id_{X^*}", sloped]
 \\ 
 & (X\ot Y)\ot  M \ot   (X\ot Y)^* &
\end{tikzcd}
}
\vsk The bottom triangle is commutative by the definition of $q_M$. The upper rectangle is commutative by the definition of $r^{(2)}_M$.
\epf
\bl \label{q2}
 One has that 
 $$
 q^{(2)}_M=\barz(q_M)\circ q_{Z(M)}=q_{\barz(M)}\circ  Z(q_M).
 $$
 \el
 \bpf 
 For the first equality one has to show that the following diagram is commutative:
 \vsk
{\Small \begin{tikzcd}[column sep=large, row sep=large]
 Z^2(M) \arrow[r,"q_{Z(M)}"] \arrow[d, "{\pi}_{Z(M);\;X}"] & \barz Z(M) \arrow[dl, dashed,  "\bar{\pi}_{Z(M);\;X}"]  \arrow[r,"\barz(q_M)"]& \bar Z^2(M) \arrow[d,"\bar{\pi}_{\barzm;\;X}"]
 \\  X\ot Z M \ot X^* \arrow[rr,"\id_X\ot q_M\ot \id _{X^*}"]  & & X\ot  \barzm \ot X^*
\end{tikzcd}
}
\vsk Note that the left triangle is commutative by the definition of $q_{Z(M)}$. The next square commutes by the definition of $\barz(q_M)$.
\vsk  For the second equality one has to show that the following diagram is commutative:
 \vsk
 {\Small \begin{tikzcd}[column sep=large, row sep=large]
 Z^2(M) \arrow[r,"Z(q_M)"] \arrow[d, "{\pi}_{Z(M);\;X}"] 
 & 
 Z\barz (M) 
 \arrow[dr, dashed,  "\bar{\pi}_{Z(M);\;X}"]  \arrow[r,"q_{\barzm}"]
 & 
 \bar Z^2(M) 
 \arrow[d,"\bar{\pi}_{\barzm;\;X}"]
 \\  
 X\ot Z M \ot X^* 
 \arrow[rr,"\id_X\ot q_M\ot \id _{X^*}"]  
 & 
 & X\ot  \barzm \ot X^*
\end{tikzcd}
}
\vsk Note that the left square commutes by the definition of $Z(q_M)$ while the right triangle commutes by the definition of $q_{\barzm}$.
 \epf
\subsection{On the map $v_M$}

\bl \label{vm}
There is a unique  natural transformation $v_M:Z(M)\ra \barz^2(M)$ such that 
$$\bar{\pi}^2_{M;\; (X,Y)}\circ v_M=  \bar{\pi}_{M;\; X\ot Y}\circ q_M$$
\vsk 
i.e. such that the following diagram is commutative:
\vsk
{\Small
\begin{tikzcd}[column sep=large,row sep=large]
 Z(M) \arrow[r,"v_M"] \arrow[d,"q_M"]&  \bar Z^2(M) \arrow[d,"\bar{\pi}^2_{M;\; (X, Y)}"]
 \\ \barz(M) \arrow[r,"\bar{\pi}_{M; X\ot Y}"]& (X\ot Y)\ot  M \ot   (X\ot Y)^* .
\end{tikzcd}
}
\vsk
One has that $$v_M=\bar{\delta}_M\circ q_M.$$\
\el
\bpf One has as above $$\barz^2(M)=\int_{(X, Y)\in \cd\times \cd }(X\ot Y)\ot M \ot (X\ot Y)^*.$$
\vsk The maps
$$Z(M)\xra{q_M}\barz(M)\xra{\bar{\pi}_{M, (X,Y)}} (X\ot Y)\ot M \ot (X\ot Y)^*$$ are dinatural with respect to $(X,Y)$. Thus by the universal property of the coend there is a unique map $v_M$ making the above diagram commutative. It remains to show that $v_M=\bar{\delta}_M\circ q_M.$
Thus we have to show that the following diagram is commutative \vsk
{\Small
\begin{tikzcd}[column sep=large,row sep=large]
 Z(M) \arrow[r,"q_M"] \arrow[drr, dashed, "{\pi}_{M; X\ot Y}", sloped]   \arrow[d,"q_M"]& \bar Z(M) \arrow[r,"\bar{\delta}_{M}"]  \arrow[dr, dashed, "\bar{\pi}_{M; X\ot Y}", sloped] & \barz^2(M)  \arrow[d,"\bar{\pi}^2_{M,\; (X, Y)}"]
 \\
 \barz(M) \arrow[rr,"\bar{\pi}_{M; X\ot Y}"]   &  & (X\ot Y)\ot  M \ot   (X\ot Y)^* 
\end{tikzcd}
}
\vsk
This follows since all three inside triangles are commutative.
\epf
\subsection{Morphism of comonads}
Let $(T, \delta,\eps) ,\;(S, \delta',\eps')$ be two comonads on a monoidal category $\cc$.
Recall that a morphism $q:T\ra S$ is a morphism of comonads if 
$$\eps'\circ q=\eps,\; \;\delta'\circ q=q_{S()}\circ T(q)\circ \delta,\;\;  \delta'\circ q=S(q)\circ q_{T()}\circ \delta$$
For any $M\in \cocc$ these relations are equivalent to the commutativity of  the following  diagrams:
{\Tiny
\begin{tikzcd}[column sep=large,row sep=large]
 T(M) \arrow[rr,"q_{M}"] \arrow[dr,"\eps_{M}", sloped] & & \bar S(M) \arrow[dl,"\eps'_{M}", sloped]
 \\ & M 
\end{tikzcd}
}
and \vsk 
{\Tiny
\begin{tikzcd}[column sep=large,row sep=large]
 T(M) \arrow[rr,"q_{M}"] \arrow[d, "{\delta}_{M}"]  & &  S(M) \arrow[d, "{\delta'}_{M}"]
\\ T^2(M) \arrow[r,"T(q_M)"]   &  TS(M)  \arrow[r,"q(S(M))"]& S^2(M)
\end{tikzcd}
\begin{tikzcd}[column sep=large,row sep=large]
 T(M) \arrow[rr,"q_{M}"] \arrow[d, "{\delta}_{M}"]  & &  S(M) \arrow[d, "{\delta'}_{M}"]
\\ T^2(M) \arrow[r,"q(T(M))"]   &  ST(M)  \arrow[r,"S(q_M)"]& S^2(M)
\end{tikzcd}
}

\bt\label{mhc} Let $\cd\subseteq \cc$ be a tensor subcategory of a finite tensor category.
With the above notations $q:Z\ra \barz$ is a morphism of Hopf comonads.
\et

\bpf
We have already shown that $q$ is a morphism of monoidal functors. It remains to show that $q$ is a morphism of comonads.
\vsk
In order to show that $q$ is a natural transformation of comonads we need to show that for any object $M \in \co(\cd)$ one has 
$$\bar{\eps}_M\circ q_M=\eps_M$$
and the following diagrams are commutative:
\vsk 
{\Tiny
\begin{tikzcd}[column sep=large,row sep=large]
 Z(M) \arrow[r,"\delta_M"] \arrow[d,"q_M"] & Z^2(M) \arrow[r,"q_{Z(M)}"] &  \barz Z(M)    \arrow[d,"\barz(q_{M}) "]
\\
\barz(M)  \arrow[rr,"\bar{\delta}_M"]  & &\bar Z^2(M) 
\end{tikzcd}
}
 {\Tiny
\begin{tikzcd}[column sep=large,row sep=large]
 Z(M) \arrow[r,"\delta_M"] \arrow[d,"q_M"] & Z^2(M) \arrow[r,"Z(q_M)"] &  Z\barz(M)   \arrow[d,"q_{\barz(M)}"]
\\
\barz(M)  \arrow[rr,"\bar{\delta}_M"]  & &\bar Z^2(M) 
\end{tikzcd}
}
\vsk  Since $\eps_M=\pi_{M;\unu}$ and $\bar{\eps}_M=\bar{\pi}_{M;\unu}$ the first formula follows from the definition of $q_M$.
\vsk We first show that the  second  diagram is commutative. By Lemma \ref{vm} it its enough to show that $$v_M=q_{\barz(M)}\circ Z(q_M)\circ \delta_M$$ i.e. the following diagram is commutative:
\vsk
{\Small
\begin{tikzcd}[column sep=large, row sep=large]
 Z(M) \arrow[r,"\delta_M"] \arrow[d,"q_M"] \arrow[drrr, dashed, "\pi_{M;\;X\ot Y}"] & Z^2(M) \arrow[r,"Z(q_M)"] \arrow[drr, dashed, "\pi^2_{M;\;X\ot Y}"]&  Z\barz(M)   \arrow[r,"q_{\barz(M)}"]
 & \barz^2(M)   \arrow[d,"\bar{\pi}^2_{M;(X,Y)}"]\\
\barz(M)  \arrow[rrr,"\bar{\pi}_{M;\;(X\ot Y)}"]  & & &(X\ot Y)\ot M \ot (X\ot Y)^* 
\end{tikzcd}
}
\vsk
The most left bottom triangle is commutative by the definition of $q_M$. The middle triangle is commutative by the definition of $\delta_M$. Thus it remains to show that the diagram:
\vsk
{\Small
\begin{tikzcd}[column sep=large, row sep=large]
 Z^2(M) 
 \arrow[r,"Z(q_M)"] \arrow[d, "\pi^2_{M;\;X\ot Y}"] 
 & 
 Z\barz(M) \arrow[r,"q_{Z(M)}"] 
 &  
 \barz^2(M)   \arrow[d,"\bar{\pi}^2_{M;\;(X,Y)}"]
  \\
(X\ot Y)\ot M \ot (X\ot Y)^*   \arrow[rr, equal]  &  &(X\ot Y)\ot M \ot (X\ot Y)^* 
\end{tikzcd}
}
\vsk is commutative. Writing down the maps $\pi^2_{M;(X,Y)}$ and $\bar{\pi}^2_{M;(X,Y)}$ one needs to show that the following diagram is commutative:
\vsk
{\SMALL
\begin{tikzcd}[column sep=large, row sep=large]
 Z^2(M) 
 \arrow[r,"Z(q_M)"] \arrow[d, "\pi_{Z(M);\;X\ot Y}"] 
 & 
 Z\barz(M) \arrow[r,"q_{\bar Z(M)}"] \arrow[d, dashed, "{\pi}_{\barz(M);\;X}"]
 &  
 \barz^2(M)  \arrow[d, "\bar{\pi}_{\barz(M);\;X\ot Y}"]
  \\
  X\ot Z(M)\ot X^* \arrow[d, "\id_X\ot \pi_{M;\;Y} \ot \id_{X^*}"] \arrow[r, dashed, "\id_X\ot q_M\ot \id_{X^*}"] &   X\ot \barz(M)\ot X^* \arrow[r, equal] &   X\ot \barz(M)\ot X^* \arrow[d, "\id_X\ot \bar{\pi}_{M;\;Y} \ot \id_{X^*}"]
  \\
(X\ot Y)\ot M \ot (X\ot Y)^*   \arrow[rr, equal]  &  &(X\ot Y)\ot M \ot (X\ot Y)^* 
\end{tikzcd}
}
\vsk
Note that the upper-left rectangle is commutative by the definition of $Z(q_M)$. The upper right rectangle is commutative by the definition of $q_{\barz(M)}$. The bottom rectangle is commutative by the definition $q_M$.
\vsk 
Similarly for the proof of the first diagram we have to show that
$$v_M=\barz(q_M)\circ q_{\barz(M)}\circ \delta_M,$$
i.e. the following diagram is commutative:
\vsk
{\Tiny
\begin{tikzcd}[column sep=large, row sep=large]
 Z(M) \arrow[r,"\delta_M"] \arrow[d,"q_M"] 
 & 
 Z^2(M) \arrow[r,"q_{\barz(M)}"] 
 &  
 Z\barz(M)   \arrow[r,"\barz(q_M)"]
 & 
 \barz^2(M)   \arrow[d,"\bar{\pi}^2_{M;(X,Y)}"]
 \\
(X\ot Y)\ot M \ot (X\ot Y)^*   \arrow[rrr,equal]  & & &(X\ot Y)\ot M \ot (X\ot Y)^* 
\end{tikzcd}
}
\vsk Writing down again the map $\bar{\pi}^2_{M;(X,Y)}$ one needs to show that the following diagram is commutative:
\vsk
{\Tiny
\begin{tikzcd}[column sep=large, row sep=large]
 Z(M) 
 \arrow[r,"\delta_{M}"] \arrow[ddd, "q_M"] \arrow[dddrrrr, dashed, "\pi_{M;\;X\ot Y}"] 
 & 
 Z^2(M) \arrow[r,"q_{Z(M)}"] \arrow[drr, dashed, "{\pi}_{Z(M);\;X}"] 
 &  
 Z\barz(M) 
  \arrow[rr, "{\barz(q_M)}"] 
  &
 & 
 \barz^2(M) 
  \arrow[d, "\bar{\pi}_{\barz(M);\;X}"]  \\
 & 
 &
 &  
 X\ot Z(M) \ot X^* \arrow[ddr, dashed, "\id_X\ot {\pi}_{M;\; Y}\ot \id_{X^*}"] \arrow[r, dashed, "\id_X\ot q_{M}\ot \id_{X^*}"]
  & 
  X\ot \barz(M)\ot X^* \arrow[dd, "\id_X\ot \bar{\pi}_{M;\;Y} \ot \id_{X^*}"] 
  \\ & & & &
  \\
\barz(M)  \arrow[rrrr, "\bar{\pi}_{M;\;X\ot Y}"] &  &  & & (X\ot Y)\ot M \ot (X\ot Y)^*  
\end{tikzcd}
}
\vsk
 Note that the most left down triangle commutes by the definition of $q_M$. The middle inside rectangle commutes by definition of $\delta_M$. The upper right corner rectangle commutes by the definition of $Z(q_M)$. The most right bottom triangle commutes by the definition of $q_M$.
\epf
\br\label{rvq}
Using Lemma \ref{q2} note that the commutativity of both diagrams is equivalent to the equality $$v_M=\delta_M\circ q^{(2)}_M.$$
\er
\subsection{Compatibility with characters and central elements}
The next lemma shows that the two maps induced by $q_\unu:Z(\unu)\ra \barzu$  are $\kk$-algebra homomorphisms.
\bl With the above notations:
\bne
\item
The epimorphism $q$ induces a surjective algebra map
\beq
\qust_\unu:\cec \ra \ced,\;\; z\mapsto q\circ z
\eeq
\item
The epimorphism $q$ induces an injective algebra map
\beq
{\qst}_\unu:\cfcd \ra \cfc,\;\;\mu \mapsto \mu \circ q
\eeq
\ene
\el
\bpf 
By \cite[Section 4.3]{scalg} one has that $q_\unu$ is an epimorphism of algebras. It follows that $\qust_\unu$ is surjective and ${\qst}_\unu$ is injective.
\bne
\item Let $a, b \in \cfcc$. Recall that $a.b=m\circ (a\ot b)$. Since $q^*(a)=a\ot q$ this is clearly a surjective map. In order to see that $q^*(a.b)=q^*(a).q^*(b)$ one has to consider the following diagram:
\vsk
\newcommand{\zu}{Z(\unu)}

{\center \tiny
\begin{tikzcd}[column sep=large,row sep=large]
\unu \simeq \unu \ot \unu \arrow[r,"a \ot b"] \arrow[d, equal] & \zu\ot \zu \arrow[r,"m"]  \arrow[d, "q\ot q"]&  \zu    \arrow[d,"q "]
\\
\unu \simeq \unu \ot \unu \arrow[r,"q^*(a) \ot q^*(b)"]  & \barzu\ot \barzu \arrow[r,"\bar m"]  &\barzu 
\end{tikzcd}
}
\vsk
The left rectangle is commutative by the naturality of the tensor product and the second one since $q_\unu$ is an algebra map.
\vsk
\item
\vsk Let $f, g \in \cfcd$. In order to show that ${\qst}_\unu(f \star g)={\qst}_\unu(f )\star {\qst}_\unu(g)$ one has to consider the following diagram.

{\tiny
\begin{tikzcd}[column sep=large,row sep=large]
\zu
  \arrow[r,"\delta_\unu"]
  \arrow[d,"q_\unu"]
&
Z(\zu)
  \arrow[r,"Z(q_\unu)"]
& 
Z(\barzu)  \arrow[r,"Z(g)"]   \arrow[d,"q_{\bar A}"]
& 
Z(\unu) \arrow[r,"q_\unu"] \arrow[dr,"\mathclap{q_\unu}"',sloped]
  & 
  \barzu  \arrow[r,"f"] \arrow[d, equal]
&
\unu \arrow[d, equal]
\\
\barzu  \arrow[rr,"\bar{\delta}_\unu"]  & &\bar Z(\barzu)   \arrow[rr,"\barz(g)"] & &  \barzu  \arrow[r,"f"] & \unu 
\end{tikzcd}
}
\vsk
Note that the left rectangle is commutative since $q$ is a natural transformations of comonads. The second rectangle is commutative by the naturality of $q$.
\ene
\epf
\bibliographystyle{alpha}
\bibliography{24nov}
\ed
\section{Appendix} Let $\cc$ be any finite tensor category and $\cd$ a full subcategory of $\cc$.
In this appendix we prove the main properties of the ends 
$$
Z_\cc(M)=\int_{X\in \cc} X\ot M\ot X^*
$$ 
and 
$$
Z_\cd(M)=\int_{X\in \cd} X\ot M\ot X^*
$$  
that were used in the previous sections of the paper. For simplicity we write $Z:=Z_\cc$ and $\barz:=Z_\cd$. We will also denote by $\pi_{M;\; X}: Z(M)\ra X\ot M\ot X^*$ and respectively $\bar{\pi}_{M;\; X}: \barz(M)\ra X\ot M\ot X^*$ their universal diantural transformations.

\vsk
By  Fubini theorem for ends one has that 
$$
Z^2(M)=\int_{(X, Y)\in \cc\times \cc} (X\ot Y)\ot M \ot (X\ot Y)^*
$$ 
and 
$$
\barz^2(M)=\int_{(X, Y)\in \cd\times \cd} (X\ot Y)\ot M \ot (X\ot Y)^*
$$ 
\vsk with universal dinatural transformations:
$$\pi^2_{M; (X, Y)}=Z^{2}(M)\xra{\pi_{Z(M); X}} X\ot Z(M) \ot X^*\xra{\id\ot \pi_{M,Y}\ot \id}(X\ot Y)\ot M \ot (X\ot Y)^*$$
respectively
$$\bar{\pi}^2_{M; (X, Y)}=\barz^{2}(M)\xra{\bar{\pi}_{\barz(M); X}} X\ot \barz(M) \ot X^*\xra{\id\ot \bar{\pi}_{M,Y}\ot \id}(X\ot Y)\ot M \ot (X\ot Y)^*.$$
\vsk
Using the universal property of the ends there are unique maps 
$$\delta_M: Z(M)\ra Z^2(M), \; \eps_M:Z(M)\ra M$$
making $Z:\cc \ra \cc$ a comonad on $\cc$. Similalry $\barz:\cc \ra \cc$ is a comonad with structure maps 
$$\bar{\delta}_M: \bar Z(M)\ra \barz^2(M), \; \bar{\eps}_M:\bar Z(M)\ra M$$
\subsection{Canonical  action of $A$ on the objects of $\cc$.}
 Any object $X\in \cc$ is canonically an $A$-module in $\cc$, via:
 \beq\label{rox}
 \ro_X:A\ot X\xra{\pi_{\unu, X}\ot \id_X}X\ot X^*\ot X\xra{\id_X\ot ev_X}X.
 \eeq
 \noindent
 By \cite[Equation 3.12]{scalg} one has that
 \beq\label{rox2}
 \ro_X=(A\ot X\xra{\sg_X}X\ot A\xra{\id_X\ot \eps_\unu}X).
 \eeq
\subsection{On the compatibility $q\ot q$ and $Z_2, \barz_2$. Morphisms of comonoidal functors.}
\bll{\Small The diagram on the last page of yellow book.}
\subsection{On the map $q^{(2)}_M$. }

By the universal property of the end $\barz^2(M)=\int_{(X, Y)\in \cd\times \cd} X\ot Y\ot M\ot (X\ot Y)^*$ there  is a unique natural transformation $q^{(2)}: Z^2\ra \bar Z^2$ such that for any two objects $X, Y\in \co(\cd)$  the following diagram
\vsk
{\Small
\begin{tikzcd}[column sep=large,row sep=large]
 Z^2(M) \arrow[rr,"q^{(2)}_{M}"] \arrow[dr,"{\pi}^2_{M;\;(X, Y)}", sloped] & & \bar Z^2(M) \arrow[dl,"\bar{\pi}^2_{M;\;(X, Y)}", sloped]
 \\ & (X\ot Y)\ot  M \ot   (X\ot Y)^* 
\end{tikzcd}
}
\vsk is commutative.

\subsubsection{Definition of $r^{(2)}_M$}

There is a dinatural map in $X \in \co(\cd)$ given by
$$Z^2(M)\xra{\pi_{Z(M);\;X}} X\ot Z(M) \ot X^*\xra{\id_X\ot q_M\ot \id_{X^*}} X\ot \barzm \ot X^*. $$
Since $\barz^2(M)=\int_{X\in \cc} X\ot \barzv \ot X^*$ it follows that there exists a unique map $$r^{(2)}_M: Z^2(M) \ra \barz^2(M)$$
\vsk such that:
{\Small
\begin{tikzcd}[column sep=large,row sep=large]
 Z^2(M) \arrow[rr,"r^{(2)}_M"] \arrow[d, "{\pi}_{Z(M);\;X}"] & & \bar Z^2(M) \arrow[d,"\bar{\pi}_{\barzm;\;X}"]
 \\  X\ot Z M \ot X^* \arrow[rr,"\id_X\ot q_M\ot \id _{X^*}"]  & & X\ot  \barzm \ot X^*
\end{tikzcd}
}
\vsk We will show that $r^{(2)}_M=q^{(2)}_M$. For this we have to show that 
\vsk
{\Small
\begin{tikzcd}[column sep=large,row sep=large]
 Z^2(M) \arrow[rr,"r^{(2)}_M"] \arrow[dr,"{\pi}^2_{M, (X, Y)}", sloped] & & \bar Z^2(M) \arrow[dl,"\bar{\pi}^2_{M, (X, Y)}", sloped]
 \\ & (X\ot Y)\ot  M \ot   (X\ot Y)^* 
\end{tikzcd}
}
\vsk Expanding $\pi^2$ and $\bar{\pi}^2$ we have to show that:
\vsk
{\Small
\begin{tikzcd}[column sep=large,row sep=large]
 Z^2(M) 
 \arrow[rr,"r^{(2)}_M"] \arrow[d,"{\pi}^2_{Z(M);\;X}"] 
 & 
 & 
 \bar Z^2(M) 
 \arrow[d,"\bar{\pi}^2_{Z(M);\; X}"]
 \\ 
 X\ot Z(M) \ot X^*  
 \arrow[dr,"\idx \ot {\pi}^2_{M;\; Y}\ot \id_{X^*}", sloped] \arrow[rr, , dashed, "\idx \ot q_M\ot \id_{X^*}"] 
 & 
 & 
 X\ot \barzv \ot  X^*  \arrow[dl,"\idx \ot \bar{\pi}^2_{M, Y}\ot \id_{X^*}", sloped]
 \\ 
 & (X\ot Y)\ot  M \ot   (X\ot Y)^* &
\end{tikzcd}
}
\vsk The bottom triangle is commutative by the definition of $q_M$. The upper rectangle is commutative by the definition of $r^{(2)}_M$.
\bl \label{q2}
 One has that 
 $$
 q^{(2)}_M=\barz(q_M)\circ q_{Z(M)}=q_{\barz(M)}\circ  Z(q_M).
 $$
 \el
 \bpf 
 For the first equality one has to show that the following diagram is commutative:
 \vsk
{\Small \begin{tikzcd}[column sep=large, row sep=large]
 Z^2(M) \arrow[r,"q_{Z(M)}"] \arrow[d, "{\pi}_{Z(M);\;X}"] & \barz Z(M) \arrow[dl, dashed,  "\bar{\pi}_{Z(M);\;X}"]  \arrow[r,"\barz(q_M)"]& \bar Z^2(M) \arrow[d,"\bar{\pi}_{\barzm;\;X}"]
 \\  X\ot Z M \ot X^* \arrow[rr,"\id_X\ot q_M\ot \id _{X^*}"]  & & X\ot  \barzm \ot X^*
\end{tikzcd}
}
\vsk Note that the left triangle is commutative by the definition of $q_{Z(M)}$. The next square commutes by the definition of $\barz(q_M)$.
\vsk  For the second equality one has to show that the following diagram is commutative:
 \vsk
 {\Small \begin{tikzcd}[column sep=large, row sep=large]
 Z^2(M) \arrow[r,"Z(q_M)"] \arrow[d, "{\pi}_{Z(M);\;X}"] 
 & 
 Z\barz (M) 
 \arrow[dr, dashed,  "\bar{\pi}_{Z(M);\;X}"]  \arrow[r,"q_{\barzm}"]
 & 
 \bar Z^2(M) 
 \arrow[d,"\bar{\pi}_{\barzm;\;X}"]
 \\  
 X\ot Z M \ot X^* 
 \arrow[rr,"\id_X\ot q_M\ot \id _{X^*}"]  
 & 
 & X\ot  \barzm \ot X^*
\end{tikzcd}
}
\vsk Note that the left square commutes by the definition of $Z(q_M)$ while the right triangle commutes by the definition of $q_{\barzm}$.
 \epf
\subsection{On the map $v_M$}

\bl \label{vm}
There is a unique  natural transformation $v_M:Z(M)\ra \barz^2(M)$ such that 
$$\bar{\pi}^2_{M;\; (X,Y)}\circ v_M=  \bar{\pi}_{M;\; X\ot Y}\circ q_M$$
\vsk 
i.e. such that the following diagram is commutative:
\vsk
{\Small
\begin{tikzcd}[column sep=large,row sep=large]
 Z(M) \arrow[r,"v_M"] \arrow[d,"q_M"]&  \bar Z^2(M) \arrow[d,"\bar{\pi}^2_{M;\; (X, Y)}"]
 \\ \barz(M) \arrow[r,"\bar{\pi}_{M; X\ot Y}"]& (X\ot Y)\ot  M \ot   (X\ot Y)^* .
\end{tikzcd}
}
\vsk
One has that $$v_M=\bar{\delta}_M\circ q_M.$$\
\el
\bpf One has as above $$\barz^2(M)=\int_{(X, Y)\in \cd\times \cd }(X\ot Y)\ot M \ot (X\ot Y)^*.$$
\vsk The maps
$$Z(M)\xra{q_M}\barz(M)\xra{\bar{\pi}_{M, (X,Y)}} (X\ot Y)\ot M \ot (X\ot Y)^*$$ are dinatural with respect to $(X,Y)$. Thus by the universal property of the coend there is a unique map $v_M$ making the above diagram commutative. It remains to show that $v_M=\bar{\delta}_M\circ q_M.$
Thus we have to show that the following diagram is commutative \vsk
{\Small
\begin{tikzcd}[column sep=large,row sep=large]
 Z(M) \arrow[r,"q_M"] \arrow[drr, dashed, "{\pi}_{M; X\ot Y}", sloped]   \arrow[d,"q_M"]& \bar Z(M) \arrow[r,"\bar{\delta}_{M}"]  \arrow[dr, dashed, "\bar{\pi}_{M; X\ot Y}", sloped] & \barz^2(M)  \arrow[d,"\bar{\pi}^2_{M,\; (X, Y)}"]
 \\
 \barz(M) \arrow[rr,"\bar{\pi}_{M; X\ot Y}"]   &  & (X\ot Y)\ot  M \ot   (X\ot Y)^* 
\end{tikzcd}
}
\vsk
This follows since all three inside triangles are commutative.
\epf
\subsection{Natural transformation of comonads.}\vsk \blue{\Small Write the definition and compatibility with the antipodes.}
\bl $q:Z\ra \barz$ is a natural transformations of Hopf comonads.\el
\bpf
Note that $\barz:=\barz_\cd$ has also a structure of a comonad and it is associative by Fubini theorem applied for a cube transformation. We denote by $\bar{\eps}:\barz \ra \unu$ and $\bar{\delta}:\barz \ra \barz^2$ the structures of this comonad.
\vsk
In order to show that $q$ is a natural transformation of comonads we need to show that for any object $M \in \co(\cd)$ one has 
$$\bar{\eps}_M\circ q_M=\eps_M$$
and the following diagrams are commutative:
\vsk 
{\Tiny
\begin{tikzcd}[column sep=large,row sep=large]
 Z(M) \arrow[r,"\delta_M"] \arrow[d,"q_M"] & Z^2(M) \arrow[r,"q_{Z(M)}"] &  \barz Z(M)    \arrow[d,"\barz(q_{M}) "]
\\
\barz(M)  \arrow[rr,"\bar{\delta}_M"]  & &\bar Z^2(M) 
\end{tikzcd}
}
 {\Tiny
\begin{tikzcd}[column sep=large,row sep=large]
 Z(M) \arrow[r,"\delta_M"] \arrow[d,"q_M"] & Z^2(M) \arrow[r,"Z(q_M)"] &  Z\barz(M)   \arrow[d,"q_{\barz(M)}"]
\\
\barz(M)  \arrow[rr,"\bar{\delta}_M"]  & &\bar Z^2(M) 
\end{tikzcd}
}
\vsk  Since $\eps_M=\pi_{M;\unu}$ and $\bar{\eps}_M=\bar{\pi}_{M;\unu}$ the first formula follows from the definition of $q_M$.
\vsk We first show that the  second  diagram is commutative. By Lemma \ref{vm} it its enough to show that $$v_M=q_{\barz(M)}\circ Z(q_M)\circ \delta_M$$ i.e. the following diagram is commutative:
\vsk
{\Small
\begin{tikzcd}[column sep=large, row sep=large]
 Z(M) \arrow[r,"\delta_M"] \arrow[d,"q_M"] \arrow[drrr, dashed, "\pi_{M;\;X\ot Y}"] & Z^2(M) \arrow[r,"Z(q_M)"] \arrow[drr, dashed, "\pi^2_{M;\;X\ot Y}"]&  Z\barz(M)   \arrow[r,"q_{\barz(M)}"]
 & \barz^2(M)   \arrow[d,"\bar{\pi}^2_{M;(X,Y)}"]\\
\barz(M)  \arrow[rrr,"\bar{\pi}_{M;\;(X\ot Y)}"]  & & &(X\ot Y)\ot M \ot (X\ot Y)^* 
\end{tikzcd}
}
\vsk
The most left bottom triangle is commutative by the definition of $q_M$. The middle triangle is commutative by the definition of $\delta_M$. Thus it remains to show that the diagram:
\vsk
{\Small
\begin{tikzcd}[column sep=large, row sep=large]
 Z^2(M) 
 \arrow[r,"Z(q_M)"] \arrow[d, "\pi^2_{M;\;X\ot Y}"] 
 & 
 Z\barz(M) \arrow[r,"q_{Z(M)}"] 
 &  
 \barz^2(M)   \arrow[d,"\bar{\pi}^2_{M;\;(X,Y)}"]
  \\
(X\ot Y)\ot M \ot (X\ot Y)^*   \arrow[rr, equal]  &  &(X\ot Y)\ot M \ot (X\ot Y)^* 
\end{tikzcd}
}
\vsk is commutative. Writing down the maps $\pi^2_{M;(X,Y)}$ and $\bar{\pi}^2_{M;(X,Y)}$ one needs to show that the following diagram is commutative:
\vsk
{\SMALL
\begin{tikzcd}[column sep=large, row sep=large]
 Z^2(M) 
 \arrow[r,"Z(q_M)"] \arrow[d, "\pi_{Z(M);\;X\ot Y}"] 
 & 
 Z\barz(M) \arrow[r,"q_{\bar Z(M)}"] \arrow[d, dashed, "{\pi}_{\barz(M);\;X}"]
 &  
 \barz^2(M)  \arrow[d, "\bar{\pi}_{\barz(M);\;X\ot Y}"]
  \\
  X\ot Z(M)\ot X^* \arrow[d, "\id_X\ot \pi_{M;\;Y} \ot \id_{X^*}"] \arrow[r, dashed, "\id_X\ot q_M\ot \id_{X^*}"] &   X\ot \barz(M)\ot X^* \arrow[r, equal] &   X\ot \barz(M)\ot X^* \arrow[d, "\id_X\ot \bar{\pi}_{M;\;Y} \ot \id_{X^*}"]
  \\
(X\ot Y)\ot M \ot (X\ot Y)^*   \arrow[rr, equal]  &  &(X\ot Y)\ot M \ot (X\ot Y)^* 
\end{tikzcd}
}
\vsk
Note that the upper-left rectangle is commutative by the definition of $Z(q_M)$. The upper right rectangle is commutative by the definition of $q_{\barz(M)}$. The bottom rectangle is commutative by the definition $q_M$.
\vsk 
Similarly for the proof of the first diagram we have to show that
$$v_M=\barz(q_M)\circ q_{\barz(M)}\circ \delta_M,$$
i.e. the following diagram is commutative:
\vsk
{\Tiny
\begin{tikzcd}[column sep=large, row sep=large]
 Z(M) \arrow[r,"\delta_M"] \arrow[d,"q_M"] 
 & 
 Z^2(M) \arrow[r,"q_{\barz(M)}"] 
 &  
 Z\barz(M)   \arrow[r,"\barz(q_M)"]
 & 
 \barz^2(M)   \arrow[d,"\bar{\pi}^2_{M;(X,Y)}"]
 \\
(X\ot Y)\ot M \ot (X\ot Y)^*   \arrow[rrr,equal]  & & &(X\ot Y)\ot M \ot (X\ot Y)^* 
\end{tikzcd}
}
\vsk Writing down again the map $\bar{\pi}^2_{M;(X,Y)}$ one needs to show that the following diagram is commutative:
\vsk
{\Tiny
\begin{tikzcd}[column sep=large, row sep=large]
 Z(M) 
 \arrow[r,"\delta_{M}"] \arrow[ddd, "q_M"] \arrow[dddrrrr, dashed, "\pi_{M;\;X\ot Y}"] 
 & 
 Z^2(M) \arrow[r,"q_{Z(M)}"] \arrow[drr, dashed, "{\pi}_{Z(M);\;X}"] 
 &  
 Z\barz(M) 
  \arrow[rr, "{\barz(q_M)}"] 
  &
 & 
 \barz^2(M) 
  \arrow[d, "\bar{\pi}_{\barz(M);\;X}"]  \\
 & 
 &
 &  
 X\ot Z(M) \ot X^* \arrow[ddr, dashed, "\id_X\ot {\pi}_{M;\; Y}\ot \id_{X^*}"] \arrow[r, dashed, "\id_X\ot q_{M}\ot \id_{X^*}"]
  & 
  X\ot \barz(M)\ot X^* \arrow[dd, "\id_X\ot \bar{\pi}_{M;\;Y} \ot \id_{X^*}"] 
  \\ & & & &
  \\
\barz(M)  \arrow[rrrr, "\bar{\pi}_{M;\;X\ot Y}"] &  &  & & (X\ot Y)\ot M \ot (X\ot Y)^*  
\end{tikzcd}
}
\vsk
 Note that the most left down triangle commutes by the definition of $q_M$. The middle inside rectangle commutes by definition of $\delta_M$. The upper right corner rectangle commutes by the definition of $Z(q_M)$. The most right bottom triangle commutes by the definition of $q_M$.
\epf
\br\label{rvq}
Using Lemma \ref{q2} note that the commutativity of both diagrams is equivalent to the equality $$v_M=\delta_M\circ q^{(2)}_M.$$
\er
\subsection{Compatibility with characters and central elements}
The next Lemma shows that these two maps  are $\kk$-algebra homomorphisms.
\bl With the above notations:
\bne
\item
The epimorphism $q$ induces a surjective algebra map
\beq
\qust_\unu:\cec \ra \ced,\;\; z\mapsto q\circ z
\eeq
\item
The epimorphism $q$ induces an injective algebra map
\beq
{\qst}_\unu:\cfcd \ra \cfc,\;\;\mu \mapsto \mu \circ q
\eeq
\ene
\el
\bpf 
\bne
\item Let $a, b \in \cfcc$. Recall that $a.b=m\circ (a\ot b)$. Since $q^*(a)=a\ot q$ this is clearly a surjective map. In order to see that $q^*(a.b)=q^*(a).q^*(b)$ one has to consider the following diagram:
\vsk
{\center \tiny
\begin{tikzcd}[column sep=large,row sep=large]
\unu \simeq \unu \ot \unu \arrow[r,"a \ot b"] \arrow[d, equal] & A\ot A \arrow[r,"m"]  \arrow[d, "q\ot q"]&  A    \arrow[d,"q "]
\\
\unu \simeq \unu \ot \unu \arrow[r,"q^*(a) \ot q^*(b)"]  & \bar A\ot \bar A \arrow[r,"\bar m"]  &\bar A 
\end{tikzcd}
}
\vsk
The left rectangle is commutative by the naturality of the tensor product and the second one since $q_\unu$ is an algebra map.
\vsk
\item
\vsk Let $f, g \in \cfcd$. In order to show that ${\qst}_\unu(f \star g)={\qst}_\unu(f )\star {\qst}_\unu(g)$ one has to consider the following diagram.

{\tiny
\begin{tikzcd}[column sep=large,row sep=large]
 A
  \arrow[r,"\delta_\unu"]
  \arrow[d,"q_\unu"]
&
Z(A)
  \arrow[r,"Z(q_\unu)"]
& 
Z(\bar A)  \arrow[r,"Z(g)"]   \arrow[d,"q_{\bar A}"]
& 
Z(\unu) \arrow[r,"q_\unu"] \arrow[dr,"\mathclap{q_\unu}"',sloped]
  & 
  \barzu  \arrow[r,"f"] \arrow[d, equal]
&
\unu \arrow[d, equal]
\\
\bar A  \arrow[rr,"\bar{\delta}_\unu"]  & &\bar Z(\bar A)   \arrow[rr,"\barz(g)"] & &  \barzu  \arrow[r,"f"] & \unu 
\end{tikzcd}
}
\vsk
Note that the left rectangle is commutative since $q$ is a natural transformations of comonads. The second rectangle is commutative by the naturality of $q$.
\ene
\epf
\subsection{Extra topics that are not needed yet in the paper}
\subsubsection{On the maps $\mtc P$ and $\mtc Q$ in a finite tensor category} In any finite tensor category $\cc$, for any three objects $M, N, X\in \co(\cc)$ there are canonical natural bijections  $$\mtc P:\hm_\cc(M\ot X, N)\ra \hm_{\cc}(M, N\ot X^*)$$ and  $$ \mtc Q:\hm_{\cc}(M, N\ot X^*) \ra \hm_\cc(M\ot X, N)$$
inverse one to another. It is easy to show that $\mtc P$ and $\mtc Q$ send commutative diagrams in commutative diagrams.
\vsk This means that for any two arrows $f: A\ra A'$ and  $g: B\ra B'$ the left commutative diagram 
\vsk
\begin{tikzcd}[column sep=large,row sep=large]
A \arrow[r,"f"]  \arrow[d,"u"] 
& 
A'   \arrow[d, "v"]
\\
 X \ot B \ot X^* \arrow[r," \id_X\ot g \ot \id _{X^*}"]  &  X \ot B'\ot X^*
 \end{tikzcd}
$\xra{\mtc Q}$
\begin{tikzcd}[column sep=large,row sep=large]
A \ot X \arrow[r,"f\ot \id_X"]  \arrow[d,"\mtc Q(u)"] 
& 
A' \ot X  \arrow[d, "\mtc Q(v)"]
\\
 X \ot B \arrow[r," \id_X\ot g"]  &  X \ot B'\end{tikzcd}
\vsk is sent by $\mtc Q$ in the right commutative diagram. 
\subsubsection{On the rhomb applications and images in the center} For any $M, X\in \co(\cc)$ there is a natural map 
$$R^\cc_{M, X}:Z(M)\ot X\ra M\ot X^*$$ corresponding by $\mtc P$ to the universal maps $\pi_{M;\;X}$.
Similarly, for any $M \in \co(\cc)$ and $X \in \co(\cd)$ there is a natural map 
$$R^\cd_{M, X}:\barz(M)\ot X\ra M\ot X^*$$ corresponding by $\mtc P$ to the universal maps $\bar{\pi}_{M;\;X}$.
\subsection{Applying $\mtc Q$ to the definition of $q_M$} It follows that the following diagram is commutative:
\vsk
\begin{tikzcd}[column sep=large,row sep=large]
Z(M) \ot X \arrow[r,"q_M\ot \id_X"]  \arrow[d,"R^\cc_{M, X}"] 
& 
\barz(M) \ot X  \arrow[d, "R^\cd_{M, X}"]
\\
 X \ot M \arrow[r, equal]  &  X \ot M\end{tikzcd}
\subsubsection{Natural transformation lemma from Sakalos paper}
\bl [Lemma 6.4 \cite{sak}] If $\cc$ is a finite tensor category and $Y, W, M \in \co(\cc)$ then there is an isomorphism
\beqn
\phi:\hm_\cc(W, Y\ot Z(M))\simeq \nat (W\ot \id_{\cc}, \;Y\ot \id_{\cc} \ot M)\eeqn
\vsk given by
\beqn 
f \mapsto (\phi(f)_T:W\ot T\xra{} Y\ot Z(M) \ot  T\xra{\id \ot \ccrm} Y\ot T\ot M)
\eeqn
\el
\vsk 
Completely similalry one can prove that there is an isomorphism
\beqn
\bar{\phi}:\hm_\cc(W, Y\ot \barzm)\simeq \nat (W\ot \id_{\cd}, \;Y\ot \id_{\cd} \ot M)\eeqn
\vsk given by
\beqn 
f \mapsto (\bar{\phi}(f)_T:W\ot T\xra{} Y\ot Z(M) \ot  T\xra{\id \ot \ccrm} Y\ot T\ot M)
\eeqn
For $Y=\unu$ in both equations above one has that
\beqn
\phi:\hm_\cc(W,\; \barzm)\simeq \nat (W\ot \id_{\cc}, \; \id_{\cc} \ot M)
\eeqn
\vsk  and 
\beqn 
\bar{\phi}:\hm_\cc(W, \; \barzm)\simeq \nat (W\ot \id_{\cd}, \; \id_{\cd} \ot M)
\eeqn
\vsk \blue{\Small The first item of the  following Lemma appears in \cite{DS}, see also \cite[shimizu-integrals], write which section. The second item is proven completely similarly. }
\bl \label{cstr}Let $\cc$ be a finite tensor category and $M \in \co(\cc)$. With the above notations one has that:
\bne
\item
A map $f \in \hm_\cc(M, \; Z(M))$ has $Z$-comodule structure if and only if the corresponding natural transformation $\phi(f)$ is a half-braiding.
\item
A map $f \in \hm_\cc(M, \; \barzm)$ has $\barz$-comodule structure if and only if the corresponding natural transformation $\bar{\phi}(f)$ is a relative half-braiding.
\ene
\el
\subsubsection{On the centrality of $\delta_M$ and relative centrality of $\bdlt_M$.} Since $\delta_M:Z(M)\ra Z^2(M)$ is a comodule structure from Lemma \ref{cstr} above one has 
{\Small  $$[\sg^{st}_{Z(M), X}:Z(M)\ot X\ra X\ot Z(M)]=[Z(M)\ot X\xra{\dlt_M} Z(M)\ot X\xra{R^{\cc}_{Z(M), X}}X\ot Z(M)]$$
}
\vsk is a half-braiding. Thus one  has that  $R(M)=(Z(M), \sg^{\mtr{st}}_{M, -})\in \czcc$. Centrality of $\delta_M$ then reduces to the commutativity of the following diagram 
\vsk
\begin{tikzcd}[column sep=large,row sep=large]
M \ot X \arrow[r,"\sgst_{M, X}"]  \arrow[d,"\delta_M\ot \id"]  & X \ot M \arrow[d,"\id \ot \delta_M"] 
\\
Z^2(M)\ot X  \arrow[r,"\sgst_{Z^2M, X}"]  &X \ot Z^2M
\end{tikzcd}
\vsk 
which in  turn becomes
\vsk
\begin{tikzcd}[column sep=large,row sep=large]
M \ot X \arrow[r,"\dltv\ot \id_X"]  \arrow[d,"\delta_M\ot \id_X"] 
& 
Z^2(M)\ot X \arrow[r,"R^\cc_{M, X}"] \arrow[d, dashed, "Z(\delta_M)\ot \id_X"]
&
X \ot M \arrow[d,"\id_X \ot \delta_M"] 
\\
Z^2(M)\ot X  \arrow[r,"\delta_{M}\ot \id_X"]  & Z^3(M)\ot X \arrow[r,"R^\cc_{Z^2M, X}"]
& 
X \ot Z^2M
\end{tikzcd}
\vsk Therefore it reduces to verify a compatibility of the rhomb application with $\delta_M$. This follows by applying the map $\mtc Q$  to the following diagram:
\vsk
\begin{tikzcd}[column sep=large,row sep=large]
Z^2(M) \arrow[r,"Z(\delta_M)"]  \arrow[d,"\pi_{M;\;X}"]  &  Z^3M \arrow[d,"\pi_{Z^2M;\;X}"] 
\\
 X\ot M \ot X^*  \arrow[r,"\id \ot \delta_M \ot \id"]  & X\ot Z^2(M)\ot X^*
\end{tikzcd}
\vsk Note that this diagram is commutative by definition of $Z(\delta_M)$.
\vsk Also, since $\bar{\delta}_M:\barz(M)\ra \barz^2(M)$ is a comodule structure from Lemma \ref{cstr} above one has that
{\Small  $$[\sg^{st}_{\barzm, X}:\barzm \ot X\ra X\ot \barzm]=[\barzm\ot X\xra{\bdlt_M} \barz^2(M)\ot X\xra{R^{\cd}_{\barzm, X}}X\ot \barzm].$$
}
\noindent is a relative half-braiding. Thus one also has $R_1(M)=(\barzm, \sgst_{\barzm, -} )\in \cz_\cc(\cd)$. It is proven completely similarly as above that the morphism $\bdlt_M:\barzm \ra \barz^2(M)$ is a map in the relative center $\cz_\cc(\cd)$.
\subsection{$q_M:Z(M)\ra \bar Z(M)$ is a relative central morphism} In this subsection we show that $q_M$ is a morphism in $\cz_\cc(\cd)$.
In order to do this one needs to verify that the diagram\vsk
\begin{tikzcd}[column sep=large,row sep=large]
Z(M) \ot X \arrow[r,"\sgst_{Z(M), X}"]  \arrow[d,"q_M\ot \id_X"]  & X \ot Z(M) \arrow[d,"\id_X \ot q_M"] 
\\
\barzm \ot X  \arrow[r,"\sgst_{\barzm, X}"]  &X \ot \barzm
\end{tikzcd}
\vsk 
is commutative for any object $X \in \co(\cd)$. In turn this diagram becomes:
\vsk
\vsk
\begin{tikzcd}[column sep=large,row sep=large]
Z(M) \ot X \arrow[r,"\dlt_M\ot \id_X"]  \arrow[d,"q_M\ot \id_X"] & 
Z^2(M)\ot X \arrow[r,"R^\cc_{M, X}"] \arrow[d, dashed, "q^{(2)}_M\ot \id_X"]
&
X \ot M \arrow[d,"\id_X \ot \delta_M"] 
\\
\barzm \ot X  \arrow[r,"\bdlt_{M}\ot \id_X"]  & \barz^2(M)\ot X \arrow[r,"R^\cd_{M, X}"]
& 
X \ot \barzm
\end{tikzcd}

\vsk
The left square is commutative by Remark \ref{rvq}. For the second square note that the diagram 
\vsk
\begin{tikzcd}[column sep=large, row sep=large]
Z^2(M) \arrow[d,"q^{(2)}_M"]  \arrow[r,"\pi_{Z(M);\;X}"] 
& 
X \ot Z(M) \ot X^* \arrow[d,  "\id_X \ot q_M \ot \id_{X^*}"]
\\
\barz^2(M)  \arrow[r,"\bar{\pi}_{\barzm;\;X}"]   & X \ot \bar Z(M) \ot X^* \end{tikzcd}
\vsk is commutative. \blue{Then applying $\mtc P$ to the above diagram one obtains the above diargam. }
\subsection{Reference for $q_\unu$ epimorphism of algebras}
{\Small  A reference probably for this epimorphism $q$-shimizu-further results on coends!} 
\vsk Since the epimorphism $q:A_\cc \ra A_\cd$ is a morphism of algebras {\Small in the relative center or at least in $\cc$!!} it follows that
\beq
u_\cd q=u_\cc
\eeq 
\subsubsection{Interpretation in the fusion case}
{\sml As a morphism in $\cc$ is given by the projections on $X_i\ot X_i^*$ with $X_i\in \cd$. Need to understand what it does on the idempotents $F^\cd_i$ or how the conjugacy classes from $\cc$ restrict as conjugacy classes in $\cd$.}
\subsubsection{Restriction of the dimension function}

Moreover as above this induces an epimorphism $\pi^*:\cecc \ra \cecd$ which is unitary morphism of $\kk$-algebras. This implies that the dimension function on $\cc$
\beq d_\cc: \cfcc \ra 
\kk,\;\;\ch\mapsto \langle \ch, u_\cc\rangle 
\eeq
restricts to the corresponding dimension function $d_\cd:\cfcc\ra \kk$.
\vsk
\subsection{Definition of an integral}
\bll{\Small The result follows from the fact that $q^*$ is a surjective algebra map onto $\cecd$.}\vsk
$\eps_\unu$ is a morphism of algebras. An integral $\blam:\unu \ra Z(\unu)$ is a  morphism of left $A$-modules $\blam:\unu \ra Z(\unu)$ in $\cc$.
\vsk
By definition this means that the diagram:
\vsk \begin{tikzcd}[column sep=large,row sep=large]
A=A\ot \unu \arrow[r,"\id_A\ot \blam"] \arrow[d,"\eps_\unu\ot \blam"] & A\ot A\arrow[d,"m"] 
\\
A=\unu \ot A   \arrow[r, equal]  & A
\end{tikzcd}
\vsk
It reduces to show that $\qust(\blam_\cc)=\blam_\cd$.
\subsection{Compatibility with Fourier transforms}

Check also the compatibility with Fourier transform:
\vsk
\begin{tikzcd}[column sep=large,row sep=large]
\cfcc \arrow[r,"\mtf"]  & \cecc \arrow[d,"\qust_\unu"] 
\\
\cfcd \arrow[u,"{\qst}_\unu"]  \arrow[r,"\bar{\mtf}"]  &\cecd 
\end{tikzcd}

\newpage

\ed
\subsection{General settings-commutative case-non braided case} 
Suppose further that $\cfcc$ is commutative and $\cd\subseteq \cc$ is as above a fusion subcategory. Without loss of generality one may suppose that $\irr(\cd)=\{\ch_0, \ch_1, \dots, \ch_{r'}\}$ and $ \irr(\cc)=\{\ch_0, \ch_1, \dots, \ch_{r}\}$ with $r'\leq r$.

Denote the primitive idempotents of $\cfcc$ by $\{F_i\}_{i=0}^{r}$. Also denote the primitive idempotents  of $\cfcd$ by $\{G_m\}_{m=0}^{r'}$. Write the decomposition of the idempotent cointegral $G_0$ of $\cd$ as
\beq\label{intd}
G_0:=\sum_{j \in \mtca_0}F_j.
\eeq Also write 
$G_m:=\sum_{j \in \mtca_m}F_j$ for the other central primitive idempotents of $\cfcd$. In this way we get a partition of the index set $\{0, \dots, r\}$ in $r'$ disjoint parts $\{\mtca_i\}_{i=0}^{r'}$. Inside $\cfcd$ we write
$
\ch_i=\sumjtoprp\tilde\alpha_{im}G_m
$
and inside $\cfcc$ we write
$
\ch_i=\sumjtor\al_{ij}F_j.
$
By expanding $G_j$ it follows that 
$
\al_{ij}=\al_{ij'}=\tilde\alpha_{im}
$
whenever $j,j'$ are in the same partition component $\mtca_m$. 
\mdn
\subsection{Class Equation for tensor categories}
For any two indices $0\leq, j, j'\leq r$ define
\beqn
\beta_{jj'}:=\sum_{i=0}^{r'}\al_{ij}\al_{i^* j'}
\eeqn
On the other hand in the identity from Equation \eqref{db} inside $\cfcd$ can be written as: 
\beq
\sum_{i=0}^{r'}\ch_i\ot \ch_{i^*}=\sum_{m=0}^{r'}G_m\ot \frac{\dimcd}{\tilde f_m} G_m
\eeq
\noindent By expanding $\ch_i$ and $G_j$ in terms of the central primitive idempotents $F_j$ this equation gives that
\beq
\beta_{jj'}=0,\;\;\text{if}\; j\in \ca_m,j'\in \ca_{m'},\;\text{are in different partition components}
\eeq
\beq
\beta_m:=\beta_{jj'}=\frac{\dimcd}{\tilde f_m}, \;\;\text{if}\; j,j',\;\text{are in the same partition components}\; \mtca_m
\eeq
\bp Let $\cc$ be a fusion category with a commutative character ring and $\cd$ be a fusion subcategory of $\cc$. With the above notations define linear functions on $\cfcd$ by
 \beq
h_m: \cfcd\ra \kk,\;\; \ch_j\mapsto \frac{\widetilde{\al}_{jm}}{\chj(1)},\;\;0\leq j \leq r'.
\eeq
\bne
\item
Then
\beq
h_j=h_l,\;\; \text{if the indices $j$ and $l$ are in the same partition component $\ca_m$.}
\eeq
\item
\beq
\langle h_j, h_l\rangle =\beta_{jl}=0, \;\;\text{if the indices $j$ and $l$ are in different  partition components.}
\eeq
\ene
\ep
Thus the number of different functions $\{h_m\}$ equals the number of partition components, \blue{\bf i.e. the rank of $\cd$.}
\mdn

\mdn{\blue{\bf Take the algebra with multiplication given by
\beq
h_m \star h_n:=\sum_pN^p_{mn}h_p
\eeq}}
\mdn
\subsection{On the braided case-first subsect}
Try to use or reprove the result of Drinfeld DGNO on cosets and centralizer. 
\bp \label{main2}
In a non-degenerate braided fusion category $\cc$ \blue{\bf with $fp=dim$} the number of $\cd'$ (cosets) components is equal to the rank of $\cd$, i.e the cardinality of $|\irr(\cd)|$. \mdn Equivalently the number of $\cd$-components equal the number of simple objects in $\cd'$. 
\ep

\bpf \blue{\bf The proof uses the result main1. It assumes for a moment $\fp=d$, quantum dimension, this is the case for pseudo-unitary categories.}

Note that
\beq
\rcd=\sum_{\{j \in \mtca_0=\irr(\cd')\}}F_j
\eeq
\blue{\bf since non-degenerate there is a bijection between idempt and simple objects!! Moreover $\cd''=\cd$ in this case.}

So equivalently, we have to show that the number of components of $\cd$ equal the number of primitive idempotents entering in the formula of $\rcd$.
\mdn  Then
\beq
X\simeq_\cd Y\iff \frac{\rcd\clsx}{d(X)}=\frac{\rcd\clsy}{d(Y)}\iff \frac{s_{Xj}}{d(X)\dxj}=\frac{s_{Yj}}{d(Y)\dxj}
\eeq
If $X\in \cd'$ maybe this values are all different.

\red{Two simple objects $\ch_m$ and $\ch_n$ are in the same $\cd'$-component iff
\beq
\frac{\rcdp \ch_m}{\ch_m(1)}=\frac{\rcdp \ch_n}{\ch_n(1)}
\eeq
i.e 
\beq
\frac{\al_{mj}}{\ch_m(1)}=\frac{\al_{nj}}{\ch_n(1)}
\eeq
 for all $j \in \irr(\cd)$.
}
\epf
\mdn Consider the inner product on $\cfcd$ defined in DGNO:
\beq
(\mu_j, \mu_l)_\cd:=\sum_{X\in \co(\cd)}|X|^2\mu_j(X)\mu_l(X^*)=
\eeq
\blue{\bf 
Since $|X|^2\rangle 0$ it follows that
\beqn
(\mu_j, \mu_j)_\cd=\sum_{X\in \co(\cd)}|X|^2\mu_j(X)\mu_j(X^*)=\sum_{X\in \co(\cd)}|X|^2|\mu_j(X)|^2\rangle 0
\eeqn
based on the fact that $\mu_j(X^*)=\bar{\mu_j(X)}$. Here $\mu_j$ are the linear characters of the Grothendieck ring!}

For any $Y\in \irr(\cc)$ define also the function 
\beq
h_Y:\cfcd\ra \kk,\;\; \text{given by}\;\;\;h_Y([X])=\tilde{s}_{XY}=\frac{s_{XY}}{d_+(X)d_-(Y)}
\eeq
If $Y'$ and $Y$ are in the same $\cd$-component then \blue{from the above equations } the functions $h_{Y'}$ and $h_Y$ are equal. 
\mdn If they are in  different $\cd$- components then the functions  are orthogonal. Thus they are linear independent and they are exact the columns of $\tilde S_\cd$. \blue{\bf This allows the interpretation as the rank of the tilde matrix.}
\mdn \blue{\bf The orthogonality of $h_X$ and $h_Y$ is probably the second orthogonality relation in $\cfcd$.}

\subsection{On the braided case}
Suppose first that $\cd$ is non-degenerate. Then from above one has:
\beq
r_{\cd'}=\sum_{j \in \irr(\cd)}F_j
\eeq
Define also
\blue{\bf \beq
\teh_m: \cfcd\ra \kk,\;\; \ch_j\mapsto \frac{\al_{mj}}{\chm(1)}
\eeq
and then extend it linearly.}
 It follows that $\teh_m=\teh_n$ if $\ch_m$ and $\chn$ are in the same $\cd'$ component.
\mdn If they are in $\cd'$-different components then
\beq
(\teh_m, \teh_n)_\cd=\sum_{j=0}^{r'}\chj(1)^2\teh_m(j)\teh_m(j^*)=\sum_{j=0}^{r'}\chj(1)^2\frac{\al_{mj}}{\chm(1)}\frac{\al_{nj^*}}{\chn(1)}
\eeq

\noindent \blue{\bf Since we are in the braided case one has:}
\beq
\al_{mj}=\frac{\smj}{\chj(1)}
\eeq
and 
\beq
\aljm=\frac{\sjm}{\chm(1)}=\frac{\smj}{\chm(1)}
\eeq
Thus the inner product becomes 
\beq
(\teh_m, \teh_n)_\cd=\sumjtorp \frac{\smj\snj}{\chm(1)\chn(1)}=\sum_{j=0}^{r'}\aljm\al_{j^*n}=\beta_{mn}.
\eeq
Thus if $m$ and $n$ are in different components then $(\teh_m, \teh_n)_\cd=0$. Thus the number of $\cd'$ components equals the rank of $\cd$.
\mdn \blue{\bf Note that in the braided case $h_m(\chj)=\frac{\sjm}{\chj(1)\chm(1)}=\frac{\smj}{\chj(1)\chm(1)}=\teh_m(\chj)$.} This shows that a the simple objects of a $\cd'$-component are in bijection of a certain partition component $\ca_m$. \blue{\bf Does there exists a canonical bijection?} 
\mdn \blue{\bf Maybe the $\cd'$-components are given by the simple objects of $\cd$ in the non-degenerate case! Indeed, note that
\beq
\hm(\dxi x, \dxj y)=\hm(x\dxi, \dxj y)=\hm(\dxi, \dxj yx^*)=\hm(\dxj^*\dxi, yx^*)
\eeq
but $\cd\cap \cd'=\mtr{Vect}.$}

Take an example, $D(G)$ then see what is $f_Q(\cd' \dxi)??$.
\mdn Suppose that $G_i:=\sum_{j \in \ca_i}F_i$ corresponds to the $\cd'$-component $\cd'X_{\phi(i)}$. Then

\beq
\irr(\cd' \xphii)=\{\ch_j\;|\;j\in \ca_i\}
\eeq
\subsection{Example $\cd=\cc'$}
One has that 
\beq
G_0=\sum_{j\in \mtca_0}F_j
\eeq
It follows that $\irr(\cc')=\{\ch_j\;|\; j\in \mtca_0\}$.
\mdn\blue{\bf Check that the two partitions, one from the formula for $\fq$ and one from $\cd=\cc'\subseteq \cc$ coincide!}

\newpage

\newpage
\section{Coends}

\subsection{{\bf Definition of the coend  $Z$ called central Hopf monad}} Let $\cc$ be a finite tensor category. Define the functor $ Z:\cc \ra \cc$ given on objects by
\beq\label{sh-3.3}
 Z(V)=\int^{X\in \cc}X^*\ot V\ot X.
\eeq
\blue{\bf This is the one also used in FRG in the braided case.}
Then 
$$Z=\text{coend} =\text{monad}\; \implies Z-\text{modules in $\cc$  gives the center} \;\czcc $$
\subsection{On the equivalence with the center} One has an equivalence
\dbd
_{Z}\cc\simeq \cz(\cc)
\dbd
defined using the  dinatural transformation $i_{V, X}$ as below. \mdn
Under the above identification the forgetful functor $F:\cz(\cc)\ra \cc$ has the \blue{left} adjoint 
the free module functor.
\subsection{On the central cocommutative colagebra $B_\cc$}

 Define also $$B_\cc=Z(\unu)=\int^{X\in \cc}X^*\ot X.$$ Then  $$B_\cc=Z(1)=\text{cocommutative  coalgebra in}\; \czcc$$ $$=\text{coadjoint coalgebra}$$
\blue{\bf $$Q1:  B_\cc-\text{comodules in\;}\czcc\simeq \cc??? $$}
\subsection{On the central Hopf monad and nat transf}

\mdn
As explained in {\bf Shmizu-14.02-unimodular finite tensor categoryories, page 12, after Eq 3.4} on has that
\beq\label{F2}
\hm_\cc(Z(V), W)\simeq \nat(V\ot (-), (-)\ot W)
\eeq
\blue{\bf TD2: Relate these with centers and relative centers!}
\mdn
In particular for $V=1$ one has that
\beq
\hm_\cc(Z(\unu), W)\simeq \nat(\id_{\cc}, \id_\cc\ot W)
\eeq
\br Following FGR one has that
$$
\hm_\cc(V, E_{\cc})\simeq \nat(\id_{\cc}, \id_\cc\ot V),$$ it agrees with the next subsection!!
\er
\subsection{Proof of the equivalence via Coends and the result of Day and Street}

One defines $Z:\cc \ra \cc$ by
$
Z(V)=\int^{X\in \cc}X^*\ot V\ot X
$. One has an equivalence
$
_{Z}\cc\simeq \cz(\cc)
$
defined using the  dinatural transformation $i_{V, X}$ as below.
Under the above identification the forgetful functor $F:\cz(\cc)\ra \cc$ has as \blue{left} adjoint the free module functor.
\mdn {\bf On the proof of the equivalence:}
  \mdn\blue{ It seems that $Z$ always have values in $\czcc$. Think of it as $FL$.Yes and it transforms morphisms from $\cc$ to $\czcc$. }
  
It is explained in Shmizu-char algebra after defn 3.1 that $\;_Z\cc\ra \czcc$ is given by
\beqn
(M,a)\mapsto (M, \sigma)
\eeqn
where
\beqn
\sg_X:M\ot X\xra{\delta_{M, X}} X\ot Z(M)\xra{\\id_X\ot a} X\ot M
\eeqn
The map $\delta_{M,X}$ is given by the image of $i_{M,X}$ under the isomorphism
\beqn
\hm_\cc(X^*\ot M\ot X, Z(M))\cong \hm_\cc(M\ot X, X\ot Z(M)).
\eeqn

\subsection{Hopf algebra case}

\mdn The universal dinatural transformation of $Z$ is given by
 $$M^*\ot V \ot M\xra{i_{V, M}} A^*\ot V,\;\; f\ot v \ot m\mapsto f(-m)\ot v$$
\mdn

The right adjoint of $F$ is $R=L^{!}$
\mdn
\subsection{Example Hopf algebras}

\subsubsection{The coend From FRZ for Hopf algebras} Here $\cc=\rep(A)$ for a Hopf algebra $A$.

The coend structure as an $A$-modules is given in Runkel's paper. It has the coadjoint action:
$$B_\cc=A^*=Z(\unu)^*=\int^{X\in \cc}X^*\ot X$$ with coadjoint given by
$$[a.f](x)=f(Sa_1xa_2)$$
\subsubsection{Dinatural maps}

The universal dinatural maps are given by
$$i_M:M^*\ot M\ra A^*,\; f\ot m\mapsto (a \mapsto f(am))$$
\mdn More generally $Z(V)$ for $V\in \cc$ is given by
$ Z(V)=\int^{X\in \cc}X^*\ot V \ot X=A^*\ot V$
 with the action
 $$
 a(f \ot v)=f(Sa_1?a_3)\ot a_2v
 $$
 \mdn{\blue{\bf
 \bne
 \item probably the central structure $Z(V)\in \czcc$ is given by the $A^*$-action 
 $g.f=gf$.
 \item It seems that the multiplication is not morphism of $A$-modules.
 \ene
 }}
 \subsubsection{Universal property of this coend}
 
The unique morphism in the universal property of the coend  is given by $g:B_\cc\xra{\id \ot \eta} B_\cc\ot A\xra{\phi_{_AA}} B$
\subsection{To do in this section}
\bne
\item Write the dinatural $i_{X, V}$ in general what property of universality satisfy.
\item FRZ versus FRG establish the paper and write the bibliography.
\ene
\newpage
\section{{\bf Definition of the end called central Hopf comanad}}
If $\barz (V)=\int_{X\in \cc} X^*\ot V \ot X$ then 
$$\barz =end=\text{comonad}\; \implies \blue{\barz -\text{comodules gives the center} \;\czcc\;}$$$$E_{\cc}=\barz (1)=\text{comm algebra in}\; \czcc$$ $$=\text{adjoint algebra}$$
$$E_{\cc}-\text{modules in}\; \czcc \simeq \cc$$
\mdn define delta and epsilon \mdn The map 
$$\delta: \barz \ra \bar  Z^2$$ \blue{\bf probably is} a morphism in the center $\czcc$. I verifed for Hopf algebras it is given by
$$A\ot V\to A\ot (A\ot V),\;\; a\ot v\mapsto a_1\ot a_2\ot v$$
This is used in defining the canonical isomorphism 
\beqn
\cfcc \xra{\phi} \enx_{\czcc}(\runu)
\eeqn
\mdn Since $E:=\barz =FR$ one may also apply the formulae of Virelizier and Turaev to see this. \mdn\blue{\bf TD1: Try to generalize their results corresponding to induction from the relative center!}
\mdn
For any finite tensor category$\cc$ define
\beqn
E_{\cc}:=\barz (\unu)=\int_{X\in \cc}X^*\ot X
\eeqn
Then $E_{\cc}=UR(\unu)$ is an etale algebra in $\czcc$ and one has an equivalence of $\cz(\cc)$-modules categoryories
\beq
\cz(\cc)_{E_{\cc}}\xra{\simeq} \cc
\eeq
whose inverse is given by $V\mapsto R(V)$.
One has that $B_\cc=E_{\cc}^*$. This is a coalgebra in $\cc$.
\mdn \br In the Hopf algebra example the two objects $Z(\unu), \barz (\unu)$ have the same underlying vector space. That is why one can define both harpoons in this case. In the general case the second harpoon is defined in terms of the $B_\cc=E_{\cc}^*$.\er
\subsubsection{Another natural transformation-for coends}
Define 
\beqn
\barz (V)=\int_{X\in \cc}X^*\ot V\ot X.
\eeqn
As explained in {\bf Shmizu-14.02-char algebra Eq 3.7} on has that
\beq\label{F1}
\hm_\cc(V, \barz (W))\simeq \nat(V\ot (-), (-)\ot W)
\eeq
\blue{\bf From the two sections it seems that $Z$ and $\barz $ are pair of adjoint functors.}\subsection{Relation between central ends and coends}

Explained on {\bf Shimizu-Ch alg-page 8} one has that $$\barz =UR\simeq Z^!$$ May assume even $\barz = Z^!$.

\blue{\bf Q3 What it works only with commutative character rings? Probably the fusion categoryresults.} 
\subsection{Unimodular tensor categoryories}
Let $U:Z(\C)\ra \C$ be the forgetful functor and let  $R$ be its right adjoint.
\mdn
There is a distinguished invertible object $D\in \C$ introduced by Etingof and Ostrik and rewritten by Shimizu.

\bn{defn}
A fusion categoryis called unimodular if this element is the identity.
\end{defn}
\noindent
\blue{\bf  define $\eps_{\unu}, \; \delta_{\unu}$.}
\subsection{Central elements}
By definition
\beq 
\cecc:= \hm_{\C}(\unu,A) 
\eeq
is called \blue{the set of central elements.}
\subsubsection{Description of $\mtr{End}(\id_\cc)$} One has that $T\in $ if for any morphism $f:X\ra Y$ one has that the following diagram
\[
\begin{tikzcd}
X \arrow{r}{f} \arrow[swap]{d}{T_X} & Y \arrow{d}{T_Y} \\
X \arrow{r}{f} & Y
\end{tikzcd}
\]
is commutative.

\subsubsection{On the canonical bijection}

\blue{\bf When this bijection is used??:}
There is a bijection:
\beqn
\cce\xra{} \enx(\id_{\C})
\eeqn

The inverse of it is defined in the book KL-with corners, via ends or coends.
\subsubsection{Multiplication of central elements} Let $a, b \in \cecc$ then
\beq
ab:=(\unu=\unu\ot \unu\xra{a\ot b}Z(\unu)\ot Z(\unu) \xra{m} Z(\unu))
\eeq

\section{The monoidal center and the character algebra}

\subsection{Coends and the result of Day and Street}
One has an equivalence of categories
\dbd
_{Z}\cc\simeq \cz(\cc)
\dbd
dinatural transformation $i_{V, X}$.
\mdn
under the above identification the forgetful functor $U:\cz(\cc)\ra \cc$ has as \blue{left} adjoint the free module functor.
\mdn
the right adjoint is $L^{!}$. \green{\bf Shimizu has given a relation between left and right adjoints of a tensor functor via left and right duality.}
\subsection{Unimodular tensor categories}
Let $U:Z(\C)\ra \C$ be the forgetful functor and let  $R$ be its right adjoint.
\mdn
There is a distinguished invertible object $D\in \C$ introduced by Etingof and Ostrik and rewritten by Shimizu.

\bn{defn}
A fusion category is called unimodular if this element is the identity.
\end{defn}
\subsection{Canonical action of $Z(\unu)$ on every $X$ object of $\cc$.}
$$Z(\unu) \ot X\xra{\pi_{\unu, X}\ot \id_X}X\ot X^*\ot X\xra{\id_X\ot \ev_X} X\ot \unu=X.$$
\noindent
\blue{\bf Define $A=U(1)$ as Hopf algebra in $\cc$??    Define $\eps_{1}, \; \delta_{1}$!!!}
\subsection{Central elements}

\mdn
Definition\beq 
\cecc := \hm_{\C}(1,A) 
\eeq
is called \blue{the set of central elements.}

\subsubsection{The isomorphism with natural transformations}

There is a bijection:
\beqn
\cce\xra{\psi} \enx(\id_{\C}),\;\;\psi(a)_X=\ro_X(a\ot \\id_X),\;\; X=\unu\ot X\xra{a\ot \\id_X}X\ot Z(\unu)\xra{\ro_X} X
\eeqn
\blue{\bf describe this bijection!!}\md \blue{\bf Where in the paper this bijection is used:??}

\subsection{Integral theory}
 An \blue{\bf integral} in $\C$ is a morphism $\Lambda: 1 \ra A$   in $  \C$   such that
 \beqn
 m \circ (\id_{A}\otimes \Lambda)=\eps_{1}\ot \Lam
 \eeqn
   
A \blue{\bf cointegral in $  \C $}  is a morphism $  \lambda : A\ra 1 $  such that
\beqn
\barz (\lam)\circ \delta_{1}=u \otimes \lam
\eeqn

\subsection{Integrals and Fourier transform}
Let $  \Lambda $\; be a non-zero cointegral in $  \C $  and $\lambda$ be the integral such that 
\beqn
\langle \Lambda, \lambda\rangle =1
\eeqn
\noindent
The Fourier transform associated to $\lambda$ is the linear map
\beqn
\mtc F_{\lambda}:\cecc\ra \cfcc
\eeqn
given by
$$
a \mapsto \lambda \lh S(a)
$$
where $S: \cecc\ra \cecc$ is the antipodal operator defined by:

and $\lh$ id the action of $\cecc$ on $\cfcc$ given by 
\beqn
f \lh b=f \circ m \circ (b\ot \id_{A})
\eeqn
for $f \in \cfcc$ and $b \in \cecc$.

\blue{\bf this is the harpoon $(f\lh b)(a)=f(ba)$. }
\subsection{Radford's trace formula}
\subsubsection{The Drinfeld isomorphism}
For any central object $(X, \sg)\in \cz(\cc)$ one defines  the canonical isomorphism:
\beqn
\psi_{X}:X\xra{} X\ot X^{*}\ot X^{**}\xra{\sg, \id} X\ot X\ot X^{**}\xra{\ev \ot id} X^{**}
\eeqn

\blue{By [Shim, sect 7] $\psi_{A}=j_{A}$ if $\cc$ has a pivotal structure $j$.}
\subsubsection{Definition of the trace by abuse of notation-using Drinfeld isom}
For any $\xi:A\ra A$ define
\beqn
tr(\xi):1\xra{\coev}A\ot A^{*}\xra{\psi_{A}\ot \id} A^{**}\ot A \xra{\ev}1
\eeqn
\subsubsection{Radord's formula}

\bt Let $\cc$ be a unimodular finite tens category. With the above notations one has:
\beqn
tr(\tilde{f})=\langle f, \Lam\rangle \langle \lam, u\rangle 
\eeqn
\et
\newpage

\subsection{Character algebra}
The internal character $\mtr{ch}(X)$ is defined as the class function
\beqn
\mtr{ch}(X):=tr^{X}_{A, \unu}(\rho_{X}).
\eeqn
\bt
The set of irreducible characters is linearly independent in $\cfcc$.
\et
\blue{\bf Define multiplication on $f, g \in \cfcc$
\beq
f \star g:=(Z(\unu)\xra{\delta_\unu}\barz ^2(\unu)\xra{\barz (g)}Z(\unu) \xra{f} \unu)
\eeq}
\bt The adjunction gives 
\beq
\mathrm{\bf{CF}}(\mtc C) \xra{\cong} \mtr{End}_{\Z(\mtc C)}(R(\unu)),\;\; \ch \mapsto \barz (\ch)\circ \delta_{\unu}
\eeq
an isomorphism of monoids. \blue{\bf Checking the maps - in the Hopf algebra case this becomes $\ch\mapsto \ch\rh$.}
\et
\blue{\bf In the braided case there is also a coproduct on $A$ which makes the product of characters.}

\subsection{Conjugacy classes for fusion categoryories}

A conjugacy class is defined as a direct summand of $R(\unu)$. Thus let $\mtc C_{0}\cdots, \C_{m}$ be the
conjugacy classes of $\C$. 
 \mdn
 Since the unit object $\unu_{\czcc }$ is always a subobject of $R(\unu)$, we can assume $\C_{0} = \unu_{\czcc }$.

\bt (Theorem 6.6) For a non-degenerate pivotal fusion categoryC, the following assertions are equivalent:

(1) $R(\unu) \in \czcc $ is multiplicity-free.

(2) The algebra $Gr_{\kk}(C)$ is commutative.
\et

\subsection{Radford's trace formula}
\subsubsection{The Drinfeld isomorphism in the Drinfeld center}
For any central object $(X, \sg)\in \cz(\cc)$ one defines  the canonical isomorphism:
\beqn
\psi_{X}:X\xra{} X\ot X^{*}\ot X^{**}\xra{\sg\ot  \id} X^*\ot X\ot X^{**}\xra{\ev \ot id} X^{**}
\eeqn

\blue{By [Shim, sect 7] $\psi_{A}=j_{A}$ if $\cc$ has a pivotal structure $j$.}
\subsubsection{Definition of the trace by abuse of notation-using Drinfeld isom}
For any $\xi:A\ra A$ define
\beqn
\tr(\xi):\unu\xra{\coev}A\ot A^{*}\xra{\psi_{A}\ot \id} A^{**}\ot A^*\xra{\ev}\unu
\eeqn
\subsubsection{Radord's formula}

\bt Let $\cc$ be a unimodular finite tens categoryory. With the above notations one has:
\beqn
\tr(\tilde{f})=\langle f, \Lam\rangle \langle \lam, u\rangle 
\eeqn
\blue{\bf where $f: A\ra A$??}
\et

\subsection{On harpoons and Fourier transforms} 
Define a right action of $\cecc$ on $\cfcc$ by
$$
f \lhau b=(A=\unu\ot A\xra{b\ot \id_A}A\ot A\xra{m}A\xra{f}\unu)
$$
for any $b \in \cecc$ and $f \in \cfcc$.
\mdn
There is also an action of $\cfcc$ on $\cecc$ given by
$$b\lhau f= (\unu\xra{b} A\xra{\blue{\Delta}}A\ot A\xra{f\ot \id_A}\unu\ot A=A)$$
{\bf \blue($\Delta$) exists only in the braided case on $A$.}
\subsubsection{Integrals and Fourier transform}
Let $  \Lambda $\; be a non-zero cointegral in $  \C $  and $\lambda$ be the integral such that 
$
\langle \Lambda, \lambda\rangle =1
$.
\mdn
The Fourier transform associated to $\lambda$ is the linear map
\beqn
\mtc F_{\lambda}:\cecc\ra \cfcc
\eeqn
given by
$
a \mapsto \lambda \lh S(a)
$
where $S:\cecc\ra\cecc$ is the antipodal operator defined by:

and $\lh$ id the action of $\cecc$ on $\cfcc$ given by 
\beqn
f \lh b=f \circ m \circ (b\ot \id_{A})
\eeqn
for $f \in \cfcc$ and $b \in\cecc$.

\blue{\bf this is the harpoon $(f\lh b)(a)=f(ba)$. }
\blue{\bf generalization to cw- result
$$(A//L)^{*}=\blam_{L}\rh A$$
at the level of the characters.
}
{\bf The inverse of the Fourier transform is given by the use of the morphisms $d_\lam$
 by the formula
 \beqn
 \mtc F'_\lam(f)=(\id_A\ot f)\circ d_\lam=(\unu\xra{d_
 \lam}A\ot A\xra{\id_A\ot f}A\ot \unu=A)
 \eeqn A formula 
 \beqn
 d_\lam=(\unu \xra{\mtr{coev}_\bfa}\bfa \ot \bfa^*
\xra{\lam\ot \id_{\bfa^*}} \bfa \ot \bfa)
\eeqn}
In the case of Hopf algebras this map is described in Sh-char algebra Equation 5.18 by
$$
d_\lam(1)=\blam_2\ot S^{-1}(\blam_1)
$$
and the inverse of $\mtc F'_\lam$ becomes $f\mapsto f(S^{-1}\blam_1)\blam_2$.
\mdn
Define the conjugacy class sum $\csu_j\in \cecc$ by
$$\csu_j:=\mtc F_\lam(F_j)$$
where $\mtf(\ch)=\blam\lhau \ch$ is the Fourier transform.
\mdn \blue{\bf Note that in the commutative case $\csu_j$ form a basis for $\cecc$.}

\newpage
\newpage
\section{On the harpoon relations and other relative generalizations} 
\mdn Denote by $\{V_0, V_1, \dots, V_m\}$ the simple objects of $\cc$ and suppose that \\$\{V_0, V_1, \dots, V_{m'}\}$ with $m'\leq m$ are the simple objects of $\cd$.
\mdn
\blue{\bf Need an interpretation of the sums of equivalence relations wrt to the restriction to normal $L$.}

\subsection{General settings for fusion subcategoryories-using DMNO 4.10-thm}Let $\cd$ be a fusion subcategoryof $\cc$. Then 
$L_\cd:=I_\cd(\unu)$ is an object of the relative center $\cz_\cd(\cc)$ corresponding to $\cd$ by Theorem 4.10 DMNO.
\mdn
Let $\cd$ be a fusion subcategoryof $\cc$. Define $L$ to be the etale subalgebra of $A=UR(1)$ corresponding to $\cd$ via Thm 4.10, DMNO. It is a connected etale algebra in $\cz(\cc)$.
\mdn More precisely $\cd$ coincides to those $A$-modules that are dyslectic with respect to $L$ in $\cz(\cc)$. \mdn \blue{\bf In the Hopf algebra setting this coincides to those $A$-modules that are trivial as $L$-modules}
\subsection{On the canonical epimorphism between ends}
There is a canonical epimorphism
$$\pi:E_{\cc} \ra A_\cd$$
{\bf which is an epimorphism of algebras in $\cz(\cc)$ \blue{\bf $A_\cd$ is not in that center}. I don't know if $A_\cd$ is an object in $\czcc$ or in the relative center.} 
\mdn\blue{\bf write a reference for this}
\mdn \blue{Needs explnd:}
Using the above connection with natural transformations this map corresponds to the usual restriction:
\beqn
\nat(\id_{\cc}, \id_{\cc} \ot W)\ra \nat(\id_{\cd}, \id_{\cd}\ot W)
\eeqn
\br Following FGR one has that
$$
\hm_\cc(V, E_{\cc})\simeq \nat(\id_{\cc}, \id_\cc\ot V)$$
\er
\subsection{On the inclusion $\cfd\hookrightarrow \cfc$}One has that the above
$$\pi: E_{\cc}\xra{\pi} A_\cd$$ is an epimorphism of algebras. This gives a monomorphism of algebras under  convolution
$$j={\pi^*}: \cfd\hookrightarrow\cfc$$ given by composition with $\pi$. \mdn  Recall that $\cfd=\hm_\cd(A_\cd, \unu_\cd)\hookrightarrow \cfc= \hm_\cc(E_{\cc}, \unu_\cd)$.\mdn {\bf I don't know  how to prove apriroi that this map is an injection.}
\subsection{On the epimorphism $\cecc\xra {\delta}\cecd$} There is also a surjective map on the centres
$$\delta:\cecc\ra \cecd.$$
Note that this map corresponds to the restriction
$$
\mtr{End}(\id_\cc)\ra \mtr{End}(\id_\cd)
$$
{\bf I don't know  how to prove apriroi that this map is a surjection.}
\subsection{On the integrals and cointegrals}
The integral $\blam\in \hm(\unu, E_{\cc})$ is defined as an $E_{\cc}$-module map in $\cc$. \blue{\bf Note that $\unu$ is an $E_{\cc}$ module via $\eps_\unu$??}

\subsection{On $\hm_\cc(\unu, L)$ interpreted as $Z(A)\cap L$}
Define $$\bf{CE}(L):=\hm_\cc(\unu, L)$$
\blue{\bf It corresponds to $Z(H)\cap L$ if $\cc=\rep(H)$.} Denote by $\blam_L$ {\bf the unique } morphism (up to a scalar) corresponding to it. It is unique since $L$ is etale. Is it an isomorphism of $L$-modules??
\mdn\blue{There is a notion of integral for module categoryories in the last paper. See if it coincides to $\blam_L$}
\subsection{On the Drinfeld maps in the braided case}

\beqn
\mtcd_\cc:\cfcc \ra \cecc, \;\;\ch_i\mapsto \sum_j\frac{s_{ij}}{d_j }E_j
\eeqn
an algebra map.

The Drinfeld map of $\cd$ is a restrction of this map, by construction!! Or check with nondegener pairing and universal property via $\pi$. \blue{See the definition from Bakalov! In the Hopf algebra this might the case since $\apl$ is \qtr with $(\pi \ot \pi)(R)$. Shimizu explained that one needs ribbon braided category to have this identity!!!, based on Bakalov Theorem 3.1.11 or Prop 3.1.12}
\beqn
\mtcd_\cc:\cfcc \ra \cecc, \;\;\ch_i\mapsto \sum_j\frac{s_{ij}}{d_j }E_j
\eeqn
Denote $\mtcd_\cc(F_i)=\sum_{j\in \mtca_i}E_j$.
\mdn  \blue{\bf Form here follows that the centralizer comes in chunks. Let $\ch \in \cd$. Write $$\mu=\sum_i\mu_iG_i=$$ It follows that $\langle \mu\rangle '=$.}
\mdn Write
\beqn
G_i=\sum_{j \in \ccb_i}F_i 
\eeqn
In particular write 
\beqn
\lamcd=G_0=\sum_{j \in \ccb_0}F_i 
\eeqn

From monatshefte $$\lker_\cc(\fpdimcd\lamcd)=\mtca_0.$$
\subsection{Generalization of the first criterion} $\ce=\repapm$ and $\cd=\repapl$. Then $$\ce\subseteq \cd' \iff\phi_R((A//L)^*)\subseteq M \iff \phi_R(\laml)\blam_M=\blam_M$$
It implies that $$\repapm \subseteq \repapl \iff \repapm\subseteq \rep(A//\lstar)'$$ Thus we have equality.
\bpf If $(1-\lam_L)\lam_M=0$ then $(1-\lam_L)=(1-\lam_L)(1-\lam_M)$ and then $AL^+\subseteq AM^+$, ie $L\subseteq M$.
\epf
\blue{\bf This shows that the first criterion holds for q-tr!!}
\noindent {\bf Conjecture generalization of the first criterion:}
$$\ce\subseteq \cd' \iff \mtcd_\cd(\lam_\cd)\blam_\ce=\blam_\ce$$ For this it is enough to show that
\bp
\beq
\mtcd_\cd(\lam_\cd)=\blam_{\cd'}
\eeq
in the factorizable case. In the other cases $\mtcd(\lam_L)=\blam_{\lstar}$.
\ep
\subsection{On the conjugacy classes in $\cc$} Suppose that the char ring is commutative.  Write  $$R(\unu)=\oplus_{j=0}^r\cc_j$$ for the decomposition of $R(\unu)$ in sum of simple objects in $\cz(\cc)$.
\mdn Denote by $E_j$ the primitive central idempotent of $\mtr{End}_{\czcc}(R(\unu))$ corresponding to the simple object $V_j=\cc_j$.
\mdn
Denote by $F_j$ the primitive central idempotents of $ \cfc$ corresponding to $E_j$ via the natural isomorphism
$\mtr{End}_{\czcc}(R(\unu))\simeq\cfcc.$\mdn 
Define the conjugacy class sum $\csu_j\in \cecc$ by
$$\csu_j:=\mtc F_\lam(F_j)$$
where \blue{\bf $\mtf(\ch)=\blam\lhau \ch$??} is the Fourier transform.
\mdn \blue{\bf Note that in the commutative case $\csu_j$ form a basis for $\cecc$.}

\subsection{The pairing between $\cfcc$ and $\cecc$}
To formulate an analogous fact
in our setting, we consider a paring
\beqn
\langle f , a \rangle : CF(\C) \times CE(\C)\ra L, \\ \langle f, a\rangle  \id_{\unu}= f \circ a.
 \eeqn
\blue{\bf This is nondegenrate by Shimizu-chars, read the result and complete it. }
\bl
To every conjugacy class sum corresponds
$p_{V_i}$ dual basis to this pairing. Show that they are the idempotents of $\cfcc$.
\el
\bpf
\blue{Is it explained in char shimizu?? Yes a multiple of the idempotent form dual bases. I have proven this by compatibility of multiplication.}
\epf
\subsection{On the decomposition of $L$ and  $\blam_L$}
Following the isomorphism 
$$\mtr{End}_{\cz(\cc)}(R(\kk))\simeq \cfc$$ one has  that

$$L=\bigoplus_{j\in {\mtc I_L}}\mathfrak C^j.$$ 

Since the class sums form a basis one has that
\beqn
\blam_L=\sum_{j\in J_L}\al_j\csu_j
\eeqn
with $\al_j \neq 0$. \blue{\bf Here use that $$\hm_\cc(\unu, L)\hookrightarrow \hm_\cc(\unu, E_{\cc})$$
}
\mdn \blue{\bf Try to prove the decomposition results from my paper in this general case!,i.e $J_L=I_L$} 
\mdn
Using the pairing one has that $\langle \ch, \blam_L\rangle \neq 0 $ if and only if $\ch \in \rep(A//L)$.
\subsection{Example: Hopf algebras: Interchange for the Drinfeld map on the character ring}
$$\phi_R(f)=(\id \ot f)(R_{21}R)$$
\mdn Class Equation gives in thus case that
$$\frac{1}{\tau_j}=\frac{d}{\dim(f_jH^*)}=\frac{d\dim(F_j)}{\dim(F_jH^*)}$$
\subsection{On the integrals and cointegrals}
\subsubsection{Integral theory}
 An \blue{\bf integral} in $\C$ is a morphism $\Lambda: \unu \ra A$   in $  \C$   such that
 \beqn
 m \circ (\id_{A}\otimes \Lambda)=\eps_{\unu}\ot \Lam
 \eeqn
   
A \blue{\bf cointegral in $  \C $}  is a morphism $  \lambda : A\ra \unu $  such that
\beqn
\barz (\lam)\circ \delta_{\unu}=u \otimes \lam
\eeqn
\noindent \blue{\bf In the ss case there is another characterization of cointegrals by multiplication as in the usual case.}
\mdn
Recall {\bf Shmizu-char of algebras} that $\lam_{\cc}$ is the unique morphisms in $\czcc$
\beqn
\lam_\cc:E_{\cc}\ra \unu
\eeqn
\bl
\bne
\item
 If $\cd\subseteq \cc$ then the composition
$$E_{\cc}\xra{\pi} A_\cd\xra{\lam_\cd} \unu$$
is a cointegral of $\cc$.
\item
One has that the composition 
$$\unu \xra{\blam_\cc} E_{\cc} \xra{\pi} A_\cd$$
is an integral in $\cd$.
\ene
\el

\subsection{On the left dual of $E_{\cc}$.}
It is shown in {\bf Shi 16} that ${\bf A}$ is a Frobenius algebra in $\czcc$ with trace $\lam$ in the sense that
$$
\phi_\lam:\bfa\xra{}\bfa\bfa\bfa^*\xra{mA^*}\bfa\bfa^*\xra{\lam A^*}\bfa^*
$$
is an isomorphism in $\cz(\cc)$.
\mdn

\subsection{Study also:}the space of maps $\hm_\cc(B_\cc, \unu)$ and the other one $\hm_\cc(\unu, B_\cc)$.
This is an object with the coadjoint action on $A^*$. $a.f(b)=f(S^{-1}a_2ba_1)$. This gives exactly the class functions if $S^2=\id$.
\bl With the above notations one has the following commutative diagram:
\[
\begin{tikzcd}
\cfd \arrow{r}{j} \arrow[swap]{d}{\mtc F_\cd} & \cfc \arrow{d}{\mtc F_\cc} \\
\cecd \arrow{r}{\delta} & \cecc
\end{tikzcd}
\]
\el
\bpf
\epf
\mdn
\blue{\bf The commutative above diagram should give a similar decomposition in the centers. Imitate the proof from ART for normal Hopf subalgebras.}
\mdn
Write (read section 6 on applications to fusion categoryories-shimizu ch algebra)
\beqn
\cecc=\kk E_0\times \kk e_1 \times \dots \times\kk e_m
\eeqn
and
\beqn
\cecd=\kk g_0\times \kk g_1 \times \dots \times\kk g_{m'}
\eeqn
By Lemma 6.1 one has that $E_0$ and $g_0$ are the integrals for $\cc$ and respectively $\cd$.
\mdn
\subsection{Fourier transfrom formulae}
 By Equation 6.10 from Shimizu-char alg one has that $$\mtc F_\lam(E_i)=\frac{\dim(V_i)}{\dim(\cc)}\ch_i^*$$
\mdn
By Theorem 6.9 of {Shimizu-chr alg} one has that
$$f_r=\frac{|\mathfrak C_r|}{\dim(\cc)}\sum_{i=0}^m\langle \ch_{i^*}, g_r\rangle \ch_i
$$
\mdn It seems that this implies that
$$
\mtc {F_\ch}^{-1}(f_r)=\csu_r
$$
In fact this is the way it is defined!!!
\mdn Equation 6.13 implies that
$$\ch_i=\sum_{r=0}^m\langle \ch_i, g_r\rangle f_r$$
\mdn
Write \beqn
\cfcc\simeq \kk F_0\times \kk F_1\times \dots \times \kk F_m
\eeqn
and 
\beqn
\cfcd\simeq \kk G_0\times \kk G_1\times \dots \times \kk G_m'
\eeqn
One has $F_0=\lam_\cc$ and $G_0=\lam_\cd$, by {\bf Shmizu -char algebra lemma ??}. Suppose that
$$j(G_0)=\sum_{i\in \mtc H_\cd} F_i$$
\bl
One has that a class function $\lam \in \cfcc$ is am integral if and only if 
\beqn
f\star \lam=\langle f,u\rangle \lam
\eeqn
for any functional $f \in \cfcc$.
\el 
\mdn
This implies that the regular char 
$$\Theta(d):=\sum_{i=0}^m\dim(V_i^*)\ch_i$$ is a cointegral for $\cc$.
\mdn
\blue{\bf Show that $$\blam_L\rh \cfcc=\cfcd$$ and $$\lam_\cd\circ \pi=\blam_L\rh \lam_\cc$$}
\mdn
Using theorem 4.10 DMNO the category$\cd$ can be recovered as those $A$-modules in $\czcc$ that are dyslectic wrt $L$, but also as the induction of those $L$-modules that are dyslecti wrt $L$.
\newpage

\section{A comment on Shimizu's formula}

\subsection{More relations with Fourier transform from Simizu}Below we write the main formulae from \cite{shimizu-char-alg.}
\blue{\Small  \br The philosophy is that he does not define the harppons of $\cfcc$ on $\cecc$, instead he uses the inverse of the Fourier map via $d_\lam$.
\er}
\beq\label{Eq (6.1)}
\langle \chii,\; E_j\rangle =d_i \delta_{i,j}
\eeq
\beq\label{Eq (6.2)}
E_iE_j=\delta_{i,j}E_i,\;\;\; S(E_i)=E_{i^*}
\eeq
\beq\label{Lemma 6.1}
E_0=\blam
\eeq 
\beq\label{Eq (6.3)}
a\blam =\eps(a)\blam,\; \;\text{for all a} \in \cecc
\eeq
where $\eps:\cecc\ra \kk$ given by $a\mapsto \eps_\unu(a)$ is a morphism of algebras.
\bl\label{Lemma 6.2}
$\mathrm{Gr}_\kk(\cc)$ is a symmetric Frobenius algebra with the trace $\tau$ given by 
\beq
[X]\mapsto \dimk\hm_\cc(\unu, X)
\eeq
\el
It follows that
\beq\label{Eq (6.4)}
\tau(\ch_i)=\delta_{i,0}=\langle \ch_i, E_0\rangle  \;\;\tau(\chii\star \chj)=\delta_{i,j}
\eeq
This shows that 
\beq\label{recoverform}
\tau(\ch)=\langle \ch, E_0\rangle 
\eeq since $\chii$ form a basis on $\cfcc$.
\mdn
Thus for
\beqn
E:=\sumitom \chi\ot \chistar
\eeqn
one has that
\beq\label{Eq (6.5)}
\Theta :\hm_\kk(\cfcc, \;\kk)\ra \cfcc,\;\;\ch\mapsto (\ch \ot \id)(E)
\eeq
is bijective.
\mdn By general Frobenius theory if $\al:\cfcc\ra \kk$ is an algebra map then
\beq\label{Eq (6.6)}
f\star \lam_{\al}=\al(f)\lam_\al, \text{for all f}\in \cfc
\eeq
where 
\beqn
\lam_\al:=\Theta(\al)
\eeqn
Moreoveor 
\beqn
\Theta(\al)\in \mtc Z(\cfcc)
\eeqn
since $\tau$ is symmetric.
\mdn Consider the map
\beq
d:\cfcc\ra \kk,\;\; f\mapsto \langle f,\;u\rangle 
\eeq
sending $\ch([X])$ into the  pivotal dimension of $X$. 
It is to see that is an algebra map.

 \green{\bf \small Compute the kernel of this morphism.}
\bl\label{Lemma 6.3}
A class function $\lam \in \cfcc$ is a cointegral iff
\beqn
f\star\lam=\langle f, \;u\rangle \lam, \;\;\text{for all f}\in \cfcc
\eeqn
Moreover, by the above argument 
\beq\label{reg char}
\theta(d)=\dimcc({V_i}^*)\ch_i\in \cfcc,\;\;\text{is a cointegral in}\;\cc
\eeq
and it belongs to the center of $\cfcc$.
\el
\noindent \blue{Read the proof of it!!}
Define \beq\label{Eq (6.9)}
g_r:=|\cc_r|\csu_r
\eeq
divided by the pivotal dimension, it is non-zero since $\cc_r$ is a simple object in the modular tens categ $\cz(\cc)$.
\subsection{Compatibility with multiplicity}Maybe in the braided case where I have braided Hopf algebra structures on coend and end. Are they self dual via the pairing?
\beq
\mu(a\lh\ch)=\ch\mu(a)
\eeq

\blue{\bf Assuming compatibility with multiplicity:}
\beq
\langle F_i, \csu_j\rangle =\langle F_i, \dimcc\blam\lh F_j\rangle =\delta_{i,j}\dimcc F_i(\blam)
\eeq
{\bf Prove that $$F_i(\blam):=\frac{1}{n_i}\neq 0.$$ It is true since the form $\tau$ is nondegenrate on a semisimple $\kk$-algebra.!!\\Interpret it using a Class Equation argument via Radford's formula for trace. }
\mdn From here obtain dual basis for evaluation
\beq
\{F_i, \frac{n_i}{\dimcc}\csu_i\}.
\eeq
\noindent Another pair of dual bases is 
\beq
\{\ch_i, \frac{1}{d_i }E_i\}
\eeq
\noindent This implies that
\beq
\sumitom F_i\ot \frac{n_i}{\dimcc}\csu_i=\sumitom \ch_i\ot  \frac{1}{d_i }E_i\in \cfcc\ot \cecc
\eeq
\noindent Applying also $\id \ot \mtf^{-1}$ one obtains
\beq\label{db}
\sumitom F_i\ot \frac{n_i}{\dimcc}F_i=\sumitom \ch_i\ot \frac{1}{d_i }\frac{d_i }{\dimcc}\ch_i^*
\eeq
\noindent One also has
\beq
\sumitom \ch_i\ot \chistar=\sumitom F_i\ot \frac{\dimcc}{f_i}F_i
\eeq
\mdn The two Equations show that
\beq
n_i=\frac{\dimcc}{f_i}
\eeq
\blue{\bf This shows that $1/n_i$ is not zero!}

\mdn
Then all equations form Higman's ideals follow:
\mdn Denote
\beq
\ch_i=\sumjtom \al_{ij}F_j
\eeq

Apply $\mtf^{-1}$ to the above equation and get
\beq
\frac{\dimcc}{d_i }E_{i^*}=\sumjtom \al_{ij}\csu_j
\eeq

\bl
\beq
\alij=\frac{d_i }{d_j }\alji
\eeq
\el
\bpf
Use Shimizu formulae $\alij=\frac{\sij}{d_j }$.
\epf
\bl
\beq
C_i=\frac{d}{n_i}\sumjtom \frac{1}{d_j }\alji E_j
\eeq
\el
\bpf
Replacing $\ch_i$ by the formula $\ch_i=\sumjtom \al_{ij}F_j$ in Equation \eqref{db} it follows that
\beq
\sumitom F_i\ot \frac{n_i}{\dimcc}F_i=\sumitom(\sumjtom \al_{ij}F_j\ot \sumjtom \al_{i^*j'}F_{j'})=\sum_{j, j'}(\sumitom \alij\al_{i^*j'})F_j\ot F_{j'}
\eeq
It follows that
\beq
\sumitom \alij\al_{i^*j'}=\frac{n_i}{\dimcc}\delta_{j, j'}
\eeq
\mdn This implies that $AB^t=\id$ where
\beq
A=(\alij)_{ij}, \;\; B=(\frac{1}{n_i}\al_{i^*,j})_{ij}
\eeq
$B^t_{ij'}=b_{j'i}$
\beq
\delta_{jj'}=\sumitom a_{ji}b_{j'i}=\sumitom \alij
\eeq
\blue{\bf Anyway, this implies that also $B^tA=\id$ and this gives by inversion the desired formula for $C_i$.}
\epf
\bt Let $\cc$ be a braided balanced nondegenrate fusion category. Then
\beq\label{fqch}
f_Q(\ch_i)=\frac{1}{\dvi}\csu_i.
\eeq
\et
\bpf One has that
\beq
f_Q(\ch_i)=\sumjtom \al_{ij} E_j=\sumjtom \frac{d_i }{d_j }\alji E_i=d_i \sumjtom \frac{1}{d_j }\alji E_i=\frac{d_i n_i}{d}\csu_i
\eeq
\epf
\blue{\bf Conjecture: In a factorizable category \beq
n_i=\frac{d}{d_i ^2}
\eeq
}
Anyway, even without this formula it follows that
\beq
\csu_i\csu_j=\frac{d^2}{n_id_i n_jd_j }f_Q(\ch_i\ch_j)=\frac{d^2}{n_id_i n_jd_j }\sum_lN^l_{ij}f_Q(\ch_l)=\eeq
\beqn=\sum_l\frac{d^2}{n_id_i n_jd_j }\frac{n_ld_l}{d}N^l_ij\csu_l=d\sum_l\frac{n_ld_l}{n_id_i n_jd_j }N^l_{ij}\csu_l
\eeqn
\subsection{Formulae from the Fourier transform and its inverse}
Here $\mtf:\cecc \ra \cfcc$  is the Fourier transform given by 
\beq
\mtf(a)=\lam\lhau S^{-1}(a)
\eeq
\mdn
 The inverse of the Fourier transform is given by the use of the morphisms $d_\lam$ by the formula
 \beqn
 \mtc F'_\lam(f)=(\id_A\ot f)\circ d_\lam=(\unu\xra{d_
 \lam}A\ot A\xra{\id_A\ot f}A\ot \unu=A)
 \eeqn A formula for $d_\lam$ is given by
 \beqn
 d_\lam=(\unu \xra{\mtr{coev}_\bfa}\bfa \ot \bfa^*
\xra{\lam\ot \id_{\bfa^*}} \bfa \ot \bfa)
\eeqn
In the case of Hopf algebras this map is described in Sh-char algebra Equation (5.18) by
$$
d_\lam(1)=\blam_2\ot S^{-1}(\blam_1)
$$
and the inverse of $\mtc F'_\lam$ becomes $f\mapsto f(S^{-1}\blam_1)\blam_2=\blam \lh S^{-1}f$.
\mdn
\mdn In the fusion case
\beq
\mtf(E_i)=\frac{d_i }{\dimcc}\chistar, \;\;\mtf^{-1}(\chii)=\frac{\dimcc}{d_i }E_{i^*}
\eeq

\subsection{Definition of the conjugacy classes}

Define the conjugacy class sum $\csu^j\in \cecc$ by
\beq
\csu^i:=\mtc F_\lam^{-1}(F_i)
\eeq
\noindent where $F_i$ is the primitive idempotent corresponding to $\mathfrak C^i$ via the canonical isomorphism $\enx_\cc(\runu)\simeq \cfcc$.  \mdn  Note that in the commutative case $\csu^j$ form a basis for $\cecc$.

\bibliographystyle{amsplain}
\bibliography{24nov}

\blue{\bf What $\rep(A//L)$ represents? It sholud be related with DMNO\\ Conjecture Module categoryories are related with central subalgebra of $A$.}

The notion of centralizer in a 
braided fusion categorywas introduced by M\"uger in \cite{proclond}. It was shown in \cite[Theorem 8.21.4]{EGNO15} that the centralizer of a nondegenerate  fusion subcategoryof a braided categoryis a categoryorical complement of the nondegenerate subcategoryory.  This principle is the  basis of many classification results of braided fusion categoryories, see for example papers  \cite{dgno2,  ENO2, DGNO} and references therein. 
\mdn
Despite its importance, there is no concrete formula for the centralizer of all fusion subcategoryories of a given fusion categoryory. Only few cases are completely known in the literature.  For instance, in the same paper \cite{proclond},  M\"uger described  the centralizer of all fusion subcategoryories of the categoryof finite dimensional representations of a Drinfeld double of a finite abelian group. More generally, for the categoryof representations of a (twisted) Drinfeld double of an arbitrary finite group a similar formula was then given in \cite{nnw}. For the braided center of Tambara-Yamagami categoryories, in \cite{gnn}, the centralizer was described by computing completely the $S$-matrix of the modular categoryory.  A different formula for braided equivariantized fusion categoryories was given by the author in \cite{bcg}.
\mdn
Also, in the paper \cite{mathz}, the author has given some general formulae for the centralizer of certain subcategoryories (normal in the sense of \cite{brn}) of the categoryof representations of a factorizable Hopf algebras. In the present paper we extend these formulae for any subcategoryof representations completing the study of the centralizer for these type of braided fusion categoryories.
\mdn
Given a quasitriangular Hopf algebra $(A,\;R)$ one can define the linear map 
\beqn
\phi_R:A^* \ra A, \;\;f \mapsto (f \ot \id)(R_{21}R)=f(Q_{1})Q_{2}
\eeqn
where $Q=R_{21}R$ is the monodromy matrix.
It is well known, see \cite{majid} that $\phi_{R}$ sends Hopf subalgebras of $A^{*}$ into left normal coideal subalgebras of $A$. We also define $\ovr{\bf A}:=\phi_{R}(A^{*})$. It was proven in \cite{rqts} that in general, $\ovr{\bf A}$ is a left normal coideal subalgebra of $A$. In this paper we will show that in the semisimple case $\ovr{\bf A}$ is a normal Hopf subalgebra of $A$.
\mdn
 Given a braided fusion category $\cc$ and $\cd$ a fusion subcategoryof $\cc$,  the notion of {\it centralizer of $\cd$ } was introduced in \cite{dgno2}. The centralizer $\cd'$ is defined as the fusion subcategory$\cd'$ of $\cc$ generated by all simple objects $X$ of $\cc$ satisfying  
\beqn c_{X,\; Y}c_{Y,\;X}=\mtr{id}_{X\ot Y}\eeqn
for all objects $Y \in \co(\cd)$ (see also \cite{proclond}). 
\mdn
We prove the following theorem which gives a description for the centralizer of a fusion subcategoryof the categoryof representations of a factorizable Hopf algebra:
 \bt \label{main1}Let $(A, R)$ be a semisimple factorizable Hopf algebra and $L$ be a  left normal  coideal subalgebras of $A$. Then $$ \rep(A//{L})'= \rep(A//{M})\;\;\text{where}\;\;{M}=\phi_{R}((A//{L})^{*})$$
Moreover, in this case one has that:
$
{L}\ovr{\bf A}=\ovr{\bf A}L=\phi_{R}((A//{M})^{*}).
$
\et
\mdn
Recall that the quasitriangular Hopf algebra $(A,R)$ is called {\it factorizable} if and only if the linear map  $\phi_{R}$  is an isomorphism of vector spaces. 
\mdn
Let $V_{0}, V_{1},\dots ,V_{r}$ a  complete set of isomorphism classes of irreducible $A$-modules of a factorizable semisimple Hopf algebra $A$. Let also $\irr(A)=\{\ch_{0}, \ch_{1}, \dots ,\ch_{r}\}$ be the set of irreducible characters afforded by these modules. For any $1\leq j \leq r$  let $E_{j}\in \mtc Z(A)$ be the associated central primitive idempotent of the irreducible character $\ch_{j}$. Then
\dbd
F_{j}:=\phi_{R}^{-1}(E_{j})
\dbd 
is a primitive central idempotent in $C(A)$ since $\phi_{R}$ is an algebra isomorphism.  Following  \cite{CW4} one can define the conjugacy classes $\cc_{j}$ of $A$ as $$\cc_{j}:=\blam \lh F_{j}A^{*},$$where $\blam$ is an idempotent integral of $A$ and $a \lh f=\langle f, a_{1}\rangle a_{2}$ for all $a \in A$ and $f \in A^{*}$. It is well known that these conjugacy classes are the simple $D(A)$-submodules of the induced $D(A)$- module $k\uparrow^{D(A)}_{A}$,  see \cite{zind}. \mdn Our second main result is the following:
\bt \label{main2}Let $(A, R)$ be a semisimple quasitriangular factorizable Hopf algebra and $L$ be a left normal coideal subalgebras of $A$. Then with the above notations one has that:
$$\irr(\rep(A//L)')=\{\ch_{j}\;|\;\cc_{j}\subseteq L\}=\{\ch_{j}\;|\;F_{j}(\blam_{L}) \neq 0\}.$$
\et
\mdn
Here $\blam_{L}$ denotes the integral of $L$ with $\eps(\blam_{L})=1$. Recall that $\blam_{L}$ is  unique up to a scalar element of $k$ and satisfies $a\blam_{L}=\eps(a)\blam_{L}$ for any $l \in L$, see e.g \cite{gmj}.

\mdn \red{Through the rest of the paper we use the following  notation. Given a left normal coideal subalgebra $L$ of $A$ we denote by $L'$ the left normal coideal subalgebra of $A$ for which:
\beq\label{not}
\rep(A//L)'=\rep(A//\lprime).
\eeq
We also use the following notation:
\beq\label{not2}
\phi_{R}((A//L)^{*})=\lstar.
\eeq
Then the first main result Theorem \ref{main1} states that $L'=L^{*}$ is $A$ is a semisimple factorizable Hopf algebra.}
\mdn
We will also prove the following factorization result for factorizable Hopf algebras:

\bt\label{main3}
Let $A$ be a factorizable semisimple Hopf algebra and $K$ be a normal Hopf subalgebra of $A$ such that $\rep(A//K)$ is a nondegenerate fusion categoryory. Then there is another normal Hopf subalgebra $K$ of $A$ such that 
\dbd
A\simeq K\otimes L
\dbd
as Hopf algebras. Moreover $\rep(A//K)'=\rep(A//L)$.
\et
\mdn  We should also notice that in Proposition \ref{intc}, see also the remark following it, we correct a possible error in the results of \cite[Theorem 4.7, Theorem 4.8]{rqts}.

\mdn This paper is organized as follows. In Section \ref{prelim} we recall the basic notions of Hopf algebras and fusion categoryories that are used throughout this paper. In Section \ref{qtr} we recall the main  properties of quasitriangular Hopf algebras and their associated Drinfeld maps. In Section \ref{fmr} we prove the first main result and its consequences.  In Section \ref{smr} we prove the second main result  of this paper. In Section \ref{factoriz} we prove the factorization result mentioned above in Theorem \ref{main3}. 
 \section{ Preliminaries}\label{prelim}
 Let $A$ be a finite dimensional semisimple Hopf
algebra over an algebraically closed field $\kk$  of characteristic zero. 
Then $A$ is also cosemisimple and $S^{2}=\id$, \cite{Lard}. 
The character ring
$C(A):=\mathrm K_{0}(\cc)\ot_{\Z}\kk$ is a semisimple subalgebra of $A^*$ and it has a vector space basis
given by the set $\mtr{Irr}(A)$ of irreducible characters of $A$, see  \cite{Z}. Moreover,
$C(A)=\mtr{Cocom}(A^*)$, the space of cocommutative elements of $A^*$. By duality, the
character ring of $A^*$ is a semisimple subalgebra of $A$ and $C(A^*)=\mtr{Cocom}(A)$. If
$M$ is an $A$-representation with character $\chi$ then $M^*$ is also an
$A$-representation with character $\chi^*=\chi \circ S$. This induces an involution
$``\;^*\;":C(A)\ra C(A)$ on $C(A)$. Let also $m_{ _A}(\ch,\;\mu)$ be the usual multiplicity form on
$C(A)$. Recall that  if $M$ and $N$ are $A$-representations affording characters $\ch, \mu \in C(A)$ respectively, then $m_{ _A}(\ch,\;\mu)$ is defined by $m_{ _A}(\ch,\;\mu):=\dim_{\kk}\hm_{A}(M, N)$. We will use the notation $G(A)$ for the set of group-like elements of $A$.\mdn
Throughout of this paper we denote by $\blam$ an idempotent integral of $A$ and by $t$ an idempotent integral of $A^{*}$. Moreover one has that $t(\blam)=\frac{1}{\dim_{\kk}(A)}$. Recall also \cite{Lar} that
\dbd
\dim_{\kk}(A)\blam=\sum_{d \in \irr(A^{*})}\eps(d)d
\dbd
is the regular character of $A^{*}$. Dually 
\dbd
\dim_{\kk}(A)t=\sum_{\ch \in \irr(A)}\ch(1)\ch
\dbd
is the regular character of $A$. 

\subsection{Left coideal subalgebras}
Let $A$ be finite dimensional Hopf algebra over $\kk$. Recall a {\it left} coideal subalgebra $\bdelta(L)\subset A\ot L$. It is called left normal coideal subalgebra if $L$ is closed under the left adjoint action of $A$, i.e, $a_{1}lS(a_{2})\in L$ for any $l \in L$ and any $a\in A$.
\mdn
Following  \cite{NR} there is a bijection between fusion subcategoryories of $\rep(A)$ and Hopf subalgebras of $A^{*}$. Moreover by Takeuchi's results, see \cite[Theorem 3.2]{Tkq},  in the case of a finite dimensional Hopf algebra $A$ one has a bijection between the set of Hopf subalgebras of $A^{*}$ and the set of  \nleftcid s of $A$. This bijection is given by $L\mapsto (A//L)^{*}$ with inverse given by 
\dbd
(B\xrightarrow{i} A^{*}) \mapsto L:=A^{\mtr{co}\;i^{*}}
\dbd
\mdn
Recall \cite{gmj} there is a unique element $\blam_{L}\in L$ such that $l\blam_{L}=\eps(l)\blam_{L}$, see also \cite{kopp}. Then the coideal subalgebra $L$ is normal if  and only if $\blam_{L}$ is cocentral.
In this case the augmentation ideal $AL^{+}$ where $L^{+}:=\ker(\eps)\cap L$ has the following form 
$$
AL^{+}=A(1-\Lam_{L})=\bann_{A}(\Lam_{L}).
$$
Moreover one has that 
\beq\label{dq}
(A//L)^{*}=\{f \in A^{*}\;|\; f(al)=\eps(l)f(a)\;\;\text{for all} \;a\in A, l \in L\;\;\}.
\eeq
It follows that $\Lam_{L}\rh A^{*}=(A//L)^{*}$. Indeed one has $(\Lam_{L}\rh f)(al)=f(al\Lam_{L})=\eps(l)f(a\blam_{L})=\eps(l)(\blam\rh f)(a)$ for any  $a \in A$. On the other hand clearly $\blam_{L}\lh f=f$ for any $f \in (A//L)^{*}$.

Moreover it can also be shown that
\beq\label{dcid}
L=\{a \in a\;|\; gf(a)=f(1)g(a)\;\text{for all}\; f\in (A//L)^*, g \in A^{*}\}
\eeq
 Since $A$ is free as left $L$-module \cite{sk} it follows that the map
\dbd
A\ot_{L}\kk\simeq A\blam_{L},\;\;a\ot_{L}1\mapsto a\blam_{L}
\dbd
is an isomorphism of $A$-modules. Moreover $A\blam_{L}$ is isomorphic to the regular $A//L$-module, see \cite{iop}.
\mdn
 Moreover, by \cite[Proposition 3.11]{iop}  it follows that the regular character of the quotient  Hopf algebra $A//L$ is isomorphic to the induced module $A\ot_{L}\kk$.
\bl \label{incliusion} Let $A$ be a semisimple Hopf algebra and $L_{1}, L_{2}$ be two left normal  coideal subalgebras of $A$. Then $A{L_{1}}^{+}\subseteq A{L_{2}}^{+}$ if and only if  ${L_{1}}\subseteq {L_{2}}$
\el
\bpf  If  ${L_{1}}\subseteq {L_{2}}$ then clearly  $A{L_{1}}^{+}\subseteq A{L_{2}}^{+}$. Conversely,
suppose that $A{L_{1}}^{+}\subseteq A{L_{2}}^{+}$. Let $\pi:A//{L_{1}}\ra A//{L_{2}}$ the canonical induced projection. By duality this shows that $(A//{L_{2}})^{*}\subseteq (A//{L_{1}})^{*}\subseteq A^{*}$. Then Equation \eqref{dcid} shows that ${L_{1}}\subseteq {L_{2}}$.
\epf

\subsection{Definition of the character ring $C(A)$ and Grothendieck group.}
We denote by $F_{0}, F_{1}, \dots, F_{s}$ the central primitive idempotents of the character ring $C(A)$ where $F_{0}=t$ is the idempotent integral of $A^{*}$. \mdn

\subsection{A result concerning factorization of normal Hopf subalgebras} \label{factnh}
In the proof of Theorem \ref{main3} we also need the following result.
If  $L,K$ are two normal Hopf subalgebras of a semisimple Hopf algebra $A$ with $L\cap K=k$ then $LK\simeq L\otimes K$ as Hopf algebras, see \cite[Theorem 3.5]{cejm1}.

\subsection{On the lattice of fusion subcategoryories} 
\bt\cite{mathz}\label{interms}
Let $L$ and $M$ be two left normal  coideal subalgebras of $A$. Then the following equalities hold in $A^*$:
\bn{enumerate}
\item
\beq \lb{int}
(A//L)^*\cap (A//M)^*=(A//LM)^*.
\eeq
\item
\beq\label{cup}
\langle (A//L)^*,\;(A//M)^*\rangle =(A//(L\cap M))^*.
\eeq
\end{enumerate}
\md
If the Grothendieck ring $\G_0(A)$ is commutative then
\beq\label{coidprod}
\dim_{\kk} (LM)=\frac{ \dim_{\kk} (L\cap M)}{(\dim_{\kk} L)(\dim_{\kk} M)}.
\eeq
\et

For two \nleftcid s $L$ and $M$ of $A$ the above relations can be written as follows:
\beq\label{capsf}
\rep(A//L)\cap\rep(A//L')=\rep(A//LL')
\eeq
\beq\label{veesf}
\rep(A//L)\vee\rep(A//L')=\rep(A//L\cap L')
\eeq
Given two \nleftcid s $L$ and $M$ note that
$
LM=ML
$
since $lm=(l_{1}mS(l_{2}))l_{3}\in ML$ for all $l \in L$ and $m\in M$.
Moreover  any inclusion $M\subseteq L$ allows us to define the quotient
\dbd
L//M:=L//LM^{+}.
\dbd
For any two left normal  coideal subalgebras $L, M$  of $A$ one also has a canonical linear epimorphism
\beq\label{cmorf}
LM//M\xra{\pi} L//L\cap M,\;\;\widehat{lm}\mapsto \eps(m)\hat{l}.
\eeq
 If $C(A)$ is commutative it follows by \cite[Theorem 3.8]{mathz} that
\beq\label{fpdim}
|{L}{M}|=\frac{|{L}||M|}{|{L}\cap M|}
\eeq
for any two \nleftcid s $L$ and $M$ of $A$. Moreover in this case the above canonical epimorphism from Equation \eqref{cmorf} is in fact an isomorphism.
\subsection{Left Lernels and a Hopf-algebraic version of  Brauer's theorem}\label{lLbb}
Let $M$ be an $A$-module and let $\mtr{LKer}_{ _A}(M)$ be the {\it left kernel of $M$}. Recall \cite{gmj}  that  $\mtr{LKer}_{ _A}(M)$ is defined by:
\begin{equation}\label{L}
\mtr{LKer}_{ _{A}}(M)=\{a \in A|\; a_1\ot a_2m=a\ot m,\;\;\text{for all}\;\; m\in M\}
\end{equation}
 Then by \cite{gmj} it follows that $\lker_A(M)$ is the largest left coideal subalgebra of $A$ that acts trivially on $M$. It is also a left normal  coideal subalgebra.
 \mdn
Next theorem generalizes a well known result of Brauer in the representation theory of finite groups. 
\bn{thm}\cite[Theorem 4.2.1]{gmj}.\label{charofim} Suppose that $M$ is a finite dimensional module over a semisimple Hopf algebra $A$. Then \beq
\langle M\rangle =\mtr{Rep}(A//\mtr{LKer}_A(M))
\eeq
where $\langle M\rangle $ is the fusion subcategoryof $\mtr{Rep}(A)$ generated by $M$.
\end{thm}

\section{Quasitriangular and factorizable Hopf algebras}\label{qtr}
Recall that a Hopf algebra $A$ is called {\it quasitriangular} if $A$ admits an $R$-matrix, i.e. an element $R \in A\ot A$ satisfying the following properties:\\
$1) \;
R \Delta(x)=\Delta^{\cop}(x)R
$ for all $x \in A$.
\\$2)\;
(\Delta \otimes \id)(R)= R^{1} \ot r^{1} \ot R^{2}r^{2}
$
\\$3)\;
( \id \otimes \Delta)(R)=R^{1}r^{1} \ot r^{2}\ot R^{2}.
$
\\$4)\;
( \id \otimes \eps)(R)=1=(\eps \ot \id)(R).
$
Here $R=r=R^{1}\ot R^{2}=r^{1}\ot r^{2}$.
If $(A, R)$ is a quasitriangular Hopf algebra then the categoryof representations is a braided fusion categorywith the braiding given by
\beq
c_{M, N}:M\ot N\ra N\ot M,\; m\ot n \mapsto R_{21}(n\ot m)=R^{2}n\ot R^{1}m
\eeq
for any two left $A$-modules $M, N\in \rep(A)$ (see \cite{Kas}). Recall that $R_{21}:= R^{2}\ot R^{1}$.
Denote $Q:=R_{21}R$.
Then the monodromy of two objects is defined as:
\beq\label{mdrmy}
c_{M,N}c_{N,M}:M\ot N \ra N\ot M,\;(m\ot n)\mapsto R^{(2)}R^{(1)}m\ot R^{(2)}R^{(1)}n=Q(m\ot n)
\eeq
\md
A quasitriangular Hopf algebra $(A,\;R)$ is called {\it factorizable} if and only if the linear map  
\beq\label{dr}
\phi_R:A^* \ra A, \;\;f \mapsto (f \ot \id)(R_{21}R)
\eeq 
is an isomorphism of vector spaces. In this situation, following \cite[Lemma 2.2]{schfact} $\phi_R$ maps the character ring  $C(A)$ onto the center of $Z(A)$ of $A$.
\mdn
One can also define the map $ _{R}\phi(f)=(\id \ot f)(Q).$. Moreover by \cite[Theorem 2.1]{schfact} one  has that
\beq\lb{eqs2}
 _R\phi(f\chi)=\;_R\phi(f)\;_R\phi(\chi)
\eeq
for all $f \in A^*$ and $\chi \in C(A)$.  Thus $_{R}\phi|_{C(A)}:C(A)\ra Z(A)$ is an isomorphism of $\kk$-algebras.
\mdn
In the case of a factorizable Hopf algebra this map is also bijective and moreover by \cite{CW4} the two maps coincides on the character ring $C(A)$ of $A$.


 \subsection{On the $S$-matrix for a factorizable Hopf algebra}
Let $A$ be a semisimple factorizable Hopf algebra. By \cite{EG} one has that the $S$-matrix $(s_{ij} )$ of the modular tensor category$\Rep(A)$ is given by \beq\label{s}s_{ij}=\chi_{i}(\phi_{A}(\ch_{j^{*}}))\eeq Moreover it is not difficult to see that for all $1\leq i,j \leq s$, one has that $s_{ij} = s_{ji}$, and $s_{ij} = s_{i^*j^*}$ and
$s_{ij^*} = s_{ji^*}$ (cf. \cite{BaKi, schfact}). In this case one has the following inequality 
\beq\label{ifs}
|s_{ij}|\leq \ch_{i}(1)\ch_{j}(1).
\eeq
\subsection{Transmutation theory} In this subsection we recall some elements from \cite{majid} and \cite[Section 4]{rqts}. We assume the reader is familiar with the notion of braided Hopf algebras in monoidal categoryories.
\mdn
By \cite[Example 9.4.9]{majid}, given a quasitriangular Hopf algebra $(H,R)$ one has a braided Hopf algebra $\underline H$ in the category$H$-mod of left $H$-modules with $\ul H=H$ as algebras (thus as vector spaces too) and the following additional structures
\beq\label{coprod}
\underline \Delta(a)=a_{1}S(R^{(2)})\ot ad_{R^{(1)}}(a_{2})
\eeq
\mdn
By [14, Example 9.4.10] one also has a braided Hopf algebra  $\ul H^{*}$ algebra in the categoryof right $H^{*}$-comodules, i.e. left $H$-modules . One has that $\ul H^{*}=H^{*}$ as coalgebras and the multiplication given by
\beq\label{starprod}
p\ul .q=(S(p_{1})p_{3}\ot q_{1})(R)p_{2}q_{2}
\eeq
Moreover, $\ul H^{*}$ is regarded as right $H$-comodule with the coadjoint structure
\beq\label{coadjstr}
\rho(p)=p_{2}\ot p_{1}S(p_{3})\in \ul H^{*}\ot H
\eeq
Note that this is equivalent with the following left  $H$- action:
\dbd
h._{\mtr{coad}}p=p(h_{1}?S(h_{2}))
\dbd
\mdn
In this context, the map $\phi_{R}:\ul H^{*}\to \ul H$ given by Equation \eqref{dr} 
becomes a morphism of braided Hopf algebras . In other words $\phi_{R}$ is a left $H$-module
map that transforms the product from Equation  \eqref{starprod} into the product of H
and the coproduct of $\ul H^{*}$ into the coproduct from Equation \eqref{coprod}. See\cite[Propositions 2.1.14 and 7.4.3]{majid} for more details.
\mdn
By \cite{rqts} it follows that any  subcoalgebra in $\ul H$ corresponds to a left normal  coideal subalgebra of $H$. Moreover, by \cite[Lemma 1.1]{rqts} one has that $\phi_{R}(C)$ is a left normal  coideal subalgebra for any subcoalgebra $C$ of $H^{*}$.

\subsection{Formulae for the Frobenius-Perron dimension and double centralizer} 
Let $\ccb$ and $\cd$ be fusion subcategoryories of a braided fusion category$\cc$. Following \cite[Proposition Theorem 3.10]{DGNO} one has that
\beq\label{both}
\fp(\ccb \cap \cd')\fp(\cd)=\fp(\ccb' \cap \cd)\fp(\ccb)
\eeq

\beq\label{double}
\cd''=\cd\vee \cc'
\eeq
In particular, for $\ccb=\cc$ the first Equation \eqref{both} becomes
\beq\label{product}
\fp(\cd)\fp(\cd')=\fp(\cc)\fp(\cd\cap\cc')
\eeq
Let now $\cc=\rep(A)$ for some quasi-triangular Hopf algebra $A$. We plug in above formulae  $\cd=\rep(A//M)$ and $\mtc B=\rep(A//L)$ for any two left normal coideal subalgebras of $A$. Using the notations from Equation \eqref{not} one has that $\ccb'=\rep(A//\lprime)$ and $\cd'=\rep(A//\mprime)$. Note that by Equation \eqref{capsf}  the  equality from Equation \eqref{both} becomes 
$
|L\mprime||M|=|M\lprime||L|.
$
Using now Equation \eqref{fpdim} it follows that:
\beq\label{proport}
\frac{|\lprime|}{|\lprime\cap M|}=\frac{|\mprime|}{|\mprime\cap L|}
\eeq
\mdn
In particular for $M=k$, since $k^{*}=\overline{\bf A}$ one has that
\beq\label{proport2}
|\lprime|=\frac{|\overline{\bf A}|}{|\overline{\bf A}\cap L|}
\eeq
for any left normal coideal subalgebra $L$ of $A$. If moreover, $A$ is factorizable, i.e. $\overline{\bf A}=A$ then
\beq\label{proportfact}
|\lprime|=\frac{|A|}{|L|}
\eeq
\mdn
As noticed above, from \cite{rqts} one has that $\cc'=\rep(A//\overline{\bf A})$.  
With our notations this means that $k^{*}=\overline{\bf A}$. Since $\cd''=\cd\vee \cc'$ it follows by Equation \eqref{int} that 

\beq\label{doublel}
\mdprime=M\cap { \overline{\bf A}}
\eeq
for any left normal  coideal subalgebra of $A$.  
\mdn
\mdn
Suppose that $\cd=\rep(A//{L})$. Then $\cd'=\rep(A//{\lprime})$ and
$$\mtc Z_{2}(\cd)=\cd'\cap \cd=\rep(A//L{\lprime})$$
Note also that the M\"uger center $\cc'$ satisfies $\cc'\subseteq \cd'$  for any fusion subcategory$\cd$ of $\cc$. This gives that $\lprime\subseteq \overline{\bf A}$ for any left normal coideal subalgebra $L$ of $A$.

\section{Proof of the first main theorem on M\"uger  centralizer}\label{fmr}
Let $A$ be a semisimple quasitriangular Hopf algebra and $\cd$ be a fusion subcategoryof $\rep(A)$.  In this section we prove the first main theorem mentioned in the introduction.
\mdn

\bl\label{firstequiv} Let $(A, R)$ be a semisimple quasitriangular Hopf algebra and $L$, $M$ be two \nleftcid s. Then
$
 \rep(A//M)\subseteq \rep(A//L)'$
 if and only if 
\beq\label{cond1}
Q(\blam_{L}\ot \blam_{M})=(\blam_{L}\ot \blam_{M})
\eeq
 \el
 \bpf Two fusion subcategoryories of $\rep(A)$ centralize each other if and only if their regular representations centralize.
Thus one needs to show that the two regular characters of $A//L$ and $A//M$ centralize each other if and only if Equation \eqref{cond1} holds. On the other hand from the definition of the braiding in $\rep(A)$ the two characters centralize each other if and only if $Q=R^{2}r^{1}\ot R^{1}r^{2}$ acts as identity on their tensor product $A//L\ot A//M$. As noticed above one has $k\uw^{A}_{L}\ot k\uw^{A}_{M}=A\blam_{L}\ot A\blam_{M}$. Since $\blam_{L}$ and $\blam_{M}$ are central elements of $A$ it is clear that $Q$ acts as identity on this subspace of $A\ot A$ if and only if Equation \eqref{cond1} holds.
\epf
\subsection{Proof of Theorem \ref{main1}}
We prove the following slightly more general result:
 \bt \label{main1'}Let $(A, R)$ be a semisimple quasitriangular Hopf algebra and $L, M$ be two left normal coideal subalgebras of $A$. Then $$\rep(A//{M})\subseteq \rep(A//{L})'\;\iff\;{M}\supseteq \phi_{R}((A//{L})^{*}).$$

Moreover, in this case one has that:
$
{L}\ovr{\bf A}\supseteq \phi_{R}((A//{M})^{*}).
$
\et
\bpf 
As in  Lemma \ref{firstequiv} the above inclusion holds if and only if  $Q$ acts trivially on $A//{L}\ot A//{M}$. Therefore  
\beq\label{eqcent1}
Q^{1}\blam_{{L}}\ot Q^{2}\blam_{{M}}=\blam_{{L}}\ot \blam_{{M}}
\eeq

As above $(A//{L})^{*}=\blam_{{L}}\rh A^{*}$ and therefore $_{R}\phi(\blam_{{L}}\rh f)=f(Q^{1}\blam_{L})Q^{2}$ for any $f \in (A//L)^{*}$. From here, applying Equation \eqref{eqcent1} it follows
\beq
\phi_{R}(\blam_{{L}}\rh f)\blam_{{M}}=f(Q^{1}\blam_{{L}})Q^{2}\blam_{{M}}=f(\blam_{{L}})\blam_{{M}}
.\eeq
On the other hand note that \dd \eps(\phi_{R}(\blam_{{L}}\rh f))=f(\blam_{{L}})\dd. If ${{L}}^{*}=\phi_{R}((A//{L})^{*})$ then  it follows that $AL^{*+} \subseteq AM^{+}$ and by Lemma \ref{incliusion} one has that ${{L}}^{*}\subseteq {M}$.
\mdn
Note that $\rep(A//{L})'\supseteq \rep(A//{M})$ also implies by centralising again that $$\rep(A//{M})'\supseteq  \rep(A//{L})''=\rep(A//{L}\overline{\bf A}).$$ By the above argument one has that
$
{L}\overline{\bf A}\supseteq \phi_{R}((A//{M})^{*}).
$
\epf
\blue{\bf 
\br
One can apply $f$ on the second tensor and give that
$$_R\phi(A//M)^*\subseteq L.$$
\er
}
\br
By \cite[Lemma 2.3]{rqts} one has that
$S_{R}\phi=\phi_{R}s$ where $S$ and $s$ are the antipodes of $A$ and $A^{*}$ respectively. Thus one can also write that $$_{R}\phi((A//{L})^{*})\subseteq {S(M)} \iff \rep(A//{M})\subseteq \rep(A//{L})'.$$
\er
\section{Conjugacy classes and M\"uger centralizer}\label{smr}
In this section we will prove the second main result.
 \subsection{Duality between the character ring and the center}
Let $A$ be a semisimple Hopf algebra over the ground field $\kk$. Let us denote by $\irr(A)$ the set of irreducible characters of $A$. We suppose that $\irr(A)=\{\ch_{0}, \ch_{1}, \dots, \ch_{r}\}$. Without loss of generality we may suppose that $\ch_{0}=\eps$. Let also $E_{0}, E_{1}, \cdots, E_{r}$ be the corresponding primitive central idempotents in $A$. 
The evaluation form
\beq
C(A)\otimes Z(A)\ra k,\;\; \ch \otimes a\mapsto \ch(a)
\eeq
is nondegenerate. A pair of dual bases for this form is given by  $\{\ch_{i},\; \frac{1}{n_{i}}E_{i}\}$
since  $\langle \ch_{i},\; \frac{1}{n_{j}}E_{j}\rangle =\delta_{i,j}$ for any $1, \leq i,j \leq r$.
\mdn  According to \cite{CW4} in the case of a commutative ring $C(A)$   there are is another  pair  of dual bases corresponding to this nondegenerate form.  This pair of dual bases is given in terms of the conjugacy class sums as defined in \cite{CW4}.
Recall that the conjugacy class $\cc_{j}$ is defined as
$
\cc_{j}=\blam \lh F_{j}A^{*}.
$
One has that $\bdelta \cc_{j}\subseteq A\otimes \cc_{j}$ and $\cc_{j}$ is closed under the  left adjoint action of $A$ on itself. Thus $\cc_{j}$ is  a simple $D(A)$-submodule  of $k\uw^{D(A)}_{A}\simeq A$, see \cite{zind} and \cite{nrm-ha-art}. Recall from \cite{nrm-ha-art} that the $D(A)$-module structure of $A$ is given by:
$(f \bwt x).a=x_{1}aS(x_{2})\lh S^{-1}f$, for any $a,x \in A$ and any $f \in A^{*}$.
\mdn
Note also that $A=\bigoplus_{j=1}^{r}\cc_{j}$ since $f\mapsto \blam \lh f$ is a bijective map and therefore it preserves direct sums. One can also define the corresponding class sum
\beq\label{clssum}
C_{j}=\blam \lh (\dim A)F_{j}
\eeq
and its normalized class sum $\eta_{j}:=\frac{C_{j}}{\eps(C_{j})}$. Note that $C_{j}\in Z(A)$ and since $\dim_{\kk}Z(A)=\dim_{\kk}C(A)$ it follows that $\cc_{j}\cap Z(A)=\kk C_{j}$.
\br\label{genmz} By the Class Equation for semisimple Hopf algebras, see \cite{leq}, one has that the value
$
n_{j}:=\frac{\dim_{\kk}A^{*}}{\dim_{\kk}(A^{*}F_{j})}
$
is an integer.  Moreover as in  \cite[Equation (11)]{CW4}
one can write that
$
E_{j}(\Lam)=\frac{1}{n_{j}}.
$
\er \mdn
Then the second pair of dual bases is given by $
\{F_{i},\; \frac{n_{i}}{\dim_{\kk}(A)}C_{i}\}$, see for instance \cite[Equation (17)]{CW4}. Thus
$
\langle F_{i},\; \frac{n_{j}}{\dim_{\kk}(A)}C_{j}\rangle =\delta_{i,j}
$

\subsection{Decomposition of the integral}
Let $L$ be a \nleftcid\; of a semisimple Hopf algebra $A$ with a commutative character ring. Since $L$ is also a $D(A)$-submodule of $A$ one can write that 
$
L=\bigoplus_{j \in \mtc I_{L}}\cc_{j}
$
for some subset $\mtc I_{L}\subset \{0,1,\dots , r\}$.
Then we write
$
\blam_{L}=\sum_{j \in \mtc I_{L}}\blam_{j}
$
with $\blam_{j}\in \cc_{j}$ for the decomposition of the idempotent integral $\blam_{L}$ of $L$. It follows that for any $a \in A$ one has that
$
\eps(a)\blam_{L}=a_{1}\blam_{L}Sa_{2}=\sum_{j\in \mtc I_{L}}a_{1}\blam_{j}Sa_{2}.
$
Thus $a_{1}\blam_{j}S(a_{2})=\eps(a)\blam_{j}$ which shows that $\blam_{j}\in \cc_{j}$ are central elements.  Moreover since $\cc_{j}\cap Z(A)$ has dimension $1$ one then can write that
\beq\label{decl1}
\blam_{L}=\sum_{j \in \mtc I_{L}}\al_{j}C_{j}
\eeq
and denote by $\mtc I'_{L}$ the set of indices $j \in \mtc I_{L}$ for which these coefficients $\al_{j}\neq 0$ are not zero. Then it follows that
$
L=\blam_{L}\lh A^{*}=\bigoplus_{j \in \mtc I'_{L}}C_{j}\lh A^{*}=\bigoplus_{j \in \mtc I'_{L}}\cc_{j}. 
$
which shows that $
\mtc I'_{L}=\mtc I_{L}.
$
\bl \label{fct}(see also \cite[Theorem 5.13]{repalg}.)
Suppose that $A$ is a semisimple Hopf algebra with a commutative character ring $C(A)$. Then
$F_{j}$ coincides to the functional $p_{\cc_{j}}\in A^{*}
$ defined as the unique functional that coincides to $\eps$ on $\cc_{j}$ and it is equal to zero on the other conjugacy classes $\mtc C_{l}$ with $l \neq j$. 
\el 
\bpf
One has the following
$
\langle F_{j}, \blam\lh F_{l}f\rangle =\langle F_{l}f, \blam_{1}\rangle \langle F_{j}\blam_{2}\rangle =\langle F_{l}f, \blam_{2}\rangle \langle F_{j}, \blam_{1}\rangle =\delta_{j,l}\eps(\blam \lh F_{l}f).
$
Note that in the above computations we used the cocommutativity of $\blam$.\blue{\bf \tiny }
\epf
 We shall use the notation $\lam_{L}\in (A//L)^{*}$ for the idempotent integral of the Hopf algebra $(A//L)^{*}$. Clearly $\lam_{L}\in C((A//L)^{*})\subset C(A^{*})$.
\bl\label{charval}
Let $L$ be a left normal  coideal subalgebra of  a semisimple Hopf algebra $A$. If we write 
$\lam_{L}=\sum_{j \in \mtj_{L}}F_{j}$
then $\mtc I_{L}=\mtj_{L}$ and  
$j \in \mtj_{L} \iff  F_{j}(\blam_{L})\neq 0$.
\el 
\bpf 
If $j \in \mtj_{L}$ then $\cc_{j}=\blam \lh F_{j}H^{*}=\blam \lh F_{j}\lam_{L}H^{*}=\blam \lh \lam_{L}F_{j}H^{*}=\blam_{L} \lh F_{j}H^{*}\subseteq L$ which shows that $j \in \mtc I_{L}$.  Thus $\mtj_{L}\subseteq \mtc I_{L}$.

Suppose now that $j \in \mtc I_{L}$ where $L=\bigoplus_{j \in \mtc I_{L}}\cc_{j}$. Then by Lemma \ref{fct} one has that
$F_{j}(\blam_{L})=\al_{j}\eps(C_{j})\neq 0.$ Thus if $j \in \mtc I_{L}$ then $F_{j}(\blam_{L})\neq 0$.

On the other hand by \cite[Lemma 1.1]{CW10} one has that $\blam_{L}=\blam\lh \lam_{L}$. Thus if  $0 \neq F_{j}(\blam_{L})=(\lam_{L}F_{j})(\blam)$ then $\lam_{L}F_{j}\neq 0$ i.e, $j \in \mtj_{L}$. 
\epf

 \br With the above notations, by Equation \eqref{decl1} one has that $\al_{j}\eps(C_{j})=F_{j}(\blam_{L})=F_{j}(\blam)=\frac{1}{n_{j}}.$
On the other hand note that Equation \eqref{clssum} gives that $\eps(C_{j})=\dim_{\kk}(A) F_{j}(\blam)=\frac{\dim_{\kk}(A)}{n_{j}}.$ Thus one has that
\beq\label{decl2}
\blam_{L}=\sum_{j\in \mtc I_{L}}\frac{1}{n_{j}\eps(C_{j})}C_{j}=\frac{1}{\dim_{\kk}(A)}\sum_{j \in \mtc I_{L}}C_{j}
\eeq
\er\mdn
Let $A$ be a semisimple Hopf algebra with commutative  character ring $C(A)$. Then $\{F_{j}\}$ form a $\kk$-linear basis for $C(A)$ and for any character $\ch \in C(A)$ one can write $\ch=\sum_{j=0}^{r}\al_{\ch,j}F_{j}$ with $\al_{\ch,j }\in \kk$.
\bp If $A$ is a semisimple Hopf algebra with commutative  character ring then one has
\dd\ch \in C(A//L)\dd \;if and only if $\al_{\ch, j}=\ch(1)$ for all $j \in \mtj_{L}$.
\ep

\bpf
If $\ch \in C(A//L)$ then $\ch\lam_{L}=\ch(1)\lam_{L}$ and this implies that $\al_{\ch,j}=\ch(1)$. Conversely if $\al_{\ch, j}=\ch(1)$ for all $j \in J_{L}$ then $\ch \lam_{L}=\lam_{L}\ch(1)$ and, \dbd\ch(\blam_{L})=\ch(\lam_{L}\rh \blam)=\ch\lam_{L}(\blam)=\ch(1)\lam_{L}(\blam).
\dbd
On the other hand note that $\lam_{L}(\blam)=1$ since $\lam_{L}=\blam\rh t$ by  \cite[Lemma 1.1]{CW10}.
\epf
\mdn

\br
Suppose as above that $L=\bigoplus_{j \in \mtc I_{L}}\cc_{j}.$ Then $S(L)=\bigoplus_{j\in \mtc I_{L}}S(\cc_{j})$ and $S(\mtc C_{j})$ are conjugacy classes for the other left  adjoint  action, $a.l=a_{2}lS(a_{1})$.
\er

 \subsection{Proof of the second main result}
  Let $(A, R)$ be a semisimple quasitriangular Hopf algebra and $\cd=\rep(A//L)$ be a fusion subcategoryof $\rep(H)$. 
\bt \label{main2}Let $(H, R)$ be a semisimple quasitriangular factorizable Hopf algebra and $L$ be a left normal coideal subalgebras of $H$. Then with the above notations one has that:
$$\irr(\rep(H//L)')=\{\ch_{j}\;|\;\cc_{j}\subseteq L\}=\{\ch_{j}\;|\;F_{j}(\blam_{L})\neq 0\}.$$
\et
\bpf
 Two irreducible characters $\ch_{i}$ and $\ch_{j}$ of $A$ centralize each other if and only if
$s_{ij}=d_{i}d_{j}$. Following \cite{CW2} one has that $s_{ij}=d_{j}\langle \ch_{i}, \eta_{j}\rangle $. Thus  $\ch_{i}$ and $\ch_{j}$ centralize each other if and only if
$
\langle \ch_{i}, \eta_{j}\rangle =d_{i}.
$

By \blue{\bf \cite[Theorem 3.14]{CW2}} this happens if and only if $
\mtc C^{j}\subseteq \lker_{V_{i}}=H^{B_{i}}
$
where $B_{i}$ is the Hopf subalgebra of $H^{*}$ generated by the irreducible character $\ch_{i}$. We deduce from here that
$
\{  \ch_{i} \}'=\{ \ch_{j} \;|\; {\mtc C}^{j}\subseteq  \lker_{V_{i}}\}.
$
If $\cd=\rep(H//L)$ is fusion subcategoryof $\rep(A)$ then by intersection one has that
\dbd
\cd' = \bigcap_{\ch_{i }\in \irr(\cd)}\{\ch_{i}\}'=\{\ch_{j}\;|\:\mtc C^{j}\subseteq
\bigcap_{\ch_{i}\in \irr(\cd)}\lker_{V_{i}}\}.
\dbd
On the other hand, using Brauer's Theorem \ref{charofim} it is easy to see that if  $\cd=\rep(H//L)$ then $\cap_{\ch_{i}\in \irr(\cd)}\lker V_{i}=L.$ This implies the required formula $\cd'=\{\ch_{j}\;|\:\mtc C^{j}\subseteq L\}.$
\epf
\br From the previous proof one can deduce  that for any two characters $\chi_{i}$ and $\chi_{j}$ one has
$
\mtc C^{j}\subseteq \lker_{V_{i}}=H^{B_{i}}
$
if and only if 
$
\mtc C^{i}\subseteq \lker_{V_{j}}=H^{B_{i}}
$
\er

\br Note that using Lemma \ref{charval} our second main result now generalizes \cite[Theorem 1.4]{mathz} from normal Hopf subalgebras to coideal subalgebras. More precisely we have shown that:
\dbd
\rep(A//L)'=\{\ch_{j}\;|\; F_{j}(\blam_{L})\neq 0\}=\{\ch_{j}\;|\; j \in \mtc I_{L}\}
\dbd
\er
\mdn  Putting together the two descriptions of the centralizer of a fusion subcategoryone obtains the following corollary:
\bc\label{puttgh} Let $(A, R)$ be a factorizable semisimple Hopf algebra and $L$ be a \nleftcid of $A$. With the above notations, if $\cc_{j}$ is an irreducible representation of $A$ we have that the following assertions are equivalent:
\bne
\item $V_{j}\in \co(\rep(A//L)')$.
\item $\cc_{j}\subseteq L $
\item $\phi_{R}((A//L)^{*})\subseteq \lker_{A}(V_{j})$
\ene
\ec
\bpf
The definition of the left kernel of an $A$-module is recalled in  subsection \ref{lLbb} . Since $A$ is factorizable one has that $\rep(A//L)'=\rep(A//L^{*})$. Then the first main result can be rephrased now as $ \irr(\rep(A//L)')=\{\ch_{j}\;|\;\lker_{A}(V_{j})\supseteq L^{*}\}.$
Thus for any $j \geq 0$ one has the equivalence 
$
V_{j}\in \co(\rep(A//L)') \iff \mtc C^{j}\subseteq L\iff \lker_{A}(V_{j})\supseteq L^{*}.$\epf
\mdn
 Define $R_{V_{j}}\subset A^{*}$ as the subcoalgebra of $A^{*}$ generated by $\ch_{j}$.
\br Let $A$ be a semisimple factorizable Hopf algebra and $L$ a left normal coideal subalgebra of $A$. By \blue{\bf \cite[Lemma 4.2 i)]{CW5}} one has that $\phi_{R}(R_{V_{j}})=\mtc C_{j}.$ Therefore if $\cd:=\rep(A//L)$ then
$
(A//L)^{*}=\bigoplus_{\ch_{j} \in \irr(\cd)}R_{V_{j}}
$
and $L^{*}=\phi_{R}((A//L)^{*})=\bigoplus_{\ch_{j} \in \irr(\cd)}\phi_{R}(R_{V_{j}})$ i.e
$
L^{*}=\bigoplus_{\ch_{j} \in \irr(\cd)}\cc_{j}
$
\er
\subsection{Formulae for the integrals}Note that for any left normal coideal subalgebra $L$ of a semisimple Hopf algebra $A$ one has that 
\dbd
\blam_{L}=\sum_{\ch_{j} \in \irr(A//L)}E_{j}.
\dbd
This implies that
\dbd
\blam_{L^{*}}=\sum_{j \in \mtc I_{L}}E_{j}=\frac{1}{\dim_{\kk} A}\sum_{j\in \mtc I_{L^{*}}}C_{j}.
\dbd
If $A$ is factorizable then $L^{**}=L$ and we also have  the following formula:
\dbd
\blam_{L}=\frac{1}{\dim_{\kk} A}\sum_{j\in  \mtc I_{L}}C_{j}=\sum_{j\in \mtc I_{L^{*}}}E_{j}.
\dbd
On the other hand clearly:
\dbd
\lam_{L}=\sum_{j\in \mtc I_{L}}F_{j}=\frac{\dim_{\kk} L}{\dim_{\kk} A}\sum_{j\in  \mtc I_{L^{*}}}\ch_{j}(1)\ch_{j}.
\dbd

\section{Proof of the Hopf algebra factorization}\label{factoriz}
In this section we will show first that the map $\phi_{R}$ sends normal Hopf subalgebras to normal Hopf subalgebras. Using \cite[Lemma 2.3]{rqts}, we have to show that
$
_{R}\phi(A//L)^{*}=\phi_{R}(A//L)^{*}
$
when $L$ is a normal Hopf subalgebra of $A$.
\mdn

\bp\label{gcw}
Let $(A, R)$ be a factorizable semisimple Hopf algebra  and $\ch=\ch_{M}$ be the character of an $A$-module $M$. Let 
$
E_{M}:=\sum_{\{\ch_{i}\;|\;(\ch_{i}, \ch)\rangle 0\}}E_{i}
$
be its associated central idempotent to $M$ in $A$. If $R_{\ch}$ denotes the subcoalgebra of $A^{*}$ then $
\phi_{R}(C_{\ch})=\blam\lh \phi_{R}^{-1}(E_{M})A^{*}
$.
\ep
\bpf
One has by \cite[Lemma 4.2.i)]{CW5} that $\phi_{R}(C_{\ch_{j}})=\cc_{j}=\blam \lh F_{j}A^{*}$. Note that 
$
C_{\ch}=\bigoplus_{\{\ch_{i}\;|\;(\ch_{i}, \ch)\rangle 0\}}C_{\ch_{i}}
$
Then the result follows since $\lh$ is an isomorphism.
\epf
\br One can write a similar result for the other Drinfeld map $_{R}\phi$. Since
$S_{R}\phi(C_{\ch_{j}})=\phi_{R}(C_{\ch_{j^{*}}})$ it follows that 
$\;_{R}\phi(C_{\ch_{j}})=S^{-1}(\blam\lh F_{j^{*}}A^{*})=F_{j}A^{*}\rh \blam$ since $A$ is semisimple and $S^{2}=\id$.
\er
\bt\label{secform}
Let $L$ be a left normal  coideal subalgebra of a factorizable semisimple Hopf algebra $H$. Then $\mtc I_{L^{*}}=\{j\;|\;\ch_{j}\in\irr(A//L)\}$, i.e.
\dbd
L^{*}=\bigoplus_{\{j\;|\; \ch_{j}\in \irr(A//L)\}}\cc_{j}
\dbd
 where by $\lam_{L}$  we denote an idempotent integral in $(A//L)^{*}$.
\et
\bpf
One has that $(A//L)^{*}=\bigoplus_{\ch_{j}
\in \irr(A//L)} C_{\ch_{j}}$. By the previous proposition $L^{*}:=\phi_{R}((A//L)^{*})=\bigoplus_{\ch_{j}
\in \irr(A//L)} \phi_{R}(C_{j})=\bigoplus_{\ch_{j}
\in \irr(A//L)} \mtc C_{j}$. This shows the first equality which also implies $\mtc I_{L^{*}}=\{j\;|\;\ch_{j}\in\irr(A//L)\}$.
\epf
\br\label{nonsssd}
As already mentioned, one can also consider the second Drinfeld map $_{R}\phi:A^{*}\ra A$. By \cite[Lemma 2.3]{rqts} one has that
$
S_{R}\phi=\phi_{R}s
$
where $S$ and $s$ are the antipodes of $A$ and $A^{*}$ respectively. This shows that 
$
\;_{R}\phi((A//L)^{*})=S^{-1}(\phi_{R}((A//L)^{*}))=S^{-1}(L^{*})
$
Thus $\;_{R}\phi $ sends  subcoalgebras of $A^{*}$ into right coideal subalgebras. 
\er
\bt\label{nrmphi} Let $A$ be a semisimple factorizable Hopf algebra. If $L$ is a normal Hopf subalgebra of $A$ then $L^{*}$ is also a normal Hopf subalgebra of $A$.
\et
\bpf
Using Remark \ref{nonsssd} one has to show that $\phi_{R}((A//L)^{*})=\;_{R}\phi((A//L)^{*})$. On the other hand by Theorem \ref{secform} and the remark following it one has that
$
\phi_{R}((A//L)^{*})=\blam \lh \lam_{L}A^{*}$ and $\;_{R}\phi((A//L)^{*})= A^{*}\lam_{L} \rh \blam.$ Since $\blam$ is cocommutative and $\lam_{L}$ is  a central element  (by \cite{Mas'}) the result follows. 
\epf
\br In particular our result implies that $\ovr{\bf A}:=\phi_{R}(A^{*})$ is a normal Hopf algebra in the case of a semisimple Hopf algebra $A$. Note that in \cite{rqts} it was shown that $\ovr{\bf A}$ is a left normal coideal subalgebra.\er
\bl Suppose that $(A, R)$ is a semisimple quasitriangular Hopf algebra and $L$ a  left normal coideal subalgebra of $A$. Then with the above notations one has that:
$
\phi_{R}(\Lam_{L}\rh t)=\lam_{\lstar}.
$
\el
\bpf
 It is proven  \cite[Lemma 1.1]{CW10} one has that $\Lam_{L}\rh t=\lam_{L}.$ Since $\Lam_{L}\rh t$ is a character and $\phi_{R}$ is multiplicative when one variable is a character  one has that
$\phi_{R}((\Lam_{L}\rh f)(\Lam_{L}\rh t))=\phi_{R}(\Lam_{L}\rh f)\phi_{R}(\Lam_{L}\rh t).$
One has then to show that $\phi_{R}(\Lam_{L}\rh f)\phi_{R}(\Lam_{L}\rh t)=\eps(\phi_{R}(\Lam_{L}\rh f))\phi_{R}(\Lam_{L}\rh t)$
This follwos since $\Lam_{L}\rh t$ is an integral in $(A//L)^{*}$ and $\eps(\phi_{R}(f))=f(1)$.
\epf
\subsection{Description of an intersection} 
Let $L$ be a left normal coideal subalgebra of a semisimple factorizable Hopf algebra $A$.
Denote $B:=(A//L)^{*}$ and $B':=(A//L^{*})^{*}$.
We have the following result concerning the description of the intersection of the two Hopf subalgebras of $A^{*}$, correcting  \cite[Theorem 4.8]{rqts}.
\bt \label{intc}Suppose that $(A, R)$ is a factorizable Hopf algebra and $L$ is a left normal  coideal subalgebra of $A$. Then 
\beq\label{inp}
\phi_{R}(B\cap B')=L\cap L^*
\eeq
\et
\bpf Let $M:=\phi_{R}(B\cap B')$.
From the first main result above we know that 
\dbd
\rep(A//M)'=\rep((B\cap B')^{*})
\dbd
\mdn
On the other hand we know from Equation \eqref{int} that
$
B\cap B'=(A//L)^{*}\cap(A//L^*)^{*}=(A//LL^*)^{*}.
$
Thus $\rep(A/M)'=\rep(A//LL^*)$ which gives that 
\dbd
\rep(A/M)''=\rep(A//LL^*)'=(\rep(A//L)\vee \rep(A//L^*))'=\dbd\dbd=\rep(A//L)'\cap \rep(A//L^*)'=\rep(A//L\cap L^*).
\dbd
Thus $M=L\cap L^{*}$ as desired.
\epf

We will prove the following factorization result, see Theorem \ref{main3}.
\bt Let $A$ be a factorizable semisimple Hopf algebra and $K$ be a normal Hopf subalgebra of $A$ such that $\rep(A//K)$ is a nondegenerate fusion categoryory. Then there is another normal Hopf subalgebra $K$ of $A$ such that 
\dbd
A\simeq K\otimes L
\dbd
as Hopf algebras. Moreover $\rep(A//K)'=\rep(A//L)$.
\et
\bpf Since $\rep(A//K)$ is a nondegenerate fusion categoryit follows that $\rep(A//K)\cap \rep(A//L)=\mtr{Vec}$ and $\rep(A//K)\vee\rep(A//L)\simeq \rep(A).$ This gives by the previous results that $KL=A$ and $K \cap L=\kk$. Thus $A\simeq K\ot L$ as Hopf algebras by the factorization {result mentioned in Subsection \ref{factnh}.}
\epf
\bt \label{main2'} Let $A$ be a factorizable semisimple Hopf algebra and $K$ be a normal Hopf subalgebra of $A$ such that $\rep(A//K)$ is a nondegenerate fusion categoryory. Then the two normal Hopf subalgebras $L$ and $K$ from Theorem \ref{main3} are also factorizable Hopf algebras.
\et
\bpf
Since $A\simeq K\ot L$ as Hopf algebras we have  canonical Hopf projections $\pi_{K}: A\ra K$ and $\pi_{L}:A\ra L$. They induces Hopf algebra isomorphisms $A//L\xra{\pi_{L}} K$ and $A//L\xra{\pi_{L}}K$.
\mdn
We construct the $R$-matrices 
\dbd
R_{K}:=(\pi_{K}\ot \pi_{K})(R)
\dbd
and
\dbd
R_{L}:=(\pi_{L}\ot \pi_{L})(R)
\dbd
Then we consider the composition
\dbd
L^{*}\xra{\pi_{L}^{*}} (A//L)^{*}\xra{\phi_{R}|_{(A//L)^{*}}}L
\dbd
Clearly this map is bijective. We will show that this map coincides with the Drinfeld's map 
\dbd
f_{Q_{L}}:L^{*}\ra L
\dbd
where $Q_{L}:=R_{L}^{21}R_{L}$.
Indeed if $g\in L^{*}$ then \dbd
(\phi_{R}\star \pi_{L}^{*})(g)=g(\pi_{L}(R^{2}r^{1}))R^{1}r^{2}
\dbd
On the other hand note that $\phi_{R}(A//L)^{*}=L\subset A$ via $x \hookrightarrow x \ot 1_{L}$. Thus $\pi_{L}(\phi_{R}(f))=\phi_{R}(f)$ for any $f \in (A//L)^{*}$. In particular for $f=\pi_{L}^{*}(g)$ we obtain $\pi_{L}(\phi_{R}((\pi_{L}^{*}(g)))=\phi_{R}((\pi_{L}^{*}(g))$ for any $g \in L^{*}$.
\epf
\br One has that $A//K$ and $A//L$ are factorizable Hopf algebras with the quotient $R$-matrices since their associated  fusion subcategoryories are nondegenerate. Then the Hopf isomorphism above gives also that the Hopf subalgebras $K$ and $L$ are also factorizable Hopf algebras with the $R$-matrices transposed by these isomorphism.
\er
\subsection{On the commutativity for normal Hopf subalgebras}\label{comm}
Let $(A, R)$ be a finite dimensional quasitriangular Hopf algebra. Recall that there are two Hopf algebra maps $f_{R}:A^{*cop}\ra A$ and $f_{R_{21}}:A^{*}\ra A^{op}$ defined by:
$$p \mapsto p(R^{1})R^{2}\;\text{and}\;p \mapsto p(R^{2})R^{1}$$ for all $p \in A^{*}$.
\mdn
For a Hopf algebra map $\pi:A\ra B$ define as usually:
$$A^{co B}=A^{co \pi}:=\{a \in A\;|\; a_{1}\ot \pi(a_{2})=a\ot 1\}$$ and
$$\;^{co B}A=\;^{co \pi}A:=\{a \in A\;|\; a_{2}\ot \pi(a_{1})=a\ot 1\}.$$

\bp Let $q:A\ra B$ be a surjective Hopf algebra map and suppose that $(A, R)$ is quasitriangular. Then
\beq
a\lh p=ad'_{F_{R}(p)}a
\eeq
for all $p \in B^{*}$ and $a \in \;^{co B}A^{}$.
where $ad'_{x}(a)=x_{2}aS^{-1}(x_{1})$.
\ep

\bpf
Note that with the above notations one has
$
a \lh p=p(1)a,
$
since $a \in ^{co B}A=\;^{B^{*}}A$. On the other hand since $R^{-1}=R^{1}\ot S^{-1}(R^{2})$ it follows that
\beqn
a_{2}\ot a_{1}=\Delta^{cop}(a)=R\Delta(a)R^{-1}=
(R^{1}a_{1}r^{1})\ot (R^{2}a_{2}S^{-1}(r^{2}))
\eeqn
Thus one has that
$$ad'_{F_{R}(p)}a=(F_{R}(p))_{2}aS^{-1}((F_{R}(p))_{1})=F_{R}(p_{1})aS^{-1}F_{R}(p_{2})=$$ $$=F_{R}(p_{1})(a\lh p_{2})S^{-1}F_{R}(p_{3})=p_{1}(R^{1})R^{2}a_{2}p_{2}(a_{1})p_{r}(r^{1})S^{-1}(r_{2})=$$ $$=p(R^{1}a_{1}r^{1})R^{2}a_{2}S^{-1}(r^{2})=p(a_{2})a_{1}=a \lh p,$$ for any $p \in  B^{*}$.
\epf
\bt Let $q:A\ra B$ be a surjective Hopf algebra map and suppose that $(A, R)$ is quasitriangular. Then the following assertions are equivalent:
\bne
\item $q:A\ra B$ is a normal map.
\item $f_{R_{21}}(B^{*})\subseteq (A^{co \pi})'$
\item $f_{R_{}}(B^{*})\subseteq (\;^{co \pi}A)'$
\ene
 where for any subalgebra $S\subseteq A$ we denote by $S'$ the centralizer of $S$ in $A$.
\et
\bpf The equivalence between $i)$ and $iii)$ is contained \cite[Proposition 3.10]{rqts}. The equivalence between $i)$ and $ii)$ follows by a similar proof which we will include below for completeness. 

$i)\implies ii)$ If $q$ is normal then $A^{co \pi}=\;^{co \pi}A$ and thus
\beqn
p(1)a=a \lh p=ad'_{F_{R}(p)}a
\eeqn
Thus $xa=(x_{3}aS^{-1}x_{2})x_{1}=ax$ for all $x \in X:=F_{R}(B^{*})$ and $a \in A^{co \pi}$.

$ii)\implies i)$ Suppose now that $X\subseteq (A^{co \pi})'$. Then  $x_2aS^{-1}x_{1}=ax_{2}S^{-1}x_{1}=\eps(x)a$ which shows that
$a\lh p=p(1)a$ for any $p\in B^{*}$ and all $a \in A^{co \pi}$ Thus $A^{co \pi}\subseteq A^{B^{*}}=\;^{co B}A=\;^{co \pi}A=$.
\epf
\bt\label{commphi} Suppose that $(A,R)$ is a semisimple quasitriangular factorizable Hopf algebra and $L$ is a normal Hopf subalgebra of $A$. Then 
\bne
\item The image $K:=\phi_{R}((A//L)^{*})$ is also a normal Hopf subalgebra of $A$.
\item One has $[K, L]=1$ which means that $bl=lb$ for any $l \in L$ and $b \in K$.
\ene
\et
\bpf
The proof of the first item is contained in Theorem \ref{nrmphi}. For the second item note that $\phi_{R}=f_{R_{21}}\star f_{R}$.
Since $\pi:A\ra A//L$ is normal it follows that for any $b=\phi_{R}(f)$ with $f \in (A//L)^{*}$ one has that
$b=f_{R_{21}}(f_{1})f_{R}(f_{2})$. Thus from the previous theorem $b$ commutes with all elements from $L$
\epf
\bibliographystyle{amsplain}
\bibliography{mac-bob}
\newpage \br Note also that $
\;_{R}\phi=f_{R}\star f_{R_{21}}.
$
\er

\section{Example: the case of a Drinfeld double}\label{dda}
\blue{\bf write te definition of the drinfeld double}
\subsection{\blue{Formulae for Drinfeld maps}}

For the Drinfeld double $D(A)$ of a semisimple Hopf algebra $A$ one has
$
R=\sum_{i,}(\eps \bwt b_{i})\ot (b_{i}^{*}\bwt 1)
$
and then
$
Q=R_{21}R=\sum_{i,j}(b_{i}^{*}\bwt b_{j})\ot (b_{i}{b_{j}}^{*}).
$
Thus in this case the Drinfeld map $\;_{R}\phi:D(A)^{*}\ra D(A)$ is given by
$
F\mapsto F(Q^{2})Q_{1}=\sum_{i,j}F(b_{i}{b_{j}}^{*})(b_{i}^{*}\bwt b_{j})$ which can also be written as
\dbd
F\mapsto \sum_{i,j}F(b_{j}^{*}(S{b_{i}}_{3}?b_{i}{_{1}})\bwt b_{i}{_{2}})(b_{i}^{*}\bwt b_{j})
\dbd
On the other hand, the other Drinfeld map $\phi\;_{R}:D(A)^{*}\ra D(A)$ is given by
\dbd
\phi\;_{R}:D(A)^{*}\ra D(A),\;\;F\mapsto F(Q^{1})Q_{2}=\sum_{i,j}F(b_{i}^{*}\bwt b_{j})(b_{i}{b_{j}}^{*}).
\dbd
\blue{\bf explain the identification}
\mdn
 In particular if $F=f \ot a\in D(A)^{*}$  then one has that 
\dbd
\phi_{R}(f\ot a)= \sum_{i}b_{i}^{*}(a)f(b_{j})(b_{i}{b_{j}}^{*})=fa
\dbd
and
$
\;_{R}\phi(f \ot a)= \sum_{i}b_{i}^{*}\bwt S{b_{i}}_{3}af(b_{i}{_{2}})b_{i}{_{1}}.
$
\subsection{\blue{On a perp subspace }}
For a subspace $C$ of $A$ we define
$
C^{{\perp\;l}_{A}}:=(A//C)_{l}^{*}=\{f \in A^{*}\;|\; f(ac)=\eps(c)f(a)\;\;\text{for all}\;\;a \in A\}.
$
\blue{Note that $C^{{\perp\;l}_{A}}:=(A//AC^{+})^{*}$. If $C$ is a coideal subalgebra then 
\beqn
\dim (A//AC^{+})^{*}=\frac{\dim A}{\dim C}
\eeqn
}
Define also
$
C^{{\perp\;r}_{A}}:=(A//C)_{r}^{*}=\{f \in A^{*}\;|\; f(ca)=\eps(c)f(a)\;\;\text{for all}\;\;a \in A\} .
$

 It is easy to verify that
$
(S^{-1}C)^{\perp\;r_{A}}=S(C^{\perp\;l_{A}})
$
and
$
(S^{-1}C)^{\perp\;l_{A}}=S(C^{\perp\;r_{A}}).
$
\blue{\bf Note that
\beqn
(D\bwt C)^{\perp}=D^{\perp}\bwt C^{\perp} 
\eeqn
if $C$ and $D$ are unitarian, i.e they contain the units.\mdn Moreover if $C$ is a subcoalgebra then $C^{\perp}$ is a subalgebra of $A^{*}$.}
\bp Let $C$ be a linear subspace of $A$. With the above notations one has  the following: 
$
C^{\perp\;l_{D(A)}}=C^{\perp\;l_{A}}\otimes A.
$
\ep
\bpf
It is easy to see that $ C^{\perp_{A}\;l}\otimes A\subseteq C^{\perp_{D(A)}\;l}$. \onh\; suppose that $\sum_{i}f_{i}\ot a_{i}\in C^{\perp \;l_{D(A)}}$ with $\{a_{i}\}_{i}$ a linear  basis for $A$. It follows that:
\dbd
\sum_{i}f_{i}(bc)g(a_{i})=\eps(c)f_{i}(b)g(a_{i})
\dbd
for all $b  \in A$, $c \in C$ and $g \in A^{*}$.  Since $\{a_{i}\}_{i}$ is a  basis and $g\in A^*$ is arbitrary it follows
$
\sum_{i}f_{i}(bc)a_{i}=\eps(c)f_{i}(b)a_{i}
$ and therefore $f_{i}\in C^{\perp_{A}\;l}$.
\epf

\subsubsection{Dual situation} Note that if $D\subseteq A^{*}$ then since $A^{**}\simeq A$ one has that
$
D^{{\perp\;l}_{A^{*}}}:=\{a \in A\;|\; fg(a)=f(1)g(a)\;\text{for all}\;\;g \in A^{*}\}
$
Note also
$
D^{{\perp\;r}_{A^{*}}}:=\{a \in A\;|\; gf(a)=f(1)g(a)
\;\;\text{for all}\;g \in A^{*}\}
$
\bp
Let $D$ be a subspace of $A^{*}$. Then
\dbd
D^{\perp\;r_{D(A)}}=A^{*}\otimes D^{\perp
\;r_{A^{*}}}
\dbd
\ep
\bpf
The proof is analogue to the previous result.
\epf
\mdn
 Note that this implies:
$
S(D^{\perp\;l_{D(A)}})=(S^{-1}D)^{\perp r_{D(A)}}=(A^{*\cop}\bwt (S^{-1}D)^{\perp\;l_{A^{*}}})
$ and therefore
$
D^{\perp\;l_{D(A)}}=S^{-1}(A^{*\cop}\ot (S^{-1}D)^{\perp\;l_{A^{*}}})=S^{-1}( (S^{-1}D)^{\perp\;l_{A^{*}}})A^{*\cop}=A^{*\cop}\ot S^{-1}( (S^{-1}D)^{\perp\;l_{A^{*}}})=A^{*\cop}\ot D^{\perp\;r_{A^{*}}}
$ Note that the antipode commutes since the identification isomorphism $D(A)^{*}\simeq A^{*}\ot A^{\op}$ is as $k$-algebras. 
\subsection{Some conjectural formulae}
\beqn
(D(A)//(A//K)^{*})^{*}=A^{*}\ot K
\eeqn
and then
\beqn
\phi_{R}((D(A)//(A//K)^{*})^{*})=\phi_{R}(A^{*}\ot K)=KA^{*}=A^{*}\bwt K
\eeqn
\blue{\bf A left normal coideal subalgebra of $D(A)$ containing $A^{*}$ is of the type $A^{*}\bwt L$.}
\blue{\bf Then one needs to show onl;y that $K$ is included.}
\subsection{\blue{\bf On the image of $\phi_{R}$-some other conjectural level}}
We will prove the following:
\red{\bp Let $D$ be a subcoalgebra of $A^{*}$ and $C$ a subcoalgebra of $A$ then
\dbd
\phi_{R}((D\bwt {C})^{\perp\;l_{D(A)}})=C^{\perp\;l_{A}}\bwt D^{\perp \;l_{A^{*}}}
\dbd
\ep
}
\bpf Note that  $C^{\perp\;l_{A}}\bwt D^{\perp \;l_{A^{*}}}=(C^{\perp\;l_{A}}\bwt A)\cap (A^{*\cop}\bwt D^{\perp \;l_{A^{*}}}).$
\mdn 
\blue{\Small  Probably one needs to show separately  that 
\beqn
\phi_{R}(D^{\perp\;l_{A}})=(D\bwt A)
\eeqn
and analogue one and then intersect.
}
 On the other hand one has
\dbd
\sum_{i}f_{i}\bwt a_{i}\in C^{\perp\;l_{A}}\bwt A \iff \sum_{i}f_{i}(ac)a_{i}=\eps(c)(\sum_{i}f_{i}(a) a_{i})
\dbd
if and only if 
\dbd
\sum_{i}f_{i}(ac)a_{i}=\eps(c)(\sum_{i}f_{i}(a) a_{i})
\dbd
it cannot be written without dual basis but nevertheless try to prove the conjecture  by evaluation. \blue{\bf Show one inclusion and then the other one follows from dimension argument.}
\epf

\subsection{Description of normal Hopf subalgebras of $D(A)$}
Let $L$ be a normal Hopf subalgebra of $A$. Recall that an irreducible character $\al$ of $L$ is called $A$-stable if there is a character $\ch \in \mtr{Rep}(A)$ such that 
$\ch\dw^{ ^A}_{ _L}=\frac{\ch(1)}{\al(1)}\al$. Such a character $\ch \in \Irr(A)$ is said to seat over the character $\al \in \Irr(L)$. The set of all irreducible $A$-characters seating over $\al$ is denoted by $\Irr(A|_{\al})$. Denote by
$G_{A}^{st}(K)$ the set of all $A$-stable linear characters of $K$. Clearly $G_{A}^{st}(K)$ is a subgroup of the group of grouplike elements $G(K^*)$ of the dual
Hopf algebra of $K$.
\mdn
Suppose that $L$ and $K$ are two normal Hopf subalgebras of $A$
and let $G$ be a finite group that can be simultaneously embedded in  $G^{st}_A(L)$ and
$ G^{st}_{A^*}((A//K)^*)$ via the embeddings $\psi_1:G \hookrightarrow G^{st}_A(L)$
and respectively $\psi_2:G \hookrightarrow G^{st}_{A^*}((A//K)^*)$. Let $B(K, L,G, \cX, \psi_{1}, \psi_{2})$
be the subcoalgebra of $D(A)$ defined by \beq B(L, K, G, \cX, \psi_{1}, \psi_{2})=\bigoplus_{x \in
G}C_{\psi_1(x)\uw^{ ^{A}}_{ _{L}}}\bowtie C_{ \psi_2(x)\uw^{ ^{A^*}}_{ _{(A//L)^*}}}.
\eeq 
\mdn Recall by \cite[Theorem 4.2]{normalda} that any normal Hopf subalgebra of $D(A)$ is of the type $B(L,L, \cX, \psi)$ where $L, L$ are normal Hopf subalgebras of $D(A)$ and $\cX$  is a finite group satisfying the above properties, for more details see \cite[Subsection 3.2 ]{normalda}. \mdn 
Let  $\cd(L, L, \cX,\psi)$ be the categoryof representations of the quotient Hopf algebra $D(A)//B(L,L,\cX,\psi)$.
\mdn

It would be an interesting question to determine M\"uger's  centralizer for any normal fusion subcategoryof $\rep(D(A))$. \blue{\bf We now know that the centralizer is also a normal fusion subcategoryory. } Note that the proof of Theorem \ref{centraliz} does not work since the Frobenius-Perron dimensions of the fusion subcategoryories do not match anymore

\subsection{\blue{\bf Description of the coalgebras}}
\bn{defn} Define
$
C_{\ch}=C_{M}=\mtr{Ann}_{A}(M)^{\perp}
$. 
\end{defn}
Thus $
C_{M}=\{f \in A^{*}|f(a)=0\;\text{for all}\;
\; a \in \bann_{A}(M)\;\}
$

\br Suppose that $\psi$ is an $A$-stable linear $L$-module. Then the associated primitive central idempotent
$e_{\psi}$ in $L$ is also a central idempotent in $A$. \green{\bf since $\psi$ is stable.}
\er
\bpf
See the Rieffel relation in terms of idempotents from ART.
\epf
\bl Let $L$ be a  left normal  coideal subalgebra of $A$. Then
\dbd
C_{\psi\uw_{L}^{A}}=\{f \in A^{*}\;|\; f(al)=f(a)\psi(l)\;\text{for all}\;a \in A, l \in L\}=e_{\psi}\rh A^{*}
\dbd
\el

\green{\bf do first the proof for $L$ a normal Hopf subalgebra.}
\bpf Let $M:=A\ot_{L}\unpsi$.  It is easy to see that
 \dbd
 l(Sa\ot_{L}1_{\psi})=Sa_{1}\ot_{L}(a_{2}lS(a_{3})1_{\psi})=\eps(l)(Sa\otL 1_{\psi}) 
 \dbd
 Thus 
 \dbd
C_{M}\subseteq  e_{\psi}\rh A^{*}
 \dbd
 Conversely if $a \in \bann_{A}(M)$ then $a\ot_{L}\unpsi=0$. But $A\ot_{L}\unpsi\simeq A\epsi$ via $a\ot_{L}1\simeq a\epsi$. It follows that $a\epsi=0$ and clearly $(\epsi \rh f)(a)=f(a\epsi)=0$.
\epf

\bp Dually one has the following:
\dbd
C_{\psi(\al)\uw^{A^{*}}_{(A//L)^{*}}}=\{a \in A\;|\;\langle fg, a\rangle =\langle g,a\rangle f(\psi(\al))\;\text{for all}\; f \in (A//L)^{*}, g \in A^{*}\}
\dbd
\ep
\subsection{\blue{\bf A conjectural formula for the Mueger centralizer}}

\bt
One has that 
\dbd
\cd(K,L,\mtc X, \psi)'=\cd(\tilde L, \tilde K, \tilde{\mtc X}, \tilde \psi)
\dbd
where $\tl\supseteq L$, $\tilde K\subseteq K$. Moreover
\dbd
\tilde \al\dw^{\tilde L}_{L}=\eps_{L}
\dbd
\dbd
\tpsi(\tilde \al)\dw^{(A//L)^{*}}_{(A//\tL)^{*}}=\eps_{(A//\tL)^{*}}.
\dbd
\blue{\bf for any $\al$??.}
\et
\bpf
One has that $\cd(L,L)\supseteq \cd(L,L, \cx, \psi)$. Thus $\cd(L,L)'=\cd(L,L)\subseteq \cd(L,L,\mtc X, \psi)'=\cd(\tl, \tL, \tilde{\mtc X}, \tilde \psi)$. This implies that $B(L,L)\supseteq B(\tl, \tL, \tilde{\mtc X}, \tilde \psi)$ and all the results follows.
\epf
\blue{\bf conjecture Prove that the description gives te largest and the smallest $\tl$ and $\tilde K$ respectively satisfying the above conditions}
\dbd
\tL=\cap_{\al \in \cx}\hker_{A}(\al)
\dbd
\dbd
\tl=\pi^{-1}(\psi(\al)L\;|\;\al \in \cx)
\dbd
and
\dbd
\tilde{\cx}=
\dbd
\mdn
\blue{\bf State exactly where to apply this? Ce fac cu formula daca o gasesc?}
\mdn
\blue{\bf Is it possible to rewrite symmetrically as in NNW?}
\blue{\bf Apply Theorem \ref{commphi} to Drinfeld doubles:
\beqn
[B(K,L, \mtc X, \psi), B(L^{*}, K^{*}, \mtc X^{*}, \psi^{*})]=1
\eeqn
which gives some commutativity on components and also in cruce.
\mdn
Need to deduce that $K^{*}\subseteq L$ and $L^{*}\supseteq K$.
\mdn In particular apply it to $\rep(D(G)$.
\beqn
[B(H, K, \lam), B(K,H, \lam^{op}]=1.
\eeqn}
\blue{\Small  Apply Theorem \ref{commphi} in particular for $B(L,K)$ commutes to $B(K,L)$.}
\mdn
\blue{\bf The results from mathz give some formulae}
\mdn
\red{\bf Need a conjugacy  class in normal
\dbd
\sum_{i}b_{i}^{*}(xl) \langle g, S{b_{i}}_{3}af(b_{i}{_{2}})b_{i}{_{1}}\rangle =
\dbd
Evaluating $\phi_{R}(f \ot a)$ by $xl\ot g$ one obtains that
\dbd
 \sum_{i}b_{i}^{*}(xl)g(S{b_{i}}_{3}af(b_{i}{_{2}})b_{i}{_{1}})=\sum_{i}b_{i}^{*}(xl)g_{2}(a)\langle b_{i}, g_{3}fSg_{1}\rangle =\dbd\dbd=\langle g_{3}fSg_{1}, xl\rangle g_{2}(a)=g_{3}(x_{1}{L_{1}})g_{1}(S(x_{3}{L_{2}}))f(x_{2})g_{2}(a)=g(S({L_{2}})S(x_{3})ax_{1}{L_{1}})
\dbd
}
\subsection{Questions on the Drinfeld maps}

\mdn
{\green{\bf Q1 Where $A^{*}\ot 1$ and the other $\eps \ot A$ are sent by the two Drinfeld maps?? in themselves as copies}}
\mdn {\green{\bf Q2
see where chars are sent by $\phi_{R}$ in the center.}
\mdn
try to work with dual basis.

\green{\bf it can be computed and it shows if true or not the conj}

\blue{\bf 
it is also known that the \dd \phi _{R}\;\dd is the twist on the character ring!!! see sommerhauser}


\mdn
for the map $\phi_{R}$ use dual basis; where $f \ot 1$ and the other one is sent; first seems something with adj action for  dual bases.
\newpage
\section{ Application-on a fusion subcategoryof $\rep(D(A)$}\label{onspec}
\subsubsection{The fusion subcategorygenerated by left normal coideals}

Let $A$ be a semisimple Hopf algebra and $D(A)$ be its Drinfeld double. Then as above $A$ is a $D(A)$-module via
$(f \bwt a).b=(a_1bS(a_2))\lh S^{-1}f$ for all $f \in A^*$ and $a, b\in A$.

The simple $D(A)$-submodules of $A$ are the minimal left normal  coideals of $A$. We will use Brauer's theorem to find the fusion subcategoryof $\mtr{Rep}(D(A))$  generated by the $D(A)$-module $A$.

\blue{\bf Recall the result from \cite{repda} that $\lker_{A}(A_{ad})=K(A)$.}
\subsection{Normal coideal subalgebras of  $D(A)$}

Let $\{x^d_{ij}\}_{ij}$ be a comatrix basis of a simple subcoalgebra  $C$ of $A$. That is $d \in \Irr(A^*)$ and $\D(x^d_{ij})=\sum_l x^d_{il}\otimes x^d_{lj}.$

Then it is easy to check the following remark:
\bn{rem} Let L be a left coideal subalgebra of $D(A)$. Suppose that $l=\sum_{dij}F_{dij}\ot x^d_{ij}\in L$ then $A\rh F_{dij}\ot x^d_{Lj} \in L$ for any $d,i,j,L$. In particular $F_{dij}\ot x^d_{ij}\in L$ for any $d,i,j$.
\end{rem}
\bt\label{main}
The left kernel of the $D(A)$-module $A$ is the Hopf center $K(A)$ of $A$. Thus the fusion subcategory$\langle A\rangle $ of $\rep(D(A))$ generated by $A$ coincides with $\rep(D(A)//K(A))$.
\et

\bn{proof}
It is enough to show that $\mtr{LKer}_{D(A)}(A)=L(A)$ and then one can apply Brauer's theorem mentioned above.  

Let $L:=\mtr{LKer}_{D(A)}(A)$. It is easy to see that $K(A)\subseteq L$. To finish the proof one has to show that $L \subseteq K(A)$.

Using the above remark it is enough to show that if $f \bowtie a \in L$ then $f=\eps$ and $a \in K(A)$. By the definition of the $D(A)$-module $A$ one has that $f \bowtie a \in \lker_{D(A)}(A)$ if and only if
$(f_2\bowtie a_1)\ot (f_1\bwt a_2).b=(f \bwt a)\ot b$ for all $b\in A$. This also can be written as $(f_2\bowtie a_1)\ot  a_2bS(a_3)\lh Sf_1=(f \bwt a)\ot b$ 
or 
\bn{equation}\label{dfLr}
(f_2\bowtie a_1)\ot  a_3b_2S(a_4)f_1(a_5S(b_1)S(a_2))=(f \bwt a)\ot b
\end{equation}

Let $b=\Lam$ the integral of $A$ in the above equation and note that

\bn{equation}
\sum_{}a\Lam_1b\ot \Lam_2=\sum \Lam_1 \ot S(a)\Lam_2S(b).
\end{equation}
for all $a, b\in A$.
Using the above identity Equation \eqref{dfLr} becomes $$f_2f_1(S(\Lam_1))\ot \Lam_2=f \ot \Lam .$$ 

Since $\{S(\Lam_1),\; \Lam_2\}$ are dual Frobenius bases for $A$ this implies that $f=f(1)\eps$. Then the Equation \eqref{dfLr} becomes
$$a_1\ot a_2bS(a_3)=a \ot b$$
 for all $b \in A$.

By the previous lemma this shows that $a \in \lker_A(A_{ad})= K(A)$.
\end{proof}

\bc
Suppose that $C(A)$ commutative, then $\Rep(D(A)//K(A))$ contains $\Rep(D(A))_{\ad}$.
\ec
\bpf
If $C(A)$ commutative then $K(A)=k\bar{G}(A)$ where $\bar{G}(A)$ is the set of all central grouplike elements of $A$. Thus $K(A)\subset K(D(A))=k\bar{G}(D(A))$.
\epf
\subsection{Centralizer rezults}
By theorem ?? \cite{mathz} one has that
\beq
\cd(A)'=\langle A\rangle 
\eeq
and therefore $\langle A\rangle '=\cd(A)$. On the other hand 
\beqn
\langle A\rangle '=\rep(D(A)//K(A))'=\cd(A, K(A))'=
\eeqn
\beqn
\cd(K(A), A)=\rep(D(A)//(A//K(A))^{*}\bwt A)
\eeqn
\bp With the following notations it follows that
\beq
N_{D}(A)=(A//K(A))^{*}\bwt A
\eeq
\ep
\bpf One has that $\cd(A)=\rep(D(A)/N_{D}(A)$.
\epf
Applying the first main result it follows that
\beq
\phi_{R}(D(A)//K(A))=(A//K(A))^{*}\bwt K(A)
\eeq
\newpage
\section{\blue{\bf A conjecture and possible generalization}}
\beq
N_{D}(L)=(A//K_{1}(L))^{*}\bwt L
\eeq
where 

\subsection{On $\lker_{D(A)}(L)$}
\bl Under the adjoint action
\dbd
\lker_{A}(L)=C_{A}(L)
\dbd
largest (normal) left coideal subalgebra that commutes with $S(L)$.
It follows that
\dbd
(A//L)^{*}\bwt C_{A}(L)\subseteq \lker_{D(A)}(L)
\dbd
\el
\subsection{Rewrite theorem 1.2 from mathz}
\beqn
\cd(L)'=\langle L\rangle 
\eeqn
can be written as $$\phi_{R}(D(A)//N_{D}(L))=\lker_{D(A)}(L)$$
\dbd
\rep(HLer_{A^{*}}(d))'=\rep(A//L_{d})'=\{\ch_{j}\;|E_{j}(d)\neq 0\}=\{\ch_{j}\;|\; E_{j}(\blam_{L_{d}})\neq 0\}
\dbd
It follows that
\dbd
L_{d}=\bigoplus_{j\;|\cc_{j}\cap C_{d}\neq 0}\cc_{j}
\dbd
\blue{\bf \small Rewrite Theorem 1.4 for those for $\cd(L)$ with $L$ not normal}
\newpage
\blue{\Small  Use the theory of byrpoducts if possible to overcome Hopf modules.}\newpage
\section{Example: The case of a group algebra $kG$ }\label{grpalg}
\blue{\bf They are almost no factorizable just in the case of an abelian group}
$$
R=\sum_{a,b \in \hat{\Gamma}}\rho(a, b) e_{a}\ot e_{b}
$$
\newpage
\section{Example: the Drinfeld double  $D(kG)$ of a group }\label{ddg}
In this section we consider the case $A=kG$, for some finite group $G$. One has the following formula for the $R$-matrix of $D(A)$:
$
R=\sum_{g}(\eps \bwt g)\ot (p_{g}\bwt 1).
$
Thus
\beqn
Q=\sum_{g,h\in G}(p_{g}\bwt h)\ot gp_{h}=\sum_{g,h \in G}(p_{g}\bwt h)\ot (p_{ghg^{-1}}\bwt g). 
\eeqn
After the identification $D(kG)^{*}\simeq k^{G}\otimes kG^{op}$ as algebras, the Drinfeld maps are given by
\dbd
_{R}\phi:D(G)^{*}\ra D(G),\;\;x\ot p_{g}\mapsto (x\ot p_{g})(Q^{1})Q^{2}=p_{xgx^{-1}}\bwt x=xp_{g}
\dbd
and 
\dbd
\;\phi_{R}:D(G)^{*}\ra D(G),\;\;x\ot p_{g}\mapsto (x\ot p_{g})(Q^{2})Q^{1}=p_{g}\bwt g^{-1}xg.
\dbd
\subsection{Presentation of normal coideal subalgebras of $D(kG)$}
\red{Left coideal subalgebras are determined by $\D(L)\subset A\ot L$ and right coideal subalgebra by $\D(R)\subseteq R\ot A$.}\mdn\red{One has the formulae
\beqn
\D_{A^{*}}(f)=\sum_{i}e_{i}^{*}\ot(f \lh e_{i})=\sum_{i}(e_{i}\rh f)\ot e_{i}^{*}
\eeqn
\mdn
Thus $L$ is left coideal subalgebra if and only if $(f \lh e_{i})\in L$ and $R$ is right coidela subalgebra if and only if $(e_{i}\rh f)\in R$.}
\mdn
Let $H$ and $L$ be two normal subgroups  of $G$ such that they commute point-wise, i.e $[H,L]=1$. Let also $\lam:M\times H\ra k^{*}$ be a  twisted bicharacter, i.e. a function satisfying the following two properties
\beq\label{tw1}\lam(mm', h)=\lam(m,h)\lam(m', h)
\eeq
\beq\label{tw2}\lam(m, hh')=\lam(m,h)\lam(h^{-1}mh, h').\eeq
Following \cite{jpaa} we define the following subspace of $D(kG)$
\dbd
C(M,H,\lam):=\bigoplus_{h \in H}C_{\lam}^{M}(h)\bwt h
\dbd
where
$
C_{\lam}^{M}(h)=\{f \in L^{G}\;|\;f\lh m=\lam(m,h)f\;\text{for all}\; m \in M\}=\{f \in L^{G}\;|\; f(mx)=\lam(m,h)f(x)\;\text{for all}\; m \in M, x \in G\}.
$\mdn
\red{By above $C_{\lam}^{M}(h)$ is a right coideal subalgebra of $k^{G}$.
}
\mdn
These are all left  coideal subalgebras of $D(kG)$  and they are normal if and only if the (twisted) bicharacter is $G$-invariant, i.e.
\beq
\lam(x^{-1}mx, h)=\lam(m, xhx^{-1})
\eeq
\mdn 
A basis for $C(M,H,\lam)$ is given by 
\beqn
f_{g}^{h}=\sum_{m \in M}\lam(m,h)p_{gm}
\eeqn
Thus
$$
\dim_{k}C(M,H,\lam)=|H||G:M|
$$
\mdn
\blue{\bf Record also the formula:
\beqn
f^{h}_{s}(g^{-1}?g)=\sum_{m \in M}\lam(g^{-1}mg, h)p_{gsg^{-1}m}=\sum_{m \in M}\lam(m, ghg^{-1})p_{gsg^{-1}m}=f^{ghg^{-1}}_{gsg^{-1}}
\eeqn
}
\mdn
By  \cite[Theorem 4.1]{jpaa} a left normal  coideal  subalgebra $C(M,H,\lam)$ is in fact a normal Hopf subalgebra if and only if the bicharacter $\lam$ is strongly $G$-invariant; i.e invariant on each side.:
\beq
\lam(xmx^{-1}, h)=\lam(m, xhx^{-1})=\lam(m,h)
\eeq
for all $m \in L$ and $h \in H$.
\bp
One can check the following multiplication formulae inside $C(M,H, \lam)$:
\beqn
(f^{(h)}_{s}\bwt h)(f^{(h')}_{t}\bwt h')=\delta_{Ms, \;\;Mhth^{-1}L}\lam(s^{-1}hth^{-1}, h')f^{(hh')}_{s}\bwt hh'
\eeqn
\ep
\blue{\bf neeed a def for $f^{h}_{s}$ and that in the same classis a scalar of the initial.
}
\bpf One has that
\beqn
(f^{(h)}_{s}\bwt h)(f^{(h')}_{t}\bwt h')=\sum_{m,m'\in M}(p_{ms}\bwt h)(p_{m't}\bwt h')\lam(m, h)\lam(m',h')=$$
$$=\sum_{m,m'\in M}(p_{ms}p_{hm'th^{-1}}\bwt hh')\lam(m, h)\lam(m',h')=\delta_{Ms, \;\;Mhth^{-1}}\sum_{m \in M}\lam(m,h)\lam(mm_{0}, h')
\eeqn
where $m_{0}$ is determined by $s=m_{0}hth^{-1}$.
\epf
\subsubsection{Integrals}
\bl One has that
\dbd
\blam_{L(M,H,\lam)}=p_{1}\bwt \frac{1}{|H|}\sum_{h \in H}h
\dbd
\el
\bpf Suppose that
\dbd
\blam_{L(M,H,\lam)}=\sum_{ h \in H}(f^{h}\bwt h)
\dbd
with $f_{h}\in C_{\lam}(h)$. Then for any $g \in H$ one has that
\dbd
(f^{x}_{[g]}\bwt x)\blam_{L(M,H,\lam)}=\sum_{h  \in H}f^{x}_{[g]}(xf^{h}x^{-1})\bwt xh=\sum_{r=xh  \in H}(f^{x}_{[g]}(xf^{rx^{-1}}x^{-1})\bwt r)
\dbd
On the other hand
\dbd
\delta_{g, M}\lam(g^{-1}, x)\blam_{L(M,H,\lam)}=\delta_{g, M}\lam(g^{-1}, x)\sum_{ h}(f^{h}\bwt h)
\dbd
Thus
\dbd
\sum_{r  \in H}f^{x}_{[g]}(xf^{rx^{-1}}x^{-1})\bwt r=\delta_{g, M}\lam(g^{-1}, x)\sum_{ r \in H}(\delta_{g, M}\lam(g^{-1}, x)f^{r})\bwt r
\dbd
which gives that
\dbd
f^{x}_{[g]}(xf^{rx^{-1}}x^{-1})=\delta_{g, M}\lam(g^{-1}, x)f^{r}
\dbd
for any $r \in H$. Thus implies that $f^{r}=0$ if $r\neq 1$
\epf
\subsection{The irreducible $D(kG)$-modules}
The associated simple $D(kG)$ module is $kG\ot_{kC_{G}(a)}M_{\ch}$. The action of $kG$ is just the usual left action.The action of $k^{G}$ is given by 
\dbd
p_{x}(g\ot_{C_{G}(a)}m)=\delta_{x, gag^{-1}}(g\ot_{C_{G}(a)}m).
\dbd
The irred modules coming from $a=1$ are $\rep(G)$ with $$
p_{x}.m=\delta_{x,1}m.
$$
\subsection{\blue{\bf Recall the parameterization from \cite{nnw}}}
The category$\cs(M, H, \lam)$ consists of the simple objects $S_{(a, \ch)}$ with the property that $a \in M$ and
$
\lam(a,h)=\frac{\ch(h)}{\ch(1)}\;\text{for all}\;h \in H
$
It was shown in \cite{nnw} that
\beqn
\fp(\cs(M,H,\lam))=|M||G:H|
\eeqn
Conversely ginen $\cd$ one obtains $H_{\cd}$ by interswctin with $\rep(G)$ and $K_{\cd}$ can be obtained as the restrcitiomn to $K^{G\;cop}$ Morever $B_{\cd}:K_{\cd}\times H_{\cd}\ra k^{*}$ is gi9ven
\beqn
B_{\cd}(g^{-1}ag, h)=\frac{\ch(ghg^{-1})}{\ch(1)}.
\eeqn
\subsection{Identification of categoryories} 
\bp 
With the above notations one has that:
\beqn
\rep(D(G)//C(M,H,\lam))=\cs(M, H, \lam^{-1})
\eeqn
\ep
\bpf Let $\cd:=\rep(D(G)//C(M,H,\lam))$.
One defines $H_{\cd}$ by the formulae
\beqn
\rep(G/H_{\cd})=\cd\cap \rep(G)
\eeqn
and $K_{\cd}$ is obtained by the restriction to $k^{G \;cop}$. the bilinear bicharacter 
$B_{\cd}:K_{\cd}\ra H_{\cd}\ra k^{*}$ is given
\beqn
B_{\cd}(g^{-1}ag, h):=\frac{\ch(ghg^{-1})}{\ch(1)}
\eeqn
Note that an object $(a,\ch)\in \rep(D(G)//C(M,H,\lam))$ if and only if $C(M,H,\lam)$ acts trivially on $(a,\ch)$. This can be written as
\beqn
(f_{s}^{h}\bwt h)(g\ot_{kC_{G}(a)}v)=f_{s}^{h}(1)(g\ot_{kC_{G}(a)}v)
\eeqn
for any $s,g \in G$ and $v \in M_{\ch}$.This can also be rewritten as
\beqn
(f_{s}^{h}\bwt h)(g\ot_{kC_{G}(a)}v)=\delta_{s, M}\lam(s^{-1}, h)(g\ot_{kC_{G}(a)}v)
\eeqn
Note that the left hand side of the above equation can be written as
\beqn
\sum_{m \in M}\lam(m, h)p_{ms}(hg\ot_{kC_{G}(a)}v)=
\sum_{m \in M}\delta_{ms, \;(hg)a(hg)^{-1}}\lam(m, h)(hg\ot_{kC_{G}(a)}v)\eeqn
Equating the two terms one obtains that
\beqn
\sum_{m \in M}\delta_{ms, \;(hg)a(hg)^{-1}}\lam(m, h)(hg\ot_{kC_{G}(a)}v)=\delta_{s, M}\lam(s^{-1}, h)(g\ot_{kC_{G}(a)}v)\eeqn
If $s \notin M$ this equation implies that $(hg)a(hg)^{-1}\notin Ms$ which can happen if amd only if $a \in M$.\mdn On the other hand if $s \in M$ then the equation can be rewritten as 
\beqn
\lam((hg)a(hg)^{-1}s^{-1}, h)(hg\ot_{kC_{G}(a)}v)=\lam(s^{-1}, h)(g\ot_{kC_{G}(a)}v)\eeqn
Using the properties of $\lam$ this shows that
\beqn
\lam(gag^{-1}, h)(hg\ot_{kC_{G}(a)}v)=(g\ot_{kC_{G}(a)}v)
\eeqn
which implies that $\lam(gag^{-1}, h)(g^{-1}hg)v=v$ for any $g \in G$ and $v \in M_{\ch}$. Thus 
\beqn
\frac{\ch(g^{-1}hg)}{\ch(1)}=\lam(ga^{-1}g^{-1}, h)=\lam(a, ghg^{1-})^{-1}\eeqn
Using the reconstruction of the triple from a fusion categorythis finishes the proof.
\epf
\subsection{Identification with the fusion subcategoryories from \cite{nnw}} One needs to show that
\beqn
\rep(D(G)//C(M,H,\lam))\subseteq \cs(M, H, \lam)
\eeqn
\blue{\bf Then by checking dimension one can identify \beq\rep(D(G)//C(M,H,\lam))=\cs(M, H, B)
\eeq din [NNW].}
\subsubsection{Intersection}
One has that\\
$\rep(D(G)//C(M,H,\lam))\cap \rep(G)=\rep(D(G)//C(M,H,\lam))\cap\rep(D(G)//k^{G})=\rep(D(G)//(k^{G}C(M,H, \lam))$. Note that
$k^{G}C(M,H,\lam)=k^{G}\bwt kH$. Thus $\rep(D(G)//C(M,H,\lam))=\rep(D(G)//(k^{G}\bwt kH))=\rep(G/H)$. Thus
\beqn
H_{\rep(D(G)//C(M,H, \lam)}=H
\eeqn
\blue{\bf Note that $C(1, H, 1)=k^{G}\bwt kH$.}
\subsection{A particular case of computation}
If $N$ is a subgroup then
$$(kN)^{\perp\;l_{kG}}=k^{(G/N)_{l}}=\{\sum_{n \in N}p_{gn}\;g \in G\}
$$
$$
(kN)^{\perp\;l_{D}}=k^{(G/N)_{l}}\ot kG
$$
$$
\phi_{R}((D//kN)_{l}^{*}))=\phi_{R}(k^{(G/N)_{l}}\ot kG)=kG k^{(G/N)_{l}}=k^{(G/N)_{l}}\bwt kG
$$
If $N$ is normal we are able to write the last equality $kG k^{(G/N)_{l}}=k^{(G/N)_{l}}\bwt kG$.
\subsection{Direct computation of perp and $\phi_{R}$-image}
\bp One has that
\beqn
(C(M,H, \lam))^{\perp\;l_{D}}=
\{F:=\sum_{a,b}f_{a,b}p_{a}\ot b\;|\;|\;f_{g,x}=0 \;\text{if}\; x \notin M\;\; f_{g, x}=\lam(x, ghg^{-1})f_{gh, x}\;\text{if }\;x \in M\}
\eeqn
\ep
\bpf One has $F\in (C(M,H, \lam))^{\perp\;l_{D}}$ if and only if
\beq\label{con1}
F((p_{x}\bwt g)(f_{s}^{h}\bwt h))=F(p_{x}\bwt g)f^{h}_{s}(1)
\eeq
for any $x,g, s \in G$ and any $h \in H$. 
\mdn Note that
\beqn
f_{s}^{h}(1)=0 \;\text{if}\; s \notin M
\eeqn
and 
\beqn
f_{s}^{h}(1)=\lam(s^{-1},h) \;\text{if}\; s \in M
\eeqn
Thus the RHS of Equation \eqref{con1} equals
\beqn
\delta_{s, M}\lam(s^{-1}, h)f_{g, x}
\eeqn
For the LHS of the same Equation note that:
$F((p_{x}\bwt g)(f_{s}^{h}\bwt h))=F(p_{x}f_{s}^{h}(g^{-1}?g)\bwt gh)=F(p_{x}f_{gsg^{-1}}^{ghg^{-1}}\bwt gh)$
\mdn On the other hand note that
\beqn
p_{x}f_{gsg^{-1}}^{ghg^{-1}}=\sum_{m \in M}\delta_{x, gsg^{-1}m}\lam(m, ghg^{-1})p_{x}
\eeqn
Then equating the two terms we obtain that:
\beqn
\sum_{m \in M}\delta_{x, gsg^{-1}m}\lam(m, ghg^{-1})f_{gh, x}=\delta_{s, M}\lam(s^{-1}, h)f_{g, x}
\eeqn
{\bf Case 1:} If $s \notin M$ then the left hand-side is zero. then the other side has also to be zero which says that if $x=gsg^{-1}m$ then $f_{gh, x}=0$. This  is equivalent to $f_{g,x}=0$ whenever $x \notin M$.
\mdn {\bf Case 2:} Suppose that $s \in M$. Then the equality becomes
\beqn
\sum_{m \in M}\delta_{x, gsg^{-1}m}\lam(m, ghg^{-1})f_{gh, x}=\lam(s^{-1}, h)f_{g, x}
\eeqn
This shows that  if $f_{g,x}\neq 0$ then also $f_{gh, x}\neq 0$ and $\lam(s^{-1}, h)f_{g,x}=\lam(m, ghg^{-1})f_{gh, x}$ for the unique $m \in M$ with $x=gsg^{-1}m$. Thus $\lam(m, ghg^{-1})=\lam(xgs^{-1}g, ghg^{-1})=\lam(x, ghg^{-1})\lam(gsg^{-1}, ghg^{-1})=\blue{\bf \lam(x, ghg^{-1})\lam(s^{-1}, h)}$. Thus in this case:
\beqn
f_{gh, x}=f_{g,x}\lam(x, ghg^{-1})
\eeqn

Denote by 
\beqn
t_{g}=\sum_{h \in H}f_{gh, x}p_{gh}
\eeqn 

Then $t_{g}\lh h'=\sum_{h \in H}f_{gh, x}p_{gh}\lh h'=\sum_{h \in H}f_{gh, x}p_{gh}(h'?)=\sum_{h \in H}F_{gh,x}p_{h'^{-1}gh}$
\epf
\blue{\bf Note that
\beqn
hf_{s}^{h}=f_{hsh^{-1}}^{h}\bwt h
\eeqn
This shows that $(kh)(C_{\lam}(h))=C_{\lam}(h)\bwt (kh)$
and therefore
\beqn
 \phi_{R}^{{-1}}(C(M,H, \lam))=\oplus_{h \in H} \phi_{R}^{{-1}}((C_{\lam}(h))\bwt (kh))=\oplus_{h \in H} \phi_{R}^{{-1}}((kh)(C_{\lam}(h)))=\oplus_{h \in H}C_{\lam}(h)\ot kh
 \eeqn
 Probably the previous proposition shows that
 \beqn
 C(H, M, \lam)^{\perp\;l_{D}}=\oplus_{h \in H}C_{\lam}(h)\ot kh
\eeqn}
\red{The whole section should be just a remark in the paper!}
\newpage
\subsection{The work on April 25-describing central idempotents and conjugacy classes in $D(G)$}
\bne
\item Following \cite{dn} noik-naidu one has that
\dbd
[L_{a}, L_{b}]=1\;\text{and}\; \ch(gbg^{-1})\ch'(g^{-1}ag)=\ch(1)\ch'(1)
\dbd
\item describing all the characters of $D(G)$
\dbd
\ch=\sum_{x,y}\al_{x,y}p_{x}\ot y
\dbd
need to satisfy 
\dbd
\al_{x,y}=\al_{axa^{-1}, aya^{-1}}\text{and}\;\; \al_{w,z}=0\;\;\text{if}\;\; wz\neq zw
\dbd
\item It follows that
\dbd
\ch_{a,b}=\sum_{y \in G/G_{a}\cap G_{b}}p_{yay^{-1}}\ot yby^{-1}
\dbd
is a virtual character.
\item
\dbd
V_{a,b}=\blam\lh \ch_{a,b}D(G)^{*}=w_{a,b}\lh D(G)^{*} 
\dbd
is a $D(D(G))$ submodule of $D(G)$
\item
It seems that
\dbd
V_{a,b}=\{\sum_{z \in G_{a}/G_{a}\cap G_{b}}p_{yzb(yz)^{-1}\bf{v^{-1}}}\bwt yay^{-1}\;|y \in G\;v \in G\}
\dbd
\item On the other hand I should be able to compute
\dbd
\ch_{a,b}\ch_{m,n}
\dbd
 which is zero if $a$ is not conjugate to $m$
 \item It is easy to see that 
 \dbd
 \ch_{xax^{-1}, xbx^{-1}}=\ch_{a,b}
 \dbd
 \item We denote by
 \dbd
 w_{b}=\sum_{z \in G_{a}/G_{a}\cap G_{b}}xbx^{-1}
 \dbd
 the class sums in $G_{a}$.
 \item Then it is easy to verify that
 \dbd
 \ch_{a,b}=\sum_{y \in G/G_{a}}p_{yay^{-1}}\ot yw_{b}y^{-1}
 \dbd
 \item We need the structure constants from $G_{a}$ in order to write the product formula
 \dbd
 w_{b}w_{n}=\sum_{x\in cl(G_{a})}C(a)^{x}_{b, n}w_{x}
 \dbd
 \item In particular
 \dbd
\ch_{a,b}\ch_{a,n}=\sum_{x\in Cl(G_{a})}C(a)_{n, b}^{x}\ch_{a, x}
\dbd
\item Therefore we have an isomorphism of Grothendieck rings
\dbd
R(D(G))\xra{\simeq}\bigoplus_{a \in Cl(G)}Z(\mathbb CG_{a})
\dbd
given by
\dbd
\ch_{a, b}\mapsto w(a)_{b}
\dbd
\item It recovers a result of Lusztig and Sarah
\item From here it follows that
\dbd
F_{a, \ch_{i}}:=\ch_{a, \ch_{i}}:=\frac{\ch_{i}(1)}{|G_{a}|}\sum_{b\in cl(G_{a})}\ch_{i}(b^{-1})\ch_{a, b}
\dbd
are the central idempotents of the Grothendieck ring.
\item The Frobenius map $\mtc F$ sends  $f \mapsto \blam \lh f$ is given in this case by
\dbd
p_{x}\ot y\mapsto p_{y^{-1}}\bwt x
\dbd
\item The unnormalized conjugacy classes satisfy
\dbd
C_{(a,\ch_{i})}=\mtc F(\frac{\ch_{i}(1)}{|G_{a}|}(\sum_{b\in cl(G_{a})}\sum_{y \in G/G_{a}}p_{yay^{-1}}\ot yw_{b}y^{-1}))=\dbd\dbd=\sum_{b\in cl(G_{a})}\sum_{y \in G/G_{a}}(\sum_{z\in G_{a}/G_{a}\cap G_{b}}p_{zb^{-1}z})\bwt yay^{-1}
\dbd
\blue{\small change the order summation; the fraction should be there}
\item \green{\bf Need to compute
\dbd
C_{(a,\ch_{i})}C_{(m,\ch_{j})}=
\dbd
}
\item \green{\bf verify that $\eps(C_{j})=m_{j}=1$}
\item \green{\bf compute also the entire conjugacy class.}
\item the characters of the irreducible modules $\widehat(a, \ch_{i})$ are computed in cejm
\item From here I can write the central idempotents of $Z(D(G))$.
\item We have again in this case
\dbd
R=\sum_{g}(\eps \bwt g)\ot (p_{g}\bwt 1)
\dbd
adn
\dbd
Q=\sum_{g,h}(p_{g}\bwt h)\ot gp_{h}=\sum_{g,h}(p_{g}\bwt h)\ot (p_{ghg^{-1}}\bwt g)
\dbd
\dbd
_{R}\phi(p)=P(Q^{1})Q^{2}
\dbd

\dbd
_{R}\phi:D(G)^{*}\ra D(G),\;\;x\ot q_{g}\mapsto p_{xgx^{-1}}\bwt x=xp_{g}
\dbd
and 

\item The first Drinfeld map is given by
 \dbd
\;\phi_{R}:D(G)^{*}\ra D(G),\;\;x\ot p_{g}\mapsto p_{g}\bwt g^{-1}xg
\dbd
\item On the character ring they are just swipe
 \dbd
\;\phi_{R}:C(D(G))\ra Z(D(G)),\;\;\sum_{x, g}\al_{x,g}(x\ot p_{g})\mapsto \sum_{x, g}\al_{x,g}(p_{g}\bwt x)
\dbd
\item
Thus
\dbd
\phi_{R}(\ch_{a,b})=\sum_{y \in G/G_{a}}p_{yay^{-1}}\bwt w_{b}
\dbd
and
\dbd
\phi_{R}(F_{a, \ch_{i}})=\frac{\ch_{i}(1)}{|G_{a}|}\sum_{b\in cl(G_{a})}\ch_{i}(b^{-1})(\sum_{y \in G/G_{a}}p_{yay^{-1}}\bwt w_{b})
\dbd
\item the formula from Lemma 4.7 from CEJM, uses Lemma 3.1
\item Probably 
\dbd
\phi_{R}(F_{a,\ch_{i}})=E_{a, \ch_{i}^{*}}
\dbd
\ene

\newpage
\green{\bf I need to write a formula for the intersection of two such coideal subalgebras}
\dbd
L(L,H, \lam)\cap L(H,L, \lam)
\dbd

\subsubsection{NNW formulae and description of the corresponding coideal subalgebra}
\dbd
\cs(L,H, \lam)'=\cs(H, L, \lam^{\mtr{op}^{-1}})
\dbd

in their definition $\lam$ is a $G$-invariant bicharacter. If $H$ and $M$ commutes mine is also a bicharacter.
\mdn
\green{\bf I verified that \dbd
\cs(L,H, \lam)=\rep(D(G)//L(L,H,\lam))
\dbd
write it as $B(L,L, \mtc X, \psi)$.}

\bl\label{claim1}
\dbd
B(M,L,\cx, \psi)=C(M, H, \lam)
\dbd
where
\dbd
H:=\pi^{-1}(\psi(\al)l\;|\;\al \in \cx, l\in L)\supseteq L
\dbd
and
\dbd
\lam: L\times H\ra L^{*}\; (m, \psi(\al)l)\mapsto \al(m)
\dbd
\el
\bl\label{claim2}
\dbd
C(L, H, \lam)=B(L,H_{1},\cx, \psi)
\dbd
where
\dbd
H_{1}:=\{h \in H\;|\; \lam(-,h)=\eps\; \text{on } M\}
\dbd
and
\dbd
\cx=\{\lam(-,h)\;|\;h \in H\},\; \cx^{*}=\{\hat h\;|\; h \in H\}
\dbd
and 
\dbd
\psi:\mtc X\ra \cx^{*},\;\lam(-,h)\mapsto \hat h
\dbd
is the isomorphism of groups of stable characters.
\el
\bt\label{fn} In the case $A=D(kG)$
One has that
\dbd
B(L,L, \cx, \psi)=B(\tl, \tL, \tmx, \tpsi)
\dbd
where
\dbd
\tl=\pi^{-1}(\psi(\al)l\;|\;\al \in \cx, l\in L)\supseteq L
\dbd
\dbd
\tL=\{m \in L\;|\; \al(m)=1\;\text{for all}\; \al \in \cx
\dbd
and
\dbd
\tilde{\cx}=\{\psi_{m}\;|\;m \in L\}\; \;\text{where}\;\;\psi_{m}: \psi(\al)l\mapsto \lam^{-1}(-,m)
\dbd
and
\dbd
\tilde{\cx}^{*}=\{\hat m\;|\; m\in L\}.
\dbd
\et
\blue{\bf Look at the coefficients of the integral use the structures constants for factorizable Hopf algebras.}
\subsection{On the structure constants}
For factorizable Hopf algebras one has that
\dbd
c^{i}_{jt}=\frac{d_{i}d_{j}}{d_{t}}m^{i}_{jt}
\dbd
\bl Let $N$ be a normal subgroup of a finite group $G$. Then
\dbd
\lker_{D(G)}(kN)=C(N, C_{G}(N), 1)
\dbd
\el
\bpf
First note that $C_{G}(N)$ is a normal subgroup. Indeed
\dbd
(hgh^{-1})n=hg(h^{-1}nh)h^{-1}=h(h^{-1}nh)gh^{-1}
\dbd
The action of $D(G)$ on $N$ is given by
\dbd
(p_{g}\bwt x).n=(xnx^{-1})\lh S^{-1}p_{g}=\delta_{g^{-1}, xn^{-1}x }xnx^{-1}
\dbd
Suppose that 
\dbd
\lker_{D}(N)=C(L,H, \lam)
\dbd
Then $f_{[g]}\bwt h$ acts as a scalar on $LN$, i.e
\dbd
f_{[g]}(hn^{-1}h^{-1})hnh^{-1}=f_{[g]}(1)n
\dbd
for any $g \in G$, $h \in H$ and $n \in N$.
If $g \notin L$ then $f_{[g]}(1)=0$ and therefore $f_{[g]}(hn^{-1}h^{-1})=0.$ Thus if $g \notin L$ then $hnh^{-1}\notin Lg$. This shows that $N\subseteq L$.
\mdn Suppose that $g \in L$ then
\dbd
f_{[g]}=\sum_{m \in L}p_{mg}\lam(m,h)=\sum_{r \in L}p_{r}\lam(rg^{-1}, h)
\dbd
Thus $hn=nh$ i.e $H\subseteq C_{G}(N)$ and
\dbd
\lam(rg^{-1}, h)=\lam(g^{-1}, h)
\dbd
i.e $\lam(r,h)=1$ fro any $r\in L$. The largest such coideal subalgebra is
\dbd
C(N, C_{G}(N), 1).
\dbd
\epf
\subsection{\blue{\bf On the commutativity conjecture}}

\red{\bf \small It seems that $C(N, C_{G}(N), 1)$ does not always commute with its prime $C(C_{G}(N), N,1)$.
One has 
\dbd
\langle LN\rangle =\rep(N, C_{G}(N))
\dbd
It follows that
\dbd
\langle LN\rangle '=\cd(C_{G}(N), N)
\dbd
The question is if
\dbd
\cd(C_{G}(N), N)=\cd(N)
\dbd
probably not; look also in the file everything.}
\bp
With the above notations one has that
\dbd
N_{D(G)}(N)=k^{G/C_{G}(N)}\bwt N
\dbd
and then the conjecture is true for $D(kG)$. Thus
\dbd
\rep(D(G)//N_{D}(N))=\langle kN\rangle .
\dbd
\ep
\bpf
By definition $N_{D}(N)$ is the smallest coideal subalgebra of $D(G)$ containing $N$. One has that $L^{G/C_{G}(N)}\bwt N$ is a \nleftcid containing $N$. It remains to show it is the smallest.
\mdn
If $C(L,H,\lam)$ contains $N$ then $C(L, H,\lam)\cap LG$ contains $N$, i.e $LH$ contains $N$. Then $L\subset C_{G}(H)\subset C_{G}(N)$.
\epf
\newpage
\section{Rest from $D(G)$}
\subsection{On the subspaces of $k^{G}$}
Let $M$ be a subgroup of $G$ and $h \in G$ an arbitrary element. We define the subspaces
\dbd
C_{\lam}(h)=\{f \in k^{G}\;|\; f\lh m=\lam(m,h)f\;\text{for all}\; m \in M\}
\dbd
\subsubsection{Non-vanishing}

Then it is easy to see that $C_{\lam}(h)\neq 0$ if and only if
\beq\label{C1}
\lam(mm',h)=\lam(m,h)\lam(m',h)\;\text{for all}\;m,m'\in M
\eeq
\subsubsection{A basis for these subspaces}

A basis on $C_{\lam}(h)$ is given by the elements $$\{f_{[g]}:=\sum_{m \in L}\lam(m,h)p_{mg}\}_{g \in G}.$$
 Note that
\beq
f_{[g]}(m)=\lam(m,h).
\eeq 
\bl
One has that 
\dbd
\bdelta_{D}(C_{\lam}(h))\subset k^{G}\ot C_{\lam}(h)
\dbd
\el
\bpf Note that
\dbd
\bdelta_{D}(f)=\sum_{g}p_{g}\ot ( g\rh f)
\dbd
and $(g \rh f)\lh m=\lam(m,h)(g \rh f)$.
\epf
\subsection{Description of certain coideals of $D(G) $: $\D(L)\subset D(kG)\ot L$}
Suppose that are given two subgroups $H$ and $L$ of $G$ and $\lam:L\times H\ra k^{*}$ is  a multiplicative function  in the first argument (on m). Consider the subspace
\dbd
C(L,H, \lam):=\bigoplus_{h \in H}C_{\lam}(h)\bwt h
\dbd
Then clearly $\bdelta(C(L,H,\lam))\subset  D(G)\otimes C(L,H,\lam)$ and
\beqn
\dim_{k}C(L,H, \lam)=\frac{|G||H|}{|L|}.
\eeqn
\subsection{Condition for these to be subalgebras}
One needs that 
\blue{\large
\beq\label{C2}
H \;\;\;\text{is a subgroup}
\eeq 
}
and 
\beq
f_{g}f'_{g'}(h^{-1}?h)\in C_{\lam}(hh')
\eeq
\blue{\large \beq\label{C3}
\lam(m, hh')=\lam(m,h)\lam(h^{-1}mh, h').
\eeq
and
\beq\label{C4}
\lam(m,1)=1
\eeq}
\subsection{When the coideal subalgebras are normal}
\subsubsection{Adjoint action of $kG$}
One has that
$$x(f_{[g]}\bowtie h)x^{-1}=f_{[g]}(x^{-1}?x)\bowtie xhx^{-1}=(\sum_{m \in L}\lam(m,h)p_{xgmx^{-1}})\boxtimes xhx^{-1}$$
First of all we need that 
\blue{\large
\beq\label{C5}
H \;\;\text{is a normal subgroup}
\eeq
}
On the other hand $xgLx^{-1}$ to be another coset $lL$. Thus 
\blue{\large
\beq\label{C6}
L\;\;\;\text{is a normal subgroup}.
\eeq
}
Then
$$xf_{[g]}x^{-1}=\sum_{m \in L}\lam(m,h)p_{xgx^{-1}xmx^{-1}}=\sum_{m \in L}\lam(x^{-1}mx,h)p_{xgx^{-1}m}$$
Therefore one should have
$$xf_{[g]}x^{-1}=\sum_{m \in L}\lam(x^{-1}mx,h)p_{xgx^{-1}m}=\al_{x,g} f_{[xgx^{-1}]}$$
for some scalar $\al_{x,g}\in L$.
This gives that
\blue{\large
\beq\label{C7}
\lam(x^{-1}Mx, h)=\lam(M. xhx^{-1})
\eeq
}
\subsubsection{Adjoint action of ${L^{G}}^{cop}$} One has the following formula:
\dbd
p_{x}.(f_{[g]}\bwt h)=\sum_{m\in G,\;uv=x}\lam(m, h)p_{v}p_{mg}\bwt hp_{u^{-1}}=\dbd \dbd=\sum_{m\in G,\;umg=x}\lam(m,h)p_{mg}\bwt hp_{mgx^{-1}}=\dbd \dbd=\sum_{m\in G}\lam(m,h)p_{mg}p_{hm(hx)^{-1}(hx)g(hx)^{-1}}\bwt h
\dbd
The result of the action is not zero  if $Mg=Mhg(hx)^{-1}$. Suppose that $g=hgx^{-1}h^{-1}m_{0}$ then the above
\dbd
p_{x}.(f_{[g]}\bwt h)=\sum_{m\in G}\lam(m,h)p_{gm}\bwt h
\dbd
\blue{\large Probably we get that
\beq\label{C8}
[H,L]=1
\eeq
}
\newpage
\section{On the rest of S8}
Note that $f \ot a \in C^{\perp_{A}\;l}\otimes A$ if and only if one has that $\langle f \ot a, (g \bwt b)c\rangle =\eps(c)\langle f \ot a, (g \bwt b)\rangle $ for all $c \in C$ and $b \in A$, $g \in A^{*}$. On the other hand
$
\langle f \ot a, (g \bwt b)c\rangle =f(bc)g(a)
$
and $\eps(c)\langle f \ot a, (g \bwt b)\rangle =\eps(c)f(b)g(a)$.
\section{On the rest of Section 6}
From \cite[Theorem 4.8]{rqts} one has that $B \cap B'=k$ there are
\dbd
A, B\subseteq H^{*}
\dbd
Hopf subalgebras
such that $\dim A\dim B=\dim H$ and $A\cap B=k$.
\vskip 0,05cm
\blue{TD: \Small If $A$ is normal then $AB=BA$ and $A\ot B\ra H$ is also factorization.}
\vskip 0,05cm
It follows form here using formulae from my math z that
\beq
H\simeq L_{A}\otimes L_{B}
\eeq
is an exact factorization. In terms of categoryories one has that
\beq
\rep(A//L_{A})\cap \rep(A//L_{B})=Vec
\eeq
\beq
\rep(A//L_{A})\vee \rep(A//L_{B})=\rep(A)
\eeq
\br With the above notations one can check directly that $$\blam_{{L_{1}}}\blam_{{L_{2}}}=\blam_{{L_{1}}{L_{2}}}$$
This shows that
\dbd
\lam_{B_{1}}\lam_{B_{2}}=t.
\dbd
From here we get the direct product of the Grothendieck rings.
\er
\newpage
\section{Link between the two results via BB}\label{link}
\bt
Let $H$ be a Hopf algebra with $C(H)$ a commutative ring (e.g. quasitriangular). Then 
\beq
\ker_{\rep(H)}(M)=\{\mu_V\;|\;V\subset \lker_H(M)\}
\eeq
\et
\blue{\bf This statemment requires the definition of $V'$s.}

\bibliographystyle{amsplain}
\bibliography{mac-bob}
\blue{\bf: adress is missing.}
\newpage
\section*{To Do:}
\subsection{On the structure constants}
For factorizable Hopf algebras one has that
\dbd
c^{i}_{jt}=\frac{d_{i}d_{j}}{d_{t}}m^{i}_{jt}
\dbd

\subsection{On the work of $D(A)$}
The braiding on $H$-mod is given by the twisted action of $R$. 
\dbd
v\ot w\mapsto R^{(2)}w\ot R^{(1)}v
\dbd
there is a formula  with antipode relating both images of both dr maps
\dbd
_{R}\phi=S \phi_{R} S
\dbd
\dbd
S^{-1}\;_{R}\phi S^{-1}= \phi_{R} S
\dbd
\blue{\large \b this is very important!}
and it is easier with the second
\dbd
f \ot h\xra{_{R}\phi} \sum_{i,j}f(b_{j})b_{i}^{*}(h)b_{i}b_{j}^{*}=hf
\dbd
\green{\bf what is the antipode on the dual of da?
\dbd
S(f\ot a):=(f\ot a)\star S :g \bwt x\mapsto  \langle f\ot a, S(x)S(g)\rangle =
\dbd
\dbd=\langle f\ot a, S(g)(S(Sx_{1})?S(x_{3}))\bwt S(x_{2})\rangle =
\dbd
\dbd=\langle f\ot a, g(x_{3}s(?)S(x_{1}))\bwt S(x_{2})=g(x_{3}S(a)S(x_{1}))f(Sx_{2})
\dbd
seems it cannot be identified as tensorand!!}
\dbd
f \ot h\xra{\phi_{R}} S^{-1}\;_{R}\phi S^{-1}(f\ot h)
\dbd

\subsection{On the conjecture for coideal subalgebra and others}

\mdn\blue{\bf Write also a version of the second main result for right coideal subalgebras.}
\mdn
\blue{\bf Try for the group case if the factorization preserves for non-degenerae. See in \cite{gnn} when a subcategoryis nondegenrate.}
\bpf
One has that $L\cap L^{*}=L$ since $\cd\vee \cd'=\mtr{Vec}$. Thus $\dim_{\kk}(LL^{*})=\dim_{\kk}(L)\dim_{\kk}(L^{*})$ and therefore $LL^{*}=A$. It follows that $M: L\ot L^{*}\ra A$ is bijective abd therefore

\dbd
\blam_{L}\blam_{L^{*}}=\blam
\dbd
On the other hand $\phi_{R}(\lam_{L})=\blam_{L^{*}}$ and $\phi_{R}(\lam_{L^{*}})=\blam_{L}$. Thus $\phi_{R}(\lam_{L}\lam_{L^{*}} )=\blam_{L}\blam_{L^{*}}=\blam$. Since $\phi_{R}(t)=\blam$ this implies that $\lam_{L}\lam_{L^{*}}=t$ which gives that $m^{*}:(A//L)^{*}\ot (A//L^{*})^{*}\ra A^{*}$ is bijective.
\epf
\blue{\bf TD4 Is in this case $A//L$ factorizable?}
\bn{defn}
We call a semisimple factorizable Hopf algebra $A$ {\it indecopmosable} if $A$ does not contain any normal Hopf subalgebra.
\end{defn}
\bc Let $L, K$ be two normal Hopf subalgebras such that $mn=nm$ for any $m \in L$ and $n \in K$.
\bne
\item
\dbd
D(A)\simeq B(L,K)\otimes B(K,L)
\dbd
\blue{\bf How do I konw that $B(K,L)$ is nondegenrate?}
\item If \dbd k \ra k^{G}\ra A \ra kF\ra k
\dbd
then
\green{\bf 
\dbd
D(A)\simeq B(k^{G}, k^{G})\otimes B(k^{F}, k^{F})
\dbd
}
\ene
\ec
\green{\bf
\bc
 Let $A$ be a quasitriangular Hopf algebra. If $\rep(A)$ contains a Lagrangian fusion subcategory$
\cd=\rep(A//L)$ then $\phi_{R}(A//L)^{*}=L$.
\ec
}
\subsection{Formulae for integrals}
\blue{\bf In the case of normal Hopf subalgebras the conjugacy classes commute each other. \\ What about the other case? of \nleftcid ?}
\bne
\item the main results are rewritten in the intro
\item In the introduction rewrite
\item The characteristic of the field. Natale works over \blue{\bf any field??}
\item \blue{\small later change to $E_{i}$.}
\item under which adjoint class are closed these $\cc_{j}$?
\item compute the muger centraliz of $\rep(A//K)$ to get a conditions for nondegenracy $(A//K)^{*}\subseteq A^{*}$.
\item 
\blue{\Small try to show that centralizer of a normal is normal in any qtriangular Hopf algebras}
\item check that the two maps coincide by CW4 left R and right R
\item \blue{\Small try something similar for non-ss factorizable Hopf algebras.}
\item \blue{as a corollary of the first main if it has distributive lattice dimension???}
\item How the proof of the second result goes?
\item Make clear that Natale was wrong!!
\item \blue{\Small Definition of the associated subcolagebra of a comodule if needed.}
\item \blue{\bf \small Moreover in this case we use the notation 
\dbd
\lam_{M}:=\lam_{\lker_{A}(M)}
\dbd
}
\ene
\mdn We also mention the following compatibilities with the antipode for the $R$-matrix:
\bne
\item 
\beq\label{compant1}
(S\ot \id)(R)=R^{-1}=(\id \ot S^{-1})(R)
\eeq
\beq\label{compant2}
(S\ot S)(R)=R
\eeq
\item The formula for the antipode-is not given in natale. If needed! I will write it down.
\ene
The following proposition appears in \cite{rqts}:
\bp\label{cmltch}
 Note that if \dbd
 \ul \Delta(a)=\ul{a_{1}}\ot \ul{a_{2}}
 \dbd
 then
 \dbd
  \Delta(a)=\ul{a_{1}}R^{2}\ot ad_{R^{1}}(\ul{a_{2}})
 \dbd
\ep

\blue{
\bne
\item
Note that a normal coideal subalgebra of $H$ corresponds to a coideal subalgebra in $\ul H$.
\item \green{\bf In the case $S^{2}=\id$ (semisimple case) one has that $(H//L)^{*}$ is a submodule under the coadjoint action of $H$ on $H^{*}$.} Indeed:
\beqn
h.p(al)=p(h_{1}alS(h_{2}))=p(h_{1}aS^{-1}(h_{2})(h_{3}lS(h_{4})))=
\eeqn
\beqn
\eps(h_{3}lS(h_{4}))p(h_{1}aS^{-1}(h_{2}))=\eps(l)p(h_{1}aS^{-1}(h_{2}))=\eps(l)h._{'}p(a)
\eeqn
\item $(H//L)^{*}$ is a subcoalgebra of $\ul H^{*}$\blue{and a subalgebra!!!}. Indeed
\dbd
p\ul .q(al)=(S(p_{1})p_{3}\ot q_{1})(R)p_{2}q_{2}(al)=(S(p_{1})p_{3}\ot q_{1})(R)p_{2}(a_{1}l_{1})q_{2}(a_{2}l_{2})=\dbd
\dbd
=(S(p_{1})p_{3}\ot q_{1})(R)p_{2}(a_{1}l_{1}\eps(l_{2}))q_{2}(a_{2})=\eps(l)(S(p_{1})p_{3}\ot q_{1})(R)p_{2}(a_{1})q_{2}(a_{2})
\dbd
On the other hand it can be checked that
\dbd
[h.p](la)=\eps(l)[h.p](a)
\dbd
\item
\green{\bf In order to coincide left and right
 in description of $(A//L)^{*}$ one needs $S^{2}=\id$. !!!, i.e semisimple} 
\item It follows that $\phi_{R}((H//L)^{*})$ is a sub-bialgebra of $\ul H$ which shows that it is a left normal coideal subalgebra of $H$.
\item \green{\bf In the semisimple case one does not need to take subalgebra generated; it is already.}
\ene
}
\mdn
\blue{\bf wrte it witheq}
\mdn
 For the rest of the paper we use the notation
\beq
\rep(A//L)'=\rep(A//L^{*})
\eeq
where
 \beq
 L^{*}=\phi_{R}((A//L)^{*}).
 \eeq
 
 \bne
\item combine with the projective centralizing formula from gelaki nikshych -nilpotent case.
\ene
\section{To do section 5}
\bne
\item 
\green{\bf It looks very much as conjugacy classes and probably they commute when the character centralizes}
\item Combining 22 and 25 it follows that
\dbd
\mtc C^{j}\subseteq L\;\;\text{iff}\;\; \lker_{A}(V_{j})\supseteq L^{*}=\bigoplus_{\ch_{i} \in \cd}\cc^{i}
\dbd
\green{\bf \item My intuition is that
\dbd
L=\bigoplus_{j \in \cd'}\cc_{j}
\dbd
and
\dbd
L'=\bigoplus_{i \in \cd}\cc^{i}
\dbd
}
\item Since $m:L\ot L'\ra H$ is a bijection maybe they commute.
\item Study as example the Drinfeld double; for $S_{3}$ it is in the paper of CW from comm algebra.
\item \blue{\bf Try to write a similar result for quasitriangular Hopf allgebra, not necessarily factorizable; and apply it for $kG$.}
\ene
\subsection{Explanation for Drinfeld's result for $\cd'$-components}
\bt(DGNO2) Suppose that $\cc$ is nondegenerate. The number of $\cd'$ components equals the number of irreducible objects of $\cd$.
\et
In the non-degenerate case of $\cd$ this is clear since if the char's have the same restriction to $L'$ then they coincide.

\bp Let $\ck$ be a fusion subcategoryof $\rep(A)$. Then
\beq
(\mtc K_{ad})'=\ck=\ck'^{co}
\eeq
\ep
\bpf
Suppose that $\ck=\rep(A//L)$. Then $\ck_{ad}=\rep(A//L_{1})$ where $L_{1}:=\pi^{-1}(L(A//L))$. Thus
\dbd
(\mtc L_{ad})'=\rep(A//L_{1})'
\dbd
On the other hand 
\dbd
\ck'^{co}=\rep(A//L')^{\co}=\rep(A//[A,L'])
\dbd
It follows that
\dbd
\phi_{R}((A//[A,L'])^{*})=L_{1}'
\dbd

\epf
\subsection{On section 5}
The following lemma will be needed throughout the paper.
\bl \cite[Lemma 1.1]{CW10} Let $T$ be a left coideal subalgebra of $A$ and $B=T^{\perp_{l}}$. Let $\blam_{T}$ and $\lam$ be an  idempotent integral of $T$ and $A^{*}$. Then
\bne
\item $\blam_{T}\rh \lam=\lam_{B}$ is a nonzero integral in $B$.
\item $\lam_{B}(\blam_{T})=\frac{1}{|T|}$ since the regular character $A\blam_{T}$ is free of rank $|A|/|T|$.
\item 
$
\blam_{T}=\lam_{B}\rh \blam
$
\item $
T=\blam_{T}\lh A^{*}=\lam_{B}\rh A
$
\item $B=\lam_{B}\lh A=\blam_{T}\rh A^{*}$
\item $S\blam_{T}=\blam_{T}$ and $ST=A^{*}\rh \blam_{T}$.
\ene
\el

\subsection{Maximal versus minimal and dimension results}

Use Etingof descr weaLly-grp tht pentru category de fuziune simple.
\\ 
Maybe in genral for braided fusion categoryories there is the fact that normal are centralized into normal.
\\ *
In the non-degenrate case the dimensions are not preserved, for $\bar A$ and $L$. There is a different formula, take it from DGNO.
\\ * $(A//A^{*})$ is symmetric, thus $A//A^{*}$ is quasitriangular
\\ *
Minimal versus maximal coideal subalgebras; $\lker_{A}(M)$ are maximal for irreducible characters: 
\dbd
L \subseteq M \iff L^{*}\supseteq \mprime
\dbd
with reciprocal equality. 
\\ *
La fel pentru subalgebrele Hopf normale.
\\
\dbd
{\bar A}^{*}=L
\dbd
\\
$\bar A$ is a maximal coideal subalgebra since
\dbd
\bar A\subset L \iff k ={\bar A}^{*}\supseteq L^{*} \iff L^{*}=k\iff L=A
\dbd
\\ 
One has that
\dbd
k \subset L \implies \bar A\supset L^{*}
\dbd
Thus all stars are in fact normal coideal subalgebras of $\bar A$.
\\ 
One has that
\dbd
(LM)^{*}=L^{*}\cap \mprime\;\;\text{and}\;\; (L\cap M)^{*}=L^{*}\mprime
\dbd
\\ Moroever
\dbd
L^{**}=L\cap \bar A
\dbd
\\ Thus
\dbd
(LL^{*})^{*}=L^{*}\cap L^{**}=L^{*}\cap L \cap A
\dbd
\\ In the quasi-triangular setting if $L$ minimal then $L^{*}$ is a maximal \red{left normal coideal} subalgebra.

\subsection{On section 4}

\green{\bf TD2: Now apply it when one has equality of centralizer M\"uger  categoryories. It should give equality in the second case at least.}
\mdn
\green{\bf TD1:Then apply it to the case from my previous paper for $L_{d}$. Does $L_{d}
$ probably is minimal?}

 Note that by Equation \eqref{fpdim} one has
 \beq\label{fpdim2}
 \dim_{\kk}(L\lprime)\dim_{\kk}(L\cap \lprime)=\dim_{\kk}(A)
\eeq
 Previous proposition implies that in the factorizable case
$
 L'= L^{*}
$
 and
$
 (\ovr AL)'\supseteq \lprime
$.
\subsection{On the first factorization}
\bc\label{first-factoriz}
Suppose that $A$ is a semisimple factorizable Hopf algebra and $\cd=\rep(A//L)$ is a non-degenerate fusion subcategorywith $L$ a normal \blue{\bf coideal} algebra.  Then  $$A^{*}\simeq (A//L)^{*} \otimes (A//\lprime)^{*}$$ as Hopf algebras. 
\ec
\bpf
In this situation we have 
$
\cd \vee \cd'=\cc
$
and
$
\cd \cap \cd'=\mtr{Vec}
$. Then by Equations \eqref{capsf} and \eqref{veesf} one has that $L\cap \lprime=k$ and $L\lprime=A$. This implies that $(\dim_{\kk}L)(\dim_{\kk}\lprime)=\dim_{kk} A$. Therefore $\langle (A//L)^{*},\;(A//\lprime)^{*}\rangle =A^{*}$ and $(A//L)^{*}\cap (A//L^{*})^{*}=k$ and \blue{\bf the map $m: (A//L)^{*}\otimes (A//\lprime)^{*}\ra k$ is surjective.}
\epf
\subsection{Closure of relation}
Let $A$ be a factorizable Hopf algebra and $L$ be a  Hopf subalgebra.
Then
\beq
\cd(L)'=\{\ch_{j}\;|\;\}
\eeq

\beq
\cd(L)'=\cd(N_{A}(L))'=\{\ch_{j}\;|\; E_{j}(\blam_{N_{A}(L)})\neq 0\}
\eeq

\dbd
\blam_{N_{A}(L)}=\prod_{d \in \irr(L^{*})}L_{d}
\dbd
 \subsection{Application on the double centralizer}
One has from \cite{nrq} ''natale r-qotients'' that
$\cc'=\rep(A//\overline{\bf A})$ where $\overline{\bf A}=\phi_{R}(A)$.

Therefore one has
\beq
\rep(A//L)''=\rep(A//L)\vee \rep(A//\overline{\bf A})\simeq \rep(A//L\overline{\bf A})
\eeq
\subsubsection{Minimality of $\overline{\bf A}$ is equivalent to that of $\bf H$ as in Natale's paper}

A consequence of this result is that: 
$$
\cap_{n \geq 0}\mtr{Ann}_A(M^{\ot\;n})=A(\lker_{ _A}(M))^+.
$$
\blue{\Small  Where is this used? In what form?}
\subsection{On the projective centralizer}
\br Following \cite[Theorem 3.14]{CW5} there is a similar criterion for projective centralizers. More precisely, $V_{i}$ projectively centralize $V_{j}$ if and only if 
\dbd
s_{ij}=\omega d_{i}d_{j}
\dbd
for some root of unity $\omega$. Thus this happens 
if and only if
$
\langle \ch_{i}, \eta_{j}\rangle =\omega d_{i}
$
 and by Theorem 3.14 of \cite{CW5} if and only if 
$
\mtc C^{j}\subseteq \lker^{\omega_{i}}_{V_{i}}
$.
\er
\subsection{Correlating with the results from GN}\red{\bf keep this}
Recall the definition of the commutator $A_{comm}$ of a Hopf algebra \cite{iop} as the smallest left normal coideal subalgebra $L$ with the property that $A//L$ is a commutative Hopf algebra.

If $\cc$ is a pseudo-unitary modular categorythen
$
(\cc_{ad})'=\cc_{pt}
$
by a result of \cite{gnn}. \red{see in DGNO if this is also the case in the braided case}
\mdn
If $\cc=\rep(A)$ for a semisimple factorizable Hopf algebra $A$ then $\cc_{ad}=\rep(A//K(A))$ where $K(A)$ is the Hopf centre of $A$, i.e. largest central Hopf subalgebra of $A$. Thus by our first main result it follows that $\phi_{R}((A//K(A))^{*})=A_{\mtr{comm}}$ the commutator coidela subalgebra of $A$.
\bpf
It is shwon, see \cite[Prop \blue{\bf ??}]{iop} that $A_{\mtr{comm}}$ is the intersection of left Lernels of all linear characters of $A$. Thus $\rep(A)_{pt}=\rep(A//A_{\mtr{comm}})$.
\epf
 In particular 
$$A\simeq A//K(A) \otimes A//L(A)$$ if $\cc_{ad}$ is a nondegenerate braided fusion categoryory. Note that $A//L(A)$ is a commutative Hopf algebra.\
\bc\label{c1} With the above notations if $(A,R)$ is a semisimple factorizable Hopf algebra then
\beq\label{c1}
\rep(A//LL^{*})'=\rep(A//L\cap L^{*})
\eeq
\ec
\bpf
Indeed, using Equation \eqref{veesf} one has that
\dbd
\rep(A//LL^*)'=(\rep(A//L)\vee \rep(A//L^*))'=\rep(A//L)'\cap \rep(A//L^*)'=\dbd \dbd=\rep(A//L^*)\cap \rep(A//L)=\rep(A//L\cap L^*)
\dbd
\epf

\green{\bf Do also the case of a q-triang Hopf algebra, probably H bar enters the details.}
\red{Recall that a fusion categoryis called nondegenerate if \green{\bf $\mtc Z_{2}(\cd)=\mtr {Vec}$.} By egno15!!!\cite[Prop \blue{\bf ??}]{proclond} if $\cd$ is nondegenerate then $\cc \simeq \cd\boxtimes \cd'$ where $\boxtimes$ denotes the Deligne tensor product of abelian categoryories.
}

 \section{ Connection with old results-M\"uger  centralizer in the categoryof representations of Drinfeld doubles.}\label{conold}
\subsection{Normal sent to normal in the factorizable case}

 By transmutation theory try to prove that any normal Hopf subalgebra is sent to a normal Hopf subalgebra. It depends on a result of CW which is true for factorizable.\blue{\bf \small It should be similar computations with the second map.}
In \cite{mathz} MathZ we proved the following results:
\bt\label{centraliz} Let $A$ be a semisimple Hopf algebra and $L,K$ be two normal Hopf subalgebras of $A$ such that $B(L,K)$ is a normal Hopf subalgebra of $D(A)$. Then
\beqn
\cd(L, K)'=\cd(K,L)
\eeqn
as fusion subcategoryories of $\rep(D(A))$.
\et
\subsection{A conjecture on coideal subalgebras}

Let $A$ be a semisimple Hopf algebra and $L$ be a Hopf subalgebra of  the Drinfeld double $D(A)$. We denote by $\cd(L)$ the fusion subcategoryof $\rep(D(A))$ whose objects  are those $D(A)$-representations that receive a trivial $L$-action.

Our second main result gives a description for the centralizer of $\cd(L)$:
\bt\label{justonecomp}
Let $A$ be a semisimple Hopf subalgebra and $L$ be a  normal Hopf subalgebra of $A$. Then 
\beqn
\cD(L)'=\langle L\rangle 
\eeqn
where $\langle L\rangle $ is the fusion subcategoryof $\Rep (D(A))$ generated by the $D(A)$-module $L$. Here $L$ is regarded as $D(A)$-module via the action
\beq
(f \bwt a)x=a_{1}xS(a_{2})\lh S^{-1}f
\eeq
where $a \in A$, $f \in A^{*}$ and $x \in L$.
\et
{\bf Conjecture:} Let $L$ be a normal coideal subalgebra.\dbd
\cd(L)'=\langle L\rangle 
\dbd
How did I prove it for normal Hopf subalgebras? I applied the result for factorizable Hopf algebras.
\mdn
\blue{\bf See in the case of $D(G)$ what one can obtain?}

\subsubsection{
\blue{\bf A second approach to compute $\langle L\rangle '=L^*$}}
\dbd
M\ot V\xra{c^{2}} M\ot V; m\ot v\mapsto Q^{1}m
\ot Q^{2}v
\dbd
\subsubsection{\blue{\bf A third approach} }Apply first main result: It follows that
\dbd
\cd(L)'=\rep(D(A)//N_{D}(L))'=\rep(D(A)//N_{D}(L)^{*})
\dbd
where 
\dbd
N_{D}(L)^{*}=\phi_{R}((D(A)//N_{D}(L))^{*})
\dbd
\blue{\bf \Large
So the conjecture is equivalent to:
\dbd
\phi_{R}((D(A)/N_{D}(L))^{*})=\lker_{D(A)}(L).
\dbd
\mdn work here with the second $_{R}\phi$ is easier $f \ot h\mapsto hf$. 
\mdn One has that $S\star\;_{R}\phi=\phi_{R}\star S$.
So the conjecture is equivalent to:
\dbd
_{R}\phi((D(A)/N_{D}(L))^{*})=S(\lker_{D(A)}(L)).
\dbd
From definition of right kernel this is equivalent to
\dbd
(h.f)_{2}\ot (h.f)_{1}l=h.f\ot l
\dbd
for any $f\ot h\in (D(A)/N_{D}(L))^{*}$ and $l \in L$.}
\bp Let $A$ be a Hopf algebra and $T$ be a coideal subalgebra of $A$. Then with the above notation
\dbd
L^{\perp_{l}}??=N_{A}(L)^{\perp_{l}}
\dbd
\blue{\bf Dimensions do not matches} 
\ep
\subsubsection{\blue{\bf A fourth approach; apply the third creiterion}}
\dbd
\rep(D(A)//N_{D}(L))'=\langle \ch_{j}\;|\; \cc_{j}\subseteq N_{D}(L)\rangle 
\dbd

\green{\bf \small guess a formula for the integral of the closure; generalize the core formula from Hopf subalgebras}

  \subsection{ Connection with old results-M\"uger  centralizer in the categoryof representations of factorizable Hopf algebras.}
\newpage
\section{A formula for the monodromy element $Q$}\label{formmonodr}
\subsection{Conjugacy classes and decomposition of the integral}
Fix idempotent integral $\blam$ and $t$. One has that $\langle t, \blam\rangle =\frac{1}{d}$.
One has
\dbd
C_{j}=\blam \lh dF_{j}
\dbd
and from here
\dbd
\sum_{j}C_{j}=d\Lam
\dbd
This implies that 
\dbd
\blam_{j}=\frac{C_{j}}{d}
\dbd
Moreover dim $\dim (\cc_{j}\cap Z(A))=1$ since $C_{j}\in \cc_{j}\cap Z(A)$ and $\dim Z(A)=r$.
\subsection{Decomposition of elements}
Note that
\dbd
E_{m}=\ch_{m}(1)(\blam \lh \ch_{m}^{*})=\sum_{j}\ch_{m}(1)(\blam \lh \ch_{m}^{*})
\dbd
This implies that
\dbd
E_{m[j}=\ch_{m}(1)(\blam \lh \ch_{m}^{*})
\dbd
\subsection{Choosing the idempotents in the character ring}
We choose the idempotents such that $E_{m}$ is the idempotent corresponding to the conjugacy class $\cc_{m}$. In other words $E_{m}(c)=\eps(c_{[m})$ for any $c \in A$.
\green{\bf Is it possible to generalise this result to fusion categoryories with commutative Grothendieck ring? In fact all results from ART?}
\subsection{Formula for $Q$}

We write 
\dbd
Q=\sum_{l,j}Q_{l,j}
\dbd
with $Q_{l,j}\in \cc_{l}\ot \cc_{j}$.
\mdn \green{\bf need the first index to give it.}
\mdn
One has that
\dbd
E_{m}=\phi_{R}(F_{m})=\sum_{j}\eps(Q^{1}_{m,j})Q^{2}_{j,m}
\dbd
and 
\dbd
E_{m}=_{R}\phi(F_{m})=\sum_{j}\eps(Q^{2}_{l,m})Q^{1}_{m,l}
\dbd
This implies that
\dbd
\ch_{m}(1)(\blam \lh \ch_{m}^{*})={E_{m}}_{[j}=(\eps \ot \id)Q_{m,j}=(\id \ot \eps)Q_{j,m}
\dbd
\subsection{Define $A_{j}:=\cc_{j}^{\perp}$}

\newpage
\section{A new list to do-Hopf algebras}
\bne
\item try to see what is the core as in DGNO for factorizable Hopf algebras.
\item read the paper of dong on q-triangular Hopf algebras ASF
\item give entire examples; $D(H_{8})$ for example.
\item Sommerhauser's rule probably is in ENO, annals.
\item Need to have clear the equation class and from here the chracters on the Grothendieck group.
\item give entire examples; $D(H_{8})$ for example.
\item \blue{\bf chiar daca fac si categoryorical, intai il pun pe acesta ptr algebre si Hopf si apoi pe celalalt. pe Hopf la ant.}
\item \green{\bf dupa ce termi articolul o intreb pe sonia de greseala din al ei.}
\item as another application maybe a use of te central characters.
\item give entire examples; $D(H_{8})$ for example.
\item try to use central chars too
\dbd
d_{\ch}=\sum_{i\in A_{l}}\ch_{i}(1)\ch_{1}
\dbd
In this case $L_{d_{\ch}}=H_{d_{\ch}}$.
\item In categoryorical settings the left coideal subalg replaced by etale-subalgebras of $R(1)$.
\item $f_{Q}$ defined only in terms of characters, only at this level.
\item ``Shimizu's conjecture''
\item Sommerhauser's rule probably is in ENO, annals.
\item Need to have clear the equation class and from here the chracters on the Groethendieck group.
\item \item \green{Look in Natale's pap which gives a left normal coideal subalgebra and apply for that one.}
\item 
\dbd 
a_{1}\ot a_{2}m=R^{-1}(a_{2}\ot a_{1})R(1\ot m)
\dbd
\item what is $\lker_{D(H)}(\cc_{j})$ where $\cc^{r}_{j}$ is made such that to be a left coideal.

\ene

\section{A new list to do-categoryorical settings}
\bne
\item read the paper of sonia and dong for ASF-modular
\item dong has another one for integral asf with li 
\item my main first result uses integrals for Hopf and coideals; basically these can be rewritten in terms of idempotents and characters.
\item In categoryorical settings the left coideal subalg replaced by etale-subalgebras of $R(1)$.
\item $f_{Q}$ defined only in terms of characters, only at this level.
\item ``Shimizu's conjecture''
\item Sommerhauser's rule probably is in ENO, annals.
\item Need to have clear the equation class for fusion categoryories and from here the characters on the Grothendieck group.
\ene
\newpage
\section{15.04-The monoidal center and the character algebra}
\blue{It works only with commutative character rings}
\subsection{Coends and the result of Day and Street}
\dbd
_{Z}\cc\simeq \cz(\cc)
\dbd
dinatural transformation $i_{V, X}$.
\mdn
under the above identification the forgetful functor $U:\cz(\cc)\ra \cc$ has as \blue{left} adjoint the free module functor.
\mdn
the right adjoint is $L^{!}$
\mdn
\green{\bf shimizu has given a relation between left and right adjoints of a tensor functor via left and right duality.}
\subsection{Unimodular tensor categories}
Let $U:Z(\C)\ra \C$ be the forgetful functor and let  $R$ be its right adjoint.
\mdn
There is a distinguished invertible object $D\in \C$ introduced by Etingof and Ostrik and rewritten by Shimizu.

\bn{defn}
A tensor category is called unimodular if this element is the identity.
\end{defn}
\noindent
\blue{\bf define $A=UR(1)$ as  algebra in $\cc$??\mdn define $\eps_{1}, \; \delta_{1}$.}
\subsection{Central elements}

\mdn
Definition\beq 
\cecc:= \hm_{\C}(1,A) 
\eeq
is called \blue{the set of central elements.}

\blue{\bf When this bijection is used:}
There is a bijection:
\beqn
\cecc\xra{} \enx(\id_{\C})
\eeqn
\subsection{Definition of the character algebra}
The internal character $\mtr{ch}(X)$ is defined as the class function
\beqn
\mtr{ch}(X):=tr^{X}_{A, 1}(\rho_{X}).
\eeqn

\subsubsection{The pairing}
In the theory of finite-dimensional Hopf algebras, it is well-known that the pairing
between cointegrals and integrals is non-degenerate. To formulate an analogous fact
in our setting, we consider a paring
\beqn
\langle f , a \rangle : CF(\C) \times CE(\C)\ra L, \\ \langle f, a\rangle  \id_{1}= f \circ a.
 \eeqn
\subsection{Conjugacy classes}

A conjugacy class is defined as a direct summand of $R(1)$.

 Thus let $\mtc C_{0}\cdots, \C_{m}$ be the
conjugacy classes of $\C$. 
 
 Since the unit object $1_{\czcc }$ is always a subobject of $R(1)$, we can assume $\C_{0} = 1_{\czcc }$.

\bt (Theorem 6.6) For a non-degenerate pivotal fusion categoryC, the following assertions are equivalent:

(1) $R(1) \in \czcc $ is multiplicity-free.

(2) The algebra $Gr_{L}(C)$ is commutative.
\et

\subsection{Integral theory}
 An \blue{\bf integral} in $\C$ is a morphism $\Lambda: 1 \ra A$   in $  \C$   such that
 \beqn
 m \circ (\id_{A}\otimes \Lambda)=\eps_{1}\ot \Lam
 \eeqn
   
A \blue{\bf cointegral in $  \C $}  is a morphism $  \lambda : A\ra 1 $  such that
\beqn
\barz (\lam)\circ \delta_{1}=u \otimes \lam
\eeqn

\subsection{Integrals and Fourier transform}
Let $  \Lambda $\; be a non-zero cointegral in $  \C $  and $\lambda$ be the integral such that 
\beqn
\langle \Lambda, \lambda\rangle =1
\eeqn
\noindent
The Fourier transform associated to $\lambda$ is the linear map
\beqn
\mtc F_{\lambda}:\cecc\ra \cfcc
\eeqn
given by
$$
a \mapsto \lambda \lh S(a)
$$
where $S:\cecc\ra\cecc$ is the antipodal operator defined by:

and $\lh$ id the action of $\cecc$ on $\cfcc$ given by 
\beqn
f \lh b=f \circ m \circ (b\ot \id_{A})
\eeqn
for $f \in \mathrm{\bf{CF}}(\cc)$ and $b \in\cecc$.

\blue{\bf this is the harpoon $(f\lh b)(a)=f(ba)$. }
\subsection{Radford's trace formula}
\subsubsection{The Drinfeld isomorphism}
For any central object $(X, \sg)\in \cz(\cc)$ one defines  the canonical isomorphism:
\beqn
\psi_{X}:X\xra{} X\ot X^{*}\ot X^{**}\xra{\sg, \id} X\ot X\ot X^{**}\xra{\ev \ot id} X^{**}
\eeqn

\blue{By [Shim, sect 7] $\psi_{A}=j_{A}$ if $\cc$ has a pivotal structure $j$.}
\subsubsection{Definition of the trace by abuse of notation-using Drinfeld isom}
For any $\xi:A\ra A$ define
\beqn
\tr(\xi):1\xra{\coev}A\ot A^{*}\xra{\psi_{A}\ot \id} A^{**}\ot A \xra{\ev}1
\eeqn
\subsubsection{Radord's formula}

\bt Let $\cc$ be a unimodular finite tens categoryory. With the above notations one has:
\beqn
\tr(\tilde{f})=\langle f, \Lam\rangle \langle \lam, u\rangle 
\eeqn
\et
\newpage
\subsection{Application-example to the modular case}Let $\cc$ be a braided fusion categoryory. The map 
\dbd
f_{Q}:\ccf\ra \cce,\;\; \ch_{i}\mapsto \sum_{j}\frac{s_{ij}}{d_{j}}e_{j}
\dbd
is morphism of $L$-algebras by  [BL01, Theorem 3.1.1] and its proof.
\mdn
It follows that
\dbd
\ch_{i}=\sum_{j}\frac{s_{ij}}{d_{j}}f_{j}
\dbd
after an appropriate renumeration of the idempotents.

\mdn \green{\bf Q1 Is there a sommerhauser rule concerning the restriction fro the center? }
\mdn \green{\bf \small need to know the multiplications on both $\ccf$ and $\cce$. One is by convolution, how it is the other one?}
\mdn \blue{\bf Here it is the conjecture: maybe it can be solved using cosets as in dgno}
\bt
Let $\cd$ be a fusion subcategoryof $\cc$ and $R_{\cd}$ be the regular character of $\cd$. Then
\dbd
\phi_{R}(R_{\cd})=\Lam_{\cd'}
\dbd
\et
\mdn
how to identify the idempotents of $\cce$ with ``those'' of the characters.
\newpage

\newpage
\section{On Higman's ideals and conjugacy classes}
\subsubsection{ Frobenius and symmetric algebras}
\bne
\item first section studies Frobenius and symmetric algebras
\item for dual bases one has from the isomorphism $A\ot A\simeq ENd_{L}(A)$ that
\dbd
\sum_{i}xa_{i}\ot b_{i}\simeq \sum_{i}a_{i}\ot b_{i}x
\dbd
\item \green{\bf only in the case of symmetric algebras when $a_{i}$ and $b_{i}$ can be switched it can be proven the other relation.}
\item define the higman map
\dbd
\tau_{A}(x)=\sum_{i}a_{i}xb_{i}
\dbd
\green{\bf it has values in the center of $A$.}
\item if $e$ is idempotent primitive and $E$ the unique idempotent primitve central such that $eE=e$ then
\dbd
\tau(e)=\al E
\dbd
with $\al\neq 0$.
\item one can check various formulae; mainly Porp 1.13
\dbd
p_{a}=t\lh \tau(a)
\dbd
\dbd
\tau(E)=\dim(Ae)\tau(e)
\dbd
\dbd
\langle p_{e}, \tau(e)\rangle =\frac{\langle p_{e},\tau(1)\rangle }{\dim (Ae)}
\dbd
\dbd
\al=\frac{\langle p_{e},\tau(e)\rangle }{\dim (Ae)}=\frac{\langle p_{e},\tau(1)\rangle }{\dim (Ae)^{2}}
\dbd
\dbd
\dbd
\green{\bf In the non-ss case I didn't understand why since $\tau(e)E=\tau(e)=\al E$ for some scalar $\al$. In the semisimple case it is clear from idempotents splitting.}
\ene
\subsubsection{Study of the chracter ring as a symmetric algebra}

\bne
\item the conjugacy classes are not closed under the left adjoint action; probably they are closed under the right adjoint action
\item from my paper normal Hopf ART, the Frobeniuss map $a \mapsto a\rh t_{H}$ is an isomorphism of $D(H)$-modules
\item however it seems that babutele showed that $E_{j}$ coincid with $\eps$ on the conjugacy class $\cc_{j}$ and are zero outside.
\item Babutele wrote a verlinde formula in general in terms of the matrix $A=\al_{ij}$ thath changes between irreducible characters and primitive central idempotents. 
\item \green{\bf it seems that what I mean by left coideal is right coideal in their paper.}
\item Radford's trace formula worLs for Frobenius Hopf algebras
\item the argument from class equations also worLs in general for symmetric Hopf algebras. It follows that in the commutative case
\dbd
\dim (H^{*}E)=\langle \ch_{H^{*}}, E\rangle =\dim (H^{*})\langle \Lam, E\rangle 
\dbd
Thus
\dbd
\langle E_{j}, \blam\rangle =\frac{1}{n_{j}}
\dbd
where $n_{j}=\frac{\dim H^{*}}{\dim(H^{*}E_{j})}=\frac{\dim H^{*}}{m_{j}\dim(H^{*}e_{j})}$
\item Moreover it can be shown that
\dbd
\al_{j}=\frac{\dim H^{*}}{\dim(H^{*}e_{j})}=m_{j}n_{j}
\dbd
is a positive integer.
\item suppose that
\dbd
\ch_{i}=\sum_{j}\al_{ij}E_{j}
\dbd
\item it easy to see then that 
\dbd
\langle \ch_{i}, C_{j}\rangle =\langle \ch_{i}, dE_{j}\lh \blam\rangle =d\langle E_{j}\ch_{i}, \blam\rangle =d\al_{ij}\langle E_{j}, \blam\rangle 
\dbd
\item
\item
\ene

\subsection{\green{\bf To do generalise to quasitriangular Hopf algebras}} If it is quasitriangular and simple then it is factorizable since $\phi_{R}(A^{*})$ is a normal Hopf subalgebra of $A$.
\subsection{\green{\bf To do; generalise to braided categoryories using the conjugacy classes introduced by Shimizu.}}
\newpage
\subsection{math z continuare-apply the previous results on their result on forgetful functor!!! laregest subalg etale din indus!!!}
\newpage
\subsection{Lax monoidal/comonoidal adjunction for the left adjoint}
suppose $F:\cd\ra \cc$ is a strong monoidal functor and $I+L$ is a left adjoint while $R$ is a right adjoint.
One has that \green{\bf for $X\in \cc$ and $Y \in \cd$:}
\beq
_{\cd}(I(X), Y)\simeq _{\cc}(X, F(Y))
\eeq
and
\beq
_{\cd}(Y, R(X))\simeq \;_{\cc}(F(Y), X)
\eeq
\subsection{For the left adjunction}
One has canical
\dbd 
X=F(Y) \implies \eps_{Y}:I(F(Y))\ra Y \;\text{corresponding to the  identity}\; \id_{F(Y)}
\dbd
\dbd
Y=I(X)\implies \eta_{X}:X\ra F(I(X)) \;\text{corresponding to the  identity}\; \id_{I(X)}
\dbd
Then 
\dbd
I^{2}(X, X'):I(X\ot X')\ra I(X)\ot I(X')
\dbd is given by
{\small \dbd
I(X\ot X')\xra{I(\eps_{X}\ot \eps_{X'})} I(FI(X)\ot FI(X'))\xra{I(F_{2}(I(X), I(X'))} IF(I(X)\ot I(X'))\xra{\eta_{I(X)\ot I(X')}} I(X)\ot I(X')
\dbd
}
\subsection{For the right adjunction}
One has canical
\dbd 
X=F(Y) \implies \eps_{Y}:Y\ra R(F(Y)) \;\text{corresponding to the  identity}\; \id_{F(Y)}
\dbd
\dbd
Y=R(X)\implies \eta_{X}:F(R(X))\ra X \;\text{corresponding to the  identity}\; \id_{I(X)}
\dbd
Then 
\dbd
R^{2}(X, X'):R(X)\ot R(X')\ra R(X\ot X')
\dbd is given by
{\small \dbd
R(X)\ot R(X')\xra{\eta_{R(X)\ot R(X')}} RF(R(X)\ot R(X'))\xra{RF_{2}(R(X), R(X'))} R(FR(X)\ot FR(X'))\xra{R(\eps_{X}\ot \eps_{X'}} R(X\ot X')
\dbd
}
\green{\bf Need to use the adjunction that left and right are isomorph in the case of fusion categoryories}
\subsection{Co-induction and induction are isomorphic} In the paper CW2 comm algebra the
Frobenius map eq 1 gives that
\dbd
p(h\Lam_{1})\Lam_{2}=p(\Lam_{1})S^{-1}(h)\Lam_{2}
\dbd
i.e 
\dbd
h\Lam_{1}\ot \Lam_{2}=\Lam_{1}\ot S^{-1}h\Lam_{2}
\dbd
apply $\id \ot S$ gives
\dbd
h\Lam_{1}\ot S(\Lam_{2})=\Lam_{1}\ot S(\Lam_{2})h
\dbd
where $\Lam$ is a left integral of $H$ and $t$ a  right integral of $H^{*}$.
\green{\bf If $\blam$ commutative then 
\dbd
h\Lam_{2}\ot S(\Lam_{1})=\Lam_{2}\ot S(\Lam_{1})h
\dbd
}
\bp
This gives that
\beq
A^{*}\xra{} A,\;\; p\mapsto \Lam \lh Sp
\eeq 
is a morphism of left $A$-modules.
Its inverse is given by $$a \mapsto a \rh t$$ with $t(\Lam)=1$
\ep
\green{\bf  need its inverse}
\bp Regard $A^{*}$ as letf $A$-module via left hit. Let $L$ be a \green{Left} coideal subalgebra of $A$ and $M$ be a left $L$-module. Then
\beq
\phi:A\ot_{L}M\simeq \hm_{L}(A^{*}, M)\; \;a\ot_{L}m\mapsto (T_{a, m}:
 (f \mapsto f(a)m))\eeq
 is an isomrphism of $A$-modules.
\ep
\green{\bf it is not well defined}
\bpf
We define the inverse by
\dbd
f \mapsto \sum_{i}a_{i}\ot_{L}f(a_{i}^{*})
\dbd
Then the composition is
\dbd
f \mapsto \sum_{i}a_{i}\ot_{L}f(a_{i}^{*})\mapsto \sum_{i}T_{a_{i}, f(a_{i}^{*})}=f
\dbd
and
\dbd
a\ot_{L}m\mapsto T_{a, m}\mapsto \sum_{i}a_{i}\ot_{L}T_{a, m}(a_{i}^{*})=\sum_{i}a_{i}\ot_{L}a_{i}^{*}(a)m=a\ot_{L}m
\dbd
Moreover
\dbd
\phi(ba\ot_{L}m)=T_{ba, m}=bT_{a,m}
\dbd
\epf
The two combined give that
\dbd
\psi:A\ot_{L}M\xra{\simeq} \hm_{L}(A,M)\;\; a\ot_{L}m\mapsto \tilde T_{a,m}(b \mapsto t(ba)m)
\dbd
\green{\bf when transp the structure of left $A$-module from $\hm_{L}(A^{*},M)$ to $\hm_{L}(A,M)$ one gets the usual structure of the coinduced module.}
\mdn
The inverse of this map is 
\dbd
\phi: \hm_{L}(A,M) \xra{\simeq}A\ot_{L}M \;\; g \mapsto \sum_{i}a_{i}\ot_{L}g(\blam \lh S(a_{i}^{*}))
\dbd
\green{\bf does not maLe sense the first part, the second, yes it does!need to consider relative bases; it is defined in terms of relative bases. but does not depend on it!}
\dbd
f\mapsto Sx_{i}\ot_{L}f(x_{i})
\dbd
\subsection{Easy to worL with coinduced modules:}
\dbd
_{A}(D,M)\ot _{A}(D, N)\xra{\phi} _{A}(D, M\ot N)\;\; f\ot g\mapsto T_{f,g}:(d \mapsto f(d_{1})\ot g(d_{2}))
\dbd
write the monodromy between two such modules
\dbd
c^{2}:f \ot g\mapsto \sum_{i,j}(b_{j}b_{i}^{*}).f\ot (b_{j}^{*}\bwt b_{i})g
\dbd
and dyslecti correspond to 
\dbd
\phi=c^{2}\phi
\dbd
for $M=\mathrm 1$

\newpage
\br
 This proposition shows that in general $B(L, L)\supseteq B(L', L, \psi')$.
 \er

\subsection{Fusion subcategoryories of Drinfeld centers $\cZ(\cA)$ and their centralizers}\lb{qst}
\blue{Two questions on the formulae Lernels la for left  Lernels and $\cd(L)'$}.\mdn
Let $\ca$ be a fusion categoryand let $\cZ(\cA)$ be \blue{its Drinfeld center} with forgetful functor $F : \cZ(\cA) \ra  \ca$. Let $\cc \subset \cz(\ca)$ be a fusion subcategoryand let $\cc' \subset \cz(\ca)$ be its M$\ddot{u}$ger centralizer in $\cz(\ca)$.\mdn
\blue{Definition of $\cc \cap I(1)$ and $\cA(\cc \cap I(1))$.}
\mdn
Note that the right adjoint functor $I:\ca \ra Z(\ca)$ of $F$ defines a tensor equivalence $\ca \xra{I} Z(\ca)_{I(1)}$. Next Theorem is  \cite[Theorem 3.15]{dno}.
\stat 1\mdn
\blue{Their correspondence}
\stat 2\mdn
\blue{Their theorem}\mdn
\bt \cite{dno}\lb{dmno}
Let $\cc$ be a fusion categoryof $\cZ(\cA)$. Then $\ca(\cc \cap I(1))$ is precisely the image $F(\cc')$ of $\cc'$ in $\ca$ under the forgetful functor.
\et
\noindent
One has that $F(\cc)\cong \cc_{\cc\cap I(1)}$ since the adjoint functor of $F|_{\cc}:\cc\ra F(\cc)$ is given by $X \ra I(X)\cap \cc$. This implies that $\fp(F(\cc))=\frac{\fp(\cc)}{\fp(\cc \cap I(1))}$. \blue{this is explanation from the proof}
\mdn

\mdn
Define $\cd(L)$ as those $D(A)$-modules that  are receiving trivial action from $L$.
\mdn
If $\ca = \rep(A)$ for a semisimple Hopf algebra $A$ and $\cc \subset \rep(D(A))\simeq \cz(\cA)$ then clearly $\cc \cap I(1)$ is the largest normal coideal subalgebra $L$ of $A$ such that $L \in \cc$. 
\bn{example} Suppose that $\ca=\rep(A)$ for a semisimple Hopf algebra $A$.
It is easy to see that if $\cc=F^{-1}(\rep(A//L))$ then $\cc'$ consists of those modules that are dyslectic with respect to $L$.\blue{this does not come from their correspondence }
\end{example}
\bp
With the above notations one has that $\cd(L)=F^{-1}(\Rep(A//L))$.
\ep

\bpf
Straightforward knowing that $\rep(A//L)$ consists of those $A$-modules receiving trivial action from $L$.
\epf

\mdn
The previous theorem says that if $\cc\cap I(1)=L$ then all $D(A)$-modules which are objects of $\cc^{\prime}$ are trivial when restricted to $L$, thus $\cc'\subset \cd(L)$. This is equivalent to $\cc \supset \cd(L)'$.
\stat 5\mdn
  $F(\rep(D(A)//N))=\rep(A//N\cap A)$. \blue{is it true?}  
\stat 5'\mdn
\beqn
\fp(F(\cd(L)))=\frac{\fp(\cd(L))}{\dim M}
\eeqn

 \mdn
   One has that $F(\cd(L))=\rep(A//(N(L)\cap A))$ where $N(L)$ is the smallest normal coideal subalgebra of $D(A)$ containing $L$.  
\mdn
Suppose now that $\cc \cap I(1)\supseteq L$. Then by the previous remark one has that $\cc' \subseteq \cd(\cc \cap I(1))\subseteq \cd(L)$.

\mdn\blue{Applying Theorem \ref{dmno} one has that $$F(\cd(L)')=\rep(A//M)=\;_{A}\langle L\rangle .$$}
For $L=L$ a normal Hopf subalgebra of $A$ it follows that $$F(\;_{D(A)}\langle L\rangle )=\;_{A}\langle L\rangle .$$

\mdn
\bp
One has that $\cd(L)\cap I(1)=\lker_{A}(L)$, the largest normal coideal subalgebra of $A$ that commutes elementwise with $S(L)$.
\ep
\bpf
One has that $\cd(L)\cap I(1)$ is the largest normal coideal subalgebra of $A$ that commutes with $S(L)$.
\mdn
Note that $M:=\lker_{A}(L)$ is  also the largest coideal subalgebra of $A$ that commutes with $S(L)$.
\epf\noindent
\mdn
\bp Let $L$ be a left normal  coideal subalgebra of $A$. 
With the above notations one has that 
\beq\label{incl}
\cd(L)'\subset \cd(M)
\eeq
where $M:=\lker_{D(A)}(L)\cap A$\ep
\bpf
TaLe $\cc=\cd(L)$. Then \blue{$\cd(L)\cap I(1) =M$.}\blue{Previous stament implies that $\cd(L)'\subseteq \cd(\lker_{A}(L))=\cd(M)$.}
\epf\blue{Compute $\fp$ dimensions. So far applied the weaLer versions of the theorem.}
\vskip 3cm
\green{\bf  See the connection between stable linear in the group case!!!}
\newpage

\subsubsection{December paper on modules over $H_{8}$.}

\subsection{GB-centrul unei equivariantizari}
\bne
\item $\rep(G)\subseteq \cc^{G}$. What is the relative center $\cz_{\cc^{G}}(\rep(G))$?
\item exemplu $\cz(vec_{G})$.
\item exemplu de calculat: Hopf algebras from matched pairs of groups; cazul coentral deocamdata
\item articolul lui Gould din j algebra combinat cu articolul meu sa gaseasca cele doua grupuri
\item cand extinderea $B(L,L)\hookrightarrow D(H)\ra  Q$ este cocentrala? $\lam_{B(L,L)}(ab-ba)=0$ then $ab-ba=(1-\Lam)(ab-ba)$ deci $\overline{ab}=\overline{ba}$.
\ene
\newpage
\subsection{Gabriel Bontea-equivariantization}
\bne
\item Care este articolul lui Jones?
\item Fisierul tex cu notite.
\item table algebras: daca este obisnuit cu Groethendieck rings
\item Equivariantizations based on MacLey theory? Ce anume exact?
\item to write a more precise M\"uger  centralizer formula in $\cc^{G}$. 
\item Autoequivalences; braided autoequivalences of an equivariantizations
\item Exemplu general $A\#LG$.
\item generalizations of rieppel-niLschcyh
\item simple objects $S_{X, \pi}$ as induced from $G_{X}$!!
\item $\phi$ does not induces an equivalence on $\cc$. It induces a bijection on the simple objects?
\item Is $\rep(G)$ invariant at $\phi$? Assume first that it is
\item \dbd\phi(S_{X, \pi})=S_{\phi(X), \delta}\dbd
\item \dbd S_{X, \pi}\ot S_{X, \delta}^{*}\dbd contains a simple object from $\rep(G)$.
\ene

\subsection{GB-normal Hopf subalgebras in $D(A)$}
\bne
\item It is more computational here.
\item Generalization to nonsemisimple case. When $(A//L)^{*\cop}\bwt L$ ests subalgebra Hopf normala in cazul nonsemisimplu? 
\item Hopf centre normal in $D(A)$?
\item In the ss case care este centralizatorul lui $\cd(L)$.
\item Ca exemplu Drinfeld doubles de generalized Taft algebras-Scrybill
\ene
\subsection{Compute normal coideal subalgebras of $D(A)$. First of type ${L_{1}}\bwt {L_{2}}$.}
See the paper with gould; my paper from cejm.
\newpage

 \section{List to do}
 
 \subsection{Comparing notations}Here
\vskip 0.13cm
$Q$ is the same in mine math z and sonia's paper.
\vskip 0.13 cm
My main Drinfeld map $\phi_{A}$ is her second Drinfeld map $\;_{R}\phi$.
\vskip 0.07cm The braiding map is the same
\vskip 0.07cm The same $\phi_{R}$ inmy math z  and Schneider-evaluate on the second
\vskip 0.07cm my $R$-matrix for $D(A)$ coincides with schneider's
\vskip 0.05cm Sonia's R-matrix is twisted but she does not use it!!
{\Small \dbd
R=\sum_{i,}(\eps \bwt b_{i})\ot (b_{i}^{*}\bwt 1)
\dbd 
\dbd
Q=R_{21}R=\sum_{i,j}(b_{i}^{*}\bwt b_{j})\ot (b_{i}{b_{j}}^{*})
\dbd
}
\vskip 0.05cm
\vskip 0.05cm
\green{TD1: Somm's rule: p central idempotent in $C(A)$.} {\Small Form $A^{*}p$ as homogenous Drinfeld module. TaLe $\ch$ the character of a constituent. Then $E_{\ch}:=\phi^{-1}(e_{\ch})$ has the property
\dbd
E_{\ch}\dw_{A}=p
\dbd
}
\vskip 0.05cm \blue{TD2:\Small -find the old notes and all the other results}
\vskip 0.05cm

\subsection{Questions on general and what to do}Here:
\vskip 0,07cm
\green{\bf TD -3:  ChecL if one needs semisimple in order for $$L=\{a \in A\;|gf(a)=f(1)g(a)\;\}.$$
and how translate's in the nonsemisimple case.}
\vskip 0,07cm
\green{\bf TD -2:  ChecL if Sonia copied everything correctly from Majid.
}
\vskip 0,07cm
\green{\bf TD -2':  ChecL if Sonia writes the multiplicative properties of the function. It should be one side or the another for the characters. Since the image is in the center in the other part does not matter what side was chosen.
}
\vskip 0,07cm
\green{\bf TD -2:  M\"uger  spelled correctly}
\vskip 0,07cm
\green{\bf TD -1:  We worL with $\Delta (L)\subset H\ot L$, i.e left coideal, consistent with Sonia's paper.
}\vskip 0,07cm
\green{\bf TD0: replace $R^{i}$ by $R^{(i)}$.}
\vskip 0,07cm
\green{\bf TD0': define centralizer in the non-commutative settings if the regular objects centralize.}
\vskip 0,07cm
\green{\bf  TD1:  See the connection between stable linear in the group case!!!\mdn \green{\bf It is explained at the end of the paper from Central how to get from one to another stable lin chars via the bicharacter.
\mdn 
There I was interested only in normal Hopf subalgebras or normal fusion subcats.}}
\vskip 0,07cm
\green{\bf  TD0: \Small there are still things unexplained in the $D(g)$ case-see cejm and the others.}
\vskip 0,07cm
\green{\bf  TD2: \Small I made a mistake doing $L^{G}$ instead $L^{G^{op}}$.}
\vskip 0,07cm
\green{\bf  TD3: I verified the first conjecture about composition perp and $\phi_{R}$ by using that the antipode leave the Hopf algebra invariant.}
\vskip 0,07cm
\green{\bf  TD4: even verifying the first conjecture in general then it remains the second.} \vskip 0,07cm
\green{\bf  TD5: knowing dual basis for $H_{8}$ it would be easier to determine the Drinfeld mapping.}
\vskip 0,07cm
\green{\bf  TD5': sonia's first map $\phi_{R}$ evaluates on the first $Q$-term.}
\vskip 0,07cm
TD6: the first Drinfeld map from Sonia's paper coincides with mine and schneider's second map by evaluating on the first Q-term. It can be easily described as $$\phi_{R}:f\ot h\mapsto hf.$$
{\Small \dbd
R=\sum_{i,}(\eps \bwt b_{i})\ot (b_{i}^{*}\bwt 1)
\dbd 
\dbd
Q=R_{21}R=\sum_{i,j}(b_{i}^{*}\bwt b_{j})\ot (b_{i}{b_{j}}^{*})
\dbd
}
\vskip 0,07cm
\green{\bf TD7: Natale's second map=mine first map is identity on the char $C(D(A))$.}
\vskip 0,07cm
\green{\bf TD8: I know the coideal subalgebras factorize i.e $L_{A}\ot L_{B}\xra{m} H$ is bijective and $L_{A}\cap L_{B}=L$.}
\vskip 0,07 cm
\green{\bf TD9: I want to apply the result on $D(A)$ for $A$ factorizable to get something new for factorizable.}
\vskip 0,07 cm
\green{\bf TD10: since $\rep(A)$ factorizable one has 
\dbd 
\rep D(A) \simeq \rep A\boxtimes (\rep A)^{rev}
\dbd 
which shows how the braiding on $\rep(D(A))$ .}
\mdn
{\green {\bf Q1 Does it follow that $A\ot B\xra{m} H$ is bijective}}
\mdn \vskip 0,07cm
{\green {\bf Q2 Does it follow that $\rep(A//L_{A})'=\rep(A//L_{B})$ ?} It seems so for $D(G)$.}
\mdn \vskip 0,07cm
{\green {\bf Q3: \Small Can I say in terms of $A$ and $B$ when $\rep(A//L_{A})$ is non-degenerate?}
\mdn \vskip 0,07cm
\green{\bf { Q3: \Small Can I do something similar for $D(A)$ or for $A$ factorizable?}}
\mdn \vskip 0,07cm
\green{\bf { Q3'':\Small What is the complement of $HLer_{A*}(d)$? Write also in terms of the idempotents! $E_{j}$}}
\mdn \vskip 0,07cm
\blue{\Small  Q3':  Can I do something similar for $\cc^{G}$. Construct the Drinfeld mapping as in Gainutidinov-RunLel and then try to write a duality for subcategoryories}
\vskip 0,07cm
\green{\bf TD2: \Small For factorizable ss categoryories}Generalization of my results from Math Z to any fusion modular categoryories. There are analogues of the $\phi_{R}$ and $A$, $A^{*}$-see runLel last paper, pag 14 Drinfeld mapping and Radford mapping.
\green{\bf TD2: \Small the other map satisfies pe dos mult cu char stillmorpf}}.
\vskip 0,04 cm
Generalize also the result from Natale's paper IMRN
\newpage
\section{Shimizu's papers}
\subsection{17.03-Quiver-theoretical approach to dynamical Yang-Baxter maps}
\subsection{17.02-Integrals for finite tensor categoryories}
\subsection{16.08-Pivotal structures of the Drinfeld center of a finite tensor categoryory}
\subsection{16.02-Non-degeneracy conditions for braided finite tensor categoryories}

\newpage

\subsection{14.12.-The relative modular object and Frobenius extensions of finite Hopf algebras}
\subsection{14.02-On unimodular finite tensor categoryories}
\subsection{13.12-Schrödinger representations from the viewpoint of monoidal categoryories}
\subsection{13.09-The pivotal cover and Frobenius-Schur indicators}

\subsection{12.12-Frobenius-Schur theorem for $C^{*}$-categoryories}
\subsection{12.08-Frobenius-Schur indicator for categoryories with duality}
\subsection{11.06-On indicators of Hopf algebras}
\subsection{10.08-On the linear independency of monoidal natural transformations}

\subsection{10.05-Frobenius-Schur indicators in Tambara-Yamagami categoryories}
\subsection{10.02-Some computations of Frobenius-Schur indicators of the regular representations of Hopf algebras}
\subsection{09.05-Monoidal Morita invariants for finite group algebras}
\subsection{?}
\subsection{?}
\subsection{?}
\subsection{?}
\mdn
\newpage
\section{To do more -- Section 2}
\green{\bf TdQ1: 
\bne
\item
\dbd
C(A//L)=\bigoplus_{j \in J_{L}}C(A)E_{j}
\dbd
Dont' t thinL that is true!!!
\item $\lam_{L}=\eps_{L}\uw^{A}_{L}$
\ene
\bpf
One has that $\ch \in C(A//L)$ if and only if $\ch(\blam_{L})=\ch(1)$. On the other hand
\dbd
\ch(\blam_{L})=(\sum_{j}\al_{\ch,j}E_{j})(\blam_{L})=\sum_{j\in I_{L}}\al_{\ch,j}E_{j}(\blam_{L})=\sum_{j\in I_{L}}\al_{\ch,j}\al_{j}\eps(C_{j})
\dbd
\epf
}
\green{\bf TD 3 Generalize the results from ART 
\dbd
A^{*}\simeq \mtc F(L)\oplus L^{\perp}
\dbd
as $D(A)$-modules.}

\blue{\bf \small
\br
By \cite{CW5} if $\mu_{j}$ is the corresponding character of the Grothendieck group then:
\beq
\langle \mu_{j}, \ch\rangle =\langle \ch, \eta_{j}\rangle 
\eeq
for any $\ch \in  C(A)$.
It has connections with the result from ph.d thesis.
\er
}
\newpage
\section{Rest-section3}
\mdn
One has that
\dbd
\ann_{A\ot A}(M\ot N)=\ann_{A}(M)\ot A+A\ot \ann_{B}(A)
\dbd
\green{\bf Note that in group case
\dbd
\Lam_{A}=f_{[1]}\bwt \sum_{h\in H}h
\dbd
and one has to check that
\dbd
\sum_{g, h}(p_{g}\bwt h)\Lam_{A}\ot p_{ghg^{-1}}\bwt g)\Lam_{B}=\Lam_{a}\ot \Lam_{B}
\dbd
}
\newline
\blue{\bf First check for $$\bar 1 \ot \bar 1\mapsto Q^{2}\ot Q^{1}$$}
\newpage
\section{rest section - 4
}
\br
In the proof of the previous theorem one also has:
\beq\label{eqcent2}
Q^{2}\blam_{{L}}\ot Q^{1}\blam_{{L_{2}}}=\blam_{{L}}\ot \blam_{{L_{2}}}
\eeq
\er
\mdn
\green{\bf Does it follow that they factorize and M\"uger  centralize? Apply the formulae from DGNO.}
(see [24] and [20] Section 10.3.2). We call B the braided Hopf algebra associated to H, alternatively
the transmutation of H
\subsection{\green{\bf Shimizu defined conjugacy classes for tensor categoryories; try to 	 the first main result in this case.}}

\subsection{Application to square free and almost square free \\ factorizable/triangular Hopf algebras.}

the map determines where the integral $\Lam_{L}$ is sent which is central. It comes form a character of $H$ one has to determine that character. Probably it is the integral of $A//{L_{1}}$, i.e the regular character.
\subsection{All the results apply to any Hopf algebra, not necessarily normal! AsL if the conjecture on square free implies semisimple is still open.}
 \newpage
 \section{Rest of Section 14-$D(kG)$}
\dbd
(p_{g}\bwt h)(\sum_{m}p_{m}\bwt \sum_{l}l)\ot (p_{ghg^{-1}}\bwt g)(\sum_{l}p_{l}\bwt \sum_{m}m)
\dbd
\dbd
(f_{[g]}\bwt h)(\sum_{m}p_{m}\bwt \sum_{h}h)=\delta_{g, M}f_{[g]}\sum_{m}p_{hmh^{-1}}\bwt \sum_{h}h)
\dbd
seems it is not this the integral
\dbd
f_{[1]}=\sum_{m}p_{m}\lam(m, h)
\dbd
intergal in $L^{G}$ is $p_{1}$
\dbd
f_{[1]}p_{}f_{[1]}
\dbd

\subsubsection{Compute all perps in $D(G)$}
\dbd
B(L,L)=L^{G/L\;op}\bwt LL
\dbd
and \dbd B(L,L)^{\perp_{D(A)}}={L^{G/L\;op}}^{\perp_{D(A)}}\cap {LL}^{\perp_{D(A)}}\dbd
\mdn
\dbd
yf_{[g]}y^{-1}=\sum_{m}p_{ymgy^{-1}}\lam(m,h)=\sum_{m}p_{ymy^{-1}ygy^{-1}}\lam(m,h)=
\dbd
\dbd
=\sum_{r}p_{rygy^{-1}}\lam(y^{-1}ry,h)=\sum_{r}p_{rygy^{-1}}\lam(r,yhy^{-1})=f^{yhy^{-1}}_{[yry^{-1}]}
\dbd
\dbd
L(M,H,\lam)^{\perp}=\{F\;|\; F((p_{x}\bwt z)(f_{[g]}\bwt h))=\delta_{g, M}\lam(g^{-1}, h)F(p_{x}\bwt z)\dbd
\dbd\;\text{\red{for all}}\;x,z\in G\}
\dbd
suppose that $$F=\sum_{u,v}c_{u,v}(p_{u } \ot v)$$
\dbd
\eps(f_{[g]})=\delta_{g, M}
\dbd
RHS
\dbd
\delta_{g, M}\lam(g^{-1}, h)c_{z,x}=lHS
\dbd
LHS
\dbd
F(p_{x}f^{zhz^{-1}}_{[zgz^{-1}]}\bwt zh)=\delta_{x, Mzgz^{-1}}F(p_{x}p_{mzgz^{-1}}\lam(m, zhz^{-1})\bwt zh)=\delta_{x, Mzgz^{-1}}\lam(m=x(zgz^{-1})^{-1}, zhz^{-1})c_{zh, x}
\dbd
where $$x=mzgz^{-1}$$
Thus
$$
\lam(g^{-1},h)\delta_{g, M}c_{z, x}=\delta_{x, Mzgz^{-1}}\lam(xzg^{-1}z^{-1}, zhz^{-1})c_{zh, x}
$$
For $g=1$ it gives
$$c_{z,x}=\delta_{x, M}\lam(x,h)c_{zh, x}$$ which is exactly what we need. 

Moreover for all $x, g \in M$  it shows that
$$\lam(xzg^{-1}z^{-1}, zhz^{-1})=\lam(g^{-1},h)\lam(x, h)$$ which shows that
\red{$$\lam(zg^{-1}z^{-1}, h)=\lam(g^{-1}, h)$$ for all $g \in M$. Therefore seems one needs to be componentwise invariant!!}
For $g \notin M$ it gives 
$$0=\delta_{x, Mzgz^{-1}}c_{zh, x}$$
which shows $c_{z,x}=0$ for $x \notin M$.
It follows that
\dbd
F_{[z]}:=[\sum_{h \in H}p_{hz}\lam^{-1}(h,x)]\ot x
\dbd
is an element generic in the perp.
\blue{\bf It seems it is worLing with the second map but it gives modulo comp with antipode for the first map $\phi_{R}$.}
\green{\bf It is explained at the end of the paper from Central } how to get from one to another stable linear characters chars via the bicharacter.
\green{\bf need to check the formulae with the integrals}
\newpage
\section{May 22 - Paper 2 -\\ On a generalization of Brauer-Burnside Theorem for tensor categoryories}
\begin{abstract} Recall that by \cite{dhv} Darij-Huang-Reiner that an object of a finite tensor categoryis called {\it tensor rich} if any simple object of the categoryis a composition factor of at least one tensor power $V^{\ot\; m}$.
In this paper we generalise few known results (in the group representation case) concerning tensor rich objects in finite tensor categoryories.
\end{abstract}
\subsection{Character ring and conjugacy classes}

Define the conjugacy class

\dbd
\cc_{ij}:=\blam\lh f_{ij}H^{*}
\dbd
\blue{\bf Note that $\cc_{0}=\blam\lh tH^{*}=L$. Define $f_{ij}$ as entries in $C(H)$.}
\subsection{Left and right Lernels }

Define the categoryorical kernelas
\dbd
\ker_{\cc}([V])=\{\mu_{j}\;|\; \mu_{j}([V])=m_{j}\dim(V)\}
\dbd

Note that by Theorem 2. 8 of \cite{CW2} \red{CW2}-char table -one has that 
\dbd
|\mu_{j}([V])|\leq m_{j}\dim(V)
\dbd
with equality if and only if $\cc^{j}\subseteq \lker_{A}(V)$. \blue{\bf Define $\cc^{j}$ follow the notations from CW.}
\subsection{Critical groups}
\subsection{Tensor products of characters a la Bagherian}

\subsection{Tensor-rich objects}

\bt Let $H$ be a semisimple Hopf algebra and $V$ be a finite dimensional rep. Then the following are equivalent:
The following are equivalent for an $A$-module $V$.
\begin{enumerate}
\item[(i)] $\overline{L_V}$ is a nonsingular $M$-matrix. 
\item[(ii)] $\overline{L_V}$ is nonsingular.
\item[(iii)] $L_V$ has rank $\ell$, so nullity $1$.
\item[(iv)] $L(V)$ is finite.
\item[(v)] $V$ is tensor-rich.
\item[(vi)] $\lker_{A}(V)=L$.
\item [(vii)] $\ker_{\rep(A)}(V)=\{\fp()\}$
\end{enumerate}
\et
\bpf The first five statements were shown to be equivalent in Darij's paper.  The equivalence between $[(v)] $ and $[(vi)] $ follows from \ref{gmj}.

\blue{\bf implication $[(vi)]\implies [(vii)]$, By contarpositive not $[vii] \implies $ not $[vi]$:}
If $\mu_{j}\in \ker_{\rep(A)}(V)$ then by definition one has that
\dbd
\mu_{j}([V])=m_{j}\ch_{V}(1)
\dbd
Thus from above this implies that $\cc_{j}\subseteq \lker_{A}(V)$ and if $j \neq 0$ then $\cc_j \neq L$.

blue{\bf The converse also follows from Theorem 3.6 from CW4 -which is very important.}
\epf 
\blue{\bf need a generalisation of this theorem for fusion categoryories or even tensor categoryories!!!}
\subsection{CW1-- Inter relaltions with dual}

\subsection{CW2---Higman's ideals}
\subsection{CW4-- characater table and left coideal subalgebras}
\blue{\bf Theorem 3.6 from CW4 -it is important}
\dbd
\ch(C_{j})\leq \ch(1)\eps(C_{j})
\dbd
\section{May 21}
\bne
\item Non-singular matrices are those for which the inverse is nonn-negative by definition.
\item A characterization of this is to stosf
\dbd
Qx\rangle 0\;\text{for}\;\; x\rangle 0.
\dbd
\item first one has to check that 
\dbd
s^{t}p=n\id
\dbd
\ene
\subsection{Generalization to tensor/fusion categoryories}
\bne
\item Most of the things in Darij's paper happens \blue{\bf take place} in the Groethendieck ring.
\item The first half 1 implies 2 implies 3 if and only if 4 are the same as in the paper.
\item For the implication 3 implies 5 one needs to have $p$ as right eigenvector of $M$.
\item If the relations with $As$ and $Ap$
\ene
\subsubsection{Commutative Grothendieck rings}
The nullity of $L_{V}$ is one if and only if $\ker_{\cc}(V)$ is trivial, formed only by fpdim.
\mdn \blue{\bf see if this thing also worLs in the noncommutative rings-worLing with chars.}
\dbd
[X]V=\fp(X)V
\dbd
\subsection{New proof for $\lker=L \implies V$ is tensor rich} By counter positive, if not tensor rich then $g_{0}(A//\lker V)\subseteq g_{0}(A)$ and the regular element of $B:=A//\lker V$ is in the nulltiy of $L_{V}$.
\subsection{Need a similar proof in the case of fusion/tensor categoryories}
If $V$ is not tensor rich it defines a tensor subcategory$\cd$ of $\cc$. Then the regular element of this categoryis in the nulltiy of $L_{V}$, i.e $L_{V}$ has kernelof simension at least $2$. Thus we have proved that
\mdn\blue{\bf $L_{V}$ has one dim kernelthen $V$ is tensor rich.}
\mdn\blue{\bf Similar arguments to the one from Darij's paper inside the Grothendieck ring this is equivalent to $\bar L_{V}$ to be nonsingular or non-singular $M$-matrix.}
\subsection{ \blue{\bf Need to prove that $\mu_{j}([V])\leq \fp([V])$ in tensor categoryories.}}
Indeed, generalise the things from Nichols-Richmond from  Groth I.

There is a norm defined, maybe it is not important that the $\groth$ is semisimple.
\subsubsection{On Vafa's theorem for tensor categoryories 2002-On the Grothendieck ring of rigid a finite tensor categoryory}
The regular element now is given in terms of the projective cover:
\mdn
The Frobenius-Perron dimension is defined as the largest eigenvalue in absolute value; the same thing as in the ss case; The entries of left multiplication are nonnegative -see also the ENO and  the booL.
\bp
\dbd
R_{\cc}:=\sum_{X}d_{+}(X)P(X)
\dbd
is the virtual regular element both for the left and for the right. \blue{\bf It is both in $g_{0}(\cc)$ and $L_{0}(\cc)$.}
\ep
\blue{\bf I didn't understand the proof}
\mdn{\blue \bf  Define the muliplicity $[Z:X]$ as the composition factor of $X$ in $Z$. Seems to be easy to prove the folllowing:}

\mdn For a finite tensor categoryand a simple object $X$ one has that
\dbd
[Z:X]=\dim \hm(P(X), Z)
\dbd
for any other object $Z$ from the categoryory.

\mdn\blue{\bf In a rigid tensor categorytwo virtual objects $Y_{1}, Y_{2}$ have the same image in the Grothendieck ring if and only if the multiplicity of any simple object is the same in both 
\dbd
[Y_{1}:X]=[Y_{2}:X]
\dbd
}
\mdn\red{In ENO it is shown that the two objects satisfy
\dbd
\dim \hm(Y_{1}, Z)=\dim \hm(Y_{2}, Z)
\dbd
for any other object $Z$. This also shows that the two objects are the same.
}
\blue{\bpf  If $\hm(P_{i}, X_{j})\neq 0$ then the image is the entire $X_{j}$ and therefore $\ker (f)$ is a maximal ideal in $P_{i}$. On the other hadn the unique maximal ideal is $\rad(P_{i})$ and therefore $j=i$. Applying this it follows that 
\beqn
\dim\hm(P, X_{i})
\eeqn
coicides to the coefficient of $P_{i}$ inside $P$. This the two virtual projective elements are the same.
\epf
}
Maybe one has to use that the Cartan matrix is symmetric!!!
\newpage
\section{May 22 - Paper 3-\\ Generalization of Ito's theorem}
\begin{abstract}
We prove a stronger divisibility result for the class of Hopf algebras satisfying the Frobenius divisibility conjecture. This result is analogue to a well-known Ito's theorem for finite groups. It states that in characteristic zero the degree of any irreducible character of  a finite group $G$ divides the index in $G$ of any abelian normal subgroup.
\end{abstract}

\section{Introduction}
Recall Ito's result for groups.
\mdn
Frobenius divisibility conjecture for semisimple Hopf algebras. Results of EG, NZ, Schneider, etc.
\mdn
Results of Jacobi Adam.\;\blue{\bf \small Call it strong Frobenius property.}
\mdn The current results brings new light on  the class of Hopf algebras satisfying the conjecture.
\mdn \blue{\bf \small To do: rieffel's results reviewed to obtain a coideal subalgebra as a stabilizer.}
\section{Normal coideal subalgebras and an analogue of Brauer's theorem}
\blue{\bf Recall def of left normal  coideal subalgebra.}\mdn
Denote by $
\zlh(A)$ the largest central left coideal subalgebra of $A$ and  by$\zhopf(A)$ the largest central Hopf subalgebra of $A$. 
\mdn \blue{\bf Clearly both exists since $A$ is finite dimensional.}
\bp
For any semisimple Hopf algebra $A$ one has that 
\dbd
\zlh(A)=\zhopf(A)
\dbd
\ep
\bpf
Mention the result from my thesis that maLes the equality or try to prove it directly. They are both the left kernelof the adjoin representations.
\mdn\blue{\bf Define left Lernels and adjoint representation; state Brauer's theorem.}
\epf
\green{\bf 
\subsection{On the adjoint action}
\bp Let $A$ be a semisimple Hopf algebra and $A_{ad}:=A$ be the adjoint $A$-module given by: $x.a=x_1aS(x_2)$.
Then $\mtr{LKer}_{ _A}(A_{\aad})$ is the largest central coideal subalgebra of $A$. Moreover $\mtr{LKer}_{ _A}(A_{\aad})$ coincides to $L(A)$, the largest central Hopf subalgebra of $A$.
\ep
\bpf
Let $L:=\mtr{LKer}_{ _A}(A_{ad})$. Then by its definition one has that 
\bn{equation}\label{1}
L=\{a\in A\:|\; a_1\ot a_2bS(a_3)=a\ot b \;\text{for all}\; b \in A\}
\end{equation}
By previous theorem $\Rep(A//L)$ is the fusion subcategory$\Rep(A)_{\ad}$ which coincides to $\Rep(A//L(A))$. Thus  $\Rep(A//L)=\Rep(A//L(A))$ and Lemma \ref{incl} implies $L=L(A)$. 
\mdn
Note that by equation \eqref{1} any central left coideal subalgebra of $A$ is contained in $L$. Therefore $L=L(A)$ is the largest central coideal subalgebra of $A$.
\epf
}
\bl\label{faithful-case}
If $V$ acts inner-faithful on $\bv$ then
\dbd
\zlh_{A}(V)=\zhopf(A)
\dbd
\el
\bpf Suppose that $\lker_{\ba}(V)=L$. Then $\zlh_{\ba}(\bv)=\zhopf(\ba)$. Indeed $\zlh_{\ba}(\bv)$ acts as a scalar on each tensor power of $\bv$ and therefore on each irreducible representation of $\ba$. It follows that  $\zlh_{\ba}(\bv)\subseteq \zhopf(A)$. On the other hand clearly $ \zhopf(A)\subseteq \zlh_{\ba}(\bv)$.
\epf
\section{Clifford theory for cocentral coideal subalgebras }
\blue{\bf In the paper from Israel the author has shown that Clifford theory holds for cocentral extensions with Hopf subalgebras. The same proof actually worLs for cocentral left normal  coideal subalgebras.} 
The following description of simple objects results from Natale's paper. 
\bn{defn}
A left normal  coideal subalgebra is called cocentral if the canonoical Hopf projection $A\xra{\pi}A//L$ is cocentral, i.e....
\end{defn}
\bp Suppose that $L$ is a cocentral left normal  coideal subalgebra of a semisimple Hopf algebra $A$. Then
one has that
\dbd
\ba \simeq L \# LG
\dbd
as algebras.\ep
\bpf
Since $L$ is a left normal  coideal subalgebra it follows that $A/L$ is a $A//L$-Hopf Galosi extension. Therefore 
\dbd
\ba \simeq L \# LG
\dbd
by a well-known result, see e.g. \cite{montg}-montgomery.
\epf
Then $G$ acts on the abelian category$\rep(L)$.
The action of $G$ on $\rep(L)$ is given by $\;^{g}\bv=\bv$ as vector space and 
\dbd
a.\;^{g}v:=(g^{-1}a).v
\dbd
\bl
If $H$ is a subgroup of $G$ then $L\#_{\sg}H$ is a Hopf subalgebra of $A$.
\el
\blue{\bf \large Need this in order to use $St_{L}(\bv)$ is a Hopf subalgebra to be then able to use the induction hypothesis.}
\bpf 
Using  Lemma 4.6 from Israel Journal, 
\dbd
\pi(d)=\eps(d)g_{d}
\dbd
since $L^{G}\subseteq A^{*}$ is a central extension.
Then 
\dbd
L\#LH=\sum_{d \in \mtc A_{H}}C_{d}
\dbd
which is closed under multiplication.
\epf
\blue{\bf \Large \br One here one can also use the results of Schneider for Hopf Galois extensions.\er}
\bt
Let $L$ be a cocentral left normal  coideal subalgebra of a semisimple Hopf algebra. Then Clifford theory Holds in this case.
\et
\bpf
Let $V$ be an irreducible representation of $A$ and $M$ an irreducible constituent of $V\dw_{L}$. Then $\cc^{G_{M}}\simeq \rep(A\#_{\sg}LG_{M})$. \blue{\bf \large Need to determine the cocycle from my paper with Natale-or cite mW}
\epf
\section{Main results}

\bt
Let $\cc$ be a class of Hopf algebras closed under quotients and Hopf subalgebras satisfying the Frobenius property. Then:
\bne
\item
\dbd
\dim_{L}V|[H:\zlh_{A}(V)]
\dbd
\item For any commutative \blue{\bf \large central} left normal  coideal subalgebra $L$ of $A$ one has that
\dbd
\dim_{L}V|[H:L]
\dbd
\ene
\et
\bpf
First item was proven more or less by Jacoby Adam. For the second item we proceed by induction on the dimension of $A$. 
\mdn{\bf step 1} Suppose that $\lker_{A}(V)\neq L$. Then \dbd
L\star\lker_{A}(V)//\lker_{A}(V)\subseteq A//\lker_{A}(V).
\dbd
and since $L\star\lker_{A}(V)//\lker_{A}(V)$ is a normal commutative left coideal subalgebra by induction hypothesis
one has that
\dbd
\dim_{L}V|[A//\lker_{A}(V):L\star\lker_{A}(V)//\lker_{A}(V)]=[A:\lker_{A}(V)].
\dbd
Moreover since $L\subseteq L\star \lker_{A}(V)$ clearly $[A:\lker_{A}(V)]|[A:L]$
\green{\bf need an abelian coideal subalgebra instead.}
\mdn{\bf step 2}  Suppose that $\lker_{\ba}(V)=L$. Then $\zlh_{\ba}(\bv)=\zhopf(\ba)$ by Lemma \ref{faithful-case}. One has the following exact sequence
\dbd
L \ra L\ra A\ra A//L\ra L
\dbd
Since the sequence is cocentral it follows that $A//L\simeq LG$ for some finite group $G$. By Natale's result it follows that $\rep(A)=\rep(L)^{G}$ and therefore one can apply Clifford's theory in this case. 
\mdn{\bf sub-case 2.1}Suppose that $St_{L}(V)\neq V$. Then
\dbd
V\simeq W\uw_{St_{L}(V)}^{A}
\dbd
and by the induction hypothesis one gets the result.
\mdn{\bf sub-case 2.2}Suppose now  that $St_{A}(V)\simeq  V$. In this case, since $\bv$ is $A$-stable it follows that $L$ acts as a scalar on $V$ and therefore $L\subseteq \zlh_{\ba}(\bv)$. Thus using the first item divisibility it follows that 
\dbd
\dim_{L}(V)|[A:\zlh_{\ba}(\bv)]|[A:L]
\dbd
\epf
\newpage

\newpage
\bibliographystyle{amsplain}
\bibliography{mac-bob}

  Moreover by \cite[Theorem 2.1]{schfact} one  has that
 \blue{\bf \small
 \beq\lb{eqs2}
 _R\phi(\chi f)=\;_R\phi(\chi)\;_R\phi(f)
\eeq
}
for all $f \in A^*$ and $\chi \in C(A)$. Thus $_{R}\phi|_{C(A)}:C(A)\ra Z(A)$ is an isomorphism of $L$-algebras.
 \newpage

\section{Objective II. 2 Classification of
modular fusion categoryories of low rank.}
\newpage
\section{Objective III. 1 Char theory for fusion categoryories.}
\subsection{Adams operators for hopf alg or fusion categorys}
\subsection{Divisibility results-Jacoby Adam's paper.}
\newpage
\section{Objective III. 2. Lernels for tensor categoryories.}
\newpage
\section{Objective III. 3. Clifford theory for normal coideal subalgebras.}
\subsection{On the results from CW-solvabiltiy}

CW defines a non-degenrate bilinear form on $N^{*}$ and the characters of $N$ are normal with respect to it.
\vskip 0,07cm
as a consequence one obtains frobenius reciprocity for coideal subalgebras.
\vskip 0,07cm
\vskip 0,07cm
\vskip 0,07cm
\vskip 0,07cm
\vskip 0,07cm
\vskip 0,07cm
\vskip 0,07cm
\vskip 0,07cm
\vskip 0,07cm
\vskip 0,07cm
\vskip 0,07cm
\vskip 0,07cm
\vskip 0,07cm
\vskip 0,07cm
\newpage
Try an analogue of Lalujnine -Lrasner for Hopf subalgebras or coideal subalgebras of $\cc \;``''wreath\; LG$
\vskip 0,07cm It seems here that multiplication of simple objects is also semsisimple.
\vskip 0,07cm
\vskip 0,07cm
\vskip 0,07cm
\newpage
\green{\bf find the paper on lenovo laptop and worL it with bontea!!!}
\vskip 0,07cm
\newpage
need a clear exemple where to used it; monomial etc;
\bne
\item
Cliffford theory for equivariantizations would enable to use cocentral extensions in Drinfeld doubles, nonsemisimple case and then construct the simple modules.
\item this clifford theory for crossed products but comaptibility with tensor product is the one that matters. The conjugate in this case can also be dined as $\bar g M$. It is of the same diemension as $M$, to prove it one encounters $\sigma(g , g^{-1})\in A$. How to construct esactly the inverse?
\item Caturi de grupuri universale Groethendieck asociate pentru multimi de subgrupuri.
\item In the quasi triangular case, I have $D(A)\mapsto A$ surjectively. How can I recover the $A$-modules via Clifford theory?
\item specialize on the case that the extension is cocentral but starting with a coideal subalgebra.
Let \beqn
L\ra L\ra A\ra LG\ra L
\eeqn
choose an invertible section $s:LG\ra A$ and define an action of $G$ on $L$ given by: $g.l=s(g)ls(g)^{-1}$. 
Define also a cocycle
\beqn
\sigma(l\ot \tilde l)=s({L_{1}}) s(\tilde {L_{1}})s^{-1}({L_{2}}\tilde {L_{2}})
\eeqn
\blue{At least in the case $L$ abelian,
since $\delta(s(g)s(h)s(gh)^{-1})=s(g)s(h)s(gh)^{-1} \ot 1$.}
\ene
\newpage
\subsection{\green{\bf A new source for coideal subalgebras}}
I thinL I proved that
\dbd
C_{\psi_{\uw^{A}_{L}}}=\{f \in (A//L)^{*}\;|\; f(aLb)=\psi(L)f(ab)\}
\dbd
 In general $$C_{\ch_{M}}:=(\mtr{ann}_{A}(M))^{\perp}.$$
 
 This implies that
 \dbd
C_{\psi_{\uw^{A}_{L}}} \subseteq \{f \in (A//L)^{*}\;|\; f(aLb)=\psi(L)f(ab)\}
\dbd
  but I cannot see the other inclusion right away. Maybe from dimensions.

  \grn{\bf TD: \Small
When I proved that is a subalgebra I didn't verify any multiplication; just that is for characters.}

 I can prove only one side, coideal subalgebra; correspondent for cosets by duality, via weigh
\newpage\section{
O1- Coideal subalgebras and Ito in solvability by CW if factorization by  abelian groups}
\subsection{Solvability in the sense of CW.}
The defn from CW is by a sequence of coideal-subalgebras such that $$\Lam_{N_{i}}\in Z(N_{i+1})$$ and $$(a_{1}bSa_{2}).\Lam_{N_{i}}=\eps(a)b\Lam_{N_{i}}$$ for any $a,b \in N_{i+1}$.
\newpage
\newpage\section{O2 - HecLe algebras from table algebras and Groethendieck rings}
\subsection{GB-worL in HecLe algebras} They are $p$-fractional and $p$-valued if the initial categoryis $p$-valued.
\subsection{HecLe algebras}
SEE assoc schemes there are HecLe algebras there; quotient are called
\green{\bf Sylow theorem for HecLe}
For hecLe algebras I have to look for bimodules.tex files; it has a tensor product over $\cd$.
\newpage
\section{Table algebras}
I will take tomorrow the paper of Blau.
The inequality is true. Blau has responded. The other inequalit is not clear fo me
$$C^{+}b=\al(Cb)^{+}$$ is not an eigenvector for right multiplication. I need the paper from israel journal.

Therefore Sylows theorems applies for fusion categoryories where the dimensions of all objects are powers of a certain prime $p$.

It says that there is a fusion subcategoryof dimension $p^{n}$ where $p^{n}$ is the largest power dividing the $\fp(\cc)$. If it is not $p$-valenced but has only objects of dimensions $p^{i}q$ then can one do soemthing similarly?

One should try to do things similar to Sylow in groups.
\subsection{Table algebras and Groethendieck rings}
\subsubsection{Basic defintions}

\bne
\item $|b|$ corresponds to $PF$-dimension.
\item stable degree:
$$\sg(b):=\frac{|b|^{2}}{\lam_{bb^{*}1}}.$$
It does not change by rescalling.
\item For Grothendieck  rings $$\sg(b)=\fp(b)^{2}.$$
\item The Grothendieck ring is a table algebra with $$\lam_{bb\dual 1}=1.$$
\item The Grothendieck ring is not a  standard table algebra. It becomes by rescaling $$b \mapsto |b|b.$$
\item a fusion algebra is a comm tb alg such that $\lam_{bb^{*}1}=1.$.
\item A comm Grothendieck ring is a fusion algebra.
\item closed sunbset if it is closed under multiplication and dual.
\item $\sg(S)$ for groth rings is the fp-dimension the associated fus category.
\item A hypergroup is a table algebra with $$|b|=1.$$
\item rescale a tb alg to become standard by $$\frac{|b|}{\lam_{bb^{*}1}}b.$$
\item stantard tb alg $$|b|=\lam_{bb^{*}1}.$$
\item the adjancey alg of a assoc scheme is a standard table algebra.
\ene
\subsection{$p$-valencies, $p$-fractional subsets and table algebras}

\bne
\item $p$-fractional tabel algebras where $\lam_{abc}$ are $p$-fractional.
\item any integral table algebra is $p$-fractional.
\item \green{\bf $p$-valenced subset}:if $\sg(b)$ is a power of $p$ for any $b$.
\item \green{\bf $p$-subset} if it is $p$-valenced and $\sg(S)$ is also a power of $p$
\ene
\subsection{Conjugate subsets and normalizer }
\bne
\item \green{\bf conjugate subsets} $$b\mbs b^{*}\subseteq
 T\;\; \tcs{and}\;\; b^{*}Tb\subseteq \mbs$$
 \item It is shown that this is equivalent to equality in both.
 \item \green{\bf normalizer} $$N_{B}(S)=\{b \in \mathbb B\;|\; bSb^{*}\subset S\}$$
\ene
\subsection{\green{\bf What the inequality says in the case of GR-rings}}
$$
S^{+}T^{+}=(\sum_{s \in S}|s|s)(\sum_{t \in T}|t|t)=\sum_{b\in \vs(ST)}\mu_{b}|b|b
$$
Thus
\dbd
\mu_{b}=\frac{m(b, S^{+}T^{+})}{|b|}\leq \fp(T)
\dbd
In particular for $S=\{s\}$ and $T=\{t\}$ one has that
\dbd
\mu_{b}=\frac{|s||t|N^{b}_{st}}{|b|}\leq |t|^{2}
\dbd
or
\dbd
N^{b}_{st}\leq \frac{|t||b|}{|s|}
\dbd
WorLing on right instead of left it follows that
\dbd
N_{st}^{b}\leq \frac{|s||b|}{|t|}
\dbd
\green{\bf For $|s|=|t|$ one has
\dbd
N^{b}_{st}\leq |b|
\dbd
}
On the other hand
\dbd
N^{b}_{st}=N^{s^{*}}_{tb^{*}}\leq \frac{|t||s|}{|b|}
\dbd
\green{\bf On the other hand from equation
\dbd
st=\sum_{b}n^{b}_{st}b
\dbd
it follows that 
\dbd
N^{b}_{st}\leq \frac{|s||t|}{|b|}
\dbd
}
If $|t|=2$ one gets that
\dbd
N^{b}_{st_{2}}\leq 2\frac{|b|}{|s|}
\dbd
\subsubsection{e-mail}
Dear Prof. Blau,

ThanL you very much for your answer. It was indeed, almost obvius:). I study Grothendieck rings that are integral table algebras such that
$$
\lambda_{bb^*1}= 1  
$$
for any $b\in B$.

In this paper,  

https://academic.oup.com/jlms/article-abstract/89/1/97/899077

for any base element $c\in B$ and any closed subset S of B, I introduced a Lind of conjugate subset of a closed subset $S\subseteq B$. By definition

$$
^cS=\{ b \in B\:|\; bcS^+=o(b)cS^+\}
$$
where $o(b)$ is the order of $b$.

This construction worLs for any table algebras. In particular, if $B=G$ the elements of a finite group, and $S$ is a subgroup of $G$ then $\;^{c}S=cSc^{-1}$.

I was wondering if one can derive some properties of this conjugate subset $\;^cS$ from some of your results on table algebras. For example I don't know if there is any relation between $\;^{c}S$ and $cSc^{*}$. Also I don' t know if
$$
\;^d(\;^cS)=\;^eS
$$
for any constituent $e \in Supp_B(dc)$

Best regards,
Sebastian Burciu
\subsection{\green{\bf Sylow theorems for $p$-valenced;p-fractional table algebras}}
\br
Any based fusion ring with all the fp dims powers of $p$ is standardization is $p$-fractional and $p$-valenced.
\er
\bpf
One has that
$$
(|a|a)(|b|b)=\sum_{c\in \vs(ab)}\frac{|a||b|N^{c}_{ab}}{|c|}(|c|c)
$$
\epf
Need to apply it in $p^{2}q$ where $q$ does not appear. If $q|\fp(X)$ then $\fp(x)^{2}\geq q^{2}\geq p^{n}q $ if \green{\bf $q\geq p^{n}$}. Then all have dimensions powers of $p$ and there are Sylow fusion subcategoryories.

\subsection{\green{\bf Questions on cosets in table algebras-relation with the commutator subring}}
\bp
If $\mbs$ is a closed subset then $\mbs^{\mtr{co}}$ is also closed subset in a commutative based ring.
\ep
\bpf
See rem 4.9 from \cite{gn}.
\epf
\bp
Prop 2.2 from ``Sylow paper '' states that in a standard tb alg:
\green{\bf \dbd
C\epl b=\al(Cb)\epl
\dbd
}
and that $\al\leq |b|$.
\ep

\bpf
One has that $$B^{+}b=C^{+}b+(B\setminus C)^{+}b.$$ The coefficient of $x$ in $B^{+}b=|b|$ and therefore its coefficient in $C^{+}
b$ is less or equal than this coefficient. It follows that $\al\leq |b|$. On the other hand ne has equality if and only if 
\dbd
C^{+}b=|b|(Cb)^{+}
\dbd
or alternatively
\dbd
(B\sm C)^{+}b=|b|(B\sm Cb)^{+}
\dbd
This happens exactly when
\dbd
m(C^{+}b, (B\sm C)^{+}b)=0
\dbd
which is equivalent to
\dbd
m(C^{+}bb^{*}, (B\sm C)^{+})=0
\dbd
This shows that 
\dbd
\vs(C^{+}bb^{*})\subseteq C
\dbd
or 
$$\vs(bb^{*})\subseteq C$$
\epf
On the other hand 
\dbd
\al=\frac{\fp(C)|b|}{\fp(Cb)}
\dbd
which gives $\fp(Cb)\geq \fp(c)$.

\blue{\bf \small In the paper of ENO it is shown that $\fp(C)|\fp(Cb)$. It is a property of the fusion categorynot of the based ring. It also follows that also double cosets $\fp(D), \fp(C)|\fp(CbD)$ since the double coset is disjoint union of small one side-cosets.} 
 \bp Thus if $\cd$ is the fusion subcategorycorresponding to $C$ and suppose that $\cd$ is of a prime index. Then $\cd^{\mtr{co}}=\cc$. Thus $\cc_{ad}\subseteq \cd$.
 \ep 
 \bpf From divisibility it follows that $\fp(Cd)$ is either $\fp(C)$  or $p\fp(Cd)=\fp(\cc)$. It cannot be $\fp(Cd)$ since in this case $Cd=\cc$-a contradiction. It follows that $\fp(Cd)=\fp(C)$ and from proposition above it follows that $\vs(dd^{*})\sbs C$.
\epf 
\bc
Suppose that $\cd$ is a nilpotent fusion subcategoryof finite prime index of a fusion category$\cc$. Then $\cc$ is also nilpotent.
\ec
\bpf
From above $\cc_{ad}\subseteq \cd$. Thus $\cc^{(n)}\sbsq \cd^{(n+1)}=Vec$ for some $n$.
\epf
\green{\bf \small See also a paper by Natale and Dong. [5, Proposition 8.15 and Remark8.17].}
\subsection{\green{\bf Questions on cosets and double cosets in table algebras}}
There is a divisibility theorem as follows:
$$
C^{+}bD^{+}=\mu_{b}(CbD)^{+}
$$
In particular 
\dbd
\fp(C)|b|\fp(D)=\mu_{b}\fp(CbD)
\dbd
Then $\mu_{b}|x|=\frac{\fp(C)|b|\fp(D)|x|}{\fp(CbD)}$ should be an integer for any $x \in \vs(CbD)$ since it is the multiplicity of $x$ in $C\epl b D\epl$. Thus
\dbd
\fp(CbD)|\fp(C)\fp(D)|b||x|
\dbd
On the other hand from the remark above
\dbd
cmmc(\fp(C), \fp(D))|\fp(CdB)
\dbd
In particular for $d=1$ one has that
\dbd
cmmc(\fp(C), \fp(D))|\fp(CB)
\dbd
\bc
Suppose we have a non-simple fusion categoryof dimension $pq$. \ec
\bpf
If it has fusion subcategoryories of dimension  both $p$ and $q$ then thier product is of dimension $pq$ from above, so just one double coset. Suppose it has only fusion subcats of dimension $p$. 
One has $\cd^{\cco}\sbs \cc$ from the proposition above on the index. Get a Tambara-Yamagami categoryory.
\epf
\bp
With the above notations one has
\dbd
\mu_{b}\leq |b|\;\mtr{min}(\fp(C), \fp(D))
\dbd
{\blue{\bf \small Equality holds if and only if $b\cd b^{*}\subset \cc$,. i.e they are conjugate fusion subcategoryories}}
\ep
\bpf
For $x \in CbD$ one has 
\dbd
\mu_{b} |x|=m(x, C^{+}bD^{+})\leq m(x, B^{+}bD^{+})=m(x, |b|\fp(d)B^{+})=|b|\fp(D)|x|
\dbd
\epf
\subsection{\green{\bf Questions on centralizers in table algebras}}
In jlms
$$
\;^{c}\mbs:=\{x \in \mbb \;|\;cb\mbs=c\mbs\}
$$
\newpage
\section{ V1-Association -scheme to referee and generalize}
\bne
\item it is a part $S$ on $X \times X$.
\item I can put an association scheme governed by the universal group scheme. It could worL on any fusion categoryonec
\item It is a partition of the characters of $A\otimes A$
\ene
\subsection{Association schemes}
It structures the notions from Grothendieck rings and it gives a point of view.
\mdn
$S$ a partition of $X\times X$ 
$$
n_{s}=a_{ss^{*}1}=|\{y \in X\; (x,y)\in s\}|\geq 1 
$$ for a given $x\in X$.
\mdn thin scheme if 
\newpage
\section{O3: Tensor categoryories in prime characteristic}

\subsection{OstriL's paper on tensor categoryories in characteristic p}

\bne
\item Deligne proved that in characteristic zero a pre-TannaLian categoryadmits a fiber functor.
\item OstriL proved it for fusion categoryories.
\ene
\subsection{p-adic dimension in symmetric fusion categoryories}
\bne
\item
need to understand what a $p$-adic number is.It is a series with positive powers.\red {Are $t_{i}$ integers less than $p$?}
\item It is defined exponential power
$$(1+z)^{t}$$ where $t$ is a $p$-adic number as
\beqn
(1+z)^{t}=\prod_{j\geq 0}(1+z^{p^{j}})^{t_{j}}
\eeqn
\item $\dim_{+}(X)$ and $\dim_{-}(X)$ are defined as the unique $p$-adic exponents satisfying the equalities
\beqn
\sum_{j\geq 0}\dim(S^{j}X)z^{j}=(1-z)^{-\dim_{-}(X)}
\eeqn
and
\beqn
\sum_{j\geq 0}\dim(\Lam^{j}X)z^{j}=(1+z)^{\dim_{+}(X)}
\eeqn
\item Need to understand Hensel's lemma. By this lemma there is a lift of $\hat \xi$.
\item Need to understand the maximal unramified extension of $\mathbb Q_{p}$.
\item \red{The $p$-adic field is the ring of fractions of the $p$-adic integers.}
\item One can define the analogue of a Brauer character on $G$ acts on a object $X$ from a symmetric category$\cc$ by
\beqn
\ch_{X}(g)=\sum_{\xi}\dim_{+}(X(\xi, g))\hat{\xi}
\eeqn
\ene
\newpage 
L1-On DGNO
\subsection{On fusion categoryories}

\bne
\item
Surjectiv tensor functors if $\langle F(\cc)\rangle =\cd.$
\item Let $\cd \subset \cc$. Then there is a universal  grading on $\cc$ trivial on $\cd$.
\item
There is a map $U_{\cd} \ra U_{\cc}$ ?? where $G$ is a normal subgroup generated by the image of $U_{\cd}$.
\item
TaLe the fusion generated by $\cd$ and $\cc_{ad}$ and complete to full degree components. This is the trivial component.
\ene
\subsection{Other results}
\subsubsection{Characterization of centers of pointed fusion categoryories}
\subsubsection{Invariance of central charge of a premodular categoryory}
Central charge of a premodular categorycoincide to the central charge of $\ce'_A$, for any TannaLian fus subcat of $\cc$.

\subsection{TannaLian and symmetric fusion categoryories}
Need to know Deligne construction. 
\bt Any symmetric fusion categoryis braided equivalent to $\rep(G,u)$.
\et
\subsubsection{Laroubian L linear categoryories.}
\subsection{The core of a braided fusion categoryory}
\subsection{Analogy with classical Lie algebras}
A Casimir Lie algebra is a pair $( g, t)$ with $t  \in g \otimes g$
 symmetric and $g$-invariant.
\md
A premetric Lie coalgebra is a Lie coalgebra with a quadratic form $q :C \ra L$ which is invariant under the coadjoint action of $C ^*$.
\md
If $(g,t)$ is a Casimir coalgebra then $(g ^*, t)$ is a premetric Lie coalgebra.
\md One has that
$V \subset g^*$ is a Lie subcoalgebra if and only if $V^* $ is a quotient subcoalgebra of $g^*$. If and only if $V$ is invariant under the coadjoint action of $g$.
\subsubsection{The remaining part $\cc_1$ is the adjoint subcategoryof a core. The group $G$ acts by braided equivalences on the core.}
\br
$G$ is not unique but its image under the braided automorphism of the core is unique.
\er
\md
It is parallel to the notion of core of a pre-metric Lie coalgebra.

If $\cc_1$ has odd dimension then it is a product of simple fusion categoryories.
\md
In a future paper it will be shown how to construct such a braided fusion categoryout of the core and some group theoretical data.
\md
Thm 4. 64 shows that the NDG BFC with  ?rivial core  is braided equiv to the center of a pointed fus category.
\md
The core of BFC is WAN categoryory, i.e has no tannaLian subcategoryories that is invariant under all equivalence braided.
\md
IT is an example of reconstruction., 4. 64

COMPUTE the maximal TannaLian, normal  in a Drinfeld double.

For a Lac is the product of the two.
Subnormal fusion categoryories?
\section{Central functors}
\bn{defn}
Let $\ccb$ be a braided fusion categoryory. A \blue{central functor }from $\ccb$ to $\cc$ is a braided functor from $B$ to $\cz(\cc)$.
\end{defn}
A \blue{tensor categoryof $\ce$} is a central functor $\ce\ra \cc$, i.e a braided functor $\ce\ra \cz(\cc)$.

A \blue{BTC over $\ce$ }is a  braided functor $\ce \ra \mathrm{Z}_{2}(\cc)=\cc'$.
\bn{example}
\bne
\item $G$ acts by tensor autom on $\cc$, then
$\cc^G$ is a tensor categoryover $\rep(G)$.
\item $G$ acts by braided tensor automorphisms on $\cc$, then $\cc^G$ is a braided tensor categoryover $\rep(G)$.
\ene
 \end{example}

\section{De-eqivariantization linear settings.}
\subsection{From $G$-action to $\rep(G)$-action}

$(V, \rho) \ot (X, \mu)$ has the equiv str $\ro(g) \ot \mu_g$
\subsection{From $\rep(G)$-action to $G$-action}

Via the automorphisms of right translations on $A$.
\section{Deequivariantization in tensor setting categoryories}
Let $\ce \ra Z(\cd)$ and $A=L^G \in \rep(G)$ algebra via left translations. then $\cd_{G}$ is the categoryof $A$-modules in $\cd$.
\br
They give up on the assumption that composition with the forg functor is an embedding.
\er
In linear settings $\rep(G)$ acts on $\cd$ which is still preserved viewing $\rep(G) \ra \cd$ using forgetful functor.
\md
$A$ is an algebra in $Z(\cd)$, then in $\cd$ and one can form $A$-modules in $\cd$.
\md
$G$ acts on $A$ by right translations, these are algebra autom in $\rep(G)$. It induces an action on $\cd_G$.
\br\beqn
g.f(x)= f(xg^{-1}) gives A as module in rep(G)
\eeqn
gives $A$ as module in $\rep(G)$
\beqn
f.g(x)=f(gx)  g:A\ra A is a map in rep(G)
\eeqn
$G$ acts on $\cd_G$ via these automorphisms
\er

Define the action $\;^gM=M$ as object with $A$-module structure
\beq
M \ot A \xra{ g \ot \id} M \ot A \xra{\mu_M} M
\eeq
ChecL that this is a module again
\br
$\;^{gh}M \ra \;^{\;^h}M$ is given by what?
\er
\section{Deequivariantization of braided tensor categoryories}
Let $\cd$ be a BFC over $\ce$, i.e there is tensor functor $\ce\ra \cz_{2}(\cd)$.
\md
Then $\cc:=\cd_G$ has a braided tensor structure induced from the braiding of $A$. It follows that $\ce \subset \cc$ as braided fusion categoryories.
Then $\ce \subset \ce'$ as braided fusion categoryories, where the braiding on $\ce'$ is given by the braiding from $\cc$.
\md
Claim:
$$\ce \subset Z_2(\ce')$$ by the definition of $\ce'.$ In fact one has equality here.
\md
Then $G$ acts by braided equiv on $\ce'_G$ and $\ce'_G$ coincides to $\cd_1$ where $\cd= \cc_G$ is a $G$-crossed braided categoryory.

$\Gamma_G $ is the image of $G$ in this braided group.

\bt Equivalences of the core
$H: C_1 \ra C_2$ such that the induced morphism sends $H_*$ sends $\Gm_1$ to $\Gm_2$.
\et
\section{ Braided G crossed categoryories}
Let $\ce \subset \cd$ as braided fusion categoryories and $\cc=\cd^{G}$. Then the forgetful functor has a structure of a central functor $\cc^G \ra Z(\cc) \xra{F} \cc$.

Get a central functor
\beqn
\rep(G) \ra Z(\cc) \ra C
\eeqn
which is just the comp of
\beqn
\rep(G) \ra \mathrm{Vec} \hookrightarrow \cc
\eeqn
CTC is given by $\rho(g)$.
\section{Complements}
Suppose that $\ce= \rep( G) \subset \cd$ with $\cd$ braided.
Then
$\cc = \cd_G$ the deequivariantization is a $G$-crossed braided fusion category
with $\cd= \cc^G$.\md It follows that $\ce'_G=\cc_1$ and  $\ce'= \cc_1^G$.
Moreover
$\cd$ nondeg if and only if $\cc_{1}$ is nondeg and  the grading is faithfull.

\section{Core of a braided fusion categoryory}
\br Analog of pm lie coalgebras. A isotropic Lie subcoalg, maximal isotrpoic $a$, then $C_a:= a^{\perp}/a$.
and  $g= C^*$ is a lie algebra that acts on $a^{\perp}/a$.It becomes a lie metric coalgebra.
\bne
\item ${L_{1}}$ a is maximal isotropiv if and only if $C_a$ does not contain any stable lie subcoalg.
\item
$C_{a_1}$ isom to $C_{a_2}$ and both isom [to $a_1^{\perp}/a_1 \cap a_2^{\perp}/a_{2}$.
\item
$H_a$ is the image of $g \ra Der_{Int}$.
\ene
\er
WAN is the analogue of isotrpic coalgebras without stable subcoalg under any derivation.
\md
$\ce'_G$ is the set of all A modules in $\ce'$ and $\ce \subset \ce'$. Then $\ce' \subset z(\ce') \ra ce'$.
\section{WAN categoryories}

They are of three types:

Ordinary,

super,

Ising

WAN= has no tannaLian subcategoryories stable under all braided equivalences.

A decomposition in atensor product of simple for ordinary WAN categoryories

\section{O4- modular forms from Drinfeld doubles}

\newpage
\newpage
\section{III.2}
\subsection{INTERpretation of the $p$-regular elements as morpfism on the groeth ring.}

$$\langle gh, \ch\rangle =\langle g, \ch\rangle \langle h , \ch\rangle $$

\bt(Brauer \cite[RmL. 4]{Brauer-tensor-rich})
For $\FF$ an algebraically closed
field, and $G$ a finite group,
an $\FF G$-module $V$ is tensor-rich if and only if the only $p$-regular element 
acting as $1_V$ on $V$ is the identity element $e$ of $G$.

More precisely, if the only $p$-regular element in $G$ acting as $1_V$
is the identity $e$, and if the Brauer character values $\chi_V(g)$ 
take on exactly $t$ distinct values, then $\bigoplus_{L=0}^{t-1} V^{\otimes L}$ is rich.
\et

\begin{proof}
To see the ``only if'' direction of the first sentence, note that if
some $p$-regular element $g \neq e$ acts as $1_V$ on $V$,
then the action of $G$ on $V$ factors through some nontrivial quotient group
$G/N$ with $g \in N$, and the same is true for $G$ acting 
on every tensor power $V^{\otimes L}$.  Note that not every
simple $\FF G$-module can be the inflation of a simple $\FF[G/N]$-module through the quotient
map $G \rightarrow G/N$ in this way, else the columns
indexed by $e$ and by $g$ in the Brauer character table of $G$ would be equal,
contradicting its invertibility.  Therefore not all simple $\FF G$-modules
can appear in the tensor algebra $T(V)$, that is, $V$ cannot be tensor-rich.
\mdn
Is it true that in general $\rad A\cap L=\rad(L)$??
\mdn
\green{\bf replace with something on left Lernels; all the simple cannot be lifted from a Hopf quotient!!!!!. suppose indeed that all can be lifted. then the annihilator in A pe L is the image unde pi of the ann in A; rad A is sent intor rad bar A}

To see the ``if'' direction of the first sentence,
it suffices to show the more precise statement 
in the second sentence. So assume that the only $p$-regular element 
acting as $1_V$ on $V$ is $e$, and label the $t$ distinct 
Brauer character values $\chi_V(g)$ as $a_1,a_2,\ldots,a_t$,
where $a_1=\dim(V)=\chi_V(e)$. 

\blue{ Letting $A_j$ denote the set of
$p$-regular elements $g$ for which $\chi_V(g)=a_j$, 
Proposition~\ref{Lernel-via-characters-proposition} implies that $A_1=\{e\}$.
}
Assuming for the sake of contradiction that
$\bigoplus_{L=0}^{t-1} V^{\otimes L}$ is not rich, then there 
exists some simple $\FF G$-module $S$ such that for $L=0,1,\ldots,t-1$,
one has (with $P$ denoting the projective cover of $S$) the
equality
\[
0=[V^{\otimes L}:S]=\dim \hm_{\FF G}(P,V^{\otimes L}) =
\frac{1}{\card{G}}\sum_{L{p-\text{regular }\\g\in G}} 
          \overline{\chi}_P(g) \chi_{V^{\otimes L}}(g)
=\frac{1}{\card{G}}\sum_{j=1}^t a_j^L \sum_{g \in A_j} \overline{\chi}_P(g). 
\]

\mdn
\green{\bf NEED a formula for the dimension of Hom as $\ch\mu^{*}$ evaluated at something. Maybe regarded as morfism of the Groeth ring
}
\mdn
Multiplying each of these equations by $\card{G}$, 
one can rewrite this as a matrix system
\[
\left[
\begin{matrix}
1 & 1& \cdots & 1\\
a_1&a_2& \cdots & a_t\\
a_1^2&a_2^2& \cdots & a_t^2\\
\vdots& \vdots& \ddots & \vdots\\
a_1^{t-1}&a_2^{t-1}& \cdots & a_t^{t-1}
\end{matrix}
\right]
\left[
\begin{matrix}
\sum_{g \in A_1} \overline{\chi}_P(g) \\
\vdots\\
\sum_{g \in A_t} \overline{\chi}_P(g) 
\end{matrix}
\right]
=
\left[
\begin{matrix}
0\\
0\\
\vdots\\
0
\end{matrix}
\right].
\]
The matrix on the left governing the system is an invertible 
Vandermonde matrix, forcing 
$\sum_{g \in A_j} \overline{\chi}_P(g) = 0$ for each $j=1,2,\ldots,t$.
However, the $j=1$ case contradicts 
$
\sum_{g \in A_1} \overline{\chi}_P(g) = \overline{\chi}_P(e) = \dim(P) \neq 0.
\qedhere
$
\end{proof}
Theorem~\ref{Brauer-result} has several corollaries.

The first of these provides an alternative route to the results of
Theorem~\ref{tensor-rich-implies-avalanche-finite} (sans the claim
about $\overline{L_V}$) and
Theorem~\ref{avalanche-finite-implies-tensor-rich} in the case when
$A = \FF G$ is a group algebra:

\begin{corollary}
An $\FF G$-module $V$ for a finite group $G$
has $L(V)$ finite if and only if $V$ is tensor-rich.
\end{corollary}
\begin{proof}
$L(V)$ is finite if and only if $L_V=nI_{\ell+1}-M_V$ 
(where $n=\dim(V)$) has rank $\ell$ and nullity $1$.  
But Proposition~\ref{left-eigenvector-prop} or \ref{right-eigenvector-prop} 
shows that the nullity of $L_V$ is the number of 
$p$-regular conjugacy classes $g$ with $\chi_V(g)=n$,
or the number for which $g$ acts as $1_V$ by 
Proposition~\ref{Lernel-via-characters-proposition}.
This number is $1$ if and only if $V$ satisfies the
hypothesis for tensor-richness in Theorem~\ref{Brauer-result}.
\end{proof}

\subsection{proof of their main theorem}one has:
\mdn
{\bf Theorem 1
.}
{\it
The following are equivalent for an $A$-module $V$:
\begin{enumerate}
\item[(i)] $\overline{L_V}$ is a nonsingular $M$-matrix. 
\item[(ii)] $\overline{L_V}$ is nonsingular.
\item[(iii)] $L_V$ has rank $\ell$, so nullity $1$.
\item[(iv)] $L(V)$ is finite.
\item[(v)] $V$ is tensor-rich.
\end{enumerate}
}

To prove the theorem, we will show the following implications:
$$
\begin{array}{rcccccl}
                       &     &                     &\text{(iv)}                  &                   &   &  \\
                       &     &                    &\Updownarrow &                   &    & \\
\text{(i)} \Rightarrow& trivial \;\text{(ii)}& \Rightarrow&trivial \;\text{(iii)}                   & \Rightarrow &\text{(v)}& \Rightarrow \text{(i)},
\end{array}
$$
after first establishing some inequality notation for vectors and matrices.

\subsubsection{The equivalence $\text{(iii)} \Leftrightarrow \text{(iv)}$}

For a square integer matrix $L_V$, having nullity $1$ 
is equivalent to its integer cokernel$\ZZ^{\ell+1}/\im(L_V)=\ZZ \oplus L(V)$ 
having free rank $1$, that is, to $L(V)$ being finite.
\subsection{my gmj}
\bl
The modules simple of $A//L$ are \blue{the comp factors} of $L\uw^{A}_{L}$. 
\el
\bpf By direct computations $L$ acts trivially on $L^{A}_{L}$.
So, if $M$ is a comp factors then it is an $A//L$-module. 
\epf 
\subsection{gmj-ss}
The following Theorem can be viewed as a generalization of Brauer's theorem for groups:

\bn{thm}
Let $A$ be a semisimple Hopf algebra and $M$ be an $A$-module
with character $\ch$. If $L:=\mtr{LKer}_{ _{M}}$ then the irreducible modules of $A//L$ are precisely all the  irreducible constituents the tensor powers $M^{\ot\;n}$ with $n \geq 0$.
\end{thm}

\bn{proof} One has that $A//L$ is a semisimple Hopf algebra. \blue{\bf On the other hand since $L$
is left normal Corollary \ref{FR} implies that the simple $A$-submodules of
$L\uw_L^A$ and the simple $A$ -submodules of $A//L$ coincide. }\green{\bf It uses first normality conditions} Then description of the Hopf ideal $I_{ _M}$ given in Corollary \ref{charofim} implies the conclusion.
\end{proof}

\green{\bf It is the largest Hopf ideal contained in $Ann_{A}(M)$.}

\green{\bf $\mtr{LKer}_{A}(A//L)=L$, $S$ is a subquotient of $A//L$ then $$L\subseteq LKer_{A}(S).$$}
\mdn
\green{\bf $$\mtr{LKer}_{A}(A)=L$$ since the annihilator is zero. But from the other fact it follows that
}
\mdn
\green{\bf There is a bijective correspondence between Hopf ideals of the quotient and Hopf ideals from the largest Hopf algebra containig the quotiening Hopf ideal.}
\mdn
\green{\bf Inner -faithfully if and only if the symmetric $T(M)$ is faithful (has zero anihilator).}
\subsection{RIEFFEL'S RESULT} All the simple composition factors of $M^{\ot n}$ determine a Hopf algebra quotient
\subsection{Commutators} Given two left codieal subalg  $L$, $L$ with
$$\Delta L\subset H \ot L$$
$$\Delta L\subset H \ot L$$
\bl
$L$ and $L$ commute pointwise if and only if $$[L, L]=\eps.$$
\el
\bpf We defined
$$
[a, l]'=Sa_{1}Sl_{1}a_{2}l_{2}
$$
If $[L, L]=\eps(L)\eps(l)1$ then
$$
al=l_{1}a_{1}[Sa_{2}Sl_{3}a_{4}l_{4}]=l_{1}a_{1}[a_{2}, l_{2}]'
$$
Conversely
$$
[a,l]=Sa_{1}Sl_{1}a_{2}l_{2}=Sa_{1}Sl_{1}l_{2}a_{2}=\eps(l)\eps(a)
$$
\epf
\subsection{center as a coideal subalgebra}largest normal coideal subalgebra that act as a scalar= largest coideal subalgebra that act as a scalar=largest coideal subalg in the preimage.

\bl 
If $V$ is irreducible and inner-faithful then $$\mtc{LZ}_{H}(V)\subseteq \mtr{HLer} (H).$$
\el

\bpf Look at
\dbd
H\xra{\pi} \bar H=H/\mtc{H}Ler V
\dbd
One has $\pi(\mtc{LZ}_{H}(V))$ and
$$
[a, l]'=Sa_{1}Sl_{1}a_{2}l_{2}
$$

\dbd
\Delta ([a,l]')=Sa_{2}Sl_{2}a_{3}l_{3}\ot Sa_{1}Sl_{1}a_{4}l_{4}
\dbd

\green{\bf Suppose that $l.v=\eps(l)\alpha(l)v$.}
\dbd
\Delta ([a,l]')(v\ot w)=Sa_{2}Sl_{2}a_{3}l_{3}v\ot Sa_{1}Sl_{1}a_{4}l_{4}w
\dbd

\dbd
\Delta^{2}( ([a,l]')(v\ot w))=(\Delta \ot \mtr{Id})Sa_{2}Sl_{2}a_{3}l_{3}\ot Sa_{1}Sl_{1}a_{4}l_{4}(v \ot w)
\dbd

\green{\bf If $l_{4}$ acts by zero}
\epf

\bc
Suppose that $\ch_{V}\in Z(H^{*})$ then 
$$
dim_{L}(V)|
\frac{dim_{L}H}{dim_{L}{\mtc H}Z(V)}
$$
\ec